\documentclass[openany]{memo-l}
\usepackage{amsmath, amsthm, amssymb, amsfonts, amscd, eucal, xspace}
\usepackage{epic}
\usepackage{eepic}

\usepackage[matrix,line,arrow]{xy}

\numberwithin{section}{chapter}
\numberwithin{equation}{chapter}

  \newcommand{\eg}{e.g.\@\xspace}
  
  \newcommand{\ie}{i.e.\@\xspace}
  
  \newcommand{\etc}{etc.\@\xspace}

  \newcommand{\role}{r\^ole\xspace}

  \newcommand{\fg}{f.g.\@\xspace}

\newcommand{\el}{l_{2q+2}(\Lambda)}

\newcommand{\kreckf}[5] {(#1 \stackrel{#2}{\longleftarrow}#3 \stackrel{#4}{\longrightarrow}{#1}^*,#5)}
\newcommand{\skreckf}[4] {(#1 \stackrel{#2}{\longleftarrow}#3 \stackrel{#4}{\longrightarrow}{#1})}
\newcommand{\skreckfa}[4] {(#1 \stackrel{#2}{\longleftarrow}#3 \stackrel{#4}{\longrightarrow}{#1}^*)}

\newcommand{\kreckfb}[6] {(#1 \stackrel{#2}{\longleftarrow}#3 \stackrel{#4}{\longrightarrow}{#5},#6)}

\newcommand{\kreckfs}{\kreckf{F}{\gamma}{G}{\mu}{\bar\theta}}
\newcommand{\skreckfs}{(F \stackrel{\gamma}{\longleftarrow}G \stackrel{\mu}{\longrightarrow}{F}^*)}

\newcommand{\bN}{\mathbb N}
\newcommand{\bZ}{\mathbb Z}

\newcommand{\bZt}{\bZ/2\bZ}
\newcommand{\bQ}{\mathbb Q}

\newcommand{\bC}{\mathbb C}

\newcommand{\svec}[1]{\left( \begin{smallmatrix} #1 \end{smallmatrix} \right)}
\newcommand{\mat}[1]{\begin{pmatrix} #1 \end{pmatrix}}

\newcommand{\mt}{\longmapsto}
\newcommand{\ra}{\longrightarrow}
\newcommand{\isora}{\stackrel{\cong}{\ra}}
\newcommand{\homra}{\stackrel{\simeq}{\ra}}

\theoremstyle{definition}
\newtheorem{defi}[subsection]{Definition}
\newtheorem{deflem}[subsection]{Definition and Lemma}
\newtheorem{rem}[subsection]{Remark}
\newtheorem{ex}[subsection]{Example}

\theoremstyle{plain}

\newtheorem{thm}[subsection]{Theorem}
\newtheorem{lem}[subsection]{Lemma}
\newtheorem{prop}[subsection]{Proposition}
\newtheorem{cor}[subsection]{Corollary}
\newtheorem*{prop*}{Proposition}
\newtheorem*{thm*}{Theorem}
\newtheorem*{cor*}{Corollary}

\DeclareMathOperator {\im} {im}
\DeclareMathOperator {\rk} {rk}

\DeclareMathOperator {\coker}{coker}
\DeclareMathOperator {\Hom}{Hom}

\DeclareMathOperator {\Tor}{Tor}

\newcommand {\cone}{\mathcal{C}}

\makeindex

\begin{document}
\frontmatter
\title{The algebraic theory of Kreck surgery}[Kreck surgery]
\author{J\"org Sixt}[J. Sixt]
\address{c/o University of Edinburgh\\School of Mathematics\\Mayfield Road\\Edinburgh EH9 3JZ\\United Kingdom}
\email{sixtj@yahoo.de}

\begin{abstract}
Surgery theory is a classification technique for manifolds of dimension bigger than $4$ which was
developed in the 1960s. The traditional Browder-Novikov-Sullivan-Wall-theory decides
whether an $(n+1)$-dimensional normal cobordism $(e,f,f')\colon (W,M,M') \ra X$ with $f$ and $f'$
homotopy equivalences is cobordant rel$\partial$ to an $s$-cobordism. There is an obstruction
in a group $L_{n+1}(\bZ[\pi_1(X)])$ which vanishes if and only if 
this is possible. Algebraic $L$-groups have been extensively studied and computed. 
For a ring $\Lambda$ with involution,
$L_{n}(\Lambda)$ is a Witt group of quadratic forms if $n$ is even and a Witt group of quadratic
formations if $n$ is odd. A formation is a quadratic form with a pair of lagrangians, \ie two
hyperbolic structures which arise from an expression of an odd-dimensional manifold as a twisted double.

In the 1980s M. Kreck generalized Wall's original approach by dealing with
cobordisms $(e,f,f')\colon (W,M,M') \ra X$ of {\it normal smoothings} in which
$f$ and $f'$ are only $\left[\frac{n+1}{2}\right]$-equivalences. 
There is an obstruction in a monoid $l_{n+1}(\bZ[\pi_1(X)])$ which is {\it elementary} 
if and only if that cobordism is cobordant rel$\partial$ to an $s$-cobordism. 
The $l$-monoids are little understood algebraically and there are no computations of them.

This memoir studies the algebraic properties of $l_{2q}(\Lambda)$ (\ie $n+1=2q$).
$l$-monoids are equivalence classes of generalized formations which we call {\it preformations}.
Preformations are algebraic models of highly-connected bordisms between highly-con\-nect\-ed 
odd-dimensional manifolds.
We introduce three obstructions to an element $z\in l_{2q}(\Lambda)$ being elementary.
Firstly, it is shown that every elementary $z \in l_{2q}(\Lambda)$ has
a stable {\it flip-isomorphism}. A flip-isomorphism of a preformation can be thought of as a kind of algebraic isomorphism
between the two ends of the bordism associated to that preformation.
In certain cases there is a close relationship between flip-isomorphisms 
and isometries of the topological linking forms of $M$ and $M'$. 
Secondly, every flip-isomorphism of $z$ determines an
asymmetric form which vanishes in the asymmetric Witt group $LAsy^0(\Lambda)$
if $z$ is elementary. At last, a quadratic signature can be defined for certain kinds of flip-isomorphisms.
$z$ is elementary if and only if one of these quadratic signatures is zero in $L_{2q}(\Lambda)$.
\end{abstract}

\maketitle

\tableofcontents

\mainmatter
\chapter*{Introduction}

{\bf All manifolds shall be compact and smooth.}

Surgery theory was pioneered in the famous paper of Kervaire and Milnor
\cite{KeMi63} on the classification of homotopy spheres. Surgery on
high-dimensional manifolds was then developed by Browder, Novikov,
Sullivan and Wall, culminating in the general theory of Wall's book
\cite{Wall99}.  A modified theory which needs weaker prerequisites has
been presented and applied by M. Kreck (see
\cite{Kreck99}).  It assigns to any cobordism of {\bf normal
smoothings} an element in the monoid $l_{2q}(\Lambda)$ which is {\bf
elementary} if and only if that cobordism is cobordant rel$\partial$ to an
$s$-cobordism.

The main aim of this treatise is to provide obstructions which can
help to decide whether an element in $l_{2q}(\Lambda)$ is elementary
or not.

The first step is to determine whether such an element has a {\bf
flip-isomorphism} which in certain simply-connected cases is the same as the existence
of an isometry of linking forms of the manifolds that one
wants to classify.

The second step is to compute {\bf asymmetric signatures} \ie certain
elements in the asymmetric Witt-group $LAsy^0(\Lambda)$ for each
flip-isomorphism.  In the case of an elementary element all these
signatures vanish.  If the element in $\el$ allows linking forms it turns out that the
asymmetric signatures only depends on a choice of isometry of those linking
forms.

Alternatively one can define {\bf quadratic signatures} for a certain
class of flip-isomorph\-isms. An element in the $l$-monoid is elementary
if it allows such a special kind of flip-isomorphism and if the
quadratic signature in the quadratic Witt-group $L_{2q+2}(\Lambda)$
vanishes for at least one of them.

The quadratic signature is technically more difficult to handle than
its asymmetric sister but they are related via the canonical map
$L_{2q+2}(\Lambda)\ra LAsy^0(\Lambda)$.

Nicer results can be obtained when one deals with all those elements
in $l_{2q}(\Lambda)$ that are represented by non-singular formations
(the objects which help to define the odd-dimensional
$L$-groups). They are the obstructions to a Kreck surgery problem where \eg
all normal smoothings are in fact normal maps and if $M$ and $M'$
are closed. Then the definition of asymmetric signatures still
requires the existence of a flip-isomorphism but they will be
independent of the particular choice. One can even tame the quadratic
signatures: they exist for all flip-isomorphisms and considerable
simplifications can be achieved.

This paper is based on the author's 2004 University of Edinburgh
doctoral thesis.

\section{Classical and Kreck's Surgery Theory}
\label{introsurgsec}

In the following the surgery theory of Browder-Novikov-Sullivan-Wall and 
its modification by Kreck are outlined.
For a more elaborate account see Section \ref{formgeosec}.

\index{Wall surgery theory|see{Surgery theory}}
In the 1960s C.T.C. Wall and others developed surgery theory as a tool
to find out when a (normal) homotopy equivalence $f\colon M
\rightarrow N$ of $n$-dimensional manifolds is homotopic to a
diffeomorphism. There is a first obstruction which decides whether $f$
can be extended to a degree 1 normal cobordism
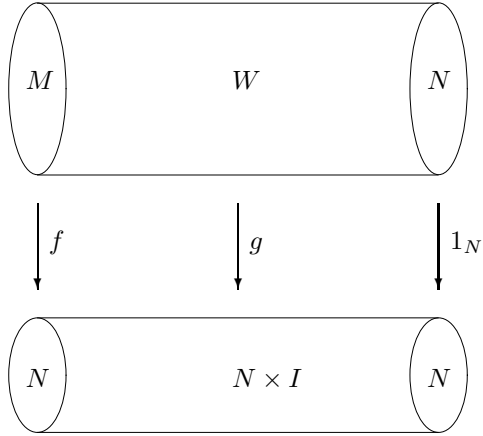
\begin{center}
\setlength{\unitlength}{0.006in}
\begin{picture}(750,400)(0,0)
\put(200,300){\ellipse{50}{150}}
\put(550,300){\ellipse{50}{150}}
\path(200,375)(550,375)
\path(200,225)(550,225)
\put(190,300){$M$}
\put(370,300){$W$}
\put(540,300){${N}$}

\put(200,200){\vector(0,-1){75}}
\put(375,200){\vector(0,-1){75}}
\put(550,200){\vector(0,-1){75}}
\put(210,160){$f$}
\put(385,160){$g$}
\put(560,160){$1_N$}

\put(200,50){\ellipse{50}{100}}
\put(550,50){\ellipse{50}{100}}
\path(200,0)(550,0)
\path(200,100)(550,100)
\put(190,40){$N$}
\put(370,40){$N\times I$}
\put(540,40){$N$}
\end{picture}
\end{center}
$$ (g,f,1_N) : (W, M, N) \rightarrow N \times (I,\{0\},\{1\})$$
\ie $f$, $g$ and $1_N$ are covered by maps of the stable normal
bundles and they map fundamental class to fundamental class.  A second
obstruction determines - in a more general setting - whether a normal
cobordism into a finite geometric Poincar\'e space $X$
\begin{eqnarray}
\label{normalmap}
    (g, f_0, f_1): (W, M_0, M_1) \rightarrow X \times (I,0,1)
\end{eqnarray}
with homotopy equivalences $f_i$ is cobordant rel$\partial$ to
a homotopy equivalence \ie an $s$-cobordism. In that case and for
$\pi_1(X)$ with vanishing Whitehead groups\footnote{Of course there is
also a version for other fundamental groups. Then we have to replace
homotopy equivalences by simple homotopy equivalences and use the simple
$L$-groups.}  the $s$-cobordism theorem tells us
that $M_0$ and $M_1$ are diffeomorphic. In the following we will
concentrate on the second obstruction.

\index{Surgery theory!Wall!even-dimensional}
Assume the dimension of $W$ is even ($\dim W=2q+2$). Surgery below the
middle dimension allows us to replace $g$ by a $(q+1)$-equivalence.
Then define $(K_{q+1}(W), \lambda, \nu)$ with 
$K_{q+1}(W)$ the {\bf kernel module}\index{Kernel modules}\index{K@$K_*(M)$} which is the homology of the
induced map $\widetilde g\colon \widetilde W \ra \widetilde X$ of
the universal covers with twisted coefficients. The form
$\lambda$ is induced by
the Poincar\'e duality on $W$ and $X$ and $\nu$ is the self-intersection
map. Together they are a non-singular $(-)^{q+1}$-quadratic form which vanishes in the Witt group
$L_{2q+2}(\bZ[\pi_1(X)])$ if and only if (\ref{normalmap}) is
cobordant rel$\partial$ to an $s$-cobordism.

\index{Surgery theory!Wall!odd-dimensional}
If the dimension of $W$ is odd ($\dim W=2q+1$), the construction of an
obstruction is slightly more complicated. One way is to extend
$g\colon W \ra Y=X\times I$ to a {\bf presentation} (see \cite{Ran01a}), that is a
$(2q+2)$-dimensional degree $1$ normal cobordism 
$(V,W, W')$
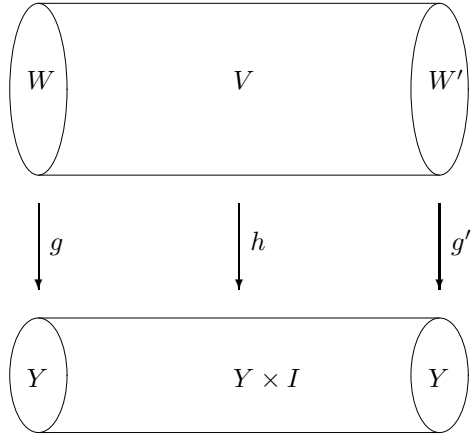
\begin{center}
\setlength{\unitlength}{0.006in}
\begin{picture}(750,400)(0,0)
\put(200,300){\ellipse{50}{150}}
\put(550,300){\ellipse{50}{150}}
\path(200,375)(550,375)
\path(200,225)(550,225)
\put(190,300){$W$}
\put(370,300){$V$}
\put(540,300){${W'}$}

\put(200,200){\vector(0,-1){75}}
\put(375,200){\vector(0,-1){75}}
\put(550,200){\vector(0,-1){75}}
\put(210,160){$g$}
\put(385,160){$h$}
\put(560,160){$g'$}

\put(200,50){\ellipse{50}{100}}
\put(550,50){\ellipse{50}{100}}
\path(200,0)(550,0)
\path(200,100)(550,100)
\put(190,40){$Y$}
\put(370,40){$Y\times I$}
\put(540,40){$Y$}
\end{picture}
\end{center}
such that also $h$ and $g'$ are highly-connected.
Define a {\bf non-singular formation $\skreckfs$} with 
\begin{eqnarray*}
\gamma\colon G=K_{q+1}(V) &\ra& F=K_{q+1}(V,W)\\
\mu \colon G=K_{q+1}(V) &\ra& K_{q+1}(V,W')\cong K^{q+1}(V,W)=F^*
\end{eqnarray*}
A non-singular formation is a tuple $\kreckfs$ such that
$\svec{\gamma\\\mu}\colon G \ra
H_{(-)^q}(F)$ is an inclusion of a lagrangian into a hyperbolic
form and $(G,\gamma^*\mu,\bar\theta)$ is a $(-)^{q+1}$-quadratic form. 
A formation determines a class of automorphisms of
the hyperbolic form sending $F$ to $G$ which represents
the obstruction in Wall's original version of 
odd-dimensional surgery theory (see \cite{Ran01a} for details).

The obstruction formation lives in some kind of Witt-group
$L_{2q+1}(\bZ[\pi_1(X)])$ of non-singular formations $\skreckfs$.
Again the obstruction vanishes if and only if surgery is
successful in producing an $s$-cobordism \ie if and only if there exists a
presentation as before with $f'\colon W'\ra Y$ a homotopy equivalence
(\ie an $s$-cobordism).

\index{Surgery theory!Kreck}
\index{Kreck surgery theory|see{Surgery theory}}
Matthias Kreck modified traditional surgery theory in the early 1980s
(see \cite{Kreck99}) for the odd- and even-dimensional case such that
it requires much weaker topological input.  In this treatise we shall
only focus on the even-dimensional case.  There are two major
differences to the classical even-dimensional surgery programme. First
of all Kreck can replace the normal maps by a considerably weaker
notion called {\bf normal smoothings}, that is, a lift of the stable
normal bundle $M \rightarrow BO$ to a fibration $B \rightarrow BO$. As
pointed out in his paper normal maps are a special case of this
concept. The second main difference is that Kreck just needs
$[n/2]$-equivalences on the boundary of a normal cobordism whereas in
Wall's theory we started with full homotopy equivalences.  Hence we look
at a $(2q+2)$-dimensional cobordism

\newpage
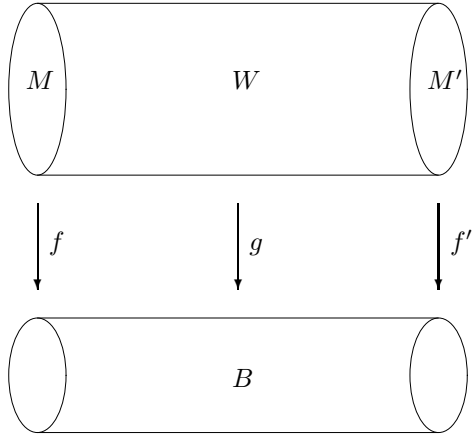
\begin{center}
\setlength{\unitlength}{0.006in}
\begin{picture}(750,400)(0,0)
\put(200,300){\ellipse{50}{150}}
\put(550,300){\ellipse{50}{150}}
\path(200,375)(550,375)
\path(200,225)(550,225)
\put(190,300){$M$}
\put(370,300){$W$}
\put(540,300){${M'}$}

\put(200,200){\vector(0,-1){75}}
\put(375,200){\vector(0,-1){75}}
\put(550,200){\vector(0,-1){75}}
\put(210,160){$f$}
\put(385,160){$g$}
\put(560,160){$f'$}

\put(200,50){\ellipse{50}{100}}
\put(550,50){\ellipse{50}{100}}
\path(200,0)(550,0)
\path(200,100)(550,100)
\put(370,40){$B$}
\end{picture}
\end{center}
with $f$ and $f'$ being $q$-equivalences and (after the usual surgery
below the middle dimension) $g$ being a $(q+1)$-equivalence.  The
obstruction here is an object
\begin{eqnarray*}
&&\kreckfs\\
&&=(H_{q+1}(W,M_0)\longleftarrow \im ( \pi_{q+2}(B,W) \rightarrow \pi_{q+1}(W))\ra H_{q+1}(W,M_1), \psi)\\
&&\in l_{2q+2}(\bZ[\pi_1(B)])
\end{eqnarray*}
$\kreckfs$ is called a {\bf preformation}. It is basically a tuple of
homomorphisms of (free) \fg modules over $\bZ[\pi_1(B)]$ such that
$(G,\gamma^*\mu,\theta)$ is a $(-)^{q+1}$-quadratic form.  The
obstruction lives in a monoid $\el$ and is {\bf elementary} if and
only if $(W,M,M')$ is cobordant rel$\partial$ to an $s$-cobordism.

One observes that on the one hand side the classical even-dimensional
case at the beginning is a special case of Kreck's surgery setting. On
the other hand the obstruction looks quite similar to the obstruction
formation in the odd-dimensional case.

The theory was successfully applied by M. Kreck and others (see also
Introduction of \cite{Kreck99}) to the classification of 4-manifolds
(see \eg\cite{Kreck01}, \cite{KreckHT94}), 7-dimensional homogeneous
spaces (see \eg\cite{KreckS88}, \cite{KreckS91}), of complete
intersections (see \cite{Kreck99}) or of classification of higher
dimensional manifolds in general (see \eg \cite{KreckT91})

Despite its successes, there are disadvantages of this programme: the
complicated algebra. The obstructions do not lie in the $L$-groups
anymore but in monoids $l_{2q}(\Lambda)$ with $\Lambda$ a ring with involution. 
The criterion that surgery is successful in creating an
$s$-cobordism is not that the obstruction vanishes but that it fulfils
certain complicated conditions (being {\bf elementary}). M. Kreck
himself writes: ``The obstruction $l$-monoids are very complicated and
algebraically - in contrast to the $L$-groups - not understood.''
(\cite{Kreck00} p.135) No-one has been able to compute $\el$ \eg for $\Lambda=\bZ$. 
Hardly any relations to the
$L$-groups are known.  The aim of this treatise is to bring some more
insight into the structure of $\el$, find relations to the quadratic
and asymmetric Witt-groups and give obstructions for elements in $\el$
to be elementary.
\pagebreak

\section{The Main Strategy}

For a moment let us forget about normal maps and smoothings and
consider a cobordism of manifolds $(W,M,M')$ with $\dim W \geq 6$ and even.  Let
$\Lambda=\bZ[\pi_1(W)]$.  The surgery theory question is: can we surger
$W$ (without affecting the boundary) in such a way that the result will be an
$s$-cobordism? Let us assume that we already know that there is a
diffeomorphism $h\colon M \stackrel{\cong}{\ra} M'$. Using that
diffeomorphism, the boundary of $W$ can be changed to a {\bf twisted
double} $\partial W = M \cup_{h|\partial M} M$ \ie two copies of $M$
glued together along its boundary by the automorphism on $\partial
M$ induced by $h$.

Twisted doubles play an important \role in topology \eg if one wants
to compute the cobordism group of diffeomorphisms or investigate
open book decompositions. Kreck computed the cobordism of automorphisms
first (see \eg \cite{Kreck84a}), followed by Quinn who
offered an alternative approach. In \cite{Quin79} he develops a theory
about open books decompositions which are strongly related to twisted doubles.
Twisted doubles were studied by Winkelnkemper in \cite{Win73} and 
both results were connected in \cite{Ran98} Chapter 30.

The main question in this context is to decide whether
manifold with a twisted double on the boundary is cobordant
rel$\partial$ to another manifold that itself carries a compatible twisted
double structure (see Section \ref{asymgeosec}).  In the
even-dimensional case, he constructed a non-singular {\bf asymmetric}
form which vanishes in the {\bf asymmetric Witt-group
$LAsy^0(\Lambda)$} if and only if the twisted double on the boundary 
can be extended (up to cobordism) to the whole manifold. If $(W,M,M')$ is
an $s$-cobordism (\ie a tube) any diffeomorphism $M\isora M'$ will
transform $W$ into a twisted double.  Hence, if $(W,M,M')$ is
cobordant rel$\partial$ to an $s$-cobordism, the asymmetric signature
will vanish for any diffeomorphism $M\isora M'$.

This approach does not lead to a practical method helping in the classification
of manifolds. After all, it starts with the assumption that the manifolds,
we want to classify, are already diffeomorphic!

The asymmetric signature becomes workable though if we use algebraic
surgery theory. This theory provides us with constructions that
imitate Quinn's asymmetric signature  for {\bf symmetric Poincar\'e pairs} (see
\cite{Ran98} Chapter 30 or Section \ref{asymalgsec} of this paper).
Symmetric pairs are purely algebraic objects but arise naturally from topology.
As an example, the symmetric Poincar\'e pair associated to a manifold $W$ with the twisted double
$M\cup_{h} M$ as a boundary consists of singular chain complexes of 
the universal covers of 
$W$, $M$ and $M$ together with chain equivalences that induce
Poincar\'e duality on those manifolds and further maps which guarantee
the symmetry properties of Poincar\'e duality.  In addition, one
needs the chain maps induced by the inclusion of the boundary into $W$
and the diffeomorphism $h\colon \partial M \isora \partial M$.

The algebraic surgery version of the asymmetric signature gives
another way of finding out whether $W$ is cobordant to a twisted
double (see \eg \cite{Ran98} Corollary 30.12). It also shows that the
answer to this question only depends on the homotopy type of the
ingredients.

But it does even more than that. It is a purely algebraic calculus and
can be used to help us testing the elementariness 
of elements in $\el$.

It will be shown that an element in $\el$ (or to be more precise a
{\bf preformation $z=\kreckfs$} representing that element) gives rise to
a {\bf quadratic Poincar\'e pair $x$}. 
A quadratic Poincar\'e pair is an algebraic model of a normal map.
Every quadratic Poincar\'e pair induces a symmetric Poincar\'e pair.
For such quadratic Poincar\'e pairs
there are notions of surgery, cobordism rel$\partial$ and $h$-cobordisms.
In particular, $x$ is cobordant rel$\partial$ to an algebraic
$h$-cobordism if and only if the element $\el$ at the beginning
was elementary. Hence $x$ behaves just like an ``algebraic Kreck surgery problem''.
This gives us some hope that a purely algebraic version of an asymmetric
signature, as presented for manifolds before, is feasible and helps
to test elementariness.

As a first step we need to transform the boundary of the 
quadratic Poincar\'e pair $x$ into some kind of algebraic equivalent of a
twisted double. Instead of a diffeomorphism $M\isora M'$ we only need
an equivalence of the quadratic complexes 
which are the two algebraic ``boundary components'' of $x$.
For the preformation $z$ this means that there exists
a (stable weak) isomorphism $t$ between $z$ and its
{\bf flip} $z'=\kreckfb{F^*}{\epsilon\mu}{G}{\gamma}{F}{-\theta}$:
a {\bf flip-isomorphism}.

There is also a more geometrical reason why flip-isomorphisms are
the correct algebraic substitute for a diffeomorphism. Let
$$(e,f,f')\colon (W,M,M') \ra X\times (I,0,1)$$ be a presentation \ie 
a degree $1$ normal cobordism such that $e$, $f$ and $f'$ are 
highly-connected and assume for simplicity that $M$ and $M'$ are closed. 
Let $z=\kreckfs$ be the Kreck obstruction formation. We encountered such a cobordism before when we 
defined the odd-dimensional Wall obstruction for $f$. By Corollary \ref{niceobstrcor}
both obstructions can be assumed to be the same. The Wall obstruction of the odd-dimensional 
highly-connected normal map $f'$ is the flip $z'=\kreckfb{F^*}{\epsilon\mu}{G}{\gamma}{F}{-\theta}$ of $z$. 
An ``algebraic'' isomorphism from the $L$-theory point of view is hence a (weak) isomorphism $t$
between $z$ and $z'$ - a flip-isomorphism:
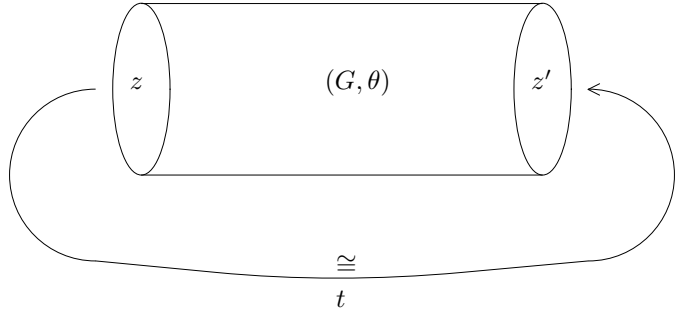
\begin{center}
\setlength{\unitlength}{0.006in}
\begin{picture}(750,300)(0,0)
\put(200,200){\ellipse{50}{150}}
\put(550,200){\ellipse{50}{150}}
\path(200,275)(550,275)
\path(200,125)(550,125)
\put(190,200){$z$}
\put(360,200){$(G,\theta)$}
\put(540,200){$z'$}

\put(160,125){\arc{150}{1.57}{4.71}}
\put(590,125){\arc{150}{4.71}{1.57}}
\path(590,200)(600,205)
\path(590,200)(600,195)
\spline(160,50)(370,30)(590,50)
\put(370,10){$t$}
\put(370,40){$\cong$}

\end{picture}
\end{center}

In any case we can prove that if $z$ is (stably) elementary 
it has indeed (stable) flip-isomorphisms.

It is easy to determine when
a preformation $z=\kreckfs$ over $\bZ$ with finite $\coker\gamma$ and
$\coker\mu$ allows flip-isomorphisms. $z$ induces {\bf linking forms} \ie
$\pm$-symmetric forms on those cokernels. If $z$ is an obstruction to
an even-dimensional Kreck surgery problem $(W,M,M')\ra B$ then those
algebraic linking forms of $z$ are induced by the topological linking
forms of $M$ and $M'$. Any isometry of those linking forms gives rise
to flip-isomorphisms and conversely any flip-isomorphism induces
an isometry of the linking forms.

But let us return to the quadratic Poincar\'e complex $x$ we constructed
out of $z$.
Every flip-isomorphism $t$ transforms the boundary of the
Poincar\'e pair $x$ into an algebraic twisted double. Now the
algebraic surgery version of the asymmetric signature yields an
asymmetric signature $\sigma^*(z,t)\in LAsy^0(\Lambda)$ depending only
on the preformation $z$ and a flip-isomorphism $t$.  (These
constructions do not only work for group rings or $\bZ$ but for any
weakly finite ring $\Lambda$ with $1$ and involution).  If $[z]\in
\el$ is elementary then the asymmetric signatures vanish for all
(stable) flip-isomorphism.

As mentioned before, isometries of linking forms are a good source for
flip-iso\-morph\-isms.  It turns out that the asymmetric signatures only
depend on those linking forms.

Obviously we are interested in whether these asymmetric signatures
will be complete obstructions to elementariness. The evidence does not
look all too optimistic. First of all asymmetric signatures completely
ignore all quadratic information simply because - unlike in other
surgery problems - their definition and application does only require
manifolds but not necessarily normal maps/smoothings. 

We again resort to the manifold world for some inspiration for a stronger
obstruction. Assume again we have an even-dimensional cobordism $(W,M,
M')$.  For a moment let us assume that $M$ and $M'$ are closed. Any
diffeomorphism $h\colon M\isora M'$ allows us to glue both ends of $W$
together. Alternatively we obtain the same manifold if we replace
$M$ by $M'$ using $h$ (then $W$ is a manifold with an (un-)twisted double $M+M$
on the boundary) and then glue the $s$-cobordism $M\times(I,0,1)$ on it.
The resulting closed manifold $V$ is null-cobordant if and
only if $(W,M,M')$ is cobordant to an $s$-cobordism. 

If $M$ and $M'$ are not closed, $h$ again turns $(W,M,M')$ into a
manifold with a twisted double as a boundary. But we have to be
careful now: not every twisted double is the boundary of an
$s$-cobordism. If, however, we demand that $h|\partial M$ is isotopic
to $1_{\partial M}$ we can glue $(W,M,M)$ onto $M\times (I,0,1)$ and
again get a closed manifold which is null-cobordant if and only if
$(W,M,M')$ is cobordant to an $s$-cobordism.  Similar constructions
also work for a normal cobordism and a compatible diffeomorphism $h$.

Again we follow our philosophy that anything that can be done for manifolds
can also be done in the algebraic world of complexes.
We imitate the procedures for the quadratic Poincar\'e pair
$x$ that we created out of a preformation $z$.
The case $\partial M=\emptyset$ corresponds to the case where
$z$ is a {\bf non-singular formation}.  As an example, the
obstruction of a Kreck surgery problem $(W,M,M')\ra X\times (I,0,1)$
is a formation when all maps involved are normal maps, $X$ is a
finite geometric Poincar\'e space and the induced map $\partial M \ra \partial X$ is a homotopy
equivalence. Just like with manifolds and normal maps, it is possible to glue
quadratic Poincar\'e pairs together and there is a notion of cobordism.
All we need is an equivalent for the diffeomorphism $h\colon M \isora M'$.
Any chain equivalence of the quadratic complexes which constitute the
two ``boundary components'' of $x$ will do the job. It is  
nothing else than the flip-isomorphism, we have encountered before.
So we use a choice of flip-isomorphism to glue the ``ends'' of the Poincar\'e pair
together and the resulting Poincar\'e complex is null-cobordant if and
only if $x$ is cobordant rel$\partial$ to an algebraic $h$-cobordism
and this is the case if and only if $[z]\in\el$ is elementary. One of
the fundamental facts of algebraic surgery theory is that the set of
cobordism classes of Poincar\'e complexes are Wall's
$L$-groups. Hence our construction leads to an obstructions in
$L_{2q+2}(\Lambda)$ - the {\bf quadratic signatures}.

Unfortunately, the generalization to arbitrary preformations is much
more unpleasant. In the case of manifolds or normal maps
not any diffeomorphism could be used for the glueing operation. The
same is true in the algebraic surgery world.  Not every
flip-isomorphism is suitable to produce a quadratic signature.
We have to introduce a new class of flip-isomorphisms
called {\bf flip-isomorphisms rel$\partial$}. 
The quadratic signatures of such a flip-isomorphism rel$\partial$
will also depend on other choices and hence are rather difficult
to handle. 

In any case, one can show that quadratic and asymmetric signatures are connected
via the canonical map
\begin{eqnarray*}
L_{2q+2}(\Lambda) &\ra& LAsy^0(\Lambda)\\
(K,\psi) &\mt& (K,\psi_0-\epsilon\psi_0^*)
\end{eqnarray*}

\section{The Results}

Let $\Lambda$ be a weakly finite ring with $1$ and an involution.  Let
$\epsilon=(-)^q$ and let $z=\kreckfs$ be an $\epsilon$-quadratic split
preformation, that is a tuple consisting of a free
\fg $\Lambda$-module $F$, a \fg $\Lambda$-module $G$, a
$\Lambda$-homomorphism $\svec{\gamma \\ \mu}\colon G \rightarrow
F\oplus F^*$ and a map $\theta\colon G \rightarrow
Q_{-\epsilon}(\Lambda)$ such that $(G, \gamma^*\mu, \theta)$ is a
$(-\epsilon)$-quadratic form.  The obstruction to a surgery problem in
Kreck's theory is such an object.

We say that two preformations are {\bf stably strongly isomorphic } if
they are isometric after one adds ``hyperbolic elements'' of the form
$$\kreckf{(P\oplus P^*)}{1}{(P\oplus P^*)}{\svec{0&1\\-\epsilon&0}}{\svec{0&1\\0&0}}$$ 
The equivalence classes form the {\bf $l$-monoid $l'_{2q+2}(\Lambda)$} 
and the equivalence classes of regular preformations (\ie $G$ is free)
define the {\bf $l$-monoid $\el$} (compare
Definition \ref{mldef} on p.\pageref{mldef}).  The main theorem in Kreck's modified surgery
theory states that surgery leads to an $s$-cobordism if the
obstruction is {\bf stably elementary}.  Hence, our aim is to find obstructions for
$z$ to be stably elementary.  By Corollary
\ref{mkregcor} (p.\pageref{mkregcor}) we can replace $z$ by a {\bf regular} preformation
\ie one with a free $G$.

In Kreck's original theory all isometries and isomorphisms were simple and all modules involved
were based. We will ignore the Whitehead obstruction in the following and hence only deal
with $h$-cobordisms. For $\Lambda=\bZ$ or $\Lambda=\bZ[\bZ^m]$ there is no difference.

There are certain obvious primitive obstruction for a preformation to be elementary.
\begin{prop*}[See Corollary \ref{simpleobstr} on p.\pageref{simpleobstr}]
If $[z]\in\el$ is elementary then $\ker\gamma \cong \ker\mu$, $\coker\gamma \cong \coker\mu$, 
$\ker\gamma^*\mu\cong \ker\gamma\oplus\ker\gamma^*$, 
$\coker\gamma^*\mu\cong \coker\gamma\oplus\coker\gamma^*$ and $\rk F$ is even
\end{prop*}

The first really important obstruction we will discuss is the existence of a {\bf
flip-isomorphism}. A flip-isomorphism is a weak isomorphism
$(\alpha,\beta,\chi)$ of $z$ with its flip
$z'=\kreckfb{F^*}{\epsilon\mu}{G}{\gamma}{F}{-\theta}$.  (Compare
Definition \ref{flipisodef} on p.\pageref{flipisodef}).  Weak
isomorphisms are a generalization of isomorphisms of formations as
they are used in the definition of the odd-dimensional $L$-groups (see
Definition \ref{weakisodef} on p.\pageref{weakisodef}). A strong isomorphism (which is used in
the definition of the even-dimensional $l$-monoids) is also a weak
isomorphism (see Remark \ref{weakisorem} on 
p.\pageref{weakisorem}).
Every elementary preformation has a flip-isomorphism:
\begin{prop*}[See Corollary \ref{elemflipcor} on p.\pageref{elemflipcor}]
Let $z$ be a regular $\epsilon$-quadratic split preformation.  If $z$
is (stably) elementary there is a (stable) strong flip-isomorphism
$(\alpha,\beta,0)$ such that $\alpha\colon F\ra F^*$ is
$\epsilon$-symmetric and zero in $L^{2q}(\Lambda)$ (and hence also in
$LAsy^0(\Lambda)$) and $\beta^2=1_G$.
\end{prop*}

Now we use the flip-isomorphism to define asymmetric and quadratic signatures.

\subsection*{Flip-Isomorphisms and Asymmetric Signatures}

For all flip-isomorph\-isms $t=(\alpha, \beta, \chi)$ (even those that do
not fulfil the stronger conditions of the previous proposition) there
is an {\bf asymmetric signature $\sigma^*(z,t)=(M,\rho) \in
LAsy^0(\Lambda)$ of a flip-isomorphism $t$ of $z$} given by
$$
    \rho=\mat{0&0&\alpha\\1&0&-\epsilon\\0&1&\epsilon\alpha(\chi-\epsilon\chi^*)^*\alpha^*}\colon M=F\oplus F^*\oplus F \ra M^*
$$
(Definition \ref{asymsignflipdef} on p.\pageref{asymsignflipdef}).
We define an abelian monoid $fl_{2q+2}(\Lambda)$ as a kind of $l$-monoid of preformations
with a choice of flip-isomorphism (Definition \ref{flmonoidef} on p.\pageref{flmonoidef}). Then
the asymmetric signatures define a map from that monoid into the asymmetric Witt-group
which vanishes for all stably elementary preformations.

\begin{thm*}[See Theorem \ref{flasymthm} and Remark \ref{pirem} on p.\pageref{flasymthm}]
The asymmetric signatures give rise to a 
well-defined homomorphism of abelian monoids
\begin{eqnarray*}
\sigma^*\colon fl_{2q+2}(\Lambda) &\ra& LAsy^0(\Lambda)
\\
{[(z,t)]} &\mt& \sigma^*(z,t)
\end{eqnarray*}

If $[z']\in l_{2q+2}(\Lambda)$ is elementary then
$\sigma^*(z,t)=0$ for all flip-isomorphisms $t$ of
all preformations $z$ with $[z]=[z']\in\el$.
\end{thm*}

We will give two proofs of this Theorem (or rather of the underlying 
Theorem \ref{elemasymsignthm} on p.\pageref{elemasymsignthm}):
one using algebraic surgery theory and a explicit stable lagrangian of
our asymmetric signature.

A test case for asymmetric signatures are boundaries of quadratic forms in particular
the submonoid of $\el$ given by the injection
\begin{eqnarray*}
L_{2q+2}(\Lambda) &\hookrightarrow& l_{2q+2}(\Lambda)\\
(K,\theta) &\mt& \partial(K,\theta)=\kreckf{K}{1_K}{K}{\theta-\epsilon\theta^*}{\theta}
\end{eqnarray*}
Here live the obstructions of traditional surgery theory interpreted as a special case
of Kreck's modified theory.

\begin{cor*}[Corollary \ref{asymbdrycor}, p.\pageref{asymbdrycor}]
Let $(K,\theta)$ be a non-singular $(-\epsilon)$-quadratic form. 
Then $z=\partial(K,\theta)$ is a non-singular $\epsilon$-quadratic split formation.
\begin{enumerate}
\item $z$ has a (stable) flip-isomorphism.
\item $[z]\in\el$ is elementary if and only if $[(K,\theta)]=0\in L_{2q+2}(\Lambda)$.
\item All asymmetric signatures equal $[(K,\theta-\epsilon\theta^*)]\in LAsy^0(\Lambda)$.
\item Assume that either $\Lambda$ is a field of characteristic different to $2$ or that $\Lambda=\bZ$ and $\epsilon=-1$.
      The preformation $z$ is elementary if and only if the asymmetric signatures vanish.
\end{enumerate}
\end{cor*}
Asymmetric signatures ignore any ``quadratic split'' information. Hence
non-sing\-ular skew-quadratic forms over $\bZ$ and $\bZt$ with non-trivial Arf-invariant
have vanishing asymmetric signatures but their boundaries aren't stably elementary
(see Example \ref{asymcountex} on p.\pageref{asymcountex}).

\subsection*{Linking Forms}

What is an easy source for flip-isomorphisms? For a certain class of
preformations the answer is linking forms.  Assume that
$S\subset\Lambda$ is a central multiplicative subset, \eg 
$\Lambda=\bZ$ and $S=\bZ\setminus\{0\}$. We call a map an
$S$-isomorphism if tensoring with $S^{-1}\Lambda$ makes it an
isomorphism. If $\mu$ is an $S$-isomorphism the preformation
determines a linking form $L_\mu$ and if $\gamma$ is an
$S$-isomorphism it defines a linking form $L_\gamma$.
\begin{prop*}[See Proposition \ref{sisoprop} on p.\pageref{sisoprop}]
Let $z=\kreckfs$ be a regular split $\epsilon$-quadratic  preformation
with either $\mu$ or $\gamma$ an $S$-isomorphism.
\begin{enumerate}
\item If $z$ allows a flip-isomorphism then both $\gamma$ and $\mu$ are 
$S$-isomorphisms. Every flip-isomorphism $t=(\alpha,\beta,\chi)$ induces
an isomorphism of split $(-\epsilon)$-quadratic linking forms 
$[\alpha^{-*}]\colon L_\mu \stackrel{\cong}{\ra} L_\gamma$.

\item Assume $\gamma$ and $\mu$ are both $S$-isomorphisms and $L_\gamma$ and $L_\mu$ are isomorphic.
Every isomorphism $l\colon L_\mu \stackrel{\cong}{\ra} L_\gamma$ induces a 
stable flip-isomorphism $(\alpha,\beta,\chi)$ of $z$
such that $[\alpha^{-*}]=l\colon L_\mu \stackrel{\cong}{\ra} L_\gamma$.
\end{enumerate}
\end{prop*}

We introduce sub-monoids of the $l$-monoids $l^{2q+2}(\Lambda)$ and $fl^{2q+2}(\Lambda)$ of preformations $\skreckfs$ 
where $\gamma$ and $\mu$ are $S$-isomorphisms and the quadratic refinement is omitted:
\begin{eqnarray*}
l_S^{2q+2}(\Lambda)  &=& \left\{[\skreckfs]\in \el | \text{$\mu$ and $\gamma$ are $S$-isomorphisms}\right\}\\
fl_S^{2q+2}(\Lambda) &=& \left\{[(z,t)] : [z]\in l_S^{2q+2}(\Lambda)\right\}
\end{eqnarray*}
Similarly we define an $l$-monoid $ll_S^{2q+2}(\Lambda)$ of preformations $\skreckfs$ with a choice
of isometry $L_\mu\cong L_\gamma$. The previous proposition can be interpreted as
the existence of a surjection $fl_S^{2q+2}(\Lambda) \ra ll_S^{2q+2}(\Lambda)$ 
(see Section \ref{linkmodsec}).

For preformations in $l_S^{2q}(\Lambda)$ the asymmetric signature $\sigma^*(z,t)$ of a flip-isomorphism $t$ 
does only depend on the preformation and the isometry of linking forms induced by $t$. 
\begin{thm*}[See Theorem \ref{elemlinkasymthm} on p.\pageref{elemlinkasymthm}]
There is a lift of the asymmetric signature map of Theorem \ref{flasymthm}
\begin{eqnarray*}
\xymatrix
{
fl_S^{2q+2}(\Lambda)    \ar@{>>}[d]_{L} 
            \ar[r]_{\sigma^*}
&
LAsy^0(\Lambda)
\\
ll_S^{2q+2}(\Lambda)    \ar[ru]_{\sigma^*}
}
\end{eqnarray*}

If $[z']\in l^{2q+2}_S(\Lambda)$ is elementary then
$\sigma^*(z,l)=0\in LAsy^0(\Lambda)$ for all
isomorphisms $l\colon L^\mu \stackrel{\cong}{\ra} L^\gamma$.
\end{thm*}
This theorem is quite an improvement. Instead of checking the
asymmetric signature for all flip-isomorphisms of all stably strongly isomorphic
preformations, we only have  to go through all
isometries of linking forms of one representative. In the case $\Lambda=\bZ$ and
$S=\bZ\setminus\{0\}$ they are only finitely many of them.
We can do even more. In certain circumstances we can read
off the linking forms from a simply-connected manifold directly.

\begin{prop*}[See Propositions \ref{mfdlinkprp} on p.\pageref{mfdlinkprp}]
Let $p\colon B \rightarrow BO$ be a fibration with $\pi_1(B)=0$.  Let $M_i$ be
$(2q+1)$-dimensional manifolds with a {\bf $(q-1)$-smoothings in $B$}  \ie a
lift of the stable normal bundle over $p$ which is $q$-connected.
Let $f\colon \partial M_0 \rightarrow \partial M_1$ be a diffeomorphism compatible
with the smoothings. Let $W$ be a cobordism of $M_0 \cup_f M_1$ with a
compatible $q$-smoothing over $B$. 
As in Corollary \ref{niceobstrcor} we define an obstruction
\begin{eqnarray*}
x(W)&=&\kreckfs\\
&=&(H_{q+1}(W,M_0)\longleftarrow H_{q+2}(B,W) \ra H_{q+1}(W,M_1), \theta)\\
  &\in& l'_{2q+2}(\bZ)
\end{eqnarray*}
Let $l_{M_i}^B$ be the linking form on $H_{q+1}(B,M_i)$ which is induced
by the topological linking form of $M_i$. 

If $\coker\gamma=H_{q+1}(B,M_0)$ is finite then $L^\gamma=-l_{M_0}^B$.

If $\coker\mu   =H_{q+1}(B,M_1)$ is finite then $L^\mu = -l_{M_1}^B$.

Assume both cokernels are finite.
If $W$ is cobordant rel$\partial$ to an $h$-cobordism then  there exist
isomorphisms $l\colon L^\mu=-l^B_{M_1}\isora L^\gamma=\epsilon -l^B_{M_0}$ 
and their asymmetric signatures $\sigma^*(x(W),l)\in LAsy^0(\bZ)$
will all vanish.
\end{prop*}

\subsection*{Non-singular Formations}

In general asymmetric signatures are not strong enough to show
elementariness, therefore we look out for a stronger obstruction - {\bf quadratic signatures}.

Quadratic signatures turn out to be rather complicated objects so we will deal
first with a special class of preformations for which they behave nicely.
Let $z=\kreckfs$ be a non-singular $\epsilon$-quadratic split formation \ie that means the map
$\svec{\gamma\\\mu}\colon G \ra H_{\epsilon}(F)$ is an inclusion of a lagrangian.
Let $(e,f,f')\colon (W, M,M') \ra X\times (I,0,1)$
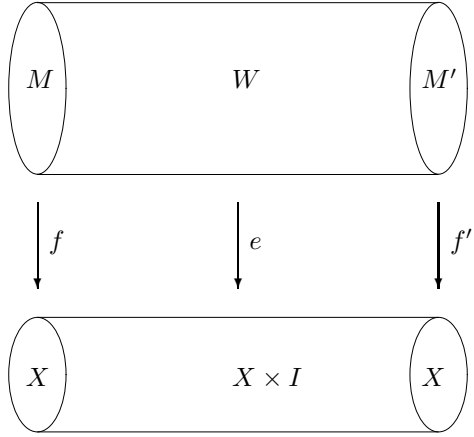
\begin{center}
\setlength{\unitlength}{0.006in}
\begin{picture}(750,400)(0,0)
\put(200,300){\ellipse{50}{150}}
\put(550,300){\ellipse{50}{150}}
\path(200,375)(550,375)
\path(200,225)(550,225)
\put(190,300){$M$}
\put(370,300){$W$}
\put(535,300){${M'}$}

\put(200,200){\vector(0,-1){75}}
\put(375,200){\vector(0,-1){75}}
\put(550,200){\vector(0,-1){75}}
\put(210,160){$f$}
\put(385,160){$e$}
\put(560,160){$f'$}

\put(200,50){\ellipse{50}{100}}
\put(550,50){\ellipse{50}{100}}
\path(200,0)(550,0)
\path(200,100)(550,100)
\put(190,40){$X$}
\put(370,40){$X\times I$}
\put(535,40){$X$}
\end{picture}
\end{center}
be a degree 1 normal cobordism (\ie  $e$ and $f$ and $f'$ are normal
maps and $(X,\partial X)$ is a finite geometric Poincar\'e pair) and
$f|\colon\partial M\ra \partial X$ is homotopy equivalence. 
Then the modified Kreck surgery obstruction of 
Corollary \ref{niceobstrcor} (p.\pageref{niceobstrcor}) 
is such a non-singular formation.

By \cite{Ran80a} Proposition 2.2, the map $\svec{\gamma\\\mu} \colon
G \ra H_\epsilon(F)$ can be extended to an isomorphism of hyperbolic
$\epsilon$-quadratic forms
$$
\left(f=\mat{\gamma & \tilde\gamma\\\mu&\tilde\mu}, 
\mat{\theta&0\\\tilde\gamma^*\mu&\tilde\theta}\right) \colon
H_\epsilon(G) \stackrel{\cong}{\ra} H_\epsilon(F)
$$

For any $\tau\colon G^*\ra G$ the maps $\tilde\gamma'=\tilde\gamma+\gamma (\tau-\epsilon\tau^*)$,
$\tilde\mu'=\tilde\mu+\mu (\tau-\epsilon\tau^*)$, $\tilde\theta'=\tilde\theta+(\tau-\epsilon\tau^*)^*\theta (\tau-\epsilon\tau^*) 
+\tilde\gamma^*\mu (\tau-\epsilon\tau^*)^* -\epsilon\tau$ define another extension to an isomorphism of hyperbolic forms.
Conversely any other extension has this form.

For any flip-isomorphism $t$ of a non-singular $\epsilon$-quadratic
split formation $z$ and a choice of extensions $\tilde\gamma$,
$\tilde\mu$, $\tilde\theta$ we can define a {\bf quadratic signature
$\tilde\rho^*(z,t,\tilde\gamma,\tilde\mu,\tilde\theta)=(M,\xi')\in
L_{2q+2}(\Lambda)$} given by
\begin{eqnarray*}
\xi'&=&\mat{    \tilde\gamma^*\tilde\mu+\tilde\gamma^*\alpha\nu\alpha^*\tilde\gamma     & -\tilde\gamma^*\alpha\gamma & 0 \\ 
        0                                   & \epsilon\theta^* & 0 \\ 
        \epsilon(\alpha^*\tilde\gamma-\tilde\mu)                & -\mu & 0}
\colon M=G^*\oplus G \oplus F^* \ra M^*
\end{eqnarray*}
(see Definition \ref{quasgndef} on p.\pageref{quasgndef}).

\begin{thm*}[See Theorem \ref{formquadthm} on p.\pageref{formquadthm}]
Let $z'$ be a non-singular $\epsilon$-quadratic split formation.
$[z']\in l_{2q+2}(\Lambda)$ is elementary if and only if there is a stably strongly isomorphic $z=\kreckfs$,
a flip-isomorphism $t$ and $\tilde\gamma$, $\tilde\mu$ and $\tilde\theta$ such that
$$\left(f=\mat{\gamma & \tilde\gamma\\\mu&\tilde\mu},\mat{\theta&0\\\tilde\gamma^*\mu&\tilde\theta}\right)\colon H_\epsilon(G) 
\stackrel{\cong}{\ra} H_\epsilon(F)$$
is an isomorphism of hyperbolic $\epsilon$-quadratic forms and $\tilde\rho^*(z,t,\tilde\gamma, \tilde\mu, \tilde\theta)\in L_{2q+2}(\Lambda)$ vanishes.
\end{thm*}

There is also a surprise about asymmetric signatures of non-singular formations:
they are independent of the choice of flip-isomorphisms.
\begin{thm*}[See Theorem \ref{asymformthm} on p.\pageref{asymformthm}]
Let $z$ be a non-singular formation. Let $t$ and $t'$ be two flip-isomorphisms. 
Then $\sigma^*(z,t)=\sigma^*(z,t')\in LAsy^0(\Lambda)$.
\end{thm*}

\subsection*{Quadratic Signatures for Arbitrary Preformations}

The general definition of a quadratic signature demands much more
preparation and we will only sketch it here.  Let $z=\kreckfs$ be a
regular $\epsilon$-quadratic split preformation and let
$t=(\alpha,\beta,\bar\nu)$ be a flip-isomorphism of $z$.  $z$ and $t$
and a choice of representatives for $\bar\nu$ and $\bar\theta$, \etc define
a self-equivalence $(h_t, \chi_t)\colon (C,\psi) \isora (C,\psi)$ of
$2q$-dimensional quadratic Poincar\'e complexes (see (\ref{htdefeqn})
on p.\pageref{htdefeqn}).  Assume there exists a homotopy
$(\Delta,\eta) \colon (1,0) \simeq (h_t,\chi_t) \colon (C,\psi) \isora (C,\psi)$ 
(see Definition \ref{hommodef} on p.\pageref{hommodef}) then
$t$ is called a {\bf flip-isomorphism rel$\partial$} (see Definition
\ref{flipreldef} on p.\pageref{flipreldef}). Those ingredients define
a quadratic signature in $L_{2q+2}(\Lambda)$ which vanishes for a choice of
$\Delta$, $\eta$, $\theta$, \etc if and only if $z$ is stably
elementary. The construction is described in Section
\ref{quadconstrsec} on p.\pageref{quadconstrsec} and the relation to
elementariness in Theorem \ref{quadelemthm} on
p.\pageref{quadelemthm}.

The quadratic signature of $z$, $t$, $\Delta$, \etc is mapped to the
asymmetric signature $\sigma^*(z,t)$ via the map
\begin{eqnarray*}
L_{2q+2}(\Lambda) &\ra& LAsy^0(\Lambda)\\
(K,\psi) &\mt& (K,\psi_0-\epsilon\psi_0^*)
\end{eqnarray*}
(see Theorem \ref{asyquadthm} on p.\pageref{asyquadthm}).
Its kernel can be computed in terms of cobordism classes of
automorphisms of quadratic Poincar\'e complexes (see
Remark \ref{asyquadrem} on p.\pageref{asyquadrem}).

In the case of $\Lambda=\bZ$ and $\epsilon=-1$ \ie $q=2m-1$
the map is an injection:

\begin{prop*}[See Corollary \ref{coolintcor} on p.\pageref{coolintcor}]
Let $q=2m-1$ \ie  $\epsilon=-1$. Let $z=\kreckfs$ be a regular skew-quadratic split preformation over $\bZ$ 
\begin{enumerate}
\item $[z]\in l_{4m}(\bZ)$ is elementary if and only if there is a flip-isomorphism rel$\partial$ $t$ such that $\sigma^*(z,t)=0\in LAsy^0(\bZ)$.
\item The quadratic signature $\rho^*(z,t,\nu,\theta,\kappa,\Delta,\eta)\in L_{4m}(\bZ)$ only depends on $z$ and $t$.
\end{enumerate}
\end{prop*}

\section{The Contents}

{\bf Chapter \ref{prefchap}} gives an introduction into the topology
and algebra of traditional and modified surgery theories.  We will
define forms, preformations and elementariness.

The next two chapters build up the foundation for the application of
algebraic surgery theory to the study of preformations. In {\bf Chapter
\ref{translchap}} we translate preformations into quadratic complexes
and pairs. We define algebraic versions of surgery and cobordism
rel$\partial$ and $h$-cobordisms for quadratic pairs.  In {\bf Chapter
\ref{flipchap}} we define the important concept of flip-isomorphisms
and discuss how they fit into the algebraic chain models we
constructed in the preceding chapter.

The following chapters discuss asymmetric and quadratic signatures in
all generality.  {\bf Chapter \ref{asymchap}} presents the theory of
asymmetric forms and complexes and how one can define asymmetric
signatures for Poincar\'e pairs that have an algebraic twisted double
as a boundary.  These general constructions are applied to the
Poincar\'e pairs defined in Chapter \ref{flipchap} and produce the
asymmetric signatures of flip-isomorphisms.  {\bf Chapter
\ref{quadchap}} deals with the definition of quadratic signatures for
general preformations.

We continue with the treatment of special classes of preformations.
{\bf Chapter \ref{nonsingchap}} covers the quadratic signatures for
the easier case of non-singular formations. It also contains a proof
for the fact their asymmetric signatures do not depend on the choice
of flip-isomorphism.  For preformations with linking forms the theory
of asymmetric signatures becomes particularly elegant as will be seen
in {\bf Chapter \ref{linkformchap}}.

The {\bf Appendix \ref{algsurchap}} contains a compilation of
constructions and formulae from algebraic surgery theory.

\chapter{Preformations}
\label{prefchap}

{\bf Section \ref{formgeosec}} presents the algebraic and
geometric concepts behind Kreck's surgery theory and its relation to
traditional surgery theories.  In {\bf Section \ref{formalgsec}} we
introduce the language of forms, formations and preformations - the
building blocks of all our various surgery obstruction groups and
monoids. Preformations are the objects that appear as obstructions in
Kreck's surgery theory; its main theorem states that surgery
succeeds in producing an $s$-cobordism if and only if that obstruction
preformation has a certain property: stable elementariness.  {\bf
Section \ref{elemgeosec}} will present a heuristic way from topology
to a definition of an elementary preformation. Then ({\bf Section
\ref{elemalgsec}}) various equivalent definitions and some simple properties of
that important concept are given, \eg simple invariants
which are obstructions to elementariness using cokernels and kernels
of $\gamma$, $\mu$ and $\gamma^*\mu$ of a preformation $\kreckfs$ (see
Corollary \ref{simpleobstr}).

An obstacle to transferring preformations
into algebraic surgery world in the next chapters is the fact that
the module $G$ in a preformation $\kreckfs$ does not need to be free.
By Corollary \ref{mkregcor} any preformation can be replaced
by a {\bf regular} preformation (\ie a preformation with $free$ G).

\section{Forms, Preformations and Formations: the Geometry}
\label{formgeosec}

In this section we want to compare even- and odd-dimensional traditional surgery theory
and Kreck's even-dimensional theory.
\index{Surgery theory!Wall!even-dimensional}
We start with the {\bf traditional even-dimensional theory} as developed by C.T.C. Wall and others.
Let $$(e,f,f')\colon (W,M,M') \rightarrow X\times (I,0,1)$$ be a $(2q+2)$-dimensional degree $1$ normal cobordism
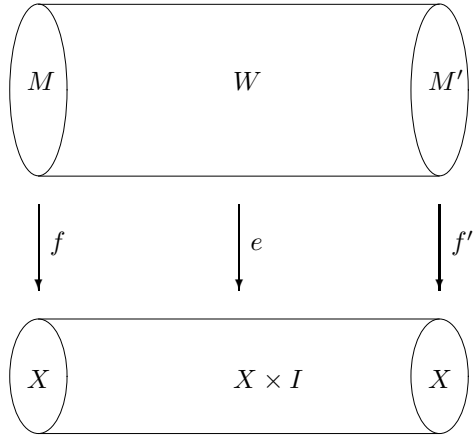
\begin{center}
\setlength{\unitlength}{0.006in}
\begin{picture}(750,400)(0,0)
\put(200,300){\ellipse{50}{150}}
\put(550,300){\ellipse{50}{150}}
\path(200,375)(550,375)
\path(200,225)(550,225)
\put(190,300){$M$}
\put(370,300){$W$}
\put(540,300){${M'}$}

\put(200,200){\vector(0,-1){75}}
\put(375,200){\vector(0,-1){75}}
\put(550,200){\vector(0,-1){75}}
\put(210,160){$f$}
\put(385,160){$e$}
\put(560,160){$f'$}

\put(200,50){\ellipse{50}{100}}
\put(550,50){\ellipse{50}{100}}
\path(200,0)(550,0)
\path(200,100)(550,100)
\put(190,40){$X$}
\put(370,40){$X\times I$}
\put(540,40){$X$}
\end{picture}
\end{center}
with a finite $(2q+2)$-dimensional geometric Poincar\'e space\footnote{
\index{Poincar\'e space} \index{Poincar\'e pair!geometric} \index{Finite Poincar\'e space}
A Poincar\'e space (or pair) is a topological space (or a pair of spaces) for which
there exists a Poincar\'e (or Poincar\'e-Lefschetz) duality.
A Poincar\'e space or pair is finite if it is a finite $CW$-complex.
All closed manifolds are finite Poincar\'e spaces.} $X$ such that $f$ and $f'$ are 
homotopy equivalences.

Our aim is to perform surgery on $W$ rel$\partial$ such that the result is an $s$-cobordism.
Then, by the $s$-cobordism theorem, $M$ and $M'$ are diffeomorphic. Surgery theory works
only in higher dimensions, hence we assume $q\ge 2$. 
 
After having made $e$ highly-connected by surgery below the middle
dimension, we can define a $(-1)^{q+1}$-dimensional quadratic form
$(K_{q+1}(W), \lambda, \mu)$ with $K_{q+1}(W)=H_{q+2}(\widetilde e)$
the homology of the induced map of the universal covers of $W$ and $X$ with local
coefficients and $\lambda$ and $\nu$ the intersection and
self-intersection numbers on $W$. This form is zero in the Witt group
$L_{2q+2}(\bZ[\pi_1(X)])$ of non-singular $(-)^{q+1}$-quadratic forms
over $\bZ[\pi_1(X)]$ if and only if $e\colon W \ra X\times I$ is
cobordant rel$\partial$ to a homotopy equivalence, \ie if and
only if we can do surgery on the inside of $W$ to obtain an
$s$-cobordism. A quadratic form vanishes in the $L$-group if (after
addition of hyperbolic forms) it has a lagrangian (\ie a free direct
summand of half dimension on which the quadratic form vanishes). If
there is a lagrangian for $(K_{q+1}(W), \lambda, \mu)$ one simply
kills its generators by surgery and the result will be an
$s$-cobordism.

\index{Surgery theory!Wall!odd-dimensional}
Before introducing Kreck's even-dimensional approach we have a look at the {\bf traditional odd-dimensional case}. 
Let $(X,\partial X)$ be a finite $(2q+1)$-dimensional geometric Poincar\'e pair.
Let $N$ and $N'$ be two $2q$-dimensional manifolds such that $\partial M= N\cup N'$
and let $f\colon (M,\partial M) \ra (X,\partial X)$ be a degree $1$ normal map such that its restriction
to the boundary $\partial M \ra \partial X$ is a homotopy equivalence. 
Surgery below the middle-dimension makes $f\ra X$ highly-connected.
We are interested in the question whether $N$ and $N'$ are diffeomorphic, that is, 
whether $f\colon M \ra X$ is cobordant rel$\partial$ to a homotopy equivalence.
One can construct an obstruction by looking at a so-called {\bf presentation of $f$}. 
\index{Presentation} A presentation is a $(2q+2)$-dimensional normal cobordism
$$
(e,f,f')\colon (W, M,M') \ra X\times (I,0,1)
$$
such that $e$ is $(q+1)$-connected and $f$ and $f'$ is $q$-connected. Such presentations exist for any such $f$ 
with the above properties: one chooses a  set of generators $\{x_1,\dots, x_k\}$ of the kernel module $K_{q}(M)$.
Because we are below the middle dimension, the generators can be realized by disjoint framed embeddings $g_i\colon S^q \times D^{q+1} \ra M$.
The trace $W$ of the surgeries performed on them will be a presentation.

There is a purely algebraic way to test whether there's any presentation that contains a homotopy equivalence $f'\colon M'\ra X$
on the other end  (see also \cite{Ran02} Chapter 12.2):
Let $U$ be the union of all the images of all $g_i$ and $M_0=\overline{M\setminus U}$.
Such a decomposition is called a \index{Heegaard-splitting}
{\bf Heegaard splitting} (see \cite{Ran02} Definition 12.6).
Then the self-intersection form on $K_{q}(\partial U)=K_q(\#_k S^q\times S^q)$ is the hyperbolic form on $\bZ[\pi_1(X)]^{2k}$.
Because $\partial M=\partial U$, the images of $K_{q+1}(U,\partial U)$ and $K_{q+1}(M_0,\partial U)$ in $K_q(\partial U)$
are lagrangians. A non-singular quadratic form with a pair of lagrangians is called a {\bf non-singular $(-)^q$-quadratic formation}.
It turns out that the formation $(K_{q}(\partial U); K_{q+1}(U,\partial U),  K_{q+1}(M_0,\partial U))$ provides enough
data to decide our surgery problem.

We can read off the same information from our presentation. We define
\begin{eqnarray*}
F&=&K_{q+1}(W,M)\cong K_{q+1}(U,\partial U)\\
G&=&K_{q+1}(W)\cong K_{q+1}(M_0,U)
\end{eqnarray*} 
Let $\gamma \colon G=K_{q+1}(W) \ra F=K_{q+1}(W,M)$ and $\mu \colon G=K_{q+1}(W) \ra F^* \cong K_{q+1}(W,M')$
be the maps induced by the long exact sequences of $(W,M)$ and $(W,M')$ and Poincar\'e-Lefschetz duality.
Then $\svec{\gamma\\\mu}$ is the inclusion of one lagrangian $G$ and $F$ is obviously another lagrangian
of the hyperbolic form $H_{(-)^q}(F)\cong K_q(\partial U)$. Hence the obstruction formation can also be written as
$(H_{(-)^q}(F); F, G)=\kreckfs$.

Different choices of presentations for $f\colon M\ra X$ change the
formations by a {\bf stable isomorphism}.  A stable isomorphism class
of formations is Wall's algebraic model for an odd-dimensional normal
map.  We will introduce this kind of isomorphism in Definition
\ref{weakisodef} as {\bf a weak isomorphism}, since we need to
distinguish it from another kind of isomorphism which we will
encounter in the discussion of the modified theory later.  A
presentation yields a homotopy equivalence $f'\colon M\ra X$ if and
only if its obstruction is a {\bf boundary} (see Definition
\ref{bdrdef}).  From this discussion equivalence relations can be
derived which are used to define the algebraic odd-dimensional surgery
obstruction groups $L_{2q+1}(\bZ[\pi_1(X)])$.  $z$ vanishes in that
$L$-group if and only if $f\colon M \ra X$ is cobordant rel$\partial$
to a homotopy equivalence.

\index{Surgery theory!Kreck}
\label{modifiedsec} In the 1980s Kreck generalized Wall's results. 
This memoir will only deal with his {\bf even-dimensional modified
theory} (compare \cite{Kreck99} p.724-732).  Surprisingly, certain
aspects of it resemble the traditional odd-dimensional theory we have
just discussed.

\index{Normal smoothings}
Kreck's theory starts off, not with normal maps, but the weaker notion
of {\bf normal smoothings}. Let $p\colon B \ra BO$ be a fibration. A
normal $B$-smoothing of a manifold is a factorization through $B$ of the classifying
map for its stable normal bundle.
(Under certain connectivity assumptions the homotopy type of $B$
actually depends only on $M$ but we will not need this fact. See
\cite{Kreck99} p. 711.)  Assume that there is a $(2q+2)$-dimensional cobordism
$(e,f,f')\colon (W,M,M') \ra B$ of normal smoothings such
that $f$ and $f'$ are $q$-equivalences.

If $B$ is a finite geometric Poincar\'e space and if $p\colon B \ra BO$ 
is its Spivak bundle, a normal $B$-smoothing is nothing but a
normal map. In addition, if $f$ and $f'$ are homotopy equivalences,
the situation is exactly the one of Wall's even-dimensional case. The
geometrical input of the modified case is considerably weaker than in
Wall's original theory.  There we started by comparing the complete
homotopy and normal bundle information of $M$ and $M'$ whereas in the
modified theory only ``half'' of that information is needed.

Surgery below the middle dimension yields a $(q+1)$-connected map $e'\colon W' \ra B$.
Now we can read off the obstruction which is a
tuple $z=\kreckfs$ with $G:=\ker (e'_*\colon \pi_{q+1}(W')\ra \pi_{q+1}(B))$, 
$\bar\theta \colon G \ra Q_{(-)^q}(\bZ[\pi_1(B)])$ the self-intersection form on $W'$, 
$F:=H_q(W',M')$ and $\gamma$ and $\mu$ are the compositions of the obvious maps in homology and the Hurewicz homomorphism.
Note that $\bar\theta +(-)^{q+1}\bar\theta^*=\gamma^*\mu$.
Such a tuple $\kreckfs$ with $(G,\gamma^*\mu, \bar\theta)$ a $\pm$-quadratic form will be called a {\bf preformation}. 

There is of course an ambiguity as there may be many ways to make
$e\colon W \ra B$ highly-connected. But the resulting manifolds $W'$
will only differ by a couple of tori $S^{q+1}\times S^{q+1}$ (see \cite{Kreck99}  p.729) and
algebraically, the obstruction preformations will differ only by
hyperbolic elements
$$\kreckf{P\oplus P^*}{1}{P\oplus P^*}{\svec{0&1\\(-)^{q+1}&0}}{\svec{0&1\\0&0}}$$
There is also a notion of isomorphism namely {\bf strong isomorphisms} (see Definition \ref{fisodef}).
If $(W,M,M')\ra B$ is changed by a diffeomorphism compatible with the normal $B$-smoothings.
the obstructions will change by such a strong isomorphism. 
The stable strong isomorphism classes (using hyperbolic preformations for stabilization) define the 
monoid $l_{2q+2}(\bZ[\pi_1(B)])$. $e\colon W\ra X$ is cobordant rel$\partial$ to an $s$-cobordism
if and only if the class of the obstruction preformation in the $l$-monoid is {\bf elementary}.
We will discuss elementariness later in Sections \ref{elemgeosec} and \ref{elemalgsec}.

The modified theory appeals by its ability to digest much simpler geometrical input than the traditional case:
\begin{enumerate}
\item The normal maps on $M$ and $M'$ are not necessarily homotopy equivalences but only $q$-equivalences.
\item The space $B$ with which we compare our cobordism does not need to be a finite Poincar\'e space.
\end{enumerate}
However, these generalizations come with a price tag: the difficult algebra. 
In fact, very little is known about the algebraic structure of those $l$-monoids compared
to the extensive literature that exists about the $L$-groups.

There is a striking similarity between the obstructions of the modified even-dimens\-ional case
and the traditional odd-dimensional case. In both cases one studies a cobordism 
$(e,f,f')\colon (W, M,M') \ra X$ of highly-connected normal maps or smoothings such that $H_{i}(W,M)=H_{i}(W,M')=0$
for $i\neq q+1$. To a certain degree the discussion about formations can be extended to preformations.
In both cases $\svec{\gamma\\\mu}\colon (G,0) \ra H_{(-)^q}(F)$ defines a map from a zero form to
a hyperbolic form.

The situation in the even-dimensional modified theory is of course more general: $X$ is not necessarily a finite Poincar\'e complex,
the maps involved don't need to be normal maps and $\partial M \ra X$ might not be a homotopy equivalence.
Therefore the map $\svec{\gamma\\\mu}$ is not always an inclusion of a lagrangian.

There are more differences if one looks at the equivalence relations in the obstruction groups/monoids.
In the modified even-dimensional case, the equivalence relations for the $l$-monoids are very strict, because
they seek to preserve all algebraic data of the whole cobordism $(e,f,f')\colon (W,M,M') \ra X$.
On the other hand, in the traditional odd-dimensional case that cobordism is just
used to write down the obstruction data of $f\colon M \ra X$. It is not important which cobordism is chosen
if only it is highly-connected. The equivalence relations for $L_{2q+1}(\Lambda)$ are designed such that
only the surgery-relevant information of $f\colon M\ra X$ is retained.

Hence, philosophically, a preformation can be interpreted as an algebraic model for
a cobordism $(e,f,f')\colon (W,M,M') \ra X$ or as a model for the map $f\colon M \ra X$ only.
In the first case, one identifies preformations via strong isomorphisms (the ones used 
to define the $l$-monoids) in the second case one uses weak isomorphisms (the ones used
to define the odd-dimensional $L$-groups). We will come back to this issue when
we introduce flip-isomorphisms in Chapter \ref{flipchap}.

In any case, the similarity of the objects in odd-dimensional $L$-theory and $l$-monoids
enable us to use algebraic surgery theory to investigate $l_{2q}(\Lambda)$. There are standard ways of translating
quadratic and symmetric complexes into forms and formations and vice versa. These
procedures extend to (regular) preformations as we will see in Chapter \ref{translchap}.

\section{Forms, Preformations and Formations: the Algebra}
\label{formalgsec}

{\bf Weakly finite} rings are a class of rings which the rank of any \fg free module is
well-defined. All the rings which we are interested in
(like fields, principal ideal domains, group rings, \etc) have this property
(\cite{Cohn89} pp.143-4 and \cite{Montgomery69}).

\begin{deflem}[(\cite{Cohn89} p.143)]
\index{Weakly finite ring}
A ring $\Lambda$ is {\bf weakly finite} if for any $n\in\bN_0$ and $\Lambda$-module $K$, 
$\Lambda^n \cong \Lambda^n\oplus K$ implies $K=0$.
\hfill\qed\end{deflem}

Let $\epsilon=(-)^q$.
(All constructions would equally work for $\epsilon \in \Lambda$ such that $\epsilon^{-1}=\bar\epsilon$.) 
Let $\Lambda$ be a weakly finite ring with $1$ and an {\bf involution} $x \mt \bar x$
(\ie an anti-automorphism 
$\Lambda \stackrel{\cong}{\ra} \Lambda^{op}, x \mt \bar x$ with $\Lambda^{op}$ the opposite ring).
\index{Involution}
All $\Lambda$-modules are left modules. 


\begin{rem}
In this section  all surgery obstruction groups and monoids of the various surgery theories 
presented before are defined. 
Strictly speaking, if we want to decide whether or whether not a cobordism can be turned
into an $s$-cobordism all the equivalence relations below must
only use \index{Simple isomorphisms}{\bf simple isomorphisms} (\ie isomorphisms for which the torsion in the Whitehead group
vanishes). We will ignore this condition in this thesis, hence the results will only deal
with $h$-cobordisms instead of $s$-cobordisms. 
A careful analysis of the proofs and constructions given in this thesis, will certainly lead
to similar results for the simple $l$-monoids.
%
\end{rem}

\subsection*{Forms and Even-Dimensional $L$-Groups}
\begin{defi}
\begin{enumerate}
\item
Let $M$ be a $\Lambda$-module.
\index{Duality involution map}\index{Tepsilon@$T_\epsilon$}
Using the canonical homomorphism $M \ra M^{**}$ 
we define the {\bf $\epsilon$-duality involution map}
\begin{eqnarray*}
T_\epsilon \colon \Hom_\Lambda(M,M^*) &\ra&\Hom_\Lambda(M,M^*)\\
\phi &\mt& (\epsilon\phi^*\colon x \mt (y\mt \epsilon \overline{\phi(y)(x)}))
\end{eqnarray*}
and the abelian groups \index{QM1@$Q^\epsilon(M)$ for a module}\index{QM2@$Q_\epsilon(M)$ for a module}
\begin{align*}
Q^\epsilon(M)&=\ker(1-T_\epsilon) &
Q_\epsilon(M)&=\coker(1-T_\epsilon)\\
Q^\epsilon(\Lambda)&=\{b\in\Lambda|b=\epsilon\overline{b}\} &
Q_\epsilon(\Lambda)&=\Lambda/\{b-\epsilon\overline{b}|b\in\Lambda\}
\end{align*}

\item \index{Symmetric form}
An {\bf $\epsilon$-symmetric form $(M,\phi)$ over $\Lambda$} is a 
$\Lambda$-module $M$ together with a $\phi \in Q^\epsilon(M)$.
It is {\bf non-singular} if $\phi\colon M \ra M^*$ is an isomorphism of $\Lambda$-modules. 
\index{Symmetric form!non-singular}

\item \index{Sublagrangian}
A {\bf sublagrangian $L$} of an $\epsilon$-symmetric form $(M,\phi)$
is a direct summand $j \colon L \hookrightarrow M$ such that $j^*\phi j=0$.
Then the {\bf annihilator} \index{Annihilator}
$L^\perp = \ker\left(j^*\phi \colon M \ra L^*\right)$
contains $L$. $L$ is a {\bf lagrangian} \index{Lagrangian} if $L=L^\perp$.
A form that allows a lagrangian is called {\bf metabolic}.\index{Metabolic}\index{Symmetric form!metabolic}

\item \index{Quadratic form}
An {\bf $\epsilon$-quadratic form $(M,\lambda,\nu)$ over $\Lambda$} is an $\epsilon$-symmetric form $(M,\lambda)$
together with a map $\nu \colon M \ra Q_\epsilon(\Lambda)$ such that for all $x,y\in M$ and $a\in\Lambda$
\begin{enumerate}
\item $\nu(x+y)-\nu(x)-\nu(y)=\lambda(x,y) \in Q_\epsilon(\Lambda)$
\item $\nu(x)+\epsilon\overline{\nu(x)} = \lambda(x,x) \in Q^\epsilon(\Lambda)$
\item $\nu(ax)=a\nu(x)\overline{a} \in Q_\epsilon(\Lambda)$
\end{enumerate}
$\nu$ is called a {\bf quadratic refinement} of the $\epsilon$-symmetric form $(M,\lambda)$.\index{Quadratic refinement}
It is {\bf non-singular} if the underlying symmetric form is non-singular.\index{Quadratic form!non-singular}

\item \index{Sublagrangian}
A {\bf sublagrangian $L$} of an $\epsilon$-quadratic form $(M,\lambda,\nu)$
is a direct summand $j \colon L \hookrightarrow M$ such that $j^*\lambda j=0$
and $\nu j=0$. Then the annihilator $L^\perp$ of the underlying $\epsilon$-symmetric
form $(M,\lambda)$ contains $L$. $L$ is a {\bf lagrangian} if $L=L^\perp$.
A form which allows a lagrangian is called {\bf metabolic}.\index{Metabolic}\index{Quadratic form!metabolic}

\item \index{Symmetric form!morphism}
A {\bf morphism $f\colon (M,\lambda) \ra (M',\lambda')$ of $\epsilon$-symmetric forms} 
is a map $f\in\Hom_\Lambda(M,M')$ such that $f^*\lambda' f=\lambda$. It is an {\bf isomorphism}
if $f\colon M \linebreak \ra M'$ is an isomorphism of $\Lambda$-modules.\index{Symmetric form!isomorphism}
\index{Quadratic form!morphism}
A {\bf morphism $f\colon (M,\lambda, \mu) \ra (M',\lambda', \mu')$ of $\epsilon$-quadratic forms} 
is a morphism $f\colon (M,\lambda) \ra (M',\lambda')$ of $\epsilon$-symmetric forms such that
$\mu'f=\mu$. It is an {\bf isomorphism} if $f\colon M \ra M'$ is an isomorphism of $\Lambda$-modules.
\index{Quadratic form!isomorphism}
\hfill\qed
\end{enumerate}
\end{defi}

\begin{rem}
\label{splitformrem} \index{Quadratic form}
For a \fg projective $\Lambda$-module $M$ there is no difference between the definition of $\epsilon$-quadratic forms and the following
alternative (see \eg \cite{Ran80a} p.117ff):

An {\bf $\epsilon$-quadratic form $(M,\psi)$ over $\Lambda$} is a tuple consisting of 
a \fg projective $\Lambda$-module $M$ together with an element $\psi\in Q_\epsilon(M)$.
It is {\bf non-singular} if $(1+T_\epsilon)\psi\colon M\ra M^*$ is an isomorphism of $\Lambda$-modules.
\index{Quadratic-form!non-singular}

A {\bf sublagrangian $L$}\index{Sublagrangian} of an $\epsilon$-quadratic form $(M,\psi)$
is a direct summand $j \colon L \hookrightarrow M$ such that $j^*\psi j=0 \in Q_\epsilon(L)$.
Then the {\bf annihilator} \index{Annihilator}
$$L^\perp = \ker\left(j^*(1+T_\epsilon)\psi \colon M \ra L^*\right)$$
contains $L$. $L$ is a {\bf lagrangian} if $L=L^\perp$.
A form which allows a lagrangian is called {\bf metabolic}.\index{Metabolic}\index{Quadratic form!metabolic}

A {\bf morphism $f\colon (M,\psi) \ra (M',\psi')$}\index{Quadratic form!morphism} is a map $f\in\Hom_\Lambda(M,M')$ 
such that $f^*\psi' f=\psi\in Q_\epsilon(M)$. It is an {\bf isomorphism}\index{Quadratic form!isomorphism}
if $f\colon M \ra M'$ is an isomorphism of $\Lambda$-modules.

An $\epsilon$-quadratic form $(M,\psi\in Q_\epsilon(M))$ defines an $\epsilon$-quadratic form $(M,(1+T_\epsilon)\psi,\nu)$
with $\nu(x)=\psi(x)(x)$. Conversely any $\epsilon$-quadratic form $(M,\lambda,\nu)$ gives rise
to an $\epsilon$-quadratic form $(M,\psi\in Q_\epsilon(M))$ (see \cite{Ran02} Proposition 11.9).
\end{rem}

If $\frac{1}{2}\in\Lambda$  then quadratic and symmetric forms are the same. More generally:
\begin{prop}
Assume there is a central $s\in\Lambda$ such that $s+\bar s=1$. Then there is a one-to-one
correspondence between $\epsilon$-quadratic and $\epsilon$-symmetric forms over \fg free (or projective) $\Lambda$-modules
given by $(G,\theta \in Q_\epsilon(G)) \mt (G,(1+T_\epsilon)\theta \in Q^\epsilon(G))$. Its inverse is
$(G,\lambda \in Q^\epsilon(G)) \mt (G,[s\lambda] \in Q_\epsilon(G))$.
\end{prop}

It can be shown (\cite{Ran80a} Proposition 2.2) that any metabolic form is isometric to
a hyperbolic form.
\begin{defi}
For any $(-\epsilon)$-symmetric form $(L^*,\phi)$ over a \fg free $\Lambda$-module $L$ we define the
non-singular {\bf hyperbolic $\epsilon$-symmetric form}\index{Hyperbolic form}
$$
    H^\epsilon(L,\phi)=\left(L \oplus L^*, \mat{0&1\\\epsilon&\phi}\in Q_\epsilon(L\oplus L^*)\right)
$$ 
We abbreviate $H^\epsilon(L)=H^\epsilon(L,0)$.

For any \fg free $\Lambda$-module $L$ we define the non-singular {\bf quadratic hyperbolic $\epsilon$-quadratic form}
$$
    H_\epsilon(L)=\left(L \oplus L^*, \mat{0&1\\0&0}\in Q_\epsilon(L\oplus L^*)\right)
$$
\hfill\qed
\end{defi}

\begin{deflem}
\label{formwittdef}\index{L-group@$L$-group!quadratic}\index{Witt-group!of forms}\index{L-group1@$L_{2q}(\Lambda)$}
{\bf The even-dimensional quadratic $L$-group\linebreak $L_{2q}(\Lambda)$} is the set of equivalence class of all 
non-singular $\epsilon$-quadratic forms on \fg free modules over $\Lambda$
where two forms are equivalent if they are isometric up to the addition of metabolic (\ie hyperbolic) forms.
It is also called the {\bf Witt-group of quadratic forms}.

Similarly we can define  {\bf the even-dimensional symmetric $L$-group $L^{2q}(\Lambda)$}.
\index{L-group@$L$-group!symmetric}\index{Witt-group!of forms}\index{L-group2@$L^{2q}(\Lambda)$}
\hfill\qed
\end{deflem}

\subsection*{Preformations and Even-Dimensional $l$-Monoids}

The building blocks of the even-dimensional $l$-monoids are {\bf preformations}.
A special case are {\bf formations} which help to define the odd-dimensional $L$-groups (see below).
We also introduce the notion of {\bf regular preformations}. They are preformations where all modules involved
are \fg free. Only they can be fed into the algebraic surgery machine which we will present in later chapters.
Corollary \ref{mkregcor} shows that the restriction to regular preformations is not a serious limitation
of the scope of our theory.
\begin{defi}
\label{formdef}
\begin{enumerate}
\item
An {\bf $\epsilon$-quadratic preformation $\skreckfs$}\index{Preformation!quadratic} is a 
tuple consisting of a free \fg $\Lambda$-module $F$, a \fg $\Lambda$-module $G$ and
$\svec{\gamma \\ \mu}\in\Hom_\Lambda(G, F\oplus F^*)$ such that $(G, \gamma^*\mu)$ is a $(-\epsilon)$-symmetric form.

An {\bf $\epsilon$-quadratic split preformation $\kreckfs$}\index{Preformation!split quadratic} is an
$\epsilon$-quadratic preformation $\skreckfs$ and a map 
$\bar\theta\colon G \rightarrow Q_{-\epsilon}(\Lambda)$ such that 
$(G, \gamma^*\mu, \bar\theta)$ is a $(-\epsilon)$-quadratic form. 

\item An $\epsilon$-quadratic preformation $\skreckfs$ is called {\bf regular}\index{Preformation!regular}
if $G$ is free. An $\epsilon$-quadratic split preformation $\kreckfs$
is regular if $G$ is a free.
In that case we interpret $\bar\theta \in Q_{-\epsilon}(G)$ as in Remark \ref{splitformrem}.

An $\epsilon$-quadratic split preformation $\kreckfs$ is called an 
{\bf $\epsilon$-quadratic split formation}\index{Formation!split quadratic} if $\svec{\gamma \\\mu}G$ is a sublagrangian of 
the $\epsilon$-quadratic hyperbolic form $H_{\epsilon}(F)$. It is called {\bf non-singular}\index{Formation!non-singular}
if the sublagrangian is indeed a lagrangian. Similar for the non-split case.\index{Formation!quadratic}
\hfill\qed
\end{enumerate}
\end{defi}

\begin{rem}
In Andrew Ranicki's work the notation for $\epsilon$-quadratic formations $\skreckfs$ is
$(H_{\epsilon}(F), F, G)$ and for $\epsilon$-quadratic split formations $\kreckfs$
it is $\left(F,\left(\svec{\gamma\\\mu},\bar\theta\right)G\right)$.
\end{rem}

In an obvious way all symmetric and quadratic forms are non-singular formations:
\begin{defi}
\label{bdrdef}
\begin{enumerate}
\item \index{Boundary!of a form}\index{Preformation!boundary}\index{$\partial(K,\psi)$}
    Let $(K,\lambda)$ be a $(-\epsilon)$-symmetric form on a free \fg $\Lambda$-module $K$. Then {\bf the boundary of $(K,\lambda)$}
      is the non-singular $\epsilon$-quadratic formation $\partial(K,\lambda)=\skreckfa{K}{1_K}{K}{\lambda}$.
\item Let $(K,\theta)$  be a $(-\epsilon)$-quadratic form on a free \fg $\Lambda$-module $K$. Then {\bf the boundary of $(K,\theta)$} 
      is the non-singular $\epsilon$-quadratic split formation $\partial(K,\theta)=\kreckf{K}{1_K}{K}{\theta-\epsilon\theta^*}{\theta}$.
\item A {\bf trivial formation}\index{Formation!trivial} is a non-singular $\epsilon$-quadratic split formation of the form $(P,P^*)=\kreckf{P}{0}{P}{1}{0}$
      with $P$ a free \fg $\Lambda$-module. Similar for the non-split case.
\hfill\qed
\end{enumerate}
\end{defi}

Now we define strong isomorphisms and stable strong isomorphism for preformations. 

\begin{defi}
\label{fisodef}
\label{mldef}
\begin{enumerate}\index{Preformation!sum}\index{Sum of preformations}
\item The {\bf sum of two $\epsilon$-quadratic split preformations $x=\kreckfs$ and $x'=\kreckf{F'}{\gamma'}{G'}{\mu'}{\bar\theta'}$} 
is the $\epsilon$-quadratic split preformations 
$
    x+x':=\kreckf{(F\oplus F')}{\gamma\oplus\gamma'}{G\oplus G'}{\mu\oplus\mu'}{\bar\theta\oplus\bar\theta'}.
$
Similar for the non-split case.

\item
\label{strisodef}\index{Preformation!strong isomorphism}
A {\bf strong isomorphism of two $\epsilon$-quadratic split preformations $\kreckfs$ and $\kreckf{F'}{\gamma'}{G'}{\mu'}{\bar\theta'}$}
is a tuple $(\alpha, \beta)$ of isomorphisms $\alpha\in\Hom_\Lambda(F,F')$ 
and $\beta\in\Hom_\Lambda(G,G')$ such that
\begin{eqnarray}
\label{isodia}
\xymatrix
{
 F \ar[d]^{\alpha}_{\cong} &G \ar[l]_{\gamma} \ar[d]^{\beta}_{\cong}
\ar[r]^{\mu} &F^* \ar[d]^{\alpha^{-*}}_{\cong} \\
  F' &G' \ar[l]_{\gamma'} \ar[r]^{\mu'} &{F'}^*
}
\end{eqnarray}
commutes and $\bar\theta = \bar\theta'\beta$.

\item \index{Preformation!stable strong isomorphism}
Two $\epsilon$-quadratic split preformations
$x$ and $x'$ are {\bf stably strongly isomorphic}
if there are boundaries of hyperbolic forms $h$ and $h'$ such that
there is a strong isomorphism between $x+h$ and $x'+h'$.
Similar for the non-split case.

\item \index{l-monoid@$l$-monoid}
{\bf The $l$-monoid $l'_{2q+2}(\Lambda)$}\index{l-monoid1@$l'_{2q+2}(\Lambda)$} 
is the set of stably strongly isomorphism classes of $\epsilon$-quadratic split preformations.
The ``simple'' version of this monoid is Kreck's original $l$-monoid.

{\bf The $l$-monoid $l_{2q+2}(\Lambda)$}\index{l-monoid2@$l_{2q+2}(\Lambda)$} 
is the set of stably strongly isomorphism classes of regular $\epsilon$-quadratic split preformations.

{\bf The $l$-monoid ${l'}^{2q+2}(\Lambda)$}\index{l-monoid3@${l'}^{2q+2}(\Lambda)$}
is the set of stably strongly isomorphism classes of $\epsilon$-quadratic preformations.

{\bf The $l$-monoid $l^{2q+2}(\Lambda)$}\index{l-monoid4@$l^{2q+2}(\Lambda)$}
is the set of stably strongly isomorphism classes of regular $\epsilon$-quadratic preformations.

All $l$-monoids are abelian monoids with zero. 
%
\hfill\qed
\end{enumerate}
\end{defi}

\subsection*{Formations and Odd-Dimensional $L$-Groups}
\index{Surgery theory!Wall!odd-dimensional}
In Section \ref{formgeosec} we explained that odd-dimensional traditional surgery theory
and the modified even-dimensional case use similar obstructions but that the equivalence
relations used in the construction of the obstruction groups/monoids are very different.
In both cases the obstruction associated to a $(2q+2)$-dimensional cobordism 
$(e,f,f')\colon (W,M,M') \ra X\times(I,0,1)$ of normal maps/smoothings with $f$ and $f'$ and $e$ highly-connected is 
some $\epsilon$-quadratic split preformation $\kreckfs$ with $F=H_{q+1}(W,M)$ and $(G,\gamma^*\mu,\bar\theta)$
containing the self-intersection form on some homology or homotopy group related to $e$.

In Kreck's surgery theory this cobordism is the very surgery problem in question. The equivalence
relations (strong isomorphism, stabilization with boundaries of hyperbolic forms) are very rigid
and preserve the important data of the whole cobordism.

In contrast, traditional odd-dimensional surgery theory uses the cobordism 
$$(e,f,f')\colon (W,M,M') \ra X\times(I,0,1)$$
just as a prop to define an obstruction to the odd-dimensional surgery problem given by $f\colon M \ra X$.
Hence the equivalence relations we present below are much more flexible - they need to filter out
ambiguities which arise by the choice of a different presentation \ie another $(2q+2)$-dimensional normal
cobordism $$(\hat e, f, \hat f')\colon (\hat W, M, \hat M') \ra X\times(I,0,1)$$ with $\hat e$ and 
$\hat f'$ highly-connected.

This leads to the unfortunate situation that there exist two notions of (stable) isomorphisms for 
preformations. Stable strong isomorphism classes of preformations are algebraic models for
diffeomorphism classes of $(2q+2)$-dimensional cobordisms $(e,f,f')\colon (W,M,M') \ra X\times(I,0,1)$ whereas stable weak isomorphism classes
are models for diffeomorphism classes of $(2q+1)$-dimensional normal maps $f\colon M \ra X$. 
Weak isomorphisms will reappear in Chapter \ref{flipchap} which
deals with {\bf flip-isomorphisms}.

We also have to be careful about stabilization. In the $l$-monoids we stabilize with ``hyperbolics''
(\ie boundaries of hyperbolics) and in $L$-theory we use trivial preformations.

\begin{defi}
\label{weakisodef}
\begin{enumerate}
\item \index{Preformation!weak isomorphism}
A {\bf (weak) isomorphism $(\alpha,\beta,\sigma)$ of two regular $\epsilon$-quadra\-tic preformations
$\skreckfs$ and $\skreckfa{F'}{\gamma'}{G'}{\mu'}$} is a triple consisting of an 
isomorphism $\alpha\in\Hom_\Lambda(F,F')$,
an isomorphism $\beta\in\Hom_\Lambda(G,G')$ and $\sigma\in Q^{-\epsilon}(F^*)$ such that
\begin{enumerate}
\item $\alpha\gamma+\alpha\sigma\mu=\gamma'\beta \in\Hom_\Lambda(G, F')$
\item $\alpha^{-*}\mu =\mu'\beta \in \Hom_\Lambda(G, {F'}^*)$.
\end{enumerate}

\item \index{Preformation!weak isomorphism}
A {\bf (weak) isomorphism $(\alpha,\beta,\nu)$ of two regular $\epsilon$-quadra\-tic split preformations
$\kreckfs$ and $\kreckf{F'}{\gamma'}{G'}{\mu'}{\bar\theta'}$} is a triple 
consisting of an isomorphism $\alpha\in\Hom_\Lambda(F,F')$,
an isomorphism $\beta\in\Hom_\Lambda(G,G')$ and $\nu\in Q_{-\epsilon}(F^*)$ such that
\begin{enumerate}
\item $\alpha\gamma+\alpha(\nu-\epsilon\nu^*)^*\mu=\gamma'\beta \in\Hom_\Lambda(G, F')$
\item $\alpha^{-*}\mu =\mu'\beta \in \Hom_\Lambda(G, {F'}^*)$
\item $\bar\theta+\mu^*\nu\mu=\beta^*\bar\theta'\beta\in Q_{-\epsilon}(G)$
\end{enumerate}

\item \index{Preformation!stable weak isomorphism}
A {\bf stable weak isomorphism of two regular $\epsilon$-quadratic split preformations $z$ and $z'$}
is a weak isomorphism $z+t \cong z'+t'$ for trivial formations $t$, $t'$.
\hfill\qed
\end{enumerate}
\end{defi}

An odd-dimensional normal map is cobordant to a homotopy-equivalence if and only if
its obstruction (pre-)formation is stably isomorphic to a boundary. It can be shown that for
any form there is another cobordant odd-dimensional map whose obstruction preformation differs from the original one
by the boundary of a form (see \cite{Ran02} Proposition 12.13, Theorem 12.29).
This motivates the definition of the odd-dimensional surgery obstruction groups.
\begin{deflem}[\cite{Ran02} Definition 12.23, Proposition 12.33]
\index{L-group@$L$-group!formations}\index{Witt-group!of formations}\index{L-group3@$L_{2q+1}(\Lambda)$}
We call two non-singular $\epsilon$-quadratic (split) formations $z$ and $z'$ equivalent if 
there is a stable weak isomorphism between $z+b$ and $z'+b'$ for some boundaries $b$ and $b'$.
In both cases (split and not-split) the equivalence classes form the {\bf odd-dimensional $L$-groups $L_{2q+1}(\Lambda)$}.
\hfill\qed\end{deflem}

\begin{rem}
There are also odd-dimensional symmetric $L$-groups which are defined as the Witt-group of {\bf $\epsilon$-symmetric formations}
(see \cite{Ran80a} Chapter 5).
\end{rem}

\begin{rem}
\label{weakisorem} 
\begin{enumerate}
\item
Any strong isomorphism between regular preformations is also a weak isomorphism.

\item
Let $z=\kreckfs$ and $z'=\kreckf{F'}{\gamma'}{G'}{\mu'}{\bar \theta'}$ be regular $\epsilon$-quadratic split preformations and
$t=(\alpha,\beta,\nu)\colon z \ra z'$ a weak isomorphism.
Then $(\alpha, \beta, (\nu-\epsilon\nu^*)^*)$ is a weak isomorphism
of the underlying regular $\epsilon$-quadratic preformation $\skreckfs$ and $\skreckfa{F'}{\gamma'}{G'}{\mu'}$.

\item 
If $z$ is a formation, weak isomorphisms 
are nothing but isomorphisms of quadratic (split) formations 
as defined in \cite{Ran80a} p.122 and p.128. For example, an isomorphism $(\alpha,\beta,\sigma)$ of 
$\epsilon$-quadratic formations $z=\skreckfs$
and $z'=\skreckf{F'}{\gamma'}{G'}{\mu'}$ is an isomorphism 
of $\epsilon$-quadratic hyperbolic forms
\begin{eqnarray*}
\mat{\alpha&\alpha\sigma\\0&\alpha^{-*}}\colon H^{\epsilon}(F) \isora H^{\epsilon}(F')
\end{eqnarray*}
which maps the (sub-)lagrangians $F$ and $G$ onto $F'$ and $G'$ respectively.

\item 
\label{isorem2}
Despite the different ways of stabilizing, every stable strong isomorphism is also a stable weak isomorphism.
That's because there is a weak isomorphism between an
even-dimensional trivial formation and a boundary of hyperbolic forms:
$$\left(1, \mat{0&1\\-\epsilon&0}, \mat{0&0\\\epsilon&0}\right)\colon \partial H_{\epsilon}(P) \ra (P\oplus P^*,(P\oplus P^*)^*)$$

The converse is not true: Let $Q$ be a free $\Lambda$-module of rank $1$. 
Let $y=(Q,Q^*)$ and $z=\partial H_{-\epsilon}(Q)$.
By the above, both preformations are stably weakly isomorphic, 
but for rank reasons they cannot be stably strongly isomorphic.
\end{enumerate}
\end{rem}

\section{Elementariness: the Geometry}
\label{elemgeosec}

\index{Surgery theory!Wall!even-dimensional}
We haven't quite explained yet how the obstruction preformation in Kreck's surgery theory 
can tell us whether surgery is able to turn our cobordism into an $h$-cobordism. 
Again let $q\ge 2$. In the case of the
{\bf traditional even-dimensional case} we only have to check that
the obstruction $(K_{q+1}(W), \lambda, \mu)$ of a highly-connected $(2q+2)$-dimensional normal cobordism 
$(W, M, M')\ra X\times (I,0,1)$ is zero in the Witt-group $L_{2q+2}(\bZ[\pi_1(X)])$. Then we know
that there is a lagrangian $L$ of $K_{q+1}(W)\oplus H_{(-)^{q+1}}(K)$. The stable
lagrangian $L$ is a recipe for successful surgery: we
perform $\rk K$ trivial surgeries on $W$ (with the result $W\#\#_{\rk K}S^{q+1}\times S^{q+1}$) 
and then kill a basis of $L$
via surgery. The result will be an $h$-cobordism.

\index{Surgery theory!Kreck}
In the modified case the criteria for success or failure are more complicated.
The starting point of Kreck's {\bf modified surgery theory} 
is the situation we described on p.\pageref{modifiedsec}f:
Let $p\colon B \rightarrow BO$ be a fibration.  Let $M_0$ and $M_1$ be
$(2q+1)$-dimensional manifolds with $(q-1)$-smoothings $f$ and $f'$ in $B$.
Let $f\colon \partial M_0 \isora \partial M_1$ be a diffeomorphism compatible
with the smoothings.
There is a normal smooth cobordism \ie
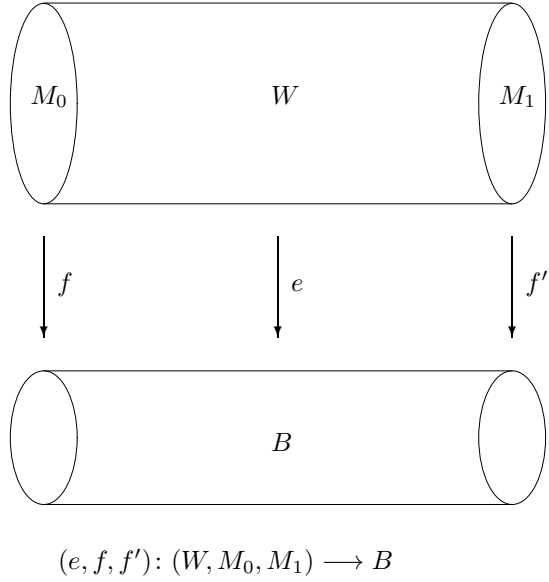
\begin{center}
\setlength{\unitlength}{0.007in}
\begin{picture}(750,400)(0,0)
\put(200,300){\ellipse{50}{150}}
\put(550,300){\ellipse{50}{150}}
\path(200,375)(550,375)
\path(200,225)(550,225)
\put(190,300){$M_0$}
\put(370,300){$W$}
\put(540,300){${M_1}$}

\put(200,200){\vector(0,-1){75}}
\put(375,200){\vector(0,-1){75}}
\put(550,200){\vector(0,-1){75}}
\put(210,160){$f$}
\put(385,160){$e$}
\put(560,160){$f'$}

\put(200,50){\ellipse{50}{100}}
\put(550,50){\ellipse{50}{100}}
\path(200,0)(550,0)
\path(200,100)(550,100)
\put(370,40){$B$}
\end{picture}
\end{center}
\begin{eqnarray*}
(e,f,f')\colon (W, M_0,M_1) \ra B
\end{eqnarray*}
Surgery below the middle dimension on $W$ is possible and yields a $q$-smoothing $e'\colon W' \ra B$.
Then we can define a $(-)^{q+1}$-quadratic split preformation.
\begin{defi}
\label{obstrdef}\index{Surgery obstruction}
\begin{eqnarray*}
y(W)&=&\kreckf{F}{\sigma}{H}{\tau}{\psi}\\
&=&(H_{q+1}(W',M_0)\leftarrow \im (\pi_{q+2}(B,W') \rightarrow \pi_{q+1}(W'))\rightarrow H_{q+1}(W',M_1), \psi)\\
  &\in& l'_{2q+2}(\bZ[\pi_1(X)])
\end{eqnarray*}
is the Kreck surgery obstruction of $W$.
\hfill\qed\end{defi}
It turns out that killing low-dimensional homotopy classes by
surgery in a different manner will not change the class $y(W) \in l'_{2q+2}(\Lambda)$.
The obstruction contains all the data to find out whether
$W$ can be made an $h$-cobordism due to the main theorem: 
\begin{thm}[\cite{Kreck99} Theorem 3 and Remark p.730ff]
\label{mainthm}
$W$ is $B$-cobordant rel$\partial$ to an $h$-cobordism if and
only if $y(W)\in l'_{2q+2}(\Lambda)$ is elementary.
\end{thm}

Before we give a strict definition of {\bf elementariness} (see the next section), 
a heuristical argument provides some geometric
motivation for this new concept. 
We will later show (Corollary \ref{niceobstrcor}) that in Definition \ref{obstrdef} 
$\im (d\colon \pi_{q+2}(B,W') \ra \pi_{q+1}(W'))$ can be replaced by $H_{q+2}(B,W')$.
Without loss of generality we assume that $W \ra B$ is $(q+1)$-connected and $W=W'$.

From the long exact sequences for $(B,W,M_i)$ we learn that $H_j(W,M_i)=0$ for $j \le q$, $i=0,1$ and
that there is an exact sequence
$$0 \ra H_{q+2}(B,M_i) \ra H_{q+2}(B,W) \ra H_{q+1}(W,M_i) \ra H_{q+1}(B, M_i)\ra 0$$
and that the canonical maps $H_j(B,M_i) \ra H_j(B,W)$ are isomorphisms for $j \ge q+3$.
Then a (possibly non-free) chain complex model of the cobordism $(W, M_0, M_1) \rightarrow B$ looks like
\begin{eqnarray*}
\xymatrix@-10pt
{
\vdots\ar[d]_0&
\vdots\ar[d]_0&
\vdots\ar[d]_0
\\
C_{j}(B,M_0)=H_j(B,M_0) \ar[r]_{\cong}\ar[d]_0&
C_{j}(B,W)=H_j(B,W) \ar[d]_0&
C_{j}(B,M_1)=H_j(B,M_1) \ar[l]^{\cong}\ar[d]_0
\\
\vdots\ar[d]_0&
\vdots\ar[d]_0&
\vdots\ar[d]_0
\\
C_{q+2}(B,M_0)=G \ar[r]_{1} \ar[d]_{\gamma}&                
C_{q+2}(B,W)=G&         
C_{q+2}(B,M_1)=G \ar[l]^{1} \ar[d]_{\mu}
\\
C_{q+1}(B,M_0)=F &
&
C_{q+1}(B,M_1)=F^*
}
\end{eqnarray*}
Assume it is possible to do simultaneous surgery on $W$ killing some homology classes
$x_1,\dots, x_k \in H_{q+1}(W)$ without changing the boundary. 
Assume further that they are the basis of a \fg free submodule 
 $j\colon U=\langle x_1,\dots, x_k \rangle\hookrightarrow H_{q+1}(W)$.
A chain complex model for the resulting cobordism $(V, M_0, M_1) \ra B$ is
\begin{eqnarray*}
\xymatrix
{
&\vdots\ar[d]_0&
\\
\vdots      \ar[d]_0&
C_{q+3}(B,V)=U  \ar[d]_{j}&
\vdots      \ar[d]_0
&
\\
C_{q+2}(B,M_0)=G \ar[r]^{1}\ar[d]_{\gamma}&
C_{q+2}(B,V)=G \ar[d]_{j^*\gamma^*\mu}&
C_{q+2}(B,M_1)=G \ar[l]_{1} \ar[d]_{\mu}
\\
C_{q+1}(B,M_0)=F \ar[r]_{\pm j^*\mu^*}&
C_{q+1}(B,V)=U^*&
C_{q+1}(B,M_1)=F^* \ar[l]^{j^*\gamma^*}
}
\end{eqnarray*}
(compare with the proof of Theorem \ref{algmainthm}). We observe that the relative 
middle-dimensional homology groups of the new cobordism are
\begin{eqnarray*}
H_{q+1}(V,M_0)= \dfrac{\im(j^*\gamma^*)}{\ker(\mu j)}\\
H_{q+1}(V,M_1)= \dfrac{\im(j^*\mu^*)}{\ker(\gamma j)}
\end{eqnarray*}

Using Poincar\'e-Lefschetz duality this means that $(V,M_0,M_1)$ is an $h$-cobordism if
and only if these homologies vanish or, equivalently, 
the mapping cones of either map $C(B,M_i) \ra C(B,V)$ (that is
$$
0 \ra U \stackrel{\gamma j}{\ra} F \stackrel{(\mu j)^*}{\ra} U^* \ra 0
$$
and its dual) are exact sequences. This is in fact one way of defining that $\skreckfs$ is elementary
(see Proposition \ref{elemprop}\ref{elemprop4}).

\section{Elementariness: the Algebra}
\label{elemalgsec}

First we will repeat the original definition of an elementary preformation before we
present alternative ways of looking at this concept.
\begin{defi}[\cite{Kreck99} p.730]
\label{elemdef}\index{Preformation!elementary}\index{Elementariness}
An $\epsilon$-quadratic split preformation $\kreckfs$ is {\bf
elementary (in respect to $U$)} if there is a \fg free submodule $j\colon U \hookrightarrow G$ with
\begin{enumerate}
\item $j^*\gamma^*\mu j=0$ and $\bar\theta j=0$,
\item $\gamma j$ and $\mu j$ are injective and their images $U_0$ and
$U_1$ are direct summands in $F$ and $F^*$ respectively, 
\item $R_1= F^*/U_1 \rightarrow U_0^*$, $f \mapsto f|U_0$ is an isomorphism.
\end{enumerate}

Such an $U$ is called an {\bf $h$-lagrangian} of the preformation.

An element in $l'_{2q+2}(\Lambda)$ is {\bf elementary} if it has an elementary representative.
All elementary elements form a submonoid ${l'}_{2q+2}^{el}(\Lambda)$ of $l'_{2q+2}(\Lambda)$.
\index{l5@${l'}_{2q+2}^{el}(\Lambda)$}

An $\epsilon$-quadratic split preformation $z=\kreckfs$ is {\bf stably elementary}
if $[z]\in l'_{2q+2}(\Lambda)$ is elementary.

Similar for non-split and regular preformations.
\hfill\qed\end{defi}


\begin{prop}
\label{elemprop} \index{Preformation!elementary}\index{Elementariness}
Let $\kreckfs$ be an $\epsilon$-quadratic split preformation.
Then the following statements are equivalent:
\begin{enumerate}

\item The preformation is elementary in respect to $U$.

\item \label{elemprop2} There is a \fg free submodule $j \colon U \hookrightarrow G$ such that
\begin{eqnarray*}
\xymatrix{
0 \ar[r]& 
U \ar[r]_{\mu j}&
F^* \ar[r]_{(\gamma j)^*}&
U^* \ar[r]&
0
}
\end{eqnarray*}
is an exact sequence and $\bar\theta|U=0$.

\item \label{elemprop4}
There is a \fg free submodule $j\colon U \hookrightarrow G$ such that the two horizontal chain maps
\begin{eqnarray*}
\xymatrix@C+30pt
{
&               U\ar[d]_{j}&\\
G\ar[d]_{\gamma}\ar[r]^{1}& G\ar[d]_{j^*\gamma^*\mu}&   G\ar[l]_{1}\ar[d]^{\mu}\\
F\ar[r]_{-\epsilon j^*\mu^*}&   U^*&                F^*\ar[l]^{j^*\gamma^*}
}
\end{eqnarray*}
are chain equivalences (\ie this is a ``chain complex model of an $h$-co\-bord\-ism'')
and $\bar\theta|U=0$.

\item \label{elemprop3}
The preformation is strongly isomorphic to a preformation of the form 
\begin{eqnarray*}
\kreckfb {U\oplus U^*}{\svec{1&0\\0&\sigma}}{U\oplus R}{\svec{0&-\epsilon\sigma\\1&\tau}}{U^*\oplus U}{\bar\theta}
\end{eqnarray*}
for some maps $\sigma\colon R \ra U^*$ and $\tau\colon R \ra U$ such that $\tau^*\sigma=-\epsilon\sigma^*\tau$
and a quadratic refinement $\theta'\colon R \ra Q_{-\epsilon}(\Lambda)$ of $\sigma^*\tau$ such that 
\begin{eqnarray*}
\bar\theta \colon U \oplus R &\ra& Q_{-\epsilon}(\Lambda)\\
(u,r) &\mt& \theta'(r) -\epsilon \sigma(r)(u)
\end{eqnarray*}
\end{enumerate}
Similar for the non-split case.
\end{prop}

\begin{proof}
The first two assertions are obviously equivalent. 
One observes that 
\begin{eqnarray*}
\xymatrix{
0 \ar[r]& 
U \ar[r]_{\mu j}&
F^* \ar[r]_{(\gamma j)^*}&
U^* \ar[r]&
0
}
\end{eqnarray*}
and its dual are the mapping cones of the chain maps in \ref{elemprop4}.
So \ref{elemprop4} is equivalent to \ref{elemprop2}.

Finally, we concentrate on \ref{elemprop3}.
Every preformation of the form described in there is elementary in respect to $U$. 
On the other hand let $\skreckfs$ be an elementary $\epsilon$-quadratic preformation. 
It is easy to show that $G=U\oplus R$ with $R=\ker(\pi\gamma)$.

Let $R_1\subset F^*$ be some complement of $U_1=\mu(U)$. We write 
\begin{eqnarray*}
\gamma=\mat{\gamma_1 & \gamma_2 \\ \gamma_3 & \gamma_4}\colon U \oplus R &\ra& U_0\oplus R_0\\
\mu'=\mat{\mu'_1&\mu'_2\\\mu'_3 & \mu'_4}\colon U \oplus R &\ra& U_1\oplus R_1\\
\Phi=\mat{\Phi_1&\Phi_2\\\Phi_3 & \Phi_4}\colon U_1\oplus R_1 &\ra& U_0^* \oplus R_0^*\\
f &\mt& (f|U_0, f|R_0)\\
\mu=\mat{\mu_1&\mu_2\\\mu_3 & \mu_4}\colon U \oplus R &\ra& U_0^*\oplus R_0^*\\
x&\mt& \Phi\mu'(x)
\end{eqnarray*}

By assumption, $\gamma_1$ and $\mu'_1$ are isomorphisms and $\gamma_3$ and $\mu'_3$
are vanishing. We can apply the strong isomorphism $\left(1_F, \svec{\gamma_1 & \gamma_2\\0&1}\right)$ to
achieve the simpler situation of $\gamma=\svec{1&0\\0&\gamma_4}$ and $U_0=U$.

We compute $\gamma^*\mu=\svec{\Phi_1\mu'_1&*\\ * & *}$ and see that $\Phi_1=0$. The last criterion of 
elementariness implies that $\Phi_2$ is an isomorphism and therefore $\Phi_3$ is bijective as well.
We use these facts to see that $$\mu=\Phi\mu'=\mat{0& \Phi_2\mu'_4\\\Phi_3\mu'_1 & \Phi_3\mu'_2+\Phi_4\mu'_4}=
\mat{0&\mu_2\\\mu_3&\mu_4}$$ Hence $\mu_3$ is an isomorphism.
Because $\gamma^*\mu=\svec{0 & \mu_2\\\gamma_4^*\mu_3 & \gamma_4^*\mu_4}$ is $(-\epsilon)$-symmetric, $\mu_2=-\epsilon\mu_3^*\gamma_4$.
We apply the strong isomorphism $(\svec{1 & 0\\0&\mu_3^*},1_G)$ and get a preformation with the properties we want.

In the case of quadratic split preformation, the same steps as before yield a strong isomorphism between
an $\epsilon$-quadratic split preformation $\kreckfs$ and 
\begin{eqnarray*}
\kreckfb {U\oplus U^*}{\svec{1&0\\0&\sigma}}{U\oplus R}{\svec{0&-\epsilon\sigma\\1&\tau}}{U^*\oplus U}{\bar\theta}
\end{eqnarray*}
Then define $\theta'=\bar\theta|R$. 
\end{proof}

The proposition allows us to derive some quite simple obstructions for elementariness.
\begin{cor}
\label{simpleobstr}
Let $z=\skreckfs$ be a regular $\epsilon$-quadratic preformation.
\begin{enumerate}
\item The isomorphism classes of kernels and cokernels of 
$\gamma$, $\mu$, $\svec{\gamma\\\mu}$, $\gamma^*\mu$ as well as $\rk G-\rk F\in\bZ$ and $\rk F \in\bZt$ 
are invariants of $[z]\in l^{2q+2}(\Lambda)$.
\item If $[z]\in l^{2q+2}(\Lambda)$ is elementary then $\ker\gamma \cong \ker\mu$, $\coker\gamma \cong \coker\mu$, 
$\ker\gamma^*\mu\cong \ker\gamma\oplus\ker\gamma^*$, 
$\coker\gamma^*\mu\cong \coker\gamma\oplus\coker\gamma^*$ and $\rk F$ is even
\end{enumerate}
Similar for the split and non-regular case.
\end{cor}


Finally we present a little lemma about elementariness which has two interesting applications.
\begin{lem}
\label{pushelemlem}
Let $x=\kreckfs$ and $y=\kreckf{F}{\sigma}{H}{\tau}{\bar\psi}$ be two $\epsilon$-quadratic split preformations
and $\pi \colon G \twoheadrightarrow H$ a surjective homomorphism such that 
\begin{eqnarray}
\label{litlemdia}
\xymatrix
{
F&
G \ar[l]_-{\gamma} \ar[r]^-{\mu} \ar@{>>}[d]_-{\pi}&
F^*\\
& H \ar[ul]^-{\sigma} \ar[ur]_-{\tau}
}
\end{eqnarray}
commutes and $\bar\theta=\bar\psi\pi$. Then $x$ is (stably) elementary if and only if $y$ is (stably) elementary.
Similar for the non-split case.
\end{lem}
%

As a first application we can slightly improve the elegance of the obstruction in Definition \ref{obstrdef}.
\begin{cor}
\label{niceobstrcor} \index{Surgery obstruction}
In the situation of Definition \ref{obstrdef} and Theorem \ref{mainthm}
we can define an alternative $(-)^q$-quadratic split preformation
\begin{eqnarray*}
x(W)&=&\kreckfs\\
&=&(H_{q+1}(W',M_0)\leftarrow H_{q+2}(B,W') \rightarrow H_{q+1}(W',M_1), \bar\theta)\\
&\in& l'_{2q+2}(\bZ[\pi_1(X)])
\end{eqnarray*}
with $\bar\theta$ being induced by the self-intersection form on $W$ and maps
$\gamma$ and $\mu$ from the long exact sequence of the triads
$(B,W,M_i)$.

Then $W$ is $B$-cobordant to an $h$-cobordism if and only if $x(W)\in  l'_{2q+2}(\bZ[\pi_1(X)])$ is 
elementary.

The long exact sequences of $(B,W,M_i)$ yields
\begin{align*}
\ker\gamma &= H_{q+2}(B,M_0)  &   \coker\gamma &=H_{q+1}(B,M_0)\\
\ker\mu    &= H_{q+2}(B,M_1)  &   \coker\mu    &=H_{q+1}(B,M_1)
\end{align*}
(Compare with Corollary~\ref{simpleobstr})
\end{cor}

The second application is a more theoretical: the decision
whether a preformation is elementary can always be replaced
by checking that a related regular preformation is
elementary
\begin{cor}
\label{mkregcor}
Let $x=\kreckf{F}{\sigma}{H}{\tau}{\psi}$ be an $\epsilon$-quadratic
split preformation and let $G$ be a free \fg module with an epimorphism
$\pi\colon G \twoheadrightarrow H$. Then there is an $\epsilon$-quadratic
split preformation $y=\kreckfs$ which makes the diagram \ref{litlemdia} commute.
$x$ is regular and it is elementary if and only if $y$ is.
\end{cor}

\chapter{Translating Kreck's Surgery into Algebraic Surgery Theory}
\label{translchap}

{\bf For the whole chapter let $q\ge 2$, $\epsilon=(-)^q$ and let $\Lambda$ be a weakly finite ring
with $1$ and involution.}

The asymmetric and quadratic signatures which will be defined in the next chapters are obstructions 
to the elementariness of a preformation. 
Constructions and proofs will use results from the vast theory of algebraic surgery. 

This section will provide the first step in the programme by translating preformations into the 
language of algebraic surgery theory:
quadratic Poincar\'e pairs and complexes (see {\bf Section \ref{constrpairsec}}).

Preformations arise as obstructions when we ask whether a cobordism $(W,M,M') \rightarrow B$ of normal 
smoothings is cobordant rel$\partial$ to an $h$-cobordism. As there is no realization result for preformations,
we cannot be sure whether they all arise from a surgery problem.
The constructions in Section \ref{constrpairsec} can be thought of as an
``algebraic realization result'': any preformation appears as an ``obstruction''
of a certain Poincar\'e pair to be cobordant rel$\partial$ to an algebraic $h$-cobordism. 
However, we will not try to develop a general $l$-obstruction theory for Poincar\'e pairs simply because we do not
need it. It suffices to create a quadratic chain complex model for a preformation
and to apply algebraic surgery theory to it.

Algebraic equivalents of concepts like cobordisms rel$\partial$ and surgery inside a manifold
will be needed to model Kreck's surgery theory. {\bf Section \ref{cobpairsec}} deals with this
rather technical issue and confirms our expectations, namely that those notions exist
and that they behave similarly to their geometric equivalents (e.g. that two Poincar\'e pairs
are cobordant if and only if one is the result of a surgery of the other).

In {\bf Section \ref{elemalgsursec}} we prove some kind of algebraic equivalent of Theorem \ref{mainthm}:
the Poincar\'e pair constructed in Section \ref{constrpairsec} is cobordant rel$\partial$
to an algebraic $h$-cobordism if and only if the preformation is (stably) elementary. This theorem is the key
to the application of algebraic surgery theory to the analysis of preformations.

\section{From Preformations to Quadratic Pairs}
\label{constrpairsec}
If we want to use the tools of algebraic surgery theory, we will need to translate preformations into
the language of quadratic chain complexes and pairs.
Readers can brush up their knowledge of algebraic surgery theory by reading
\cite{Ran80a} or the appendix 
(Chapter \ref{algsurchap}, p.~\pageref{algsurchap}).

The translation is easier for non-singular formations.
They can always be realized (\cite{Ran02} Proposition 12.17) as an
obstruction of a {\bf presentation}\index{Presentation} \ie a $(2q+2)$-dimensional cobordism of degree 1 normal maps
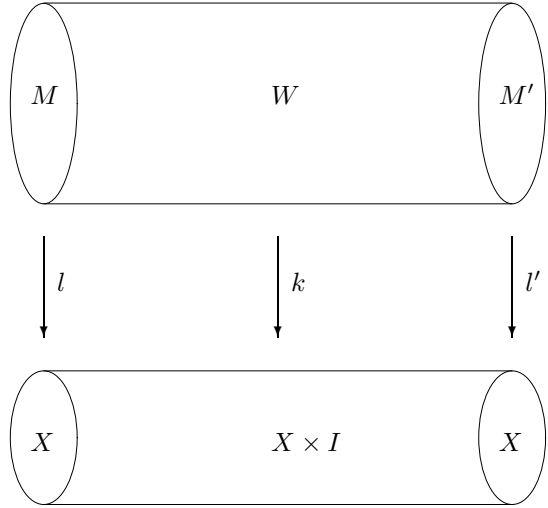
\begin{center}
\setlength{\unitlength}{0.007in}
\begin{picture}(750,400)(0,0)
\put(200,300){\ellipse{50}{150}}
\put(550,300){\ellipse{50}{150}}
\path(200,375)(550,375)
\path(200,225)(550,225)
\put(190,300){$M$}
\put(370,300){$W$}
\put(540,300){${M'}$}

\put(200,200){\vector(0,-1){75}}
\put(375,200){\vector(0,-1){75}}
\put(550,200){\vector(0,-1){75}}
\put(210,160){$l$}
\put(385,160){$k$}
\put(560,160){$l'$}

\put(200,50){\ellipse{50}{100}}
\put(550,50){\ellipse{50}{100}}
\path(200,0)(550,0)
\path(200,100)(550,100)
\put(190,40){$X$}
\put(370,40){$X\times I$}
\put(540,40){${X}$}
\end{picture}
\end{center}
\begin{eqnarray*}
(k,l,l')\colon (W, M,M') \ra X\times(I,0,1)
\end{eqnarray*}
into a finite geometric Poincar\'e pair $(X,\partial X)$ 
such that $l$ and $l'$ are $q$-connected, $k$ is $(q+1)$-connected and 
$l|\colon\partial M \ra \partial X$ is a homotopy equivalence. 
We note that a presentation is a special case of a Kreck surgery situation
but also a way to find the $L$-obstruction of the odd-dimensional normal map
$l\colon (M,\partial M)\ra (X,\partial X)$ (see also Section \ref {formgeosec}).
In both cases (see Corollary \ref{niceobstrcor} and \cite{Ran02} Chapter 12)
the obstruction is the non-singular $\epsilon$-quadratic
split formation $\kreckfs$ with $F=K_{q+1}(W,M)$, $G=K_{q+1}(W)$, \etc
It is elementary in $l_{2q+2}(\bZ[\pi_1(X)])$ if and only if $k\colon W \ra X\times I$ is cobordant rel$\partial$
to an $h$-cobordism and it vanishes in
$L_{2q+1}(\bZ[\pi_1(X)])$ if and only if $l\colon (M,\partial M)\ra (X,\partial X)$ is
cobordant rel$\partial$ to a homotopy equivalence \ie an $h$-cobordism (see also Section \ref {formgeosec}).

Algebraic surgery theory presents an alternative surgery obstruction for the normal map 
$l\colon (M,\partial M)\ra (X,\partial X)$: the {\bf quadratic kernel $(D,\nu)$ of $l$}\index{Quadratic kernel}.
It is a $(2q+1)$-dimensional {\bf quadratic Poincar\'e complex} over $\bZ[\pi_1(X)]$ where
$D=\cone(l^!)$ is the mapping cone of the so-called {\bf Umkehr chain map}\index{Umkehr chain map}
$$
l^!\colon C(\widetilde X) \homra C(\widetilde X, \widetilde{\partial X})^{2q+1-*} \stackrel{\widetilde l^*}{\ra} 
                 C(\widetilde M, \widetilde{\partial M})^{2q+1-*} \homra C(\widetilde M)
$$
with $\widetilde X$ and $\widetilde M$ the universal covers.
The homology modules of $D$ are the kernel modules $K_{*}(M)$.  The
quadratic structure $\nu\in W_\%(D)_{2q+1}$ is a family of maps $\nu_s \in \Hom(D^{2q+1-r-s}, D_r)$
which generalizes the self-intersection number.
It contains a chain equivalence $(1+T)\nu_0\colon D^{2q+1-*} \homra D$ 
inducing the Poincar\'e duality $K^{2q+1-*}(M) \isora K_*(M)$.
(For the details of the construction see \cite {Ran80b} Chapter 1 and 4.)

The algebraic surgery approach has two main advantages to
the traditional obstruction theory:
\begin{enumerate}
\item it works for normal maps $l\colon (M,\partial M)\ra (X,\partial X)$ which are {\bf not} highly-connected,
\item it provides a uniform obstruction theory for the odd- {\bf and} even-dimensional case. 
\end{enumerate}
There are notions of algebraic surgery and cobordism for Poincar\'e complexes. In the case of quadratic kernels they correspond
to geometric surgery and normal cobordism of the normal maps for which they were defined. The set of cobordism classes
of $n$-dimensional quadratic Poincar\'e complexes over a ring $\Lambda$ with involution are isomorphic to Wall's $L_n(\Lambda)$. 
The instant surgery obstruction provides
an easy formula to distill the traditional surgery obstruction form or formation out of a quadratic kernel.
(See \cite{Ran80a} Chapter 4 and \cite{Ran80b} Chapter 7 for details.) 

The quadratic kernel construction can be generalized to (odd- or even-dimensional) normal maps
which are not a homotopy equivalence on the boundary (nor on the whole manifold). The
result will be a {\bf quadratic Poincar\'e pair} of the same dimension.
Assume for example that $\partial l=l|\colon\partial M \ra \partial X$
is no longer necessarily a homotopy equivalence. Then there is chain homotopy 
commutative diagram
\begin{eqnarray*}
\xymatrix{
C(\widetilde{\partial X}) \ar[d]_{\partial l^!}  \ar[r]_{i_X} &
C(\widetilde{X}) \ar[d]^{l^!}\\
C(\widetilde{\partial M}) \ar[r]_{i_M} &
C(\widetilde{M})
}
\end{eqnarray*}
with $i_X$ and $i_M$ the inclusions of the boundary of $W$ and $X$.
It induces a map of the mapping cones
$$
    f=i_l \colon C=\cone(\partial l^!) \ra D=\cone(l^!)
$$
The quadratic kernel of $l$ is the
$(2q+1)$-dimensional quadratic Poincar\'e pair $c=(f\colon C \ra D, (\delta\psi, \psi)\in W_\%(f)_{2q+1})$.
The quadratic structure contains again self-intersection information and the maps
$$
    \mat{\delta\phi_0, f\phi_0 } \colon \cone(f)^{2q+1-*} \ra D
$$
induce the Poincar\'e-Lefschetz duality maps $K^{2q+1-*}(M,\partial M)\isora K_*(M)$.
Its boundary $(C,\psi)$ is a $2q$-dimensional quadratic Poincar\'e complex and it is by construction
the quadratic kernel of the normal map $l|\colon \partial M \ra \partial X$. (See also \cite{Ran80b} Proposition 6.5.)  

By Proposition \ref{thomthickprop}, there is a one-to-one
correspondence between homotopy classes of quadratic Poincar\'e pairs
and quadratic complexes (the latter are not necessarily
Poincar\'e). It is induced by the {\bf Thom construction} which
assigns to every quadratic Poincar\'e pair 
$c=(f\colon C \ra D, (\delta\psi, \psi))$ a quadratic complex
$(N=\cone(f),\zeta=\delta\psi/\psi)$ of the same dimension.  The
homology of $N$ are the relative kernel modules $K_*(M,\partial M)$. 
The chain map $(1+T)\zeta_0\colon N^{2q+1-*} \homra N$ induces the
maps $$K^{2q+1-*}(M,\partial M) \ra K^{2q+1-*}(M) \isora K_*(M,\partial M)$$
If $l|\colon \partial M \homra \partial X$ is a homotopy equivalence as 
in the beginning, then $C\simeq 0$, $(N,\zeta)\simeq (D,\nu)$ 
and the Poincar\'e pair $c$ is homotopy equivalent to $(0\ra D, (0,\nu))$.

In the same fashion we can translate the $(2q+2)$-dimensional normal map 
$k\colon (W, \partial W)\linebreak \ra (X\times I, X\cup_{\partial X} X)$
into a $(2q+2)$-dimensional quadratic Poincar\'e pair 
$x=(g \colon \partial E \ra E, (\delta\omega, \omega))$. The boundary $\partial W$
is the union of $M$ and $M'$ glued together along their common boundary. Similarly, the 
quadratic kernel $(\partial E, \omega)$ of the normal map $k|\partial W$ 
is the algebraic union of the quadratic kernels of $l$ and $l'$ in the sense
of Definition \ref{uniondef}.

Maps/diffeomorphisms of manifolds that are compatible with normal maps on them give rise to morphisms/isomorphisms
between the quadratic kernels (like \eg $\partial M \isora \partial M'$).

All in all, the constructions above yield a translation of normal maps (on manifolds) to quadratic complexes
as the table below illustrates:

\medskip
\begin{tabular}{|c|c|}
\hline
Topology & Algebraic surgery theory\\
\hline\hline
$\partial M$                    & $2q$-dim. quad. Poincar\'e complex $(C,\psi)$\\
\hline
$\partial M \hookrightarrow M$  & $(2q+1)$-dim. quad. Poincar\'e pair $(f\colon C \ra D, (\delta\psi,\psi))$\\
\hline
$M/\partial M$                  & $(2q+1)$-dim. quad. complex $(N=\cone(f),\zeta)$\\
\hline
$\partial M'$           & $2q$-dim. quad Poincar\'e complex $(C',\psi')$\\
\hline
$\partial M' \hookrightarrow M'$& $(2q+1)$-dim. quad. Poincar\'e pair $(f'\colon C' \ra D', (\delta\psi',\psi'))$\\
\hline
$M'/\partial M'$            & $(2q+1)$-dim. quad. complex $(N'=\cone(f'),\zeta')$\\
\hline
$\partial M \stackrel{\cong}{\ra} \partial M'$ & equivalence $(h,\chi) \colon (C,\psi) \stackrel{\simeq}{\ra} (C', \psi')$\\
\hline
$M \cup_{\partial M} M'$    & $(2q+1)$-dim. quad. Poincar\'e complex $(\partial E, \omega)$\\
\hline
$M \cup_{\partial M} M' \hookrightarrow W$ & $(2q+2)$-dim. quad. Poincar\'e pair $(g \colon \partial E \ra E, (\delta\omega, \omega))$\\
\hline
\end{tabular}
\medskip

Unfortunately, a generalization of this procedure to all preformations
will not work for two simple reasons: firstly,
there is no generalization of quadratic kernels to normal
smoothings, secondly, there is no geometric realization result known
for general preformations.

There is however a purely algebraic translation method which enables
us to construct the quadratic pairs and complexes by just using the
data given by the formation. It turns out that this method extends
without a problem to $\epsilon$-quadratic (split) preformations as
long as they are regular \ie all modules in it are \fg free.

One has to be cautious. For arbitrary preformations and arbitrary Kreck surgery problems
the relationship between geometry and algebra is not as straightforward as
for presentations.
It can happen \eg that there is a non-contractible algebraic boundary $(C,\psi)$ although $M$ is closed.
Nevertheless the philosophy remains the same. We can think of the complexes and pairs as vague algebraic models
of the manifolds or normal smoothings as in the table above - but only to boost our intuition!
If we want to prove statements about those quadratic complexes and pairs we will not be able
to use geometry but we have to resort to the methods of algebraic surgery theory alone.

One example for this strategy is Theorem \ref{algmainthm}.  It states
that a preformation is (stably) elementary if and only if the
quadratic Poincar\'e pair $x$ associated to it is
cobordant rel$\partial$ to an algebraic $h$-cobordism.
It is an almost word by word translation of the proof of 
Kreck's Theorem \ref{mainthm}. Nevertheless, it is not an automatic consequence
because there is no mathematically rigid connection between the algebraic model $x$
and the original (geometric) Kreck surgery problem.

By \cite{Ran01a} Proposition 9.4 there is a one-to-one correspondence between certain equivalence classes
of non-singular formations and {\bf short odd complexes}. A similar result can be found in \cite{Ran80a} Proposition 2.3 and 2.5.
We will not need those theorems or a generalization in detail.
We just use it as a motivation for translating a regular 
$\epsilon$-quadratic split preformation
$\kreckfs$ into a connected $(2q+1)$-dimensional quadratic chain complex $(N,\zeta)$
\begin{eqnarray}
d_N     &=& \mu^*\colon N_{q+1} = F \ra N_q=G^* \label{Ndefeqn}\\
\zeta_0 &=& \gamma \colon N^q=G \ra N_{q+1}=F \nonumber\\
\zeta_1 &=& \epsilon\theta \colon N^q=G \ra N_q=G^*\nonumber
\end{eqnarray}
with $\theta$ a representative of $\bar\theta\in Q_{-\epsilon}(G)$.
Obviously $(N,\zeta)$ depends on the choice of $\theta$. We will deal with this issue
the end of this section (see Remark \ref{welldef1rem}).

An obstruction preformation is an algebraic model for a  $(2q+2)$-dimensional cobordism $(W,M,M') \ra X$
of highly-connected normal smoothings. But the results from algebraic surgery theory
that we are using were proven with odd-dimensional traditional surgery theory in
mind. As explained in Section \ref{formgeosec}, in that context an obstruction formation 
for a $(2q+2)$-dimensional normal cobordism $(W,M,M') \ra X$ is thought of as an obstruction for
$M \ra X$ only. Hence $(N,\zeta)$ is a quadratic complex model for the normal map
$(M,\partial M) \ra (X,\partial X)$.

Now we turn around the cobordism to derive the quadratic chain complex
given by the normal map $(M',\partial M') \ra (X,\partial X)$. 
The obstruction of the ``new'' cobordism can easily be constructed out of $z$
and is called the {\bf flip of $z$}. Again we take a relationship which holds in the world
of formations and presentation and generalize it to all preformations:
\begin{defi}
\label{flipdefi}\index{Flip}
Let $z=\kreckfs$ be an $\epsilon$-quadratic split preformation. 
Its {\bf Flip} is the preformation
$$z'=\kreckfb{F^*}{\epsilon\mu}{G}{\gamma}{F}{-\bar\theta}$$
Similar for the non-split case.
\hfill\qed\end{defi}
As in (\ref{Ndefeqn}) we use the flip preformation to define a 
connected $(N',\zeta')$ be the $(2q+1)$-dimensional complex
\begin{eqnarray}
\label{Npdefeqn}
d_{N'}&=& \gamma^*\colon N'_{q+1}=F^* \ra N'_q=G^*\\
\zeta'_0 &=& \epsilon\mu \colon {N'}^q \ra N'_{q+1} \nonumber\\
\zeta'_1 &=& -\epsilon\theta \colon {N'}^q \ra N'_q \nonumber
\end{eqnarray}

The next step is to thicken up $(N,\eta)$ and $(N',\eta')$ to $(2q+1)$-dimensional quadratic
Poincar\'e pairs 
\begin{eqnarray}
c  &=& (f\colon C=\partial N \ra D=N^{2q+1-*}, (\delta\psi=0,\psi=\partial\zeta))\label{defceqn}\\
c' &=& (f\colon C'=\partial N' \ra D'={N'}^{2q+1-*}, (\delta\psi'=0,\psi'=\partial\zeta')) \label{defc2eqn}
\end{eqnarray}
(using the constructions in Definition \ref{thickthomdef}).

In the geometric situation we obviously find that $\partial M\cong\partial M'$. 
We expect a chain complex analogue and indeed
\begin{eqnarray}
\label{hdefeqn}
\xymatrix@+20pt
{
C_{q+1}=G   \ar[d]_{\svec{-\epsilon\gamma\\-\epsilon\mu}}   \ar[r]_{h_{q+1}=1_G}&
C'_{q+1}=G  \ar[d]^{\svec{-\mu\\-\epsilon\gamma}}
\\
C_q=F\oplus F^* \ar[d]_{\svec{\mu^* & \epsilon\gamma^*}}    \ar[r]_{h_q=\svec{0&\epsilon\\1&0}}&
C'_q=F^*\oplus F\ar[d]^{\svec{\gamma^* & \mu^*}}
\\
C_{q-1}=G^*                         \ar[r]_{h_{q-1}=1_{G^*}}&
C'_{q-1}=G^*
}
\end{eqnarray}
\begin{eqnarray*}
\chi_1 &=& \mat{0&-\epsilon\\0&0} \colon {C'}^q=F\oplus F^* \ra C'_q=F^*\oplus F\\
\chi_2 &=& \mat{-\mu\\0} \colon {C'}^{q-1}=G \ra C'_{q}=F^*\oplus F\\
\chi_3 &=& \theta \colon {C'}^{q-1}=G \ra C'_{q-1}=G^*
\end{eqnarray*}
defines an isomorphism $(h,\chi)\colon (C,\psi) \stackrel{\cong}{\ra} (C',\psi')$.
We glue $c$ and $c'$ together along $(h,\chi)$ \ie by
Definition \ref{uniondef} and Lemma \ref{deltapairlem} we compute the union
\begin{eqnarray*}
(\partial E, \omega)=(f'h \colon C \ra D', ((-)^{2q} f'\chi {f'}^*=0, \psi)) \cup
(f\colon C \ra D, (0, -\psi))
\end{eqnarray*}

We will try to simplify the quadratic Poincar\'e complex $(\partial E, \omega)$ and it is already clear what the result will be
if we look at the special case of formations and presentations. In this case
the long exact sequence of $(W,\partial W)$ shows that $\gamma^*\mu\colon G=K_{q+1}(W) \ra G^*=K_{q+1}(W,\partial W)$
is a chain complex model for $K_{q+1}(\partial W)$. So it makes sense to expect 
that the $(2q+1)$-dimensional quadratic Poincar\'e complex $(A,\tau)$ 
(arising from the regular $\epsilon$-quadratic split preformation $\partial(G,\theta)$ by the same process as in (\ref{Ndefeqn}))
\begin{eqnarray}
\label{Adefeqn}
d_A&=& ((1+T_{-\epsilon})\theta)^* \colon A_{q+1}=G\ra A_q=G^*\\
\tau_0&=&1 \colon A^q=G \ra A_{q+1}=G \nonumber\\
\tau_1&=&\epsilon\theta \colon A^q=G \ra A_q=G^* \nonumber
\end{eqnarray}
will be isomorphic to $(\partial E,\omega)$.
But if $(\partial E,\omega)$ looks so simple, why did we go through all the complicated procedures of thickening and glueing
in the first place? Well, the aim of Kreck's surgery theory is to decide whether $(W,M,M')$ is 
cobordant to an $h$-cobordism \ie whether the inclusions of $M$ and $M'$ into some $W'$ cobordant to $W$
are homotopy equivalences.
In our algebraic model, we will have to check whether the chain maps of $D$ and $D'$ into some 
algebraic cobordism are chain equivalences. 
Hence, we have to keep track where exactly $D$ and $D'$ are hidden in the boundary $\partial E$.

There is an equivalence $(a,\kappa) \colon (\partial E,\omega) \ra (A,\tau)$
given by
\begin{eqnarray}
\label{kdefeqn}
\xymatrix@C+45pt
{
\partial E_{q+1}\ar[d]\ar[r]_-{a_{q+1}=\svec{-1&0&0&1}}&
A_{q+1}=G   \ar[d]^{(\gamma^*\mu)^*}
\\
\partial E_q    \ar[r]_-{a_q=\svec{\epsilon\mu^*&-1& \gamma^*}}&
A_q=G^*
}
\end{eqnarray}
\begin{eqnarray*}
\kappa_2&=&\epsilon\theta \colon A^q=G \ra A_q=G^*
\end{eqnarray*}

Every boundary of a form can easily be expressed as a Poincar\'e pair as the following lemma suggests:
\begin{lem}
\label{simplepairlem}
Let $(G,\theta)$ be an $(-\epsilon)$-quadratic form. Then the $(2q+1)$-dimensional
quadratic Poincar\'e complex $(A,\tau)$ defined in (\ref{Adefeqn})
is the boundary of the
$(2q+2)$-dimensional quadratic Poincar\'e pair $y=(p\colon A \ra E, (\delta\tau=0,\tau))$ given by
$p=1\colon A_{q+1}=G \ra E_{q+1}=G$.
\end{lem}

Using Lemma \ref{deltapairlem} we find that
$$x=(g \colon \partial E \ra E, (\delta\omega=0, \omega))$$
with $g=\mat{1&0&0&-1} \colon \partial E_{q+1} \ra E_{q+1}=G$
is a $(2q+2)$-dimensional quadratic Poincar\'e pair. 

\begin{rem}
\label{welldef1rem}
It is time to investigate the effect of a choice of representative $\theta$ for $\bar\theta \in Q_{-\epsilon}(G)$
on the construction of $x$. Let $\widehat\theta=\theta+\widetilde\theta+\epsilon\widetilde\theta^*$ be another representative.
Let $(C, \widehat\psi)$ and $(C', \widehat\psi')$ be the $2q$-dimensional quadratic Poincar\'e complex
given by (\ref{defceqn}) and (\ref{defc2eqn}) using the representative $\widehat\theta$.

Then $\widehat\psi-\psi=d(\tilde\psi)$ and $\widehat\psi'-\psi'=d(-\tilde\psi)$
with $\tilde\psi\in W_\%(C)_{2q+1}$ given by $\tilde\psi_3=-\epsilon\tilde\theta^*\colon C^{q-1}\ra C_{q-1}$.

Another choice of representative for $\theta$ does not affect $[\psi]\in Q_{2q}(C)$,
$[(0,\psi)]\in Q_{2q+1}(f)$, $[(0,\psi')]\in Q_{2q+1}(f'h)$, $[\omega]\in Q_{2q+1}(\partial E)$
and $[(0,\omega)]\in Q_{2q+2}(g)$.
\end{rem}

\section{Algebraic Surgery and Cobordisms of Pairs}
\label{cobpairsec}

In Kreck's surgery theory we look at a cobordism $(e,f,f')\colon (W, M,M') \ra B$ of normal smoothings
and wonder whether it is cobordant rel$\partial$ to an $h$-cobordism or equivalently
whether surgery inside of $W$ will produce an $h$-cobordism. This is the case if and only
if the obstruction in $l'_{2q+2}(\Lambda)$ is elementary.
This section introduces algebraic versions of cobordisms and surgery for Poincar\'e pairs.

The first step will be a purely technical namely to define the notion of algebraic cobordism
of quadratic pairs rel$\partial$ and algebraic surgery inside a quadratic pair.

In geometry a cobordism rel$\partial$ between two manifolds $M$ and $M'$ with the same boundary $N$
is often thought of as a manifold with corners $(W, \partial W=M \cup N\times I \cup M')$.

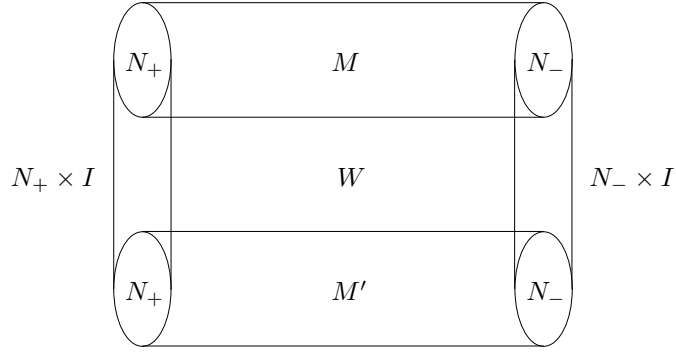
\begin{center}
\setlength{\unitlength}{0.006in}
\begin{picture}(750,300)(0,0)
\put(200,250){\ellipse{50}{100}}
\put(550,250){\ellipse{50}{100}}
\path(200,300)(550,300)
\path(200,200)(550,200)
\put(185,240){$N_+$}
\put(365,240){$M$}
\put(535,240){${N_-}$}

\path(175,250)(175,50)
\path(225,250)(225,50)
\path(525,250)(525,50)
\path(575,250)(575,50)

\put(85,140){$N_+\times I$}
\put(370,140){$W$}
\put(590,140){$N_-\times I$}

\put(200,50){\ellipse{50}{100}}
\put(550,50){\ellipse{50}{100}}
\path(200,0)(550,0)
\path(200,100)(550,100)
\put(185,40){$N_+$}
\put(365,40){$M'$}
\put(535,40){$N_-$}
\end{picture}
\end{center}

By ``collapsing'' the $N\times I$-part of $\partial W$ and glueing together the boundaries of $M$ and $M'$
we produce a new manifold $W'$ with the boundary $\partial W'= M\cup_N M'$. 
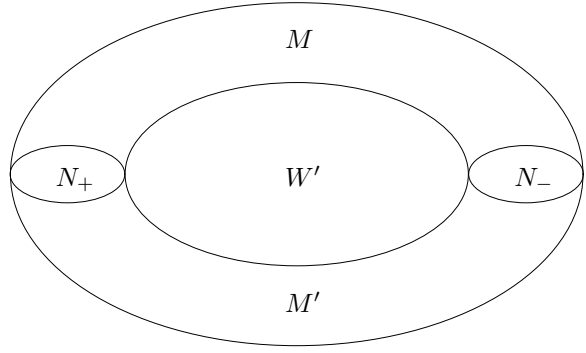
\begin{center}
\setlength{\unitlength}{0.006in}
\begin{picture}(750,300)(0,0)
\put(380,150){\ellipse{500}{300}}
\put(380,150){\ellipse{300}{160}}
\put(180,150){\ellipse{100}{50}}
\put(580,150){\ellipse{100}{50}}

\put(370,260){$M$}
\put(370,140){$W'$}
\put(370, 30){$M'$}

\put(170,140){$N_+$}
\put(570,140){$N_-$}
\end{picture}
\end{center}
Differential topology
shows us that the existence of a null-cobordism of $M\cup_N M'$ is in fact equivalent to $(M,N)$ and $(M',N)$
being cobordant rel$\partial$. We will use this picture in order to define algebraic cobordisms rel$\partial$ of quadratic pairs.
\begin{defi}
\label{cobrelbddef}
Two $(n+1)$-dimensional $\epsilon$-quadratic Poincar\'e pairs $c=(f\colon C \rightarrow D, (\delta\psi,\psi))$ 
and $c'=(f'\colon C \rightarrow D', (\delta\psi',\psi))$ are {\bf cobordant rel$\partial$}\index{Cobordism rel$\partial$} if there is an
$(n+2)$-dimensional $\epsilon$-quadratic Poincar\'e pair $(h\colon D'\cup_C D \ra E, (\delta\omega, \delta\psi'\cup_\psi -\delta\psi))$.
\hfill\qed\end{defi}

An easy example for such cobordisms are homotopy equivalences.
\begin{lem}
\label{cobhomlem}
Let $c=(f\colon C \ra D, (\delta\psi,\psi))$ 
and $c'=(f'\colon C \ra D', (\delta\psi',\psi))$ be $(n+1)$-dimensional $\epsilon$-quadratic Poincar\'e pairs.
Let $(1,h;k) \colon c \ra c'$ be a homotopy equivalence. Then $c$ and $c'$ are cobordant rel$\partial$.
\end{lem}
\begin{proof}
There is a $(\delta\chi,\chi)\in W_\%(f',\epsilon)_{n+2}$ such that 
$$(1,h;k)_\%(\delta\psi,\psi)-(\delta\psi',\psi')=d(\delta\chi,\chi).$$
Define the $(n+2)$-dimensional $\epsilon$-quadratic Poincar\'e pair 
\begin{eqnarray}
\label{bdefeqn}
(b\colon D\cup_C D' \ra D', ((-)^n\delta\chi, \delta\psi \cup_\psi -\delta\psi'))
\end{eqnarray}
by $b=(h, (-)^{r-1}k, -1) \colon (D\cup_C D')_r=D_r\oplus C_{r-1}\oplus D'_{r-1} \ra D'_r$.
\end{proof}

The next lemma proves the useful fact that changing the common boundary of two Poincar\'e pairs $c$ and $c'$
doesn't change anything about their cobordism relationship.
\begin{lem}
\label{cobbdrylem}
Let $c =(f  \colon C \ra D,  (\delta\psi,  \psi))$ and $c'=(f' \colon C \ra D', (\delta\psi', \psi))$
be two $\epsilon$-quadratic $(n+1)$-dimensional Poincar\'e pairs.
Let $(h,\chi) \colon (\widehat{C}, \widehat{\psi}) \stackrel{\simeq}{\ra} (C, \psi)$ be an equivalence.
Define the $(n+1)$-dimensional $\epsilon$-quadratic Poincar\'e pairs (using Lemma \ref{deltapairlem})
\begin{eqnarray*}
\widehat{c}&=& (\widehat{f} = fh \colon \widehat{C}\ra D, (\widehat{\delta\psi}=\delta\psi+(-)^n f\chi f^*, \widehat{\psi}))\\
\widehat{c}'&=& (\widehat{f}' = f'h \colon \widehat{C}\ra D', (\widehat{\delta\psi}'=\delta\psi'+(-)^n f'\chi {f'}^*, \widehat{\psi}))
\end{eqnarray*}
Then $c$ and $c'$ are cobordant rel$\partial$ if and only if $\widehat{c}$ and ${\widehat{c}}'$ are.
\end{lem}

\begin{proof}
If $c$ is cobordant rel$\partial$ to $c'$ then there is an $(n+2)$-dimensional $\epsilon$-quadratic Poincar\'e pair
$$(e\colon D \cup_C D \ra E, (\delta\omega, \omega=\delta\psi \cup_\psi -\delta\psi'))$$
By Lemma \ref{unionbdrylem} there is an equivalence 
$$(a,\kappa) \colon \widehat{c}\cup -\widehat{c}'=(D \cup_{\widehat{C}} D',  \widehat{\delta\psi}\cup_{\widehat{\psi}} -\widehat{\delta\psi}')
\stackrel{\simeq}{\ra} c\cup -c = (D\cup_C D', \delta\psi \cup_\psi -\delta\psi')$$
Hence, by Lemma \ref{deltapairlem} there is an $(n+2)$-dimensional $\epsilon$-quadratic Poincar\'e pair
$$(ea\colon D \cup_{\widehat{C}} D \ra E, (\delta\omega+(-)^{n+1}e\kappa e^*, \widehat{\delta\psi} \cup_{\widehat{\psi}} -\widehat{\delta\psi'}))$$
\end{proof}

It is a well-known fact that two manifolds are cobordant if and only if one manifold is derived
from the other by a finite sequence of surgeries and diffeomorphisms. There is an algebraic equivalent
for Poincar\'e complexes (Proposition \ref{cobordprop}). We will establish the same relationship
in the case of Poincar\'e pairs. First we need to define a surgery on the inside of a pair:
 
\begin{defi}
\label{surgpairdef}\index{Surgery on Poincar\'e pairs}\index{Quadratic pair!surgery}
Let $c=(f\colon C \ra D, (\delta\psi,\psi))$ be an $(n+1)$-dimensional $\epsilon$-quadratic Poincar\'e pair
and $d=(g\colon \cone(f) \ra B, (\delta\sigma, \delta\psi/\psi))$ a connected $(n+2)$-dimensional $\epsilon$-quadratic
pair. Write $g=(a,b) \colon \cone(f)_r=D_r\oplus C_{r-1} \ra B_r$.
{\bf The result of the surgery $d$ on the inside of $c$} is the $(n+1)$-dimensional $\epsilon$-quadratic Poincar\'e pair
$c'=(f'\colon C \ra D', (\delta\psi',\psi))$ given by
\begin{eqnarray*}
d_{D'}&=&\mat{  d_D & 0 & (-)^n(1+T_\epsilon)\delta\psi_0 a^* + (-)^{n} f(1+T_\epsilon)\psi_0 b^*\\
        (-)^r a & d_B   & (-)^r(1+T_\epsilon)\delta\sigma_0     + (-)^{n+1} b \psi_0 b^*\\
        0   & 0 & (-)^r d_B^*} \colon
\\
&&D'_r=D_r \oplus B_{r+1} \oplus B^{n+2-r} \ra D'_{r-1}=D_{r-1} \oplus B_{r} \oplus B^{n+3-r}
\\
f' &=&\mat{f\\-b\\0} \colon C_r \ra D'_r=D_r \oplus B_{r+1} \oplus B^{n+2-r}
\\
\delta\psi'_0&=&\mat{\delta\psi_0&0&0\\0&0&0\\0&1&0} \colon
\\
&&{D'}^{n+1-r}=D^{n+1-r}\oplus B^{n+2-r} \oplus B_{r+1}
    \ra D'_r=D_r \oplus B_{r+1} \oplus B^{n+2-r}
\\
\delta\psi'_s&=&\mat{\delta\psi_s   & (-)^s T_\epsilon\delta\psi_{s-1}a^*-f T_\epsilon\psi_{s-1}b^* &0\\
             0          & (-)^{n-r-s+1}T_\epsilon\delta\sigma_{s-1}         &0\\
             0          & 0                             &0}
\colon 
\\
&&{D'}^{n+1-r-s}=D^{n+1-r-s}\oplus B^{n+2-r-s} \oplus B_{r+s+1}
    \\&&\quad\ra D'_r=D_r \oplus B_{r+1} \oplus B^{n+2-r} 
\quad (s>0)
\end{eqnarray*}
\hfill\qed\end{defi}

The following proposition will justify the formulae above by showing that surgery inside of a pair is
nothing but the composition of the following standard procedures of algebraic surgery theory: 
Thom complex, algebraic surgery and thickening.
(The latter is the inverse operation to the Thom complex. See Proposition \ref{thomthickprop}.)

\begin{prop}
\label{surgpairlem}
We use the terminology of the previous definition.
\begin{enumerate}
\item If $C=0$ then $(D',\delta\psi')$ is the result of the surgery $(a \colon D \rightarrow B, (\delta\sigma, \delta\psi))$
as in Definition \ref{surgerydef}.
    
\item The result of the surgery $d=(g\colon \cone(f) \rightarrow B, (\delta\sigma, \delta\psi/\psi))$ on the 
Thom complex $(\cone(f), \delta\psi/\psi)$ of $c$ is isomorphic to the Thom complex $(\cone(f'), \delta\psi'/\psi)$ of $c'$.
\end{enumerate}
\end{prop}

\begin{proof}
The first part is trivial. So we turn our attention to the second claim.
The isomorphisms

\begin{eqnarray*}
u_r&=&\mat {1&0&0&0\\0&0&1&0\\0&0&0&1\\0&1&0&(-)^{n-r}\psi_0 b^*}\colon\\
&& M_r=(D_r \oplus C_{r-1})\oplus B_{r+1}\oplus B^{n+2-r}\\
&&\quad\ra \cone(f')_r=(D_r \oplus B_{r+1} \oplus B^{n+2-r})\oplus C_{r-1}
\end{eqnarray*}
define an isomorphism $(u,0)\colon  (M,\tau) \stackrel{\cong}{\ra} (\cone(f'),\delta\psi'/\psi)$ 
between the result $(M,\tau)$ of the surgery $d$ and the Thom-complex of $c'$. 
\end{proof}

At last we prove the expected relationship between cobordisms and surgery.
\begin{prop}
\label{surgcobprop}
Let $c=(f\colon C \ra D, (\delta\psi,\psi))$ and $c'=(f'\colon C \ra D', (\delta\psi',\psi))$
be $(n+1)$-dimensional $\epsilon$-quadratic Poincar\'e pairs.
They are cobordant rel$\partial$ if and only if one can be obtained from the other
by surgeries and homotopy equivalences of the type $(1,h;k)$.
\end{prop}

One direction of the proof is covered by the following lemma.
\begin{lem}
\label{cobsurlem}
Let $c=(f\colon C \ra D, (\delta\psi,\psi))$ be an $(n+1)$-dimensional $\epsilon$-quadratic Poincar\'e pair
and $d=(g=\mat{a&b}\colon \cone(f) \ra B, (\delta\sigma, \delta\psi/\psi))$ a
connected $(n+2)$-dimensional $\epsilon$-quadratic pair. 
Let $c'=(f'\colon C \ra D', (\delta\psi',\psi))$ be the result of the surgery $d$ on the inside of $c$.
\begin{enumerate}
\item Let\footnote{$Define -c=(f\colon C \ra D, (-\delta\psi,-\psi))$ for a pair $c=(f\colon C \ra D, (\delta\psi,\psi))$.}
$(D\cup_C D, \delta\psi\cup_\psi -\delta\psi)= c \cup -c$
be the union of $c$ with itself along its boundary $C$.
Then $$\tilde d = (\tilde g\colon  D\cup_C D\ra B, (\delta\sigma, \delta\psi\cup_\psi -\delta\psi))$$ given by
$\tilde g=\mat{a&b&0} \colon (D\cup_C D)_r=D_r\oplus C_{r-1} \oplus D_r \ra B_r$ is a connected 
$(n+2)$-dimensional $\epsilon$-quadratic complex. The result of the surgery $\tilde d$ is isomorphic to
$(D'\cup_C D, \delta\psi'\cup_\psi -\delta\psi)= c' \cup -c$.

\item $(h \colon D\cup_C D \ra D, (0, \delta\psi\cup_\psi -\delta\psi))$ with 
$h=\mat{1&0&-1}\colon D_r\oplus C_{r-1} \oplus D_r \ra D_r$ is an $(n+2)$-dimensional 
$\epsilon$-quadratic Poincar\'e pair. \label{cobsurit}

\item $(D'\cup_C D, \delta\psi'\cup_\psi -\delta\psi)= c' \cup -c$ is null-cobordant
\ie $c$ and $c'$ are cobordant rel$\partial$.
\end{enumerate}
\end{lem}

\begin{proof}
\begin{enumerate}
\item The philosophy of this proof is that in some sense we can transfer everything we did for the Thom complexes
in the proof of Lemma \ref{surgpairlem} to the union $D\cup_C D'$ using the morphism
\begin{eqnarray*}
\left( \mat{1&0&0\\0&1&0}, 0 \right) \colon (D\cup_C D, \delta\psi\cup_\psi -\delta\psi) \ra (\cone(f), \delta\psi/\psi)
\end{eqnarray*}
In particular we can apply it to Lemma \ref{deltapairlem} and show that $\tilde d$ 
is an $(n+2)$-dimensional $\epsilon$-quadratic pair.
The isomorphisms

\begin{eqnarray*}
u_r&=&\mat{1&0&0&0&0\\
     0&0&0&1&0\\
     0&0&0&0&1\\
     0&1&0&0&(-)^{n-r}\psi_0 b^*\\
     0&0&1&0&0} \colon
\\
&&\tilde V_r =  (D_r \oplus C_{r-1} \oplus D_{r}) \oplus B_{r+1}\oplus B^{n+2-r}
\\&&\quad\ra V'_r =  (D_r \oplus B_{r+1} \oplus B^{n+2-r}) \oplus C_{r-1} \oplus D_r
\end{eqnarray*}
define an isomorphism $(u,0) \colon (\tilde V, \tilde \sigma) \stackrel{\cong}{\ra} (V', \sigma')$
between the result $(\tilde V, \tilde \sigma)$ of the surgery $\tilde d$ and the union
$(V', \sigma')= (D'\cup_C D, \delta\psi'\cup_\psi -\delta\psi)= c' \cup -c$.

\item Exercise.

\item Follows from Propositions \ref{surgbordprop} and \ref{cobordprop}.

\end{enumerate}
\end{proof}

\begin{proof}[Proof of Proposition \ref{surgcobprop}]
It remains to show that cobordant pairs can be obtained from each other by
surgery and homotopy equivalences which leave the boundary untouched.

Let $c=(f\colon C \ra D, (\delta\psi,\psi))$ and $c'=(f'\colon C \ra D', (\delta\psi',\psi))$
be $(n+1)$-dimensional $\epsilon$-quadratic Poincar\'e pairs which are cobordant rel$\partial$
so that there exists an $(n+2)$-dimensional $\epsilon$-quadratic Poincar\'e pair
$$
    e=(h\colon D\cup_C D' \ra E, (\delta\omega, \omega=\delta\psi \cup_\psi -\delta\psi'))
$$
with $h=\mat{j_0&k&j_1} \colon D_r \oplus C_{r-1} \oplus D_r \ra E_r$.
We define the connected $(n+2)$-dimensional $\epsilon$-quadratic pair
\begin{eqnarray*}
d&=&(g \colon \cone(f) \ra B=\cone(j_1), (\delta\sigma,\sigma=\delta\psi/\psi))\\
g&=&\mat{j_0&k\\0&-f}\colon \cone(f)_r=D_r\oplus C_{r-1} \ra B_r =E_r\oplus D'_{r-1}\\
\sigma_s&=&\mat{\delta\omega_s & 0 \\ 
        (-)^{n-r-1}(\delta\psi'_s j_1^*  + (-)^{s}f'\psi_s k^*) & (-)^{n-r-s}T_\epsilon \delta\psi'_{s+1}}\colon
\\
&&B^{n+2-r-s}=E^{n+2-r-s}\oplus {D'}^{n+1-r-s} \ra B_r=E_r\oplus D'_{r-1}
\end{eqnarray*}

The result of the surgery $d$ inside of $c$ is the $(n+1)$-dimensional $\epsilon$-quadratic Poincar\'e pair
$c''=(f''\colon C \ra D'', (\delta\psi'',\psi))$. There is a homotopy equivalence $m\colon D'' \ra D'$ given by
\begin{eqnarray*}
m&=&\mat{0&0&1&0&0&\delta\psi'_0} \colon\\ 
&&D''_{r+1}=D_r\oplus(E_{r+1}\oplus D'_r)\oplus (E^{n+2-r}\oplus {D'}^{n+1-r}) \ra D'_r
\end{eqnarray*}
such that $m\delta\psi''m^*=\delta\psi'$. Hence $(1,m;0)$ defines a homotopy equivalence from $c''$ to $c'$.
\end{proof}

\section{Elementariness in Algebraic Surgery Theory}
\label{elemalgsursec}

In this section we want to understand elementariness in the context of algebraic surgery theory
by reproving Theorem \ref{mainthm} for Poincar\'e pairs: a preformation is stably elementary if
and only if the Poincar\'e pair $x$ defined in Section \ref{constrpairsec} is cobordant rel$\partial$
to an {\bf algebraic $h$-cobordism}. This theorem is the key to applying algebraic surgery theory
to preformations in this treatise.

\begin{defi}\index{Algebraic $h$-cobordism}
Let $c=(f\colon C \ra D, (\delta\psi,\psi))$ and $c'=(f\colon C \ra D', (\delta\psi',\psi))$ be 
$\epsilon$-quadratic $n$-dimensional Poincar\'e pairs whose union is the boundary of an 
$(n+1)$-dimensional $\epsilon$-quadratic Poincar\'e pair 
$d=(g\colon D\cup_C D' \ra E, (\delta\omega, \delta\psi\cup_\psi -\delta\psi'))$.
Write $g=\mat{j_0&k&j_1} \colon (D\cup_C D')_r = D_r \oplus C_{r-1} \oplus D_r \ra E_r$.
$d$ is an {\bf algebraic $h$-cobordism} if the chain maps $j_0$ and $j_1$ are chain equivalences.
\hfill\qed\end{defi}

\begin{thm}
\label{algmainthm}
Let $z=\kreckfs$ be a regular $\epsilon$-quadratic split preformation.
Let $x=(g \colon \partial E \ra E, (\delta\omega=0,\omega))$ be the $(2q+2)$-dimensional quadratic
Poincar\'e pair constructed in Section \ref{constrpairsec} for an arbitrary representative
$\theta$ of $\bar\theta\in Q_{-\epsilon}(G)$.

\begin{enumerate}
\item If $z$ is elementary then $x$ is cobordant rel$\partial$ to an algebraic $h$-cobordism.

\item If $x$ is cobordant rel$\partial$ to an algebraic $h$-cobordism then $[z] \in l_{2q+2}(\Lambda)$ is elementary.
\end{enumerate}
\end{thm}

We need some technical results before we can move on to the proof of this theorem.

\begin{lem}
\label{alghcoblem1}
Let $c=(f\colon C \ra D, (\delta\psi,\psi))$ and $c'=(f\colon C \ra D', (\delta\psi',\psi))$ be $\epsilon$-quadratic
$n$-dimensional Poincar\'e pairs.
Let $d=(g\colon D\cup_C D' \ra E, (\delta\omega, \delta\psi\cup_\psi -\delta\psi'))$
and $d'=(g'\colon D\cup_C D' \ra E', (\delta\omega', \delta\psi\cup_\psi -\delta\psi'))$
be two $(n+1)$-dimensional $\epsilon$-quadratic Poincar\'e pairs. Assume
there is a homotopy equivalence\footnote{See Definition \ref{homeqpairdef}.} 
between them that is the identity on the boundary.
If $d$ is an algebraic $h$-cobordism then so is $d'$.
\end{lem}

\begin{lem}
Let $c=(f\colon C \ra D, (\delta\psi,\psi))$ and $c'=(f\colon C \ra D', (\delta\psi',\psi))$ be $\epsilon$-quadratic
$n$-dimensional Poincar\'e pairs.
Let $d=(g\colon D\cup_C D' \ra E, (\delta\omega, \delta\psi\cup_\psi -\delta\psi'))$
be an $(n+1)$-dimensional algebraic $h$-cobordism. Then $d$ is homotopy equivalent to
an algebraic $h$-cobordism $$d'=(g'\colon D\cup_C D' \ra D, (\delta\omega', \delta\psi\cup_\psi -\delta\psi'))$$
such that $g'=\mat{1&l&h}\colon (D\cup_C D')_r=D_r\oplus C_{r-1} \oplus D'_r \ra D_r$ with $h\colon D' \stackrel{\simeq}{\ra} D$
a chain equivalence.
\end{lem}
\begin{proof}
Write $g=\mat{j_0&k&j_1} \colon (D\cup_C D')_r = D_r \oplus C_{r-1} \oplus D_r \ra E_r$.
Let $i_0 \colon E \ra D$ the chain homotopy inverse of $j_0$ and let $\Delta\colon i_0 j_0 \simeq 1 \colon D \ra D $ 
be a chain homotopy. Then there is a homotopy equivalence $(1,i_0; \mat{-\Delta&0&0}) \colon d \ra d'$ 
with $l=i_0 k+(-)^{r-1}\Delta f$ and $h=i_0 j_1$.
\end{proof}

\begin{lem}
\label{exactlem}
Let $V_n \stackrel{d_n}{\ra} V_{n-1} \stackrel{d_{n-1}}{\ra} \cdots V_0 \ra 0$ be an exact
sequence of free \fg $\Lambda$-modules. Then $\im d_i\subset V_{i-1}$ and $\ker d_i\subset V_i$
are stably \fg free direct summands for all $i\in\{1,\dots,n\}$.
\end{lem}
\begin{proof}
For $i=1$ observe that $0 \ra \ker d_1 \ra V_1 \stackrel{d_1}{\ra} V_0 \ra 0$ is an exact sequence and
since $V_0$ is free it splits. Hence $\im d_1=V_0$ and $\ker d_1$ are stably \fg free and a direct summand in the
respective $V_i$. Now assume the claim is true for $i \in \{1,\dots, n-1\}$. Then we look at the
exact sequence $0 \ra \ker d_{i+1} \ra V_{i+1} \stackrel{d_{i+1}}{\ra} \ker d_i \ra 0$. Again the 
sequence is exact and splits because by assumption $\ker d_i$ is projective. Now the claim follows for $i+1$.
\end{proof}

\begin{proof}[Proof of Theorem \ref{algmainthm}]
\begin{enumerate}

\item
Let $\kreckfs$ be elementary and $i\colon U \hookrightarrow G$ be the inclusion
of an $h$-lagrangian.
We define the chain map $m\colon \cone(g) \ra B$ (with $\cone(g)_r=E_r\oplus \partial E_{r-1}$)
\begin{eqnarray*}
\xymatrix@+25pt
{
\cone(g)_{q+3}=0\oplus G 
\ar[d]_{-\epsilon\svec{1\\\gamma\\\mu\\1}}
\\
\cone(g)_{q+2}=0\oplus (G\oplus F \oplus F^* \oplus G)
\ar[d]_{\svec{  -\epsilon   & 0     & 0         & \epsilon\\
        -\epsilon\gamma & \epsilon  & 0         & 0\\                           
        0       & \mu^*     & \epsilon\gamma^*  & 0\\
        0       & 0     & \epsilon      & -\epsilon\mu}}
\\
\cone(g)_{q+1}=G\oplus (F\oplus G^* \oplus F^*)
\ar[r]_-{m=\svec{a&b}}
&
B_{q+1}=U^*
}
\end{eqnarray*}
with $a=-i^*\gamma^*\mu$ and $b=\mat{-\epsilon i^*\mu^* & i^* & -i^*\gamma^*}$.

Because of $i^*\theta i=0 \in Q_{-\epsilon}(U)$ there is a $\delta\chi\in \Hom_\Lambda(U, U^*)$
such that $i^*\theta i = \delta\chi +\epsilon\delta\chi^* \in \Hom_\Lambda(U, U^*)$.
We can check that $(\delta\sigma,\sigma=\delta\omega/\omega)\in W_\%(m)_{2q+3}$ with $\delta\sigma_1=\epsilon\delta\chi \colon B^{q+1} \ra B_{q+1}$
is a cycle. Hence we have a connected $(2q+3)$-dimensional quadratic pair 
$d=(m \colon \cone(g) \ra B, (\delta\sigma,\sigma))$. The result of the surgery $d$ on the inside
of $x$ is the $(2q+2)$-dimensional quadratic Poincar\'e pair $x'=(g' \colon \partial E \ra E', (\delta\omega',\omega))$
given by ($\partial E_r=D'_r\oplus C_{r-1} \oplus D_r$, $E'_r=E_r\oplus B_{r+1} \oplus B^{2q+3-r}$):
\begin{eqnarray*}
\xymatrix@+30pt
{
\partial E_{q+2}=0\oplus G \oplus 0
\ar[d]_{-\epsilon\svec{1\\\gamma\\\mu\\1}}
\ar[r]_-{0} 
&
E'_{q+2}=0\oplus 0\oplus U
\ar[d]^{-i}
\\
\partial E_{q+1}=G\oplus (F\oplus F^*)\oplus G  
\ar[d]_{\svec{  -\epsilon\gamma & \epsilon  & 0         & 0\\                           
        0       & \mu^*     & \epsilon\gamma^*  & 0\\
        0       & 0     & \epsilon      & -\epsilon\mu}}
\ar[r]_-{\svec{1&0&0-1}}
&
E'_{q+1}=G\oplus 0\oplus 0
\ar[d]^{\epsilon i^*\gamma^*\mu}
\\
\partial E_q=F\oplus G^* \oplus F^*
\ar[r]_-{\svec{\epsilon i^*\mu^* & -i^* & i^*\gamma^*}}
&
E'_q=0\oplus U^*\oplus 0
}
\end{eqnarray*}
\begin{eqnarray*}
\delta\omega'_0&=& 1_U \colon {E'}^{q}=U \ra {E'}_{q+2}=U
\end{eqnarray*}
Using the inclusion $\svec{1&0\\0&0\\0&1}\colon D'\oplus D \ra \partial E$ and with the help of 
Proposition \ref{elemprop}\ref{elemprop4} we see that $x'$ is an algebraic $h$-cobordism.

\item
Now we assume that $x$ is cobordant rel$\partial$ to an algebraic $h$-cobordism 
$x'=(g'\colon \partial E \ra E', (\delta\omega', \omega))$.
In order to simplify our calculations we remember that the boundary of $x$ and $x'$ 
can be reduced in size by using the equivalence
$(a,\kappa) \colon (\partial E,\omega) \ra (A,\tau)$
defined in (\ref{kdefeqn}) on page~\pageref{kdefeqn}.
Let $y'=(p'\colon A \ra E', (\delta\tau', \tau))$ be the
$(2q+2)$-dimensional quadratic Poincar\'e pair induced by $x'$ and the inverse of $(a,\kappa)$.
Let $y=(p\colon A \ra E, (\delta\tau=0,\tau))$  be the $(2q+2)$-dimensional
quadratic Poincar\'e pair from Lemma \ref{simplepairlem}. (It is also induced by $x$ and the inverse of $(a,\kappa)$).

By Lemma \ref{cobbdrylem}, $y$ is cobordant rel$\partial$ to $y'$ and
by (the proof of) Proposition \ref{surgcobprop} and Lemma \ref{alghcoblem1} we can assume that $y'$ is the result of 
a surgery $d=(m\colon \cone(p) \ra B, (\delta\sigma, \sigma=\partial\tau/\tau))$ inside of $y$ with
\begin{eqnarray*}
\xymatrix@C+60pt
{
\cone(p)_{q+2}=G        \ar[d]^-{\svec{-\epsilon\\-\epsilon\gamma^*\mu}} \ar[r]^-{m_{q+2}=b_{q+2}}
&
B_{q+2}             \ar[d]^d
\\
\cone(p)_{q+1}=G\oplus G^*  \ar[r]^-{m_{q+1}=\svec{a_{q+1}&b_{q+1}}}
&
B_{q+1}             
}
\end{eqnarray*}
$$
\sigma_0=\mat{0&0\\0&\epsilon\theta^*} \colon \cone(p)^{q+1}=E^{q+1}\oplus A^q \ra \cone(p)_{q+1}=E_{q+1}\oplus A_q
$$
Our next step will be the analysis of the complex $E'$.
If $r \ge q+3$ or $r \le q$ the differential is given by
\begin{eqnarray*}
    d'_r=\mat{d&(-)^r(1+T)\delta\sigma_0\\0&(-)^rd^*} &\colon&
    E'_r=B_{r+1}\oplus B^{2q+3-r}\\ 
    &&\quad\ra E'_{r-1}=B_{r}\oplus B^{2q+4-r}
\end{eqnarray*}

The top differentials are dual to the bottom ones, \ie 
\begin{eqnarray*}
\mat{0&(-)^r\\1&0}d_r^*\mat{0&(-)^{r-1}\\1&0} = d_{2q+3-r}
\end{eqnarray*}
for $r \ge q+3$ and $r \le q$.

Because $E' \simeq D$ the homology groups $H_r(E')$ vanish for $r\neq q+1,q$. Hence by Lemma \ref{exactlem}
there is a \fg free submodule $X\subset E'_q$ such that $\ker d'_q\oplus X=E'_q$.
Therefore $E'_{q+2}/\ker d'_{q+2}=\coker d'_{q+3} = U$ is stably \fg free
and $U^*=\ker d'_q=\im d'_{q+1}$. This observation gives us the chance to cut away the higher and lower parts of $E'$
and establish a chain equivalence $E' \stackrel{\simeq}{\ra} E''$ using the maps
\begin{eqnarray}
\label{simplifyEeqn}
\xymatrix@C+30pt
{
E'_{q+2}=B_{q+3}\oplus B^{q+1}      \ar[r]_-{\text{proj}_X}     
                    \ar[d]_{d'_{q+2}}
&
E''_{q+2}=U             \ar[d]^{i}
\\
E'_{q+1}=G\oplus B_{q+2}\oplus B^{q+2}  \ar[d]_{d'_{q+1}}
                    \ar[r]_-{\svec{ 1   & 0 & 0\\
                                b_{q+2} & 1 & 0\\
                            0   & 0 & 1}}
&
E''_{q+1}=G\oplus B_{q+2}\oplus B^{q+2} \ar[d]^{p}
\\
E'_q= B_{q+1}\oplus B^{q+3}     \ar[r]_-{\left[\svec{0&-1\\-\epsilon&0}\right]}
&
E''_q=U^*
}
\end{eqnarray}

with
\begin{eqnarray*}
d'_{q+2}&=&\mat{0 & -b_{q+1}^*\\
        d & \epsilon(1+T)\delta\sigma_0 + b_{q+2} b_{q+1}^*\\
        0 & \epsilon d^*}\\
d'_{q+1}&=&\mat{-\epsilon a_{q+1} & d & -\epsilon(1+T)\delta\sigma_0\\
        0   & 0 & -\epsilon d^*}\\
i&=&\left[\mat{0 & -b_{q+1}^*\\
        d & \epsilon(1+T)\delta\sigma_0\\
        0 & \epsilon d^*}\right]
\\
p&=&\left[\mat{0 & 0 & \epsilon d^*\\
        -b_{q+1}\gamma^*\mu & -\epsilon d & (1+T)\delta\sigma_0}\right]
\end{eqnarray*}

Let's define a regular $\epsilon$-quadratic split preformation
\begin{eqnarray*}
z'&=&\kreckf{F'}{\gamma'}{G'}{\mu'}{\bar\theta'}\\
&=&\kreckfs \oplus \partial\left(B_{q+2}\oplus B^{q+2}, \mat{0&\epsilon\\0&0}\right)
\end{eqnarray*}
which is another representative of $[z]\in\el$. One can easily compute that $p= i^*{\gamma'}^*\mu'$.

Now we have a look at the boundary of $y'$. The map $g'\colon A \ra E'$ is given by
\begin{eqnarray*}
\xymatrix@C+30pt
{
A_{q+1}=G   \ar[d]_{(\gamma^*\mu)^*}
        \ar[r]^-{\svec{1\\-b_{q+2}\\0}}
&
E'_{q+1}=G\oplus B_{q+2}\oplus B^{q+2}  \ar[d]^{d'_{q+1}}
\\
A_{q}=G^*   \ar[r]_-{\svec{-b_{q+1}\\0}}
&
E'_q= B_{q+1}\oplus B^{q+3}
}
\end{eqnarray*}
Applying the chain equivalence (\ref{simplifyEeqn}) and the map $D'\oplus D \ra \partial E \stackrel{a}{\ra} A$
to $g'\colon A \ra E'$ we find two chain maps $D' \ra E'' \longleftarrow D$
\begin{eqnarray*}
\xymatrix@C+40pt
{
&
U       \ar[d]^{i}
\\
G       \ar[r]^{\svec{-1\\0\\0}} 
        \ar[d]_{-\epsilon\gamma}
&
G'      \ar[d]^{i^*{\gamma'}^*\mu'}
&
G       \ar[l]_{\svec{1\\0\\0}}
        \ar[d]^{-\epsilon\mu}
\\
F       \ar[r]_{\svec{0\\b_{q+1}\mu^*}}
&
U^*
&
F^*     \ar[l]^{\svec{0\\\epsilon b_{q+1}\gamma^*}}
}
\end{eqnarray*}
which by assumption are chain equivalences.
From the fact that $(\delta\sigma, \sigma)\in Q_{2q+2}(m)$ one can deduce that
$i^*\bar\theta i=0 \in Q_{-\epsilon}(U)$.
Now it is not very difficult to verify that the preformation $z'$
fulfils the assumption of Proposition \ref{elemprop}~\ref{elemprop4} in respect to the stably \fg free submodule $U$. 
Further stabilization of $z'$ by boundaries of hyperbolic forms helps to replace $U$ by a \fg free submodule. 
Hence $[z]=[z']\in l_{2q+2}(\Lambda)$ is elementary.
\end{enumerate}
\end{proof}

\chapter{Flip-Isomorphisms}
\label{flipchap}

{\bf For the whole chapter let $q\ge 2$, $\epsilon=(-)^q$ and let $\Lambda$ be a weakly finite ring
with $1$ and involution.}

Obviously, a $(2q+2)$-dimensional normal cobordism $(W,M,M') \ra X$
only stands a chance to be cobordant rel$\partial$ to an $s$-cobordism if
there is a compatible diffeomorphism between $M$ and $M'$. Some
kind of ``algebraic isomorphism'' between $M$ and $M'$ can be produced
by just using preformations. 
Let $z=\kreckfs$ be the obstruction preformation. We can interpret $z$ also as
an algebraic model for the normal map $M\ra X$ and the flip $z'=\kreckfb{F^*}{\epsilon\mu}{G}{\gamma}{F}{-\bar\theta}$
of $z$ as  a model for $M'\ra X$. 
Following that philosophy, we hope that $z$ and $z'$ are 
weakly isomorphic if $z$ is elementary. We shall call such an isomorphism {\bf flip-isomorphism}.
In {\bf Section \ref{flipisosec}} we motivate and define flip-isomorphisms
and show that, indeed, any elementary preformation - even those that aren't obstructions of
the above or any surgery problem - has at least one.

Just like we translated preformations into quadratic pairs and complexes in Section \ref{constrpairsec},
we translate flip-isomorphisms into isomorphisms of those quadratic complexes in {\bf Section  \ref{fliptranslsec}}.
Those isomorphisms can be applied to the quadratic Poincar\'e pair $x$ from Section \ref{constrpairsec}.
They transform $x$ into a Poincar\'e pair with an algebraic {\bf twisted double} on the boundary
(see {\bf Sections \ref{flipqutwdblsec} and \ref{flipsymtwdblsec}}).
This is a necessary preparation before we can define 
asymmetric and quadratic signatures in the following chapters.

\section{Flip-Isomorphism}
\label{flipisosec}
\index{Flip-isomorphism!philosophy}
The idea behind flip-isomorphisms is inspired by an observation in geometry:
a cobordism $(W,M,M')\ra X$ of highly-connected normal maps/normal smoothings can only be cobordant to an 
$s$-cobordism if $M$ and $M'$ are diffeomorphic. This hardly seems to be a very revealing 
insight. After all it is the aim of any surgery theory to establish existence or 
non-existence of such a diffeomorphisms.

But in the context of preformations we can produce a notion of some kind of ``algebraic isomorphism''
of $M$ and $M'$: {\bf flip-isomorphisms}.
First of all we remember from Section \ref{formgeosec} 
that there are actually two ways of looking at preformations.
We can think of them as algebraic vehicles for surgery-relevant data of
an even-dimensional cobordism $(W,M,M') \ra X\times (I,0,1)$ of normal smoothings/maps. In that case we identify
preformations by the very rigid equivalence relation of (stable) strong isomorphisms 
in order to preserve the essential information of the whole cobordism. This is the view of $l$-theory.
But we have also learnt that in traditional odd-dimensional surgery theory, formations encode 
the information of the normal map $M \ra X$ only. That is why the odd-dimensional obstruction groups
$L_{2q+1}(\Lambda)$ have a much more flexible equivalence relation which includes 
the use of weak isomorphisms (Definition \ref{weakisodef}).

So, philosophically, if we have a cobordism $(W,M,M') \ra X$ of normal smoothings/maps
and define its obstruction $z$ and we think of it as a description of the whole
cobordism we use strong isomorphisms. If we treat $z$ only as a description
of the map $M\ra X$ we use weak isomorphisms.
As alluded to in Section \ref{constrpairsec}, turning around the new cobordism leads to
the flip of $z$ as the obstruction for the cobordism $(W,M',M) \ra X$. The heuristics
so far suggests that the preformation-equivalent of
a diffeomorphism between $M\ra X$ and $M'\ra X$ is a weak isomorphism between $z$ and its flip. Such an isomorphism
will be called a {\bf flip-isomorphism}.

For fans of algebraic surgery theory the importance of flip-isomorphisms as an obstruction
to elementariness is even more evident:  in Section \ref{constrpairsec},
a regular $\epsilon$-quadratic split preformation $z$ defined a $(2q+2)$-dimensional quadratic Poincar\'e pair 
$x=(g \colon \partial E \ra E, (\delta\omega=0, \omega))$.
By Theorem \ref{algmainthm} we know that if $z$ is elementary, $x$ is cobordant rel$\partial$
to an $h$-cobordism. Then the $(2q+1)$-dimensional Poincar\'e pairs 
$c=(f\colon C \ra D, (0,\psi))$ and $c''=(f'h \colon C \ra D', ((-)^{2q} f\chi f^*=0, \psi))$,
whose union constitutes the boundary of $x$ (that is $(\partial E,\omega)=c''\cup-c$), have to be homotopy equivalent.
Both pairs are thickening-ups of the quadratic complexes $(N,\zeta)$ and $(N',\zeta')$
from (\ref{Ndefeqn}) and (\ref{Npdefeqn}) which themselves stem from the preformation
$z$ and its flip $z'$. Hence a homotopy equivalence of the pairs $c$ and $c''$ will 
eventually lead to some kind of equivalence between $z$ and $z'$. \cite{Ran80a} Proposition 2.3
and 2.5, which provided a recipe for translating the preformations $z$, $z'$ to quadratic complexes 
$(N,\zeta)$, $(N',\zeta')$ respectively,
states that there is a natural bijection between equivalences of those quadratic complexes 
and (stable) weak isomorphism classes of formations. A generalization of those Propositions
shows that $z$ and $z'$ are (stably) weakly isomorphic. The existence of a flip-isomorphism
for elementary preformations can also be shown quite easily without algebraic surgery theory 
(see the proof of the Proposition \ref{elemflipprp} below).

\begin{defi}
\label{flipisodef}
\begin{enumerate}
\item\index{Flip-isomorphism}
A {\bf flip-isomorphism of a regular $\epsilon$-quadratic  preformations $z=\skreckfs$}
is a weak isomorphism of $z$ with its flip \ie a triple $(\alpha,\beta,\sigma)$ consisting of
isomorphisms $\alpha\in\Hom_\Lambda(F,F^*)$ 
and $\beta\in\Hom_\Lambda(G,G)$ and an element $\sigma\in Q^{-\epsilon}(F^*)$ such that
\begin{enumerate}
\item $\alpha\gamma +\alpha\sigma\mu=\epsilon\mu\beta \in \Hom_\Lambda(G, F^*)$
\item $\alpha^{-*}\mu = \gamma\beta \in \Hom_\Lambda(G, F)$
\end{enumerate}
It is {\bf strong} if $\sigma=0$.\index{Flip-isomorphism!strong}

\item
\label{isosplititm}\index{Flip-isomorphism}
A {\bf flip-isomorphism of a regular $\epsilon$-quadratic split preformation $z=\kreckfs$}
is a 
weak isomorphism of $z$ with its flip \ie
a triple $(\alpha,\beta,\bar\nu)$ consisting of isomorphisms $\alpha\in\Hom_\Lambda(F,F^*)$ 
and $\beta\in\Hom_\Lambda(G,G)$ and an element
$\bar\nu \in Q_{-\epsilon}(F^*)$ 
such that
\begin{enumerate}
\item $\alpha\gamma +\alpha(\bar\nu - \epsilon\bar\nu^*)^*\mu=\epsilon\mu\beta\in \Hom_\Lambda(G, F^*)$ 
\item $\alpha^{-*}\mu = \gamma\beta\in \Hom_\Lambda(G, F)$
\item $\bar\theta + \mu^*\bar\nu\mu + \beta^*\bar\theta\beta=0 \in Q_{-\epsilon}(G)$
\end{enumerate}
It is {\bf strong} if $\bar\nu=0$.\index{Flip-isomorphism!strong}

\item\index{Flip-isomorphism!stable}
A {\bf stable flip-isomorphism of a regular $\epsilon$-quadratic (split) preformation $z$}
is a flip-isomorphism of $z+h$ with $h$ a boundary of a hyperbolic form. (Note that $[z]=[z+h]$ in
the $l$-monoid).\hfill\qed
\end{enumerate}
\end{defi}

\begin{prop}
\label{elemflipprp}
Let $z$ be a regular $\epsilon$-quadratic (split) preformation.
If $z$ is elementary then $z$ has a flip-isomorphism.
\end{prop}
\begin{proof}
Let $z=\kreckfs$ be an elementary regular $\epsilon$-quadratic split preformation.
We assume that our preformation has the form described in Proposition \ref{elemprop} \ref{elemprop3}.
There is a flip isomorphism $(\alpha,\beta,\bar\nu)$ of $z$ given by
\begin{eqnarray*}
\alpha  &=& \mat{0&-1\\-\epsilon&0}\colon F=U\oplus U^* \ra F^*=U^*\oplus U\\
\beta   &=& \mat{-1&-\tau\\0&1} \colon G= U \oplus R \ra G=U\oplus R\\
\bar\nu &=& 0\colon F^* \ra F
\end{eqnarray*}

In the non-split case the flip-isomorphism is $(\alpha,\beta,0)$.
\end{proof}

\begin{cor}
\label{elemflipcor}
Let $z$ be a regular $\epsilon$-quadratic split preformation.  If $z$
is (stably) elementary there is a (stable) strong flip-isomorphism
$(\alpha,\beta,0)$ such that $\alpha\colon F\ra F^*$ is
$\epsilon$-symmetric and zero in $L^{2q}(\Lambda)$ (and hence also in
$LAsy^0(\Lambda)$) and $\beta^2=1_G$.  Similar for the non-split case.
\end{cor}

\begin{cor}
Let $z$ be a regular $\epsilon$-quadratic split preformation and $z'$ its flip.
If $z$ is stably elementary then it has a stable flip-isomorphism
and $[z]=[z']\in\el$.
\end{cor}

\section{Translating Flip-Isomorphisms into Algebraic Surgery Theory}
\label{fliptranslsec}

Let $z=\kreckfs$ be a regular $\epsilon$-quadratic split preformation.
Let $t=(\alpha,\beta,\bar\nu)$ be a flip-isomorphism of $z$.

In the following we fix representatives $\theta$ and $\nu$ for $\bar\theta$ and $\bar\nu$
and let $\kappa\in\Hom(G,G^*)$ such that
$\beta^*\theta\beta+\theta+\mu^*\nu\mu=\kappa+\epsilon\kappa^* \in\Hom_\Lambda(G, G^*)$.
We must of course be aware that once we leave the realm of chain complexes and 
define concepts about preformations and flip-isomorphisms
(that is asymmetric signatures,
flip-isomorphism rel$\partial$ and quadratic signatures) we have to check
to what extent they depend on the choice of representatives. 
(Compare Remark \ref{welldef2rem})

In Section \ref{constrpairsec} we translated $z$ into a Poincar\'e pair $x$. The first step
of this construction was to use the proofs of \cite{Ran80a} Proposition 
2.3. and 2.5 to create $(2q+1)$-dimensional quadratic 
complexes $(N,\zeta)$ and $(N',\zeta')$ out of $z$ and its flip $z'$.
Those proofs also suggest that an isomorphism of two preformations gives rise
to an isomorphism of those quadratic complexes. Both propositions just cover formations but without
problems we can generalize the construction for all regular preformations.
Hence $t$ induces an isomorphism 
$(e_t,\rho_t)\colon (N,\zeta) \ra (N',\zeta')$ of $(2q+1)$-dimensional quadratic complexes given by
\begin{eqnarray}
\xymatrix@C+20pt
{
N_{q+1} = F \ar[r]_-{e_{t,q+1}=\alpha}  \ar[d]_{\mu^*}&
N'_{q+1}=F^*                    \ar[d]_{\gamma^*}
\\
N_{q} = G^* \ar[r]_-{e_{t,q}=\beta^{-*}}&
N'_{q}=G^*
}
\end{eqnarray}
\begin{eqnarray*}
\zeta_0 &=& \gamma \colon N^q \ra N_{q+1}\\
\zeta_1 &=& \epsilon\theta \colon N^q \ra N_q
\\
\zeta'_0 &=& \epsilon\mu \colon {N'}^q \ra N'_{q+1}\\
\zeta'_1 &=& -\epsilon\theta \colon {N'}^q \ra N'_q
\\
\rho_{t,0} &=& \alpha\nu\alpha^* \colon {N'}^{q+1} \ra N'_{q+1}\\
\rho_{t,1} &=& {\gamma}^*\alpha\nu\alpha^* \colon {N'}^{q+1} \ra N'_{q}\\
\rho_{t,2} &=& -\beta^{-*}\kappa^*\beta^{-1} \colon {N'}^{q} \ra N'_{q}
\end{eqnarray*}

The Poincar\'e pairs $c$ and $c'$ defined in (\ref{defceqn}) and (\ref{defc2eqn}) 
are thickenings of $(N,\zeta)$ and $(N',\zeta')$ and the isomorphism $(e_t, \rho_t)$ leads to
a homotopy equivalence of those two pairs in particular
to an isomorphism of the $2q$-dimensional quadratic Poincar\'e complexes 
$(\partial e_t,\partial\rho_t)\colon (C,\psi)=\partial(N,\zeta)$
and $(C',\psi')=\partial(N',\zeta')$ (see Lemma \ref{bdrmaplem}). 
Composing this isomorphism with the inverse of the canonical isomorphism
$$(h,\chi)\colon (C,\psi) \isora (C',\psi')$$ defined in (\ref{hdefeqn}) yields
a self-equivalence $(h_t,\chi_t)=(h,\chi)^{-1}(\partial e_t,\partial\rho_t)$ of $(C,\psi)$ given by
\begin{eqnarray}
\label{htdefeqn}
\xymatrix@+30pt
{
C_{q+1}=G \ar[r]_{\beta} \ar[d]_{\svec{-\epsilon\gamma\\-\epsilon\mu}} & 
C_{q+1}=G \ar[d]^{\svec{-\epsilon\gamma\\-\epsilon\mu}}\\
C_q    =F\oplus F^* \ar[r]_{\svec{0 & \alpha^{-*}\\\epsilon\alpha & \epsilon\alpha(\nu^*-\epsilon\nu)}} \ar[d]_{\svec{\mu^* & \epsilon\gamma^*}} & 
C_q    =F\oplus F^* \ar[d]^{\svec{\mu^* & \epsilon\gamma^*}}\\
C_{q-1}=G^* \ar[r]_{\beta^{-*}} & 
C_{q-1}=G^*
}
\end{eqnarray}
\begin{eqnarray*}
\chi_{t,1}&=& \mat{0&0\\1&-\alpha\nu^*\alpha^*}\colon C^q=F^*\oplus F \ra C_q=F\oplus F^*\\
\chi_{t,2}&=& \mat{0\\\epsilon\mu-\alpha\nu^*\alpha^*\gamma} \colon C^{q-1}=G \ra C_q=F\oplus F^*\\
\chi_{t,3}&=& -\theta+\beta^{-*}\kappa\beta^{-1} \colon C^{q-1}=G \ra C_{q-1}=G^*
\end{eqnarray*}
(see Definition \ref{compinvdef} for inverses and compositions of isomorphisms of quadratic complexes).

\begin{rem}
\label{welldef2rem}
Let $\widehat\nu$,$\widehat\theta$ be other representatives for $\nu\in Q_{-\epsilon}(z,\beta)$ and $\theta\in Q_{-\epsilon}(G)$.
Then there are $\tilde\theta\in\Hom_\Lambda(G,G^*)$, $\tilde\nu\in\Hom_\Lambda(F^*, F)$ such that
$\widehat\nu-\nu=\tilde\nu+\epsilon\tilde\nu^*\in\Hom_\Lambda(F^*, F)$ and 
$\widehat\theta-\theta =\tilde\theta+\epsilon\tilde\theta^*\in\Hom_\Lambda(G,G^*)$.
Define $\widehat\kappa=\kappa+\tilde\kappa-\epsilon\tilde\kappa^* + \beta^*\tilde\theta\beta + \tilde\theta + \mu^*\tilde\nu\mu \in\Hom_\Lambda(G,G^*)$
for some $\tilde\kappa \in \Hom_\Lambda(G,G^*)$.

$\widehat\theta$, $\widehat\theta$ and $\widehat\kappa$ induce an isomorphism $(h_t, \widehat\chi_t) \colon (C,\widehat\psi) \isora (C,\widehat\psi)$.
(Note that $h_t$ and the chain complex $C$ are not affected by the choice of representative.)
From Remark \ref{welldef1rem} we know that there is a $\tilde\psi\in W_\%(C)_{2q+1}$ such that
$\widehat\psi-\psi=d(\tilde\psi)$. Then $\widehat\chi_t-\chi_t=\tilde\psi-h_t\tilde\psi h_t^* +d\tilde\chi$
with $\tilde\chi\in W_\%(C)_{2q+2}$ given by
\begin{eqnarray*}
\tilde\chi_2 &=& \mat{0&0\\0&-\alpha\tilde\nu^*\alpha^*}\colon C^q=F^*\oplus F \ra C_q=F\oplus F^*\\
\tilde\chi_3 &=& -\tilde\chi_2^* d^* \colon C^{q-1}=G \ra C_q=F\oplus F^*\\
\tilde\chi_4 &=& \beta^{-*}(\tilde\theta+\tilde\kappa)\beta^{-1} \colon C^{q-1}=G \ra C_{q-1}=G^*
\end{eqnarray*}
See also Remarks \ref{welldef1rem}, \ref{welldef3rem} and Lemma \ref{quadamblem}.
\end{rem}

\section{Quadratic Twisted Doubles}
\label{flipqutwdblsec}

A preformation with a flip-isomorphism is some kind of algebraic equivalent
to a $(2q+2)$-dimensional normal cobordism $(W,M,M')\ra X$ with a diffeomorphism $h\colon M\isora M'$
compatible with the highly-connected normal maps $M\ra X$ and $M'\ra X$.
The boundary of $W\ra X$ can be replaced by a normal map on the 
{\bf twisted double $M\cup_{h|\partial M} M$}\index{Twisted double!geometric} 
into $X$. As sketched in the introduction this is the starting point
for constructions which yield obstructions for $(W,M,M')$ to be an $h$-cobordism:
asymmetric signatures and, if $h|\partial M \simeq 1_{\partial M}$, quadratic signatures.

These constructions can be imitated for quadratic Poincar\'e pairs 
as we will find out in Chapters \ref{asymchap} and \ref{quadchap}.
A prerequisite for those constructions is to turn the boundary of the
$(2q+2)$-dimensional $\epsilon$-quadratic Poincar\'e pair $x$ from 
Section \ref{constrpairsec} into an algebraic twisted double.

\begin{defi}[\cite{Ran98} 30.8(ii)]
\label{qutwdbldef}\index{Twisted double!quadratic}
The {\bf twisted double} of an $n$-dimensional $\epsilon$-quadratic Poincar\'e pair
$c=(f \colon C \ra D, (\delta\psi,\psi))$ over $\Lambda$ with respect to a 
self-equivalence $(h,\chi)\colon (C,\psi) \homra (C,\psi)$ is the 
$n$-dimensional $\epsilon$-quadratic Poincar\'e complex over $\Lambda$
\begin{eqnarray*}
c\cup_{(h,\chi)} -c &=& (D \cup_h D, \delta\psi \cup_\chi -\delta\psi)\\
            &=& (fh \colon C \ra D, (\delta\psi+(-)^{n-1}f\chi f^*, \psi)) \cup
            (f  \colon C \ra D, (-\delta\psi, -\psi))
\end{eqnarray*}
with
\begin{eqnarray*}
D \cup_h D &=& \cone\left(\mat{fh\\f}\colon C \ra D\oplus D\right),\\
(\delta\psi \cup_\chi \delta\psi)_s &=& 
\mat{
\delta\psi_s + (-)^{n-1}f\chi_s f^*     &0      &0\\
0                   &-\delta\psi_s  &(-)^{s-1} f\psi_s\\
(-)^{n-r}\psi_s h^*f^*          &0      &(-)^{n-r+s+1}T_\epsilon\psi_{s-1}}\colon
\\
&&(D\cup_h D)^{n-r-s} = D^{n-r-s}\oplus D^{n-r-s} \oplus C^{n-r-s-1}\\
&&\quad\ra(D\cup_h D)_r=D_r\oplus D_r \oplus C_{r-1}
\end{eqnarray*}
\hfill\qed\end{defi}

Let $z=\kreckfs$ be a regular $\epsilon$-quadratic split preformation.
Let $t=(\alpha,\beta,\bar\nu)$ be a flip-isomorphism of $z$. Again we pick
representatives $\theta$ and $\nu$ for $\bar\theta$ and $\bar\nu$
and let $\kappa\in\Hom(G,G^*)$ such that
$\beta^*\theta\beta+\theta+\mu^*\nu\mu=\kappa+\epsilon\kappa^* \in\Hom_\Lambda(G, G^*)$.

Let $(\partial E_t, \omega_t)$ be the
twisted double of $c=(f \colon C \ra D, (\partial\psi, \psi))$ of (\ref{defceqn})
with respect to the self-equivalence $(h_t,\chi_t)\colon (C,\psi) \isora (C,\psi)$ from (\ref{htdefeqn}).

There is an equivalence $(a_t,\sigma_t) \colon (\partial E_t,-\omega_t) \isora (\partial E, \omega)$
of $(2q+1)$-dimensional quadratic Poincar\'e complexes given by

\begin{eqnarray*}
a_{t,q+2}&=&\beta 
 \colon \partial E_{t,q+2}=G \ra \partial E_{q+2}=G
\\
a_{t,q+1}&=&\mat{   0&\beta &0      &0\\
            0&0 &0      &\alpha^{-*}\\
            0&0 &\epsilon\alpha &\epsilon\alpha(\nu^*-\epsilon\nu)\\
            1&0 &0      &0}\colon 
\\&&\partial E_{t,q+1}=G\oplus G\oplus (F\oplus F^*) \ra \partial E_{q+1}=G\oplus (F\oplus F^*)\oplus G
\\
a_{t,q}&=&\mat{ 0&\alpha^{-*}   &0\\
        0&0     &\beta^{-*}\\
        1&0     &0}
 \colon \partial E_{t,q}=F^*\oplus F^*\oplus G^* \ra \partial E_{q}=F\oplus G^*\oplus F^*
\\
\sigma_{t,0}&=&\mat{    0&0&0           &0\\
            0&0&0           &0\\
            0&0&\alpha\nu\alpha^*   &0\\
            0&0&0           &0}\colon 
\\&&\partial E^{q+1}=G^*\oplus (F^*\oplus F)\oplus G^* \ra \partial E_{q+1}=G\oplus (F\oplus F^*)\oplus G 
\\
\sigma_{t,0}&=&\mat{0\\1\\0}
 \colon \partial E^{q+2}=G^* \ra \partial E_{q}=F\oplus G^*\oplus F^*
\\
\sigma_{t,1}&=&\mat{    0  & 0      & \epsilon                          &0\\
            -1 & -\gamma^*  & -\mu^*+\epsilon\gamma^*\alpha\nu\alpha^*  &-1\\
            0  & 1      & \epsilon\alpha\nu\alpha^*                 &0}
 \colon \partial E^{q+1}=G^*\oplus (F^*\oplus F)\oplus G^* 
\\&&\hfill\ra \partial E_{q}=F\oplus G^*\oplus F^*
\\
\sigma_{t,2}&=&\mat{    0       & 0                     & 0\\
            -\gamma^*   & -\epsilon\beta^{-*}\kappa\beta^{-1}       & 0\\
            -1      & -\epsilon\mu-\epsilon\alpha\nu\alpha^*\gamma  &-\epsilon\alpha\nu\alpha^*}
\colon 
\\&&\partial E^{q}=F^*\oplus G\oplus F \ra \partial E_{q}=F\oplus G^*\oplus F^*
\end{eqnarray*}

Applying Lemma \ref{deltapairlem} to the Poincar\'e pair  $x=(g \colon \partial E \ra E, (\delta\omega=0, \omega))$
from Section \ref{constrpairsec} yields a $(2q+2)$-dimensional quadratic Poincar\'e pair 
given by
\begin{eqnarray*}
x_t     &=& (g_t=g a_t \colon \partial E_t \ra E, (\delta\omega_t, \omega_t))\\
g_{t,q+1}   &=& \mat{1&-\beta&0&0}\colon \partial E_{t,q+1}=G\oplus G\oplus F\oplus F^* \ra E_{q+1}=G\\
\delta\omega_t  &=& -(\delta\omega+(-)^{2q+1}g\sigma g^*)=0
\end{eqnarray*}
and a homotopy equivalence of pairs $(a_t, -1;0) \colon x_t \ra -x$. 

The equivalence $a_t$ maps each copy of $D$ in $\partial E_t$ isomorphically
onto a copy of $D$ and $D'$ in $\partial E$. Hence $x$ is cobordant rel$\partial$ to
an $h$-cobordism if and only if $x_t$ is.

\section{Symmetric Twisted Doubles}
\label{flipsymtwdblsec}

The computation of asymmetric signature of a flip-isomorphism (which we will
accomplish in Chapter \ref{asymchap}) demands a symmetric version of
the constructions of the previous section. 

\begin{defi}[\cite{Ran98} 30.8(ii)]
\label{twdbldef}\index{Twisted double!symmetric}
The {\bf twisted double} of an $n$-dimensional $\epsilon$-symmetric Poincar\'e pair
$c=(f \colon C \ra D, (\delta\phi,\phi))$ over $\Lambda$ with respect to a 
self-equivalence $(h,\chi)\colon (C,\phi) \homra (C,\phi)$ is the 
$n$-dimensional $\epsilon$-symmetric Poincar\'e complex over $\Lambda$
\begin{eqnarray*}
c\cup_{(h,\chi)} -c &=& (D \cup_h D, \delta\phi \cup_\chi -\delta\phi)\\
            &=& (fh \colon C \ra D, (\delta\phi+(-)^{n-1}f\chi f^*, \phi)) \cup
            (f  \colon C \ra D, (-\delta\phi, -\phi))
\end{eqnarray*}
with
\begin{eqnarray*}
D \cup_h D &=& \cone\left(\mat{fh\\f}\colon C \ra D\oplus D\right),\\
(\delta\phi \cup_\chi -\delta\phi)_s &=& 
\mat{
\delta\phi_s + (-)^{n-1}f\chi_s f^*     &0      &0\\
0                   &-\delta\phi_s  &(-)^{s-1} f\phi_s\\
(-)^{n-r}\phi_s h^*f^*          &0      &(-)^{n-r+s}T_\epsilon\phi_{s-1}}\colon
\\
&&(D\cup_h D)^{n-r+s} = D^{n-r+s}\oplus D^{n-r+s} \oplus C^{n-r+s-1}\\
&&\quad\ra(D\cup_h D)_r=D_r\oplus D_r \oplus C_{r-1}
\end{eqnarray*}
\hfill\qed\end{defi}

Let $z=\kreckfs$ be a regular $\epsilon$-quadratic split preformation.
Let $t=(\alpha,\beta,\bar\nu)$ be a flip-isomorphism of $z$
Let $\nu$, $\theta$ and $\kappa$ chosen as in the previous section.
We symmetrize our ingredients $(h_t,\chi_t)$ and $c$ from the previous section:

The $(2q+1)$-dimensional symmetric Poincar\'e pair
which is the product of the symmetrization of $c$ defined in (\ref{defceqn}) is given by
\begin{eqnarray}
\label{defsymceqn}
(1+T)c&=&(f\colon C \ra D, (\delta\phi=0,\phi=(1+T)\psi))
\end{eqnarray}

The symmetrization of $(h_t,\chi_t)$ defined in (\ref{htdefeqn}) leads to a
self-equivalence $$(h_t,0)\colon (C,\phi) \isora (C,\phi)$$ 
of the $2q$-dimensional symmetric Poincar\'e complex $(C,\phi)$.
Let $(\partial E_t, \theta_t)$ be the symmetric twisted double of $(1+T)c$ in respect
to $(h_t,0)$.

The twisted double construction and symmetrization are commutative operations up to
an equivalence.
\begin{lem}
Let $c=(f \colon C \ra D, (\delta\psi,\psi))$ be an $n$-dimensional quadratic Poin\-car\'e pair and
let $(h,\chi)\colon (C,\psi) \ra (C,\psi)$ be a self-equivalence. Then there is a
chain equivalence $$(1,\sigma)\colon (1+T)c \cup_{(h,(1+T)\chi)} -(1+T)c \isora (1+T)(c\cup_{(h,\chi)}-c)$$
with $(1+T)c = (f\colon C\ra D, (1+T)(\delta\psi,\psi))$ and
\begin{eqnarray*}
\sigma_0&=&\mat{0&0&0\\0&0&0\\0&0&(-)^{r-1}T\psi_0}\colon\\ 
    &&(D\cup_h D)^{n+1-r}=D^{n+1-r}\oplus D^{n+1-r}\oplus C^{n-r} \\
    &&\quad\ra (D\cup_h D)_r=D_r\oplus D_r \oplus C_{r-1}
\end{eqnarray*}

\end{lem}

We apply this lemma and Lemma \ref{deltapairlem} to the symmetrization
of the $(2q+2)$-dimensional quadratic Poincar\'e pair $x_t$ defined in
the previous section and get a $(2q+2)$-dimensional symmetric
Poincar\'e pair 
$$x^t=(g_t\colon \partial E_t \ra E,(\delta\theta_t=0, \theta_t)$$ 
If $x$ is cobordant rel$\partial$ to an $h$-cobordism then so is $x^t$.

\begin{rem}
The construction of $x^t$ and $(h_t,0)\colon (C,\psi) \ra (C,\psi)$ just depends on the ``symmetrization''
of $z$ and $t$ \ie $\skreckfs$ and $(\alpha,\beta, (\nu^*-\epsilon\nu)^*)$ (compare Remark \ref{weakisorem}).
Every choice of representative for $\bar\nu$ leads to the same Poincar\'e pair $x^t$.

For an $\epsilon$-quadratic preformation $y=\skreckfs$ and a flip-isomorphism $s=(\alpha,\beta,\sigma)$,
$x^s$ can be constructed in the same way if, in the definition of $\theta_t$ and $h_t$,
$(\nu^*-\epsilon\nu)^*$ is replaced by $\sigma$. 
\end{rem}

\chapter{Asymmetric Signatures of Flip-Isomorphisms}
\label{asymchap}

{\bf For the whole chapter let $\epsilon=(-)^q$ and let $\Lambda$ be a weakly finite ring
with $1$ and involution.}

Let $(W,M,M)$ be a cobordism with boundary $M\cup_h M$ for some automorphism $h$ of $\partial M$.
Then an asymmetric signature can be defined which vanishes in the asymmetric Witt-group
if and only if $W$ is cobordant rel$\partial$ to a twisted double. 
An $s$-cobordism $(W,M,M)$ is a twisted double, so the asymmetric signature provides an obstruction
for $W$ to be cobordant rel$\partial$ to an $s$-cobordism.

The previous chapter explained how a flip-isomorphism of a preformation 
defines a symmetric Poincar\'e pair with an algebraic twisted double as the boundary.
For such a Poincar\'e pair we can define an asymmetric signature, too. Like in geometry it vanishes
if the pair is cobordant rel$\partial$ to an $h$-cobordism \ie if the preformation is stably elementary.

Asymmetric signatures of manifolds do not require normal maps or smoothings (see \eg \cite{Ran98} Corollary 30.12.).
Similarly symmetric complexes and pairs suffice as input for the algebraic asymmetric signature,
\eg the symmetric Poincar\'e pair $x^t$ defined in Section \ref{flipsymtwdblsec}. 

In {\bf Section \ref{asymgeosec}} we will give a short introduction into the origins of
asymmetric forms and twisted doubles in geometry.
Then {\bf Section \ref{asymalgsec}}
will present the algebraic chain complex analogues. We define asymmetric forms,
complexes and pairs and show how a Poincar\'e pair with a twisted double on the boundary
determines an asymmetric complex.
We have already seen  in Section \ref{flipsymtwdblsec} how a preformation and a flip-isomorphism
can be used to construct a symmetric Poincar\'e pair with an algebraic twisted double as the boundary. In {\bf Section \ref{asymsignsec}}
we compute its asymmetric signature and show in {\bf Section \ref{asymsignelemsec}} that it vanishes
for every flip-isomorphism of an elementary preformation.

The asymmetric signatures depend on the choice of flip-isomorphism. Hence in {\bf Section \ref{asymflmonoidelemsec}} 
we define an $l$-monoid of flip-isomorphism $fl^{2q}(\Lambda)$. The asymmetric signatures define
a monoid homomorphism $fl^{2q}(\Lambda) \ra LAsy^0(\Lambda)$. It turns out that
the asymmetric signatures vanish for any flip-isomorphism of any stably elementary preformation.

\section{Geometric Twisted Doubles and Asymmetric Forms}
\label{asymgeosec}

One of the successes of asymmetric forms in surgery theory was the open book obstruction theory of Quinn 
(\cite{Quin79}) which led to a second computation of the cobordism ring of automorphism. (The first
computation was achieved by Kreck (see \eg \cite{Kreck84a}) using an early version of
his $l$-surgery theory.)

We will outline Quinn's approach to the cobordism of automorphism. For simplicity
we assume that all manifolds are simply-connected.
Let $h\colon N \stackrel{\cong}{\ra} N$ be a diffeomorphism
of a closed $2q$-dimensional manifold ($q>1$). The first obvious obstruction for
$h$ to be null-cobordant is the cobordism class of the mapping torus 
$T(h)$.

So let us assume that $T(h)$ has a $(2q+2)$-dimensional null-cobordism $V$
and after surgery on $V$ we can make $H_i(V,N)$ vanish for $i<q$.
The isomorphism
\begin{eqnarray}
\label{Quinnform}
    K=H_{q+1}(V,N) &\ra& H^{q+1}(V, \partial V - N)\\\nonumber
    &\ra& H^{q+1}(V, N \times I) \ra H^{q+1}(V, N)
\end{eqnarray}
gives us a non-singular asymmetric form $\lambda \colon K \ra K^*$ which is zero
in the asymmetric Witt group $LAsy^0(\bZ)$ if and only if there
is a diffeomorphism $H\colon Q \ra Q$ with $V=T(H)$ and hence $\partial H=h$.

More generally, an exact sequence describes the connection between the cobordism ring
of automorphisms and the asymmetric Witt group. Let $\Omega_i(X)$ be the cobordism
group of continuous maps from $i$-dimensional manifolds to $X$ and $\Delta_i(X)$ the
group of cobordism classes of triples $(F,g,h)$ with $F$ a closed $i$-dimensional manifold,
$g\colon F \ra X$ a map, $h\colon F\ra F$ an automorphism together with a homotopy
$g \simeq gh$ such that there is an induced map 
$$T(g)\colon T(h) \ra T(1\colon X \ra X)=X\times S^1$$
Then for any $k>2$ and topological space $X$ there is an exact sequence
\begin{eqnarray*}
0 
&\ra& \Delta_{2k+1}(X) 
\stackrel{T}{\ra} \Omega_{2k+2}(X\times S^1)
\ra LAsy^0(\bZ[\pi_1(X)])\\
&\ra& \Delta_{2k}(X)
\stackrel{T}{\ra} \Omega_{2k+1}(X \times S^1)
\ra 0
\end{eqnarray*}
with $T\colon \Delta_i(X) \ra \Omega_{i+1}(X \times S^1), (F,g,h) \mt (T(g), T(h))$
(see also \cite{Ran98} 30.6 (iv) or \cite{Quin79}).

In our case we start with an $(n+1)$-dimensional cobordism $(W,M,M')$ such that $M'\cong M$ and $M$ may have a boundary. 
Then there is  an isomorphism $h\colon \partial M \ra \partial M$ such that the boundary of
$W$ is the {\bf twisted double}\index{Twisted double!geometric} $M \cup_h M$ (see \cite{Win73} and \cite{Ran98} Chapter 30).
By glueing the ends of the cobordism together we obtain a manifold $V$ with boundary $T(h)$.
One can do surgery below the middle dimension to make $(V,N)$ highly-connected and read off an asymmetric
form as before. It vanishes in the asymmetric Witt group if and only if $V$ is cobordant to a mapping torus of an automorphism
and that is the case if and only if $W$ is cobordant rel$\partial$ to a twisted double.

It is also possible to define a chain complex version of that construction: an asymmetric Poincar\'e complex consisting of 
the singular chain complex $C(\widetilde V,\widetilde{ \partial M})$ together with a chain equivalence inducing the isomorphisms
$\lambda\colon H^{n+1-*}(\widetilde V,\widetilde{\partial M}) \isora H_*(\widetilde V,\widetilde{\partial M})$. 
The maps fit into a diagram of exact sequences
\begin{eqnarray}
\label{asygeoeqn}
\\\nonumber
\xymatrix
{
H^{n+1-r}(\widetilde V,\widetilde W)\cong H^{n-r}(\widetilde M) \ar[d]^{\cong} \ar[r]&
H^{n+1-r}(\widetilde V, \widetilde{\partial M}) \ar[d]^{\cong}_{\lambda} \ar[r]&
H^{n+1-r}(\widetilde W,\widetilde{\partial M}) \ar[d]^{\cong}
\\
H_r(\widetilde M,\widetilde{\partial M}) \ar[r]&
H_r(\widetilde V,\widetilde{\partial M}) \ar[r]&
H_r(\widetilde W,\widetilde M+\widetilde M)
}
\end{eqnarray}
In particular, if $\partial M=\emptyset$, the asymmetric complex is $C_*(\widetilde V)$ together with the Poincar\'e duality on $V$.
In particular, if $\partial M=\emptyset$, the asymmetric complex is $C_*(\widetilde V)$ together with the Poincar\'e duality on $V$.
One can find a twisted double cobordant rel$\partial$ to $W$ if and only if
that asymmetric complex is zero in the asymmetric $L$-group $LAsy^{n+1} (\bZ[\pi_1(V)])$
(see \eg \cite{Ran98} 30.12).

\section{Asymmetric Forms, Complexes and Pairs}
\label{asymalgsec}

We present the algebraic equivalents of the geometric constructions of the previous section.
Note that the asymmetric signatures of manifolds do not require normal maps. 
Therefore it is not surprising that the asymmetric signatures only require symmetric
complexes and pairs and not quadratic ones. (Obviously we can always symmetrize
any quadratic complexes, \etc and feed that information into the asymmetric signature
construction. Compare Section \ref{flipsymtwdblsec}.)

For the following compare with \cite{Ran98} Chapter 28F. 

\begin{defi}\index{Asymmetric form}
An {\bf asymmetric form $(M,\lambda)$ over $\Lambda$} is a \fg free $\Lambda$-module 
$M$ and a $\lambda \in \Hom_\Lambda(M, M^*)$. It is {\bf non-singular} if and only if
$\lambda$ is an isomorphism of $\Lambda$-modules.\index{Asymmetric form!non-singular}

A {\bf lagrangian $L$ of an asymmetric form $(M,\lambda)$} is a direct summand $L \subset M$
such that $L=L^\perp$ with $L^\perp = \{x \in L | \lambda(x)(K)=0\}$. If an asymmetric form
has a lagrangian we call it {\bf metabolic}.\index{Asymmetric form!metabolic}\index{Metabolic}

An {\bf isomorphism $f \colon (M,\lambda) \stackrel{\cong}{\ra} (M',\lambda')$ of asymmetric forms} is
an isomorphism of $\Lambda$-modules $f\colon M \stackrel{\cong}{\ra} M'$ such that $\lambda' = f^* \lambda f$.
\index{Asymmetric form!isomorphism}

The {\bf asymmetric Witt-group $LAsy^0(\Lambda)$} is the abelian group of equivalence classes of
non-singular asymmetric forms where
$$(N_1,\lambda_1) \sim (N_2,\lambda_2) \Leftrightarrow \exists (N_1,\lambda_1)\oplus (M_1,\kappa_1) \isora (N_2,\lambda_2)\oplus (M_2,\kappa_2)$$
for some non-singular metabolic forms $(M_i,\kappa_i)$.
\index{Asymmetric form!Witt-group}\index{Witt-group!of asymmetric forms}
\index{L-group@$L$-group!asymmetric}\index{Witt-group!of asymmetric forms}\index{LAsy1@$LAsy^{0}(\Lambda)$}
\hfill\qed\end{defi}

\begin{defi}\index{Asymmetric complex}
An {\bf $n$-dimensional asymmetric complex $(C,\lambda)$ over $\Lambda$} is a chain complex $C$
together with a chain map $\lambda \colon C^{n-*} \ra C$. $(C,\lambda)$ is {\bf Poincar\'e} if
$\lambda$ is a chain equivalence.\index{Asymmetric complex!Poincar\'e}

A {\bf morphism $f\colon (C,\lambda) \ra (C',\lambda')$ of $n$-dimensional asymmetric complexes}
\index{Asymmetric complex!morphism}
is a chain map $f\colon C \ra C'$ such that there is a chain homotopy 
$\lambda' \simeq f\lambda f^*$. The morphism is an {\bf equivalence} if $f\colon C \ra C'$ is
a chain equivalence.\index{Asymmetric complex!equivalence}

An {\bf $(n+1)$-dimensional asymmetric pair $(f\colon C \ra D, (\delta\lambda, \lambda))$}
\index{Asymmetric pair}
is an $n$-dimension\-al asymmetric complex $(C,\lambda)$, a chain map $f \colon C\ra D$
and a chain homotopy $\delta\lambda \colon f\lambda f^* \simeq 0 \colon D^{n-*} \ra D$.
It is {\bf Poincar\'e}\index{Asymmetric pair!Poincar\'e} if the chain maps given by
\begin{eqnarray*}
\mat{\delta\lambda \\ (-)^{r+1}\lambda f^*} &\colon& D^{n+1-r}        \ra \cone(f)_r=D_r\oplus C_{r-1}\\
\mat{\delta\lambda &  (-)^n f  \lambda    } &\colon& \cone(f)^{n+1-r}=D^{n+1-r}\oplus C^{n-r} \ra D_r
\end{eqnarray*}
are chain equivalences, in which case $(C,\lambda)$ is Poincar\'e as well.

Asymmetric Poincar\'e complexes $(C,\lambda)$ and $(C',\lambda')$ are {\bf cobordant}\index{Asymmetric complex!cobordism}
if $(C,\lambda) \oplus (C',-\lambda')$ is the boundary of an asymmetric Poincar\'e pair.

The {\bf asymmetric $L$-groups $LAsy^n(\Lambda)$} is the cobordism group of $n$-dimensional
Poincar\'e complexes.
\index{L-group@$L$-group!asymmetric}\index{Witt-group!of asymmetric forms}\index{LAsy2@$LAsy^{n}(\Lambda)$}
\hfill\qed\end{defi}

\begin{rem}
A $0$-dimensional asymmetric complex is an asymmetric form. It is Poincar\'e if the
form is non-singular. For details see the errata to \cite{Ran98}.
\cite{Ran98} (errata) Proposition 28.34 shows that any $2m$-dimensional asymmetric
Poincar\'e complex is cobordant to an $m$-connected $2m$-dimensional asymmetric
Poincar\'e complex which again is nothing but a $0$-dimensional asymmetric
Poincar\'e complex \ie $LAsy^{2n}(\Lambda) \cong LAsy^0(\Lambda)$. (The odd-dimensional
asymmetric Witt groups are all trivial.)
Hence we will identify asymmetric Poincar\'e complexes with asymmetric forms.
\end{rem}

We explained before that there is a geometric construction to assign an asymmetric form
to manifold with a twisted double on a boundary. We will state the algebraic analogue.
For that reason we need to define a chain equivalence
of a Poincar\'e pair with a twisted double on the boundary
(which we shall call {\bf b-duality map})
modelling the Lefschetz-duality map $H^{n+1-*}(W, \partial M) 
\stackrel{\cong}{\ra} H_*(W,M+M)$. It is mimicking the diagram of 
exact sequences with the ordinary Poincar\'e dualities of our various
manifolds
\begin{eqnarray}
\label{bdualgeoeqn}
\xymatrix
{&
\ar[r] &
H^{n-r}(\partial M)     \ar[d]_{\text{P.D.}}^{\cong}        \ar[r] \ar@{} [dr]|{\pm}&
H^{n+1-r}(W, \partial M)\ar[d]_{\text{P.D.}}^{\cong}        \ar[r]&
H^{n+1-r}(W)        \ar[d]_{h\cdot\text{P.D.}}^{\cong}  \ar[r]&
\\&
\ar[r] &
H_{r-1}(\partial M)                         \ar[r] &
H_r(W, M+M)                             \ar[r] &
H_r(W, \partial W)                      \ar[r] &
}
\end{eqnarray}
The rules for the cap product show that the first square commutes up to an alternating sign.

\begin{deflem}[\cite{Ran98} 30.10]
\label{cplxdef}
Let $x=(g\colon \partial E \ra E,(\theta,\partial\theta))$ be an $(n+1)$-dimensional symmetric
Poincar\'e pair such that the boundary $(\partial E, \partial\theta)$ is a twisted double of an $n$-dimensional symmetric
Poincar\'e pair  with respect to a self-equivalence
$(h,\chi)\colon (C,\phi) \homra (C,\phi)$ (compare Definition \ref{twdbldef}).
We write
$$
    g=\mat{j_0 & j_1 & k}\colon \partial E_r = D_r\oplus D_r \oplus C_{r-1} \ra E_r
$$

\begin{enumerate}

\item The {\bf b-duality map $\kappa \colon \cone(j_0 f)^{n+1-r} \ra \cone(j_0, j_1)$ of $x$}
    \index{b-duality map}
      is defined (up to chain homotopy) as the chain equivalence which fits
      into the chain homotopy commutative diagram of exact sequences
      \begin{eqnarray}
    \label{kappadiag}
    \xymatrix@+15pt
    {&
    0           \ar[r]&
    C^{n-*}         \ar[r]^-{\svec{0\\1}}   \ar[d]_{\xi}^{\cong} \ar@{} [dr]|{\pm}&
    \cone(j_0 f)^{n+1-*}    \ar[r]^-{\svec{1\\0}}   \ar[d]_{\kappa}^{\cong}&
    E^{n+1-r}       \ar[r]          \ar[d]_{\nu}^{\cong}&
    0
    \\&
    0           \ar[r]&
    C_{*-1}         \ar[r]_{\beta}&
    \cone(j_0, j_1)     \ar[r]_{\alpha}&
    \cone(g)        \ar[r]&
    0
    }
      \end{eqnarray}
      (compare (\ref{bdualgeoeqn})) such that
      \begin{enumerate}
      \item the first square commutes up to  an alternating sign
      \item $\alpha_r=\svec{1&0&0\\0&1&0\\0&0&1\\0&0&0}\colon\\
        \cone(j_0,j_1)_r = E_r\oplus D_{r-1}\oplus D_{r-1} \ra \cone(g)_r=E_r\oplus D_{r-1}\oplus D_{r-1}\oplus C_{r-2}$
      \item $\beta_r=\mat{k \\ -fh \\ -f}\colon C_{r-1} \ra \cone(j_0,j_1)_r$
      \item $\nu_r= \mat{\theta_0\\(-)^{n+1-r}\partial\theta_0 g^*} \colon E^{n+1-r} \ra \cone(g)_r=E_r\oplus\partial E_{r-1}$ 
        is the Poincar\'e duality map of $x$ (see Definition \ref{pairdef})
      \item $\xi_r= \phi_0 h \colon C^{n-r} \ra C_{r-1}$.
      \end{enumerate}

\item 
The {\bf asymmetric complex $(B,\lambda)$ of $x$}\index{Asymmetric complex!of a Poincar\'e pair} is 
(up to chain homotopy) the $(n+1)$-dimensional
asymmetric Poincar\'e complex with
$$
    B=\cone(j_0-j_1\colon D \ra \cone(j_0 f \colon C \ra E))
$$
and $\lambda \colon B^{n+1-*} \ra B$ a chain equivalence which fits into the chain homotopy
commutative diagram of exact sequences (compare (\ref{asygeoeqn}))
\begin{eqnarray}
\label{defdiag}
\xymatrix@+15pt
{&
0           \ar[r]&
D^{n-*}         \ar[r]^{\pi}        \ar[d]_{\zeta}^{\cong}&
B^{n+1-*}       \ar[r]^{\iota}      \ar[d]_{T\lambda}^{\cong}&
\cone(j_0 f)^{n+1-*}    \ar[r]          \ar[d]_{\kappa}^{\cong}&
0
\\&
0           \ar[r]&
\cone(f)        \ar[r]^{\tau}&
B           \ar[r]^{\sigma}&
\cone(j_0, j_1)     \ar[r]&
0
}
\end{eqnarray}
with
\begin{enumerate}
\item $\pi_r = \mat{0\\0\\1}\colon      D^{n-r}   \ra B^{n+1-r}= E^{n+1-r}\oplus C^{n-r}\oplus D^{n-r}$
\item $\iota_r=\mat{1&0&0\\0&1&0}\colon     B^{n+1-r} \ra \cone(j_0 f)^{n+1-r}= E^{n+1-r}\oplus C^{n-r}$
\item $\tau_r =\mat{j_0 & 0\\0&1\\0&0}\colon    \cone(f)_r=D_r\oplus C_{r-1} \ra B_r=E_r \oplus C_{r-1} \oplus D_{r-1}$
\item $\sigma_r=\mat{1&0&0\\0&f&1\\0&0&-1}\colon B_r \ra \cone(j_0, j_1)_r=E_r\oplus D_{r-1}\oplus D_{r-1}$
\item $\zeta_r=\mat{\delta\phi_0\\ (-)^{n-r}\phi_0 f^*}\colon D^{n-r} \ra \cone(f)_r$ the Poincar\'e 
      duality map of $(f\colon C\ra D, (\delta\phi, \phi))$ (see Definition \ref{pairdef})
\item $T \lambda_r=(-)^{rn}\lambda^* \colon B^{n+1-r} \ra B_r$ the duality involution of $\lambda$
\item $\kappa$ the b-duality map.
\end{enumerate}

\item 
{\bf The asymmetric signature of $x$}\index{Asymmetric signature!of a Poincar\'e pair} is the asymmetric cobordism class
$$\sigma^*(x)=[(B,\lambda)]\in LAsy^{n+1}(\Lambda)$$
\end{enumerate}
\end{deflem}

\begin{proof}
All we need to do is to show that the vertical sequences of Diagram (\ref{kappadiag}) and (\ref{defdiag})
are exact. First we notice that $\cone(\beta \colon C_{*-1} \ra \cone(j_0, j_1)) = \cone(g)$. Hence
the bottom sequence of Diagram (\ref{kappadiag}) is the mapping cone sequence of $\beta$ and therefore
exact. 

In the case of the bottom sequence of Diagram (\ref{defdiag}) we note that
every element of $\ker \tau_r$ is in the image of $\sigma_r$. On
the other side $\sigma_r \circ \tau_r$ is null-homotopic:
\begin{eqnarray*}
\sigma_r \circ \tau_r=
d_{\cone(j_0,j_1)} \Delta_{r+1} + \Delta_{r} d_{\cone(f)} 
&\colon&
\cone(f)_r \ra \cone(j_0,j_1)_r
\\
\Delta_r=\svec{0 & 0 \\ (-)^{r-1} & 0 \\ 0 & 0} 
&\colon& \cone(f)_{r-1}=D_{r-1}\oplus C_{r-2}\\
&&\quad\ra 
\cone(j_0, j_1)_r=E_r\oplus D_{r-1}\oplus D_{r-1}
\end{eqnarray*}
\end{proof}

\begin{rem}
The asymmetric signature is vanishing if and only if one can extend the twisted double structure
on the boundary to the whole Poincar\'e pair $x$
and it is invariant under cobordism (see \cite{Ran98} Proposition 30.11).
We will only need certain properties:
Proposition \ref{asymgenprop} states that two Poincar\'e pairs which are cobordant rel$\partial$ 
have the same asymmetric signature and Proposition \ref{hcobasymzeroprp} shows that
the asymmetric signature of an algebraic $h$-cobordism is zero.
\end{rem}

Now we will present an explicit formula for asymmetric complex $(B,\lambda)$.
\begin{prop}
\label{cplxprop}
The chain maps $\kappa$ and $\lambda$ in the previous definition are given
(up to chain homotopy) by
\begin{eqnarray*}
\kappa_r &=& \mat{  \theta_0                         & (-)^{n-r} k\phi_0 h^*\\
            (-)^{n-r+1}(\delta\phi_0+ (-)^{n-1}f \chi_0 f^*) j_0^* & (-)^{n-r+1} fh\phi_0 h^*\\
            (-)^{n-r}(\delta\phi_0 j_1^* + f\phi_0 k^*)      & (-)^{n-r+1} f \phi_0 h^*}
\\
&& \cone(j_0 f)^{n+1-r}= E^{n+1-r}\oplus C^{n-r} \ra \cone(j_0, j_1)_r=E_r\oplus D_{r-1}\oplus D_{r-1}
\\\\
T \lambda_r &=& \mat{   \theta_0                    & (-)^{n-1} j_0 f\chi_0 + (-)^{n-r}k\phi_0 h^*  & j_1\delta\phi_0\\
            (-)^{n-r}\phi_0 k^*             & (-)^{n-r+1}\phi_0 (1+h^*)         & (-)^{n-r} \phi_0 f^*\\
            (-)^{n-r+1}(\delta\phi_0 j_0^*+f\phi_0 k^*) & (-)^{n-r}f\phi_0 h^*              & (-)^{n-r+1}f \phi_0 f^*}
\\
&& B^{n+1-r}=E^{n+1-r}\oplus C^{n-r}\oplus D^{n-r} \ra B_r=E_r \oplus C_{r-1} \oplus D_{r-1}
\end{eqnarray*}
\end{prop}

\begin{proof}
First we check that our definition of $\kappa$ makes the right square of Diagram (\ref{kappadiag})
commute \ie whether $Z = \alpha\kappa - \nu \mat{1 & 0}\colon \cone(j_0 f)^{n+1-*} \ra \cone(g)$ is null-homotopic.
We define 
\begin{eqnarray*}
\Delta_r=\mat{0&0\\0&0\\0&0\\0&(-)^n \phi_0 h^*} &\colon& 
\cone(j_0 f)^{n+2-r}=E^{n+2-r}\oplus C^{n+1-r} 
\\&&\quad\ra\cone(g)_r = E_r \oplus D_{r-1}\oplus D_{r-1}\oplus C_{r-2}
\end{eqnarray*}

and find that $Z_r= d_{\cone(g)}\Delta_{r+1}+\Delta_r d_{\cone(j_0 f)^{n+1-*}} \cone(j_0 f)^{n+1-r} \ra \cone(g)_r$.

The left square commutes up to an alternating sign because
$$
(-)^{n-r}\beta_r\xi_r=\kappa_r\mat{0\\1}\colon C^{n-r} \ra \cone(j_0,j_1)_r = E_r\oplus D_{r-1}\oplus D_{r-1}.
$$

This choice of $\kappa$ helps us to confirm that our formula of $\lambda$ fits into Diagram (\ref{defdiag}).
$Y=T\lambda\pi - \tau\zeta \colon D^{n-*} \ra B$ is null-homotopic with chain homotopy
$$\Delta'_r \colon \mat{0\\0\\(-)^r\delta\phi_0} \colon D^{n-r+1} \ra B_r=E_r \oplus C_{r-1} \oplus D_{r-1}$$
which fulfils
\begin{eqnarray*}
Y_r = d_{B} \Delta'_{r+1} + \Delta'_r d_{D^{n-*}} \colon D^{n-r}  \ra B_r
\end{eqnarray*}
Hence the left square of Diagram (\ref{defdiag}) commutes up to homotopy.
For the other square define $Z=\kappa\iota-\sigma T\lambda \colon B^{n+1-*} \ra \cone(j_0, j_1)$.
Then 
\begin{eqnarray*}
\Delta''_r = \mat{0 & 0 & 0\\ 0 & (-)^{n-r+1}f\chi_0 & 0\\0&0&(-)^{r}\delta\phi_0}&\colon&
B^{n+1-r}=E^{n+1-r}\oplus C^{n-r}\oplus D^{n-r}
\\&&\quad\ra \cone(j_0, j_1)_r=E_r\oplus D_{r-1}\oplus D_{r-1}
\end{eqnarray*}

defines a chain null homotopy of $Z$ \ie
$$
Z_r= d_{\cone(j_0,j_1)} \Delta''_{r+1}+\Delta''_{r} d_{B^{n+1-*}}
$$

Finally, one has to check that $\lambda$ is a chain map. Because $Z$ is obviously a chain map and
$\iota$ consists of injective module homomorphisms it follows easily that $\kappa$ is a chain map as well.
By the five-lemma the chain maps $\lambda$ and $\kappa$ have to be chain equivalences.
\end{proof}

In the case of $\partial E=0$ or $C=0$ the asymmetric complex $(B,\lambda)$ is
equivalent to the obvious symmetric complexes.

\begin{lem}
Let $(C,\phi)$ be an $\epsilon$-symmetric $n$-dimensional Poincar\'e complex.
Then the identity induces an equivalence $1 \colon (C,\phi_0)\stackrel{\simeq}{\ra}(C, T_\epsilon\phi_0)$ 
of $n$-dimensional asymmetric complexes.
(\ie $(C,\phi_0)=(C, T_\epsilon\phi_0) \in LAsy^n(\Lambda)$)
\end{lem}

\begin{cor}
\label{asymcor}
We use the notation of Definition \ref{cplxdef}.
\begin{enumerate}
\item Let $\partial E=0$. Then $(B,\lambda)$ and $(E,\theta_0)$ are  equivalent.
\item Let $C=0$. Let $(V, \sigma)$ be the union of the fundamental $n$-dimensional symmetric
Poincar\'e pair $(g=(j_0,j_1)\colon D\oplus D \ra E, (\theta, \partial\theta=\delta\phi\oplus -\delta\phi))$
(see Definition \ref{unionfunddefi}).
Then $(B,\lambda)=(V,\sigma_0) \in LAsy^{n+1}(\Lambda)$.
\end{enumerate}
\end{cor}

\section{The Asymmetric Signature of a Flip-Isomorphism}
\label{asymsignsec}

Let $z=\kreckfs$ be a regular $\epsilon$-quadratic split preformation. Let $t=(\alpha,\beta,\bar\nu)$ 
be a flip-isomorphism of $z$. Let $\nu$ be a representative of $\bar\nu$. As an abbreviation, define 
$\sigma=(\nu-\epsilon\nu^*)^*$. $\sigma$ is independent of the choice of representative for $\bar\nu$.

The $(2q+2)$-dimensional symmetric Poincar\'e pair $x^t=(g^t \colon \partial E_t \ra E, (0, \theta_t))$
of Section \ref{flipsymtwdblsec}
has a twisted double structure on its boundary which enables us to apply the asymmetric signature
construction from  Proposition \ref{cplxprop}. The result is the $(2q+2)$-dimensional
asymmetric complex $(B,\lambda)$.

We can reduce this complex to a smaller $(2q+2)$-dimensional
asymmetric complex $(B',\lambda')$ via the chain equivalence
\begin{eqnarray*}
\xymatrix@+1em
{
B_{q+2}=G \oplus G              \ar[r]_{\Psi_2}^{\cong}             \ar[d]_{d_1}&
G \oplus G                  \ar[r]_{\svec{0&1}}             \ar[d]_{\svec{1&0\\0&d_1'}}&
B'_{q+2}=G                                          \ar[d]_{d'_1}
\\
B_{q+1}=G \oplus (F\oplus F^*) \oplus F^*   \ar[r]_-{\Psi_1}^-{\cong}           \ar[d]_{d_0}&
G \oplus (F\oplus F^*) \oplus F^*       \ar[r]_-{\svec{0&1&0&0\\0&0&1&0\\0&0&0&1}}  \ar[d]_{\svec{0& d_0'}}&
B'_{q+1}= F\oplus F^*  \oplus F^*                                   \ar[d]_{d'_0}
\\
B_q=G^*                     \ar[r]_{\Psi_0}^{\cong}&
G^*                     \ar[r]_{1}&
B'_q=G^*
}
\end{eqnarray*}
with
\begin{eqnarray*}
\Psi_2 &=& \mat{1&1+\beta\\0&1}\\
\Psi_1 &=& -\epsilon\mat{1 & 0 &0 &0\\\gamma&-1&0&-\alpha^{-*}\\0&0&0&1\\
             \mu-\epsilon\alpha\gamma&\epsilon\alpha&-1&\epsilon\alpha\alpha^{-*}-1-\epsilon\alpha\sigma}\\
\Psi_0&=&1\\
d_1'&=&\mat{\gamma\\\mu\\0}\\
d_0'&=&\mat{\epsilon(1+\beta^{-*})\mu^* & (1+\beta^{-*})\gamma^* & \gamma^*}\\
\lambda_{q+2}'&=&\epsilon\colon {B'}^{q} \ra {B'}_{q+2}\\
\lambda_{q+1}'&=&\mat{0&0&1\\\epsilon&0&-\alpha^*\\-\epsilon&\epsilon\alpha&\alpha^*-\epsilon\alpha+\epsilon\alpha\sigma\alpha^*}
        \colon {B'}^{q+1}\ra{B'}_{q+1}\\
\lambda_{q}'  &=&-\beta^{-*} \colon {B'}^{q+2}\ra{B'}_{q}
\end{eqnarray*}
(All $\lambda_r$ are in fact isomorphisms of $\Lambda$-modules.)

With the help of \cite{Ran98} (errata) 28.34, we compute a highly-connected $(2q+2)$-dimensional
asymmetric complex $(B'',\lambda'')$
which is cobordant to the asymmetric complex $(B',\lambda')$:
\begin{eqnarray*}
\xymatrix@+2em
{
{B''}^{q}\ar[d]_{\epsilon {d_0''}^*}\ar[r]_{\lambda_{q+2}''}&
B''_{q+2}=G \ar[d]_{d_1''}
\\
{B''}^{q+1}     \ar[d]_{-\epsilon{d''_1}^*}     \ar[r]_-{\lambda_{q+1}''}&
B''_{q+1}=F\oplus F^*\oplus F^*\oplus G^* \oplus G \ar[d]_{d_0''}
\\
{B''}^{q+2}                     \ar[r]_{\lambda_q''}&
B_q''=G^*
}
\end{eqnarray*}
with
\begin{eqnarray*}
d_1''&=&\mat{\gamma\\\mu\\0\\0\\\beta^{-1}}\\
d_0''&=&\mat{\epsilon(1+\beta^{-*})\mu^* & (1+\beta^{-*})\gamma^* & \gamma^* & -\epsilon\beta^{-*} & 0}\\
\lambda''_{q+1}&=&\mat{ 0&0&1&0&0\\\
            \epsilon&0&-\alpha^*&0&0\\
            -\epsilon&\epsilon\alpha&\alpha^*-\epsilon\alpha+\epsilon\alpha\sigma\alpha^*&0&0\\
            0&0&0&0&-\epsilon\beta^{-*}\\
            0&0&0&1&0}
\end{eqnarray*}

We can simplify this asymmetric complex to gain the asymmetric form we were looking for using the
isomorphism of chain complexes
\begin{eqnarray*}
\xymatrix@R+1em
{
B_{q+2}''=G \ar[r]_{\Phi_2} \ar[d]_{d''_1}&
G \ar[d]^{\svec{0\\0\\0\\0\\1}}
\\
B_{q+1}''=F\oplus F^*\oplus F^*\oplus G^* \oplus G \ar[r]_-{\Phi_1}\ar[d]_{d_0''}&
F^*\oplus F \oplus F^* \oplus G^* \oplus G \ar[d]^{\svec{0&0&0&1&0}}
\\
B_q''=G^* \ar[r]_{\Phi_0}&
G^*
}
\end{eqnarray*}
with
\begin{eqnarray*}
\Phi_2&=& 1\\
\Phi_1&=& \mat{ \alpha              &0              &0      &0          &-\alpha\gamma\beta\\
        0               &\epsilon\alpha^{-*}        &0      &0          &-\epsilon\alpha^{-*}\mu\beta\\
        0               &1              &1      &0          &-\mu\beta\\
        \epsilon(1+\beta^{-*})\mu^* &(1+\beta^{-*})\gamma^*     &\gamma^*   & -\epsilon\beta^{-*}   & 0\\
        0&0&0&0&\beta}\\
\Phi_0&=& 1
\end{eqnarray*}

Thus the asymmetric signature is given by the asymmetric form 
$$
    \rho=\mat{0&0&\alpha\\1&0&-\epsilon\\0&1&\epsilon\alpha\sigma\alpha^*}\colon M=F\oplus F^*\oplus F \ra M^*
$$

It is clear that $\rho$ does not depend on the choice of representative for $\bar\nu$. Hence we can define
\begin{defi}
\label{asymsignflipdef}\index{Asymmetric signature!of a flip-isomorphism}
Let $z=\kreckfs$ be an $\epsilon$-quadratic split preformation.
The {\bf asymmetric signature $\sigma^*(z,t)$ of a flip-isomorphism $t=(\alpha,\beta,\bar\nu)$ of $z$} 
is an element $(M,\rho) \in LAsy^0(\Lambda)$ given by
$$
    \rho=\mat{0&0&\alpha\\1&0&-\epsilon\\0&1&\epsilon\alpha(\nu^*-\epsilon\nu)\alpha^*}\colon M=F\oplus F^*\oplus F \ra M^*
$$
Similar for the non-split case.
\hfill\qed\end{defi}

\section{Asymmetric Signatures and Elementariness}
\label{asymsignelemsec}

In this section we show that the asymmetric signatures are an obstruction to elementariness
\begin{thm}
\label{elemasymsignthm}
Let $z=\kreckfs$ be a regular $\epsilon$-quadratic split preformation which allows flip-isomorphisms.
If $z$ is elementary then the asymmetric signature $\sigma^*(z,t)\in LAsy^0(\Lambda)$ vanishes for all
flip-isomorphisms $t$.
\end{thm}

\begin{rem}
Theorem \ref{flasymthm} will present a version for stably elementary preformations.
The converse is not true in general. Counterexamples are presented in
Example \ref{asymcountex}.

The asymmetric signatures are generally not-trivial as Corollary \ref{asymbdrycor}
will show.
\end{rem}

We will give two proofs for this theorem. The first one is based on algebraic surgery theory whereas
the second proof is a low-level calculation of asymmetric forms. 

The Definition \ref{asymsignflipdef} shows that the asymmetric signatures only depend on
the underlying non-split preformation. It ignores the quadratic structures of both preformation
and flip-isomorphism. A generalization of Theorem \ref{elemasymsignthm} for non-split preformations
comes without surprise:
\begin{cor}
\label{elemasymsign2thm}
Let $z=\skreckfs$ be a regular $\epsilon$-quadratic preformation which allows flip-isomorphisms.
If $z$ is elementary then the asymmetric signature $\sigma^*(z,t)\in LAsy^0(\Lambda)$ vanishes for all
flip-isomorphisms $t$.
\end{cor}
\begin{proof}
The second proof for Theorem \ref{elemasymsignthm} works also for the non-split case.
There should be no problem in using algebraic surgery again - one ``only'' needs to 
prove symmetric versions of the previous two chapters. We leave this as an exercise to the reader.
\end{proof}

The first proof of Theorem \ref{elemasymsignthm} needs some preparation. In the next two propositions
we show algebraic equivalents of the following facts from the world of manifolds:
\begin{enumerate}
\item  Two manifolds with a twisted double on their boundary 
have the same asymmetric signature if they are cobordant rel$\partial$.
\item An $s$-cobordism $(W,M,M)$ is in fact a twisted double and hence its asymmetric signature must vanish.
\end{enumerate}
\begin{prop}
\label{asymgenprop}
Let $x=(g\colon \partial E \ra E,(\theta,\partial\theta))$ and $x'=(g'\colon \partial E \ra E',(\theta',\partial\theta))$
be two $(n+1)$-dimensional symmetric
Poincar\'e pairs such that the boundary $(\partial E, \partial\theta)$ is a twisted double of an $n$-dimensional symmetric
Poincar\'e pair $(f\colon C \rightarrow D,(\delta\phi,\phi))$ with respect to a homotopy self-equivalence
$(h,\chi)\colon (C,\phi) \homra (C,\phi)$ (compare Definition \ref{twdbldef}).

\begin{enumerate}
\item If $x$ and $x'$ are cobordant rel$\partial$, then $\sigma^*(x)=\sigma^*(x')\in LAsy^{n+1}(\Lambda)$.
\item $\sigma^*(x)-\sigma^*(x')=\sigma^*(x\cup-x')\in LAsy^{n+1}(\Lambda)$.
\end{enumerate}
\end{prop}
\begin{proof}
\begin{enumerate}
\item This is a special case of \cite{Ran98} Proposition 30.11(iii).
\item By using the union construction one can easily verify that $((x \cup -x') + x') \cup x=(x\cup -x')+ (-(x\cup -x'))$.
In this formula the sum $(C,\phi) + (f\colon D\ra E, (\partial\theta, \theta))$ of an $(n+1)$-dimensional $\epsilon$-symmetric
(Poincar\'e) complex with an $(n+1)$-dimensional $\epsilon$-symmetric (Poincar\'e) pair is the $(n+1)$-dimensional $\epsilon$-symmetric
(Poincar\'e) pair $(\svec{f\\0}\colon D\ra E\oplus C, (\partial\theta\oplus \phi, \theta)$).

For any Poincar\'e complex $(C,\phi)$ (and in particular for $x\cup -x'$) $$(\mat{1&-1}\colon C\oplus C \ra C, (0,\phi\oplus-\phi))$$
defines a null-cobordism of $(C,\phi)\oplus (C,-\phi)$. Hence the pairs $x$ and $(x \cup -x')+ x'$ are cobordant rel$\partial$.
Therefore $\sigma^*(x)=\sigma^*((x \cup -x')+ x')\in LAsy^{n+1}(\Lambda)$. It is not hard to see that the latter expression is 
the same as $\sigma^*(x\cup x') +\sigma^*(x')$.
\end{enumerate}
\end{proof}
\begin{prop}
\label{hcobasymzeroprp}
Let $x=(g\colon \partial E \ra E,(\theta,\partial\theta))$ be an $(n+1)$-dimensional symmetric
Poincar\'e pair such that the boundary $(\partial E, \partial\theta)$ is a twisted double of an $n$-dimensional symmetric
Poincar\'e pair $(f\colon C \ra D,(\delta\phi,\phi))$ with respect to a self-equivalence
$(h,\chi)\colon (C,\phi) \ra (C,\phi)$ (compare Definition \ref{twdbldef}).
We write 
$$
    g=\mat{j_0 & j_1 & k}\colon \partial E_r = D_r\oplus D_r \oplus C_{r-1} \ra E_r
$$

Additionally assume that $x$ is an $h$-cobordism \ie that $j_0, j_1 \colon D \ra E$ are chain equivalences.
Then $\sigma^*(x)=0\in LAsy^{n+1}(\Lambda)$.
\end{prop}
\begin{proof}
We could refer to \cite{Ran98} Proposition 30.11(ii) but instead we give a quick and direct proof of the claim.
Obviously it is enough to construct an asymmetric null-cobordism for the asymmetric Poincar\'e complex $(B, T\lambda)$
given in Proposition \ref{cplxprop}.\footnote{In general, if $(f\colon C \ra D, (\delta\lambda, \lambda))$
is an $n$-dimensional asymmetric (Poincar\'e) pair, so is $(f\colon C \ra D, (T_\epsilon \delta\lambda, T_\epsilon \lambda))$.}

We define the $(n+2)$-dimensional asymmetric Poincar\'e pair
\begin{eqnarray*}
&&(s \colon B \ra D_{*-1}, (\delta\lambda, \lambda))\\
s&=&\mat{0&0&1} \colon B_r=E_r\oplus C_{r-1}\oplus D_{r-1} \ra D_{r-1}\\
\delta\lambda &=& (-)^{r+1}\delta\phi_0 \colon D^{(n+2)-r-1} \ra D_{r-1}
\end{eqnarray*}
In order to proof that it is Poincar\'e one observes that there is a chain equivalence
$\cone(s) \simeq \cone(j_0 f)_{*-1}$ given by
\begin{eqnarray*}
\mat{j_0-j_1&1&0&0\\0&0&1&0}
&\colon& \cone(s)_r = D_{r-1}\oplus E_{r-1}\oplus C_{r-2} \oplus D_{r-2}\\
&&\quad \ra \cone(j_0 f)_{r-1}=E_{r-1}\oplus C_{r-2}
\end{eqnarray*}
\end{proof}

\begin{proof}[First proof of Theorem \ref{elemasymsignthm}]
By Theorem \ref{algmainthm} the Poincar\'e pair $x=(g \colon \partial E\\ \ra E, (\delta\omega=0, \omega))$
from Section \ref{constrpairsec} is cobordant rel$\partial$ to an algebraic $h$-cobordism.
By Section \ref{flipsymtwdblsec} the Poincar\'e pair $x_t=(g \colon \partial E_t \ra E, (0,\partial\theta_t))$
is cobordant rel$\partial$ to an algebraic $h$-cobordism. By Proposition \ref{hcobasymzeroprp}
its asymmetric signature (which is $\sigma^*(z,t)\in LAsy^0(\Lambda)$ by Section \ref{asymsignsec}) 
is vanishing.
\end{proof}

\begin{proof}[Second proof of Theorem \ref{elemasymsignthm}]
We can also give a proof of the theorem without algebraic surgery theory.
For simplicity we assume that $\kreckfs$ has the nice form presented in Proposition~\ref{elemprop}\ref{elemprop3}
and that $i \colon U \hookrightarrow G$ is the inclusion of the $h$-lagrangian.
We define a metabolic asymmetric form $(M', \rho')$:
\begin{eqnarray*}
a &=&\mat{0&0\\1&0} \colon G^*=U^*\oplus R^* \ra F  = U  \oplus U^*\\
b &=&\mat{1&0\\0&0} \colon G^*=U^*\oplus R^* \ra F^*= U^*\oplus U\\
t &=&\mat{1&0\\\tau^*&1} \colon G^*=U^*\oplus R^* \ra G^*=U^*\oplus R^*\\
k &=&\mat{\alpha a\\0\\-\epsilon \alpha a - \epsilon b} \colon G^* \ra F^*\oplus F\oplus F^*\\
\rho'&=&\mat{\rho & k & 0\\ 0 & 0 & -\epsilon\beta\\ 0& t&0} \colon
    M'=M\oplus G^* \oplus G \ra {M'}^* = M^* \oplus G \oplus G^*
\end{eqnarray*}
A lagrangian for $(M',\rho')$ is given by
\begin{eqnarray*}
j&=&\mat{   \epsilon\gamma\beta i   & \epsilon\gamma(1+\beta)       & 0\\
        \epsilon\alpha\gamma i  & \epsilon\alpha\gamma\beta^{-1}    & \epsilon\alpha\\
        0           & \gamma                & 0\\
        0           & 0                 & \gamma^*\alpha\\
        -\epsilon i     & -\epsilon             & 0}
\colon U \oplus G \oplus F \ra M' = (F\oplus F^* \oplus F) \oplus G^* \oplus G
\end{eqnarray*}
Then $(M,\rho)\oplus (M',-\rho')$ has a lagrangian
\begin{eqnarray*}
l&=&\mat{1&0\\1&0\\0&0\\0&1} \colon M \oplus G \ra M\oplus M' = M\oplus M \oplus G^* \oplus G
\end{eqnarray*}
\end{proof}

\section{The Flip-$l$-Monoids}
\label{asymflmonoidelemsec}

The asymmetric signature of a flip-isomorphism suggests the definition of an extension of 
Kreck's $l$-monoids which includes a choice of a stable flip-isomorphisms.

\begin{defi}
\label{flmonoidef}
Let $z=\skreckfs$ and $\skreckfa{F'}{\gamma'}{G'}{\mu'}$ be regular $\epsilon$-quadratic preformations
and let $t=(\alpha,\beta,\sigma)$ and $t'=(\alpha',\beta',\sigma')$ be flip-isomorphisms of $z$ and $z'$ respectively.

An isomorphism $(\eta,\zeta)$ of the tuples $(z,t)$ and $(z',t')$ is a strong isomorphism 
$(\eta,\zeta)\colon z \isora z'$ of preformations such that $\alpha'=\eta^{-*}\alpha\eta^{-1}$,
$\beta'=\zeta\beta\zeta^{-1}$ and $\sigma'=\eta\sigma\eta^*$.

The sum $(z,t)+(z',t')$ is the well-defined tuple $(z\oplus z', t \oplus t')$.

Let $y^k=\partial H^{-\epsilon}(\Lambda^{k})$ be a hyperbolic preformation and 
$t^k=\left(\svec{0&-\epsilon\\-1&0},\svec{1&0\\0&-1},0\right)$ a (strong) flip-isomorphism of $y^k$.
\footnote{Compare Corollary \ref{elemflipcor}.}

A stable isomorphism of the tuples $(z,t)$ and $(z',t')$ is an isomorphism of $(z,t)+(y^k,t^k)$ with $(z',t')+(y^l,t^l)$ for
some $k,l\in\bN_0$.
The stable isomorphism classes form an abelian monoid namely the {\bf flip-$l$-monoid $fl^{2q+2}(\Lambda)$}.
\index{Flip-$l$-monoid}\index{fl1@$fl^{2q+2}(\Lambda)$}

Similarly we can define the {\bf flip-$l$-monoid $fl_{2q+2}(\Lambda)$} in the split case.
\index{Flip-$l$-monoid}\index{fl2@$fl_{2q+2}(\Lambda)$}
\hfill\qed
\end{defi}

\begin{rem}
\label{pirem}
\begin{enumerate}
\item There are well-defined morphisms of abelian monoids
\begin{align*}
\pi\colon fl^{2q+2}(\Lambda) &\ra l^{2q+2}(\Lambda), & [(z,t)] &\mt [z]\\
\pi\colon fl_{2q+2}(\Lambda) &\ra l_{2q+2}(\Lambda), & [(z,t)] &\mt [z]
\end{align*}
\item
There is a well-defined morphism of abelian monoids
\begin{eqnarray*}
fl_{2q+2}(\Lambda) &\ra& fl^{2q+2}(\Lambda)\\
(\kreckfs, (\alpha,\beta,\bar\nu)) &\mt& (\skreckfs, (\alpha,\beta,\bar\nu^*-\epsilon\bar\nu^*))
\end{eqnarray*}
\end{enumerate}
\end{rem}

\begin{thm}
\label{flasymthm}
The asymmetric signature of Definition \ref{asymsignflipdef} gives rise to a 
well-defined homomorphism of abelian monoids
\begin{eqnarray*}
\sigma^*\colon fl^{2q+2}(\Lambda) &\ra& LAsy^0(\Lambda)
\\
{[(z,t)]} &\mt& \sigma^*(z,t)
\end{eqnarray*}

If $[z']\in l^{2q+2}(\Lambda)$ is elementary then
$\sigma^*(\pi^{-1}[z'])=\{0\}$ (with $\pi$ as in Remark \ref{pirem})
\ie 
$\sigma^*(z,t)=0$ for all flip-isomorphisms $t$ of
all preformations $z$ with $[z]=[z']\in\el$ (\ie for all stable flip-isomorphisms $t$ of $z$).
\end{thm}

\begin{proof}
Using the notation of Definition \ref{flmonoidef} assume that there is an isomorphism 
$(\eta,\zeta)$ of $(z,t)$ and $(z',t')$. 
Let 
\begin{eqnarray*}
\rho=\svec{0&0&\alpha\\1&0&-\epsilon\\0&1&\epsilon\alpha(\nu^*-\epsilon\nu)\alpha^*}\colon M=F\oplus F^*\oplus F \ra M^*\\
\rho'=\svec{0&0&\alpha'\\1&0&-\epsilon\\0&1&\epsilon\alpha'({\nu'}^*-\epsilon\nu'){\alpha'}^*}\colon M'=F'\oplus {F'}^*\oplus {F}' \ra {M'}^*
\end{eqnarray*}
be the asymmetric forms whose image in $LAsy^0(\Lambda)$ are by definition the asymmetric signatures of $(z,t)$ and $(z',t')$.
Then $f=\svec{\eta^{-1}&0&0\\0&\eta^*&0\\0&0&\eta^{-1}}\colon M \isora M'$ is an isometry of the asymmetric forms $(M,\rho)$ and $(M',\rho')$.
So $\sigma^*(z,t)=\sigma^*(z',t')$.

Clearly $\sigma^*(y_k,t_k)=0$ 
and it is obvious that $\sigma^*(z,t)+\sigma^*(z',t')=\sigma^*(z\oplus z', t\oplus t')$.
Hence the asymmetric signature doesn't change under stable isomorphisms of tuples $(z,t)$ and it is compatible with
the actions of both monoids.
The rest follows from Corollary \ref{elemasymsign2thm}.
\end{proof}

\chapter{Quadratic Signatures of Flip-Isomorphisms}
\label{quadchap}

{\bf For the whole chapter, let $q\ge 2$, $\epsilon=(-)^q$ and let $\Lambda$ be a weakly finite ring
with $1$ and involution.}

Let $(W,M,M')$ be a cobordism such that $\partial M=\partial M'=\emptyset$.
If there is an automorphism $H\colon M \isora M'$, we can glue $M$ on $M'$
along $H$ in order to obtain a closed manifold $V_H$. If $V_H$ is null-cobordant
then $(W,M,M')$ is cobordant rel$\partial$ to the $h$-cobordism
$(M\times I \cup_{H^{-1}} M', M, M')$.

As usual we try to transfer the above into the world of algebraic surgery theory.
Let $x=(g\colon \partial E=D'\cup_C D \ra E, (\delta\omega,\omega))$ be a
$(2q+2)$-dimensional quadratic Poincar\'e pair \eg the one we constructed
out of a preformation $z$ in Section \ref{constrpairsec}. Assume that $C$ is
zero or at least contractible. For the Poincar\'e pair $x$ from Section \ref{constrpairsec} 
that is only the case if $z$ is
a non-singular formation (see also Chapter \ref{nonsingchap}).
So we deal in fact with a Poincar\'e pair of the form $(g\colon D\oplus D' \ra E, (\delta\omega', \nu\oplus\nu'))$.
A flip-isomorphism induces an isomorphism of $D$ and $D'$, so
that $x$ transforms to a fundamental pair $x'=(D\oplus D \ra E, (\delta\omega',\nu\oplus-\nu))$
which we can glue together along $D$ (see Definition \ref{unionfunddefi}). The result 
is a $(2q+2)$-dimensional Poincar\'e complex. It is (algebraically) null-cobordant if and only if
$x'$ is cobordant rel$\partial$ to $((1,1)\colon D\oplus D \ra D, (0,\nu\oplus-\nu))$ which
is the case if and only if $x$ is cobordant to an $h$-cobordism.
Using standard surgery theory (\eg Lemma \ref{trcplxlem}) this Poincar\'e complex
corresponds to a non-singular quadratic form and that form vanishes in the even-dimensional
$L$-group if and only if the Poincar\'e complex is null-cobordant. Hence we expect to
be able to define an element in $L_{2q+2}(\Lambda)$ for each flip-isomorphism of $z$
such that $z$ is elementary if and only if such a quadratic signature vanishes for
a flip-isomorphism.

The manifold case requires more care  
if the boundary of $M$ and $M'$ is non-empty. Again we go 
through all automorphisms $H\colon M \isora M'$ and replace $M'$ by $M$ using $H$.
The original cobordism becomes $(W,M,M)$ and the boundary of $W$ turns into a twisted double 
$M \cup_h M$ with $h=H|\colon \partial M \ra \partial M'\cong \partial M$. But not every 
twisted double $M\cup_h M$ is a boundary of an $h$-cobordism. 
If we want to follow the strategy of the closed case in the beginning we have to assume that
\eg $h$ is isotopic to the identity. Then $(W,M,M)$ can be glued onto 
$M\times (I,0,1)$. If the result, a closed manifold, is null-cobordant then $(W,M,M)$ 
and hence $(W,M,M')$ is cobordant rel$\partial$ to an $h$-cobordism.

Similarly, for general preformations the situation is more complicated. In Section \ref{flipqutwdblsec}
a flip-isomorphism $t$ of $z$ replaces the boundary $\partial E$ of $x$ by an algebraic twisted double.
This yields a $(2q+2)$-dimensional quadratic Poincar\'e pair
$x_t=(g_t\colon D'\cup_{h_t}D \ra E, (\delta\omega_t,\omega_t=\delta\psi\cup_{\chi_t}-\delta\psi))$
with a twisted double at the boundary. It is not always possible to find an algebraic $h$-cobordism
with that boundary except \eg if $(h_t, \chi_t)$ is homotopic to the identity. This involves a concept
of homotopies of morphisms between quadratic complexes which we develop in the rather technical
{\bf Section \ref{quadhomsec}}. We define {\bf flip-isomorphisms rel$\partial$} in 
{\bf Section \ref{quadfliprelsec}} as flip-isomorphisms for which $(h_t,\chi_t)$ is 
homotopic to $(1,0)\colon (C,\psi) \isora (C,\psi)$.

Then we deviate slightly from the example in geometry. We use the homotopy to change $x_t$ such that 
it looks like $(D\cup_C D \ra E, (*, \delta\psi\cup_{\psi}-\delta\psi))$
and then stick the standard algebraic $h$-cobordism
$(D\cup_C D \ra D, (0, \delta\psi\cup_{\psi}-\delta\psi))$ 
on it. As before, the result will be a $(2q+2)$-dimensional Poincar\'e complex which corresponds
to a non-singular quadratic form. This is the quadratic signature which will be constructed
in {\bf Section \ref{quadconstrsec}}. In {\bf Section \ref{quadelemsec}} it is proven
that a preformation $z$ is stably elementary if and only if one of its quadratic signatures
is vanishing. 

The disadvantage of the quadratic signatures is that they not only depend
on the preformation and the flip-isomorphism, but \eg also the explicit
homotopy of $(h_t, \chi_t)\simeq (1,0)$. 
Hence we do not have something like a map $fl_{2q}(\Lambda) \ra L_{2q}(\Lambda)$.
In certain cases, though, we can restrict the effect of those choices  on the
quadratic signature (see Lemma \ref{quaddeplem}).

Curiously, the quadratic and asymmetric signatures are related by the canonical map
$L_{2q}(\Lambda) \ra LAsy^0(\Lambda)$, $(K,\psi) \ra (K,\psi-\epsilon\psi^*)$
as we will show in {\bf Section \ref{quadasymsec}}.

\section{Homotopy and Twisted Doubles}
\label{quadhomsec}

This section deals with a very technical issue, the extension of the
concept of chain homotopies of chain maps 
$$\Delta \colon f \simeq f' \colon C \ra C'$$ 
to a homotopy of morphisms of quadratic (or symmetric) complexes 
$$(\Delta,\eta)\colon (f,\chi)\simeq (f',\chi')\colon (C,\psi) \ra (C',\psi')$$ 
Obviously the chain homotopy $\Delta$
will affect the quadratic structures $\chi$ and $\chi'$. Their
difference is determined by an operation $\Delta_\%\psi$ (and a
boundary $d\eta$).  In Lemma \ref{equivrellem} we show that the
homotopy of morphisms is an equivalence relation.

We will only be able to define the quadratic signature if the chain map
$(h_t,\chi_t) \colon (C,\psi)\\ \isora (C,\psi)$ defined in (\ref{htdefeqn})
is homotopic to $(1,0)$. This makes it necessary to keep track how
a twisted double changes if the self-equivalence used is changed by
a homotopy. Lemma \ref{homqutwdblem} deals with this case.

Section \ref{quadasymsec} discusses the relationship between quadratic and asymmetric
signatures. We will need Lemma \ref{asysgnhomlem} which shows that
changing the self-equivalence involved by a homotopy will not affect the
asymmetric signature.

\begin{defi}
Let $\Delta \colon f \simeq f' \colon C \ra C'$ be a chain homotopy of two
chain maps. 

Let $\phi\in W^\%(C,\epsilon)_n$. Define $\Delta^\%\phi\in W^\%(C',\epsilon)_{n+1}$ by
$$
(\Delta^\%\phi)_s=-\Delta\phi_s f^* + (-)^{r+1}(f'\phi_s+(-)^{n+1}\Delta T_\epsilon \phi_{s-1})\Delta^*
\colon {C'}^{n+1-r+s} \ra C'_r
$$

Let $\psi\in W_\%(C',\epsilon)_n$. Define $\Delta_\%\psi\in W_\%(C',\epsilon)_{n+1}$ by
$$
(\Delta_\%\psi)_s=-\Delta\psi_s f^* + (-)^{r+1}(f'\psi_s+(-)^{n}\Delta T_\epsilon \psi_{s+1})\Delta^*
\colon {C'}^{n+1-r-s} \ra C'_r
$$
\hfill\qed\end{defi}

\begin{lem}
\label{chhomlem}
Let $\Delta \colon f \simeq f' \colon C \ra C'$ be a chain homotopy of two chain maps.
\begin{enumerate}
\item Let $\phi\in W^\%(C,\epsilon)_n$.\\Then $d(\Delta^\%\phi)=-\Delta^\%(d\phi)+f\phi f^*-f'\phi {f'}^* \in W^\%(C,\epsilon)_n$.
\item Let $\psi\in W_\%(C,\epsilon)_n$.\\Then $d(\Delta_\%\psi)=-\Delta_\%(d\psi)+f\psi f^*-f'\psi {f'}^* \in W_\%(C,\epsilon)_n$.
\item Let $\psi\in W_\%(C,\epsilon)_n$. Then
\begin{eqnarray*}
&&(1+T_\epsilon)(\Delta_\%\psi)-\Delta^\%((1+T_\epsilon)\psi)\\
&&= d\xi + 
\begin{cases} 
    (-)^r\Delta d(T_\epsilon \psi)_0 \Delta^* \colon {C'}^{n+1-r} \ra {C'}_r &:s=0\\
    0 &: s\neq 0
\end{cases}\end{eqnarray*}
with $\xi\in W^\%(C',\epsilon)_{n+2}$ given by 
$\xi_0=(-)^{r+1}\Delta T_\epsilon \psi_0\Delta^* \colon {C'}^{n+2-r} \ra {C'}_r$.

\item Let $g\colon C' \ra D$ be a chain map. Then $g\Delta\colon gf \simeq gf' \colon C \ra D$ is a 
chain homotopy. 

Let $\phi\in W^\%(C,\epsilon)_n$. Then $(g\Delta)^\%\phi=g(\Delta^\%\phi)g^* \in W^\%(D,\epsilon)_{n+1}$.

Let $\psi\in W_\%(C,\epsilon)_n$. Then $(g\Delta)_\%\psi=g(\Delta_\%\psi)g^* \in W_\%(D,\epsilon)_{n+1}$.
\end{enumerate}
\end{lem}

\begin{lem}
Let $(f,\chi) \colon (C,\psi) \ra (C',\psi')$ be a morphism of $n$-dimensional $\epsilon$-quadratic 
complexes. Let $\Delta \colon f \simeq f' \colon C \ra C$ be a chain homotopy.

Then $(f',\chi+\Delta_\%\psi) \colon (C,\psi) \ra (C',\psi')$ is also a
morphism of $n$-dimensional $\epsilon$-quadratic 
complexes.
Similar in the symmetric case.
\end{lem}

\begin{defi}
\label{hommodef}\index{Homotopy of morphisms}
A {\bf homotopy $(\Delta,\eta)$ of two morphisms of $n$-dimensional $\epsilon$-quadratic 
complexes $(f,\chi),(f',\chi') \colon (C,\psi) \ra (C',\psi')$} is a chain homotopy 
$\Delta \colon f \simeq f' \colon C \ra C'$ and an element $\eta\in W_\%(C')_{n+2}$ such that 
$$\chi'-\chi=\Delta_\%\psi+d(\eta) \in W_\%(C')_{n+1}$$

Similar in the symmetric case.
\hfill\qed\end{defi}

\begin{lem}
\label{equivrellem}
Let $(\Delta,\eta)\colon (f,\chi)\simeq (f',\chi') \colon (C,\psi) \ra (C',\psi')$
be a homotopy of two morphisms of $n$-dimensional $\epsilon$-quadratic 
complexes.
\begin{enumerate}
\item Then there is a homotopy $(\Delta'=-\Delta,\eta')\colon (f',\chi')\simeq (f,\chi) \colon (C,\psi) \ra (C',\psi')$.
\item Let $(\Delta',\eta')\colon (f',\chi')\simeq (f'',\chi'') \colon (C,\psi) \ra (C',\psi')$ be another homotopy.
      Then also $(f,\chi)$ and $(f'',\chi'')$ are homotopic.
\item Homotopy induces an equivalence relation on all morphisms $(C,\psi)\ra (C',\psi')$.
\item Let $(g,\rho)\colon (C',\psi') \ra (D,\theta)$ be a morphism. Then there is a homotopy
      $$(g\Delta,g\eta g^*) \colon (g,\rho)(f,\chi) \ra (g,\rho)(f',\chi')\colon (C,\psi) \ra (D,\theta)$$
      with $(g,\rho)(f,\chi)=(gf, \rho+g\chi g^*)$ as in Definition \ref{compinvdef}.
\end{enumerate}
Similar in the symmetric case.
\end{lem}
\begin{proof}
\begin{enumerate}
\item Use $\eta'_s=-\eta_s+(-)^{r+1}\Delta\psi_s\Delta^*\colon {C'}^{n+2-r-s}\ra C'_r$.
\item Define $\Delta''=\Delta+\Delta'$ and 
$\eta''_s=\eta_s+\eta'_s+(-)^r\Delta'\psi_s\Delta^*\colon {C'}^{n+2-r-s}\ra C'_r$.
Then  $(\Delta'',\eta'')\colon (f,\chi)\simeq (f'',\chi'') \colon (C,\psi) \ra (C',\psi')$
is a homotopy.
\item Clear from the previous two previous statements.
\item Obvious.
\end{enumerate}
\end{proof}

\begin{lem}
\label{homqutwdblem}
Let $c=(f \colon C \ra D, (\delta\psi,\psi))$ be an $n$-dimensional $\epsilon$-quadratic Poincar\'e pair.
Let $(\Delta,\eta)\colon (h,\chi) \simeq (h',\chi') \colon (C,\psi) \homra (C,\psi)$ 
be a homotopy of self-equivalen\-ces.
Then there is an isomorphism $(a,\sigma)\colon c\cup_{(h,\chi)}-c \isora c\cup_{(h',\chi')}-c$
of the corresponding quadratic twisted doubles given by
\begin{eqnarray*}
a_r&=& \mat{1&0&(-)^r f\Delta\\0&1&0\\0&0&1}\colon
\\&&(D\cup_h D)_r=D_r\oplus C_{r-1}\oplus D_r \ra (D\cup_{h'} D)_r=D_r\oplus C_{r-1}\oplus D_r\\
\sigma_s &=& \mat{  (-)^{n-1} f\eta_s f^*   &0&0\\
            0           &0&0\\
            (-)^n\psi_s\Delta^* f^* &0&0}\colon
\\&&(D\cup_{h'} D)^{n+1-r-s}=D^{n+1-r-s}\oplus C^{n-r-s}\oplus D^{n+1-r-s} 
\\&&\quad\ra (D\cup_{h'} D)_r=D_r\oplus C_{r-1}\oplus D_r\\
\end{eqnarray*}
\end{lem}

\begin{lem}
\label{homsymtwdblem}
Let $c=(f \colon C \ra D, (\delta\phi,\phi))$ be an $n$-dimensional $\epsilon$-symmetric Poincar\'e pair.
Let $(\Delta,\eta)\colon (h,\chi) \simeq (h',\chi') \colon (C,\phi) \homra (C,\phi)$ 
be a homotopy of self-equivalen\-ces.
Then there is an isomorphism $(a,\sigma)\colon c\cup_{(h,\chi)}-c \isora c\cup_{(h',\chi')}-c$
of the corresponding symmetric twisted doubles given by
\begin{eqnarray*}
a_r&=& \mat{1&0&(-)^r f\Delta\\0&1&0\\0&0&1}\colon 
\\&&(D\cup_h D)_r=D_r\oplus C_{r-1}\oplus D_r \ra (D\cup_{h'} D)_r=D_r\oplus C_{r-1}\oplus D_r\\
\sigma_s &=& \mat{  (-)^{n-1} f\eta_s f^*   &0&0\\
            0           &0&0\\
            (-)^n\phi_s\Delta^* f^* &0&0}\colon
\\&&(D\cup_{h'} D)^{n+1-r+s}=D^{n+1-r+s}\oplus C^{n-r+s}\oplus D^{n+1-r+s} 
\\&&\quad\ra (D\cup_{h'} D)_r=D_r\oplus C_{r-1}\oplus D_r\\
\end{eqnarray*}
\end{lem}

\begin{lem}
\label{asysgnhomlem}
Let $c=(f \colon C \ra D, (\delta\phi,\phi))$ be an $n$-dimensional $\epsilon$-symmetric Poin\-car\'e pair.
Let $(\Delta,\eta)\colon (h,\chi) \simeq (h',\chi') \colon (C,\phi) \homra (C,\phi)$ 
be a homotopy of self-equival\-en\-ces.
Then there is an isomorphism 
$$(a,\sigma)\colon (\partial E, \theta)=c\cup_{(h,\chi)}-c \isora  (\partial E', \theta')=c\cup_{(h',\chi')}-c$$
of the corresponding symmetric twisted doubles given in Lemma \ref{homsymtwdblem}.
Let  $x'=(g'    \colon\partial E'\ra E, (\delta\theta', \theta'))$ be an $(n+1)$-dimensional $\epsilon$-symmetric Poincar\'e pair.

Then $x =(g=g'a \colon\partial E \ra E, (\delta\theta=\delta\theta'+(-)^n g'\sigma{g'}^*, \theta))$ is also 
an $(n+1)$-dimensional $\epsilon$-symmetric Poincar\'e pair (by Lemma \ref{deltapairlem}) and
the asymmetric signatures $\sigma^*(x)=\sigma^*(x')\in LAsy^0(\Lambda)$.
\end{lem}
\begin{proof}
Use Definition \ref{cplxdef} and show that the b-duality maps of $x$ and $x'$ are homotopic.
Then it follows that the asymmetric complexes of $x$ and $x'$ are equivalent.
\end{proof}

\section{Flip-Isomorphisms rel$\partial$}
\label{quadfliprelsec}

As explained in the introduction we can only produce a quadratic signature if $(h_t,\chi_t) \simeq (1,0)$.
Inconveniently, $(h_t,\chi_t)$ depends on the choice of representatives for
$\bar\theta$ and $\bar\nu$ and a map $\kappa$. This is the reason for the next, rather awkward,
definition.

\begin{defi}
\label{flipreldef}\index{Flip-isomorphism!rel$\partial$}
A {\bf flip-isomorphism $t$ rel$\partial$ of a regular $\epsilon$-quadratic split preformation
$z=\kreckfs$} is a flip-isomorphism $t=(\alpha,\beta,\bar\nu)$ of $z$ such that 
there is a representative $\theta$ of $\bar\theta\in Q_{-\epsilon}(G)$ and
a representative $\nu$ for $\bar\nu \in Q_{-\epsilon}(F^*)$
and $\kappa\in \Hom_\Lambda(G,G^*)$ such that 
$\beta^*\theta\beta+\theta+\mu^*\nu\mu=\kappa+\epsilon\kappa^*$
and such that the isomorphism
$(h_t,\chi_t) \colon (C,\psi) \isora (C,\psi)$ defined in (\ref{htdefeqn}) and (\ref{defceqn}) 
is homotopic to the identity $(1,0)\colon (C,\psi) \isora (C,\psi)$.
\hfill\qed\end{defi}

\begin{rem}
\label{welldef3rem}
To a certain extent the ``rel$\partial$''-property is independent of the choices
of $\theta$ and $\nu$ (but not necessarily of the choice of $\kappa$):
Assume that for $\nu$, $\theta$ and $\kappa$ as before, there exists a homotopy 
$$(\Delta, \eta) \colon (1,0) \simeq (h_t,\chi_t) \colon (C,\psi) \isora (C,\psi)$$
Let $\widehat\nu$, $\widehat\theta$, $\widehat\kappa$, $\tilde\nu$, $\tilde\theta$, $\tilde\kappa$
as in Remark \ref{welldef2rem}.
They induce a new isomorphism $$(h_t, \widehat\chi_t) \colon (C,\widehat\psi)\\ \isora (C,\widehat\psi)$$
which is also homotopic to the identity by
$$(\Delta, \eta+\tilde\chi) \colon (1,0) \simeq (h_t, \widehat\chi_t) \colon (C,\widehat\psi) \isora (C,\widehat\psi)$$
with $\tilde\chi$ defined in Remark \ref{welldef2rem}.
See also Remarks \ref{welldef1rem}, \ref{welldef2rem} and Lemma \ref{quadamblem}.
\end{rem}

The ``rel$\partial$''-property is invariant under the equivalence relations of $fl_{2q+2}(\Lambda)$
and any elementary preformation has such a flip-isomorphism.
\begin{prop}
\label{fliprelbprp}
\begin{enumerate}
\item \label{fliprelbit} Every elementary preformation has a flip-isomorph\-ism rel$\partial$.
\item Let $t$ and $t'$ be flip-isomorphisms of $\epsilon$-quadratic split preformations $z$ and $z'$ respectively.
If $t$ and $t'$ are flip-isomorphisms rel$\partial$ then so is $t\oplus t'$.
\item Let $[(z,t)]=[(z',t')]\in fl_{2q+2}(\Lambda)$. If $t$ is a flip-isomorphism rel$\partial$ then so is $t'$.
\end{enumerate}
\end{prop}

\begin{proof}
\begin{enumerate}
\item Let $z$ be of the form describe in Proposition \ref{elemprop}\ref{elemprop3}.
Then the flip-isomorphism defined in Proposition \ref{elemflipprp}
is a flip-isomorphism rel$\partial$ with representatives
\begin{eqnarray*}
\theta &=& \mat{0&-\epsilon\sigma\\0&\theta'} \colon G=U\oplus R \ra G^*=U^*\oplus R^*\\
\nu &=&0\\
\kappa  &=& \mat{0&0\\0&\theta'}\colon G=U\oplus R \ra G^*=U^*\oplus R^*
\end{eqnarray*}
and with a homotopy 
$(\Delta, \eta) \colon (1,0) \simeq (h_t,\chi_t) \colon (C,\psi) \isora (C,\psi)$
given by
\begin{eqnarray*}
\Delta_{q+1}=   \mat{\epsilon & 0 & 0 & \epsilon\\
             0 & 0 & 0 & 0}
        &\colon& C_q=(U\oplus U^*) \oplus (U^*\oplus U) \ra C_{q+1}=U\oplus R
\\
\Delta_q=       \mat{0 & 0 \\ -1 & 0 \\ -\epsilon & 0 \\ 0 & 0}
        &\colon& C_{q-1}=U^*\oplus R^* \ra C_q=(U\oplus U^*) \oplus (U^*\oplus U)
\\
\eta_1= \mat{0 & 0 & 1 & 0\\
             0 & 0 & 0 & 0} 
        &\colon& C^q = (U^*\oplus U)\oplus(U\oplus U^*) \ra C_{q+1}=U\oplus R
\\
\eta_1= \mat{0 & 0\\\epsilon & 0\\1& 0\\0&0}
        &\colon& C^{q+1}=U^*\oplus R^* \ra C_q=(U\oplus U^*) \oplus (U^*\oplus U)
\\
\eta_2= \mat{1&0\\0&0}
        &\colon& C^{q-1}=U\oplus R \ra C_{q+1}=U\oplus R
\\
\eta_2= \mat{0&0&0&0\\
             0&0&0&0\\
             1&0&0&0\\
             0&1&\epsilon&0}
        &\colon& C^q = (U^*\oplus U)\oplus(U\oplus U^*)\\&&\quad\quad \ra C_q=(U\oplus U^*) \oplus (U^*\oplus U)
\\
\eta_3= \mat{-\epsilon&0\\0&0\\0&0\\0&0}
        &\colon& C^{q-1}=U\oplus R \ra C_q=(U\oplus U^*) \oplus (U^*\oplus U)
\\
\eta_4=     \mat{0&0\\-\sigma^*&0}
        &\colon& C^{q-1}=U\oplus R \ra C_{q-1}=U^*\oplus R^*
\end{eqnarray*}

\item Obvious.

\item
We only sketch the proof.
Let $t$ and $t'$ be flip-isomorphisms of $\epsilon$-quadratic split preformations $z$ and $z'$ respectively.
$t$ and $t'$ induce self-equivalences $(h_t,\chi_t)\colon (C,\psi) \isora (C,\psi)$ and
$(h'_t,\chi'_t)\colon (C',\psi') \isora (C',\psi')$.
Let $t$ and $t'$ be flip-isomorphisms of $\epsilon$-quadratic split preformations $z$ and $z'$ respectively.
Let $(\eta,\zeta)$ be an isomorphism between $(z,t)$ and $(z',t')$ in the sense of Definition \ref{flmonoidef}.
Then define an isomorphism $(\tilde h,0) \colon (C,\psi) \ra (C',\psi')$ by
\begin{eqnarray*}
\tilde h_{q+1}  &=& \zeta   \colon C_{q+1}=G \isora C'_{q+1}=G'\\
\tilde h_{q}    &=& \mat{\eta&0\\0&\eta^{-*}}   \colon C_{q}=F\oplus F^* \isora C'_{q}=F'\oplus{F'}^*\\
\tilde h_{q-1}  &=& \zeta^{-*}  \colon C_{q-1}=G \isora C'_{q-1}=G'
\end{eqnarray*}
Then $\tilde h h_t \tilde h^{-1}=h'_t \colon C' \isora C'$ and $\tilde h \chi_t \tilde h^*=\chi'_t$.
Assume $$(\Delta, \eta) \colon (1,0) \simeq (h_t,\chi_t) \colon (C,\psi) \ra (C,\psi)$$ is a homotopy,
then $$(\tilde h\Delta\tilde h^{-1}, \tilde h \eta \tilde h^*)\colon (1,0)\simeq (h'_t,\chi'_t) \colon (C',\psi') \ra (C',\psi')$$
is a homotopy as well. 

Assume now that there is a $k\in\bN_0$ such that the flip-isomorphism $t+t_k$ of $z+y_k$ is 
a flip-isomorphism rel$\partial$ ($t_k$ and $y_k$ are defined in Definition \ref{flmonoidef}.)
Now observe that the $2q$-dimensional quadratic complex $C$ induced by $y_k$ is contractible.
It follows that $t$ is a flip-isomorphism rel$\partial$.
\end{enumerate}
\end{proof}

\section{Construction of the Quadratic Signature}
\label{quadconstrsec}

Let $t=(\alpha,\beta,\bar\nu)$ be a flip-isomorphism rel$\partial$ of a regular $\epsilon$-quadratic split preformation 
$z=\kreckfs$. Choose $\theta$, $\nu$ and $\kappa$ as in Definition \ref{flipreldef}
Then there exists a homotopy $(\Delta,\eta) \colon (1,0) \simeq (h_t,\chi_t) \colon (C,\psi) \isora (C,\psi)$.
Write $\Delta_{q+1}=\mat{R&S}\colon C_q=F\oplus F^* \ra C_{q+1}=G$ and 
$\Delta_q=\mat{U\\V}\colon C_{q-1}=G^* \ra C_q=F\oplus F^*$.

We use the homotopy to change the boundary of 
$x_t= (g_t \colon \partial E_t=D\cup_{h_t} D \ra E, (\delta\omega_t, \omega_t))$
to the simpler quadratic Poincar\'e complex
$$(\partial E',\omega')= c \cup -c=(D\cup_C D, \delta\psi\cup_\psi -\delta\psi)=(D\cup_{1} D, \delta\psi\cup_0 -\delta\psi)$$
by applying Lemma \ref{homqutwdblem}. We obtain an isomorphism
$(a,\sigma) \colon (\partial E',\omega') \isora (\partial E_t,\omega_t)$ given by
($\partial E'_r=D_r\oplus C_{r-1} \oplus D_r$, $\partial E_{t,r}=D_r\oplus D_r\oplus C_{r-1}$)
\begin{eqnarray*}
\xymatrix@+30pt
{
\partial E'_{q+2}=0\oplus G \oplus 0            \ar[d]_-{-\epsilon\svec{1\\\gamma\\\mu\\1}}
                            \ar[r]_-{1}
&
\partial E_{t,q+2}=0\oplus 0\oplus G            \ar[d]^-{-\epsilon\svec{\beta\\1\\\gamma\\\mu}}
\\
\partial E'_{q+1}=G\oplus (F\oplus F^*) \oplus G    \ar[d]_-{\svec{ -\epsilon\mu    & 0 & \epsilon      & 0\\
                                    0       & \mu^* & \epsilon\gamma^*  & 0\\
                                    0       & 0 & \epsilon      &-\epsilon\mu}}
                            \ar[r]_-{\svec{ 1   & -\epsilon R   & -\epsilon S   & 0\\
                                    0   & 0     & 0     & 1\\
                                    0   & 1     & 0     & 0\\
                                    0   & 0     & 1     & 0}}
&
\partial E_{t,q+1}=G\oplus G\oplus (F\oplus F^*)    \ar[d]^-{\svec{ -\epsilon\mu    & 0     & \alpha    & \alpha(\nu^*-\epsilon\nu)\\
                                    0       & -\epsilon\mu  & 0     & \epsilon\\
                                    0       & 0     & \mu^*     & \epsilon\gamma^*}}
\\
\partial E'_q= F^* \oplus G^* \oplus F^*        \ar[r]_-{\svec{ 1 & \epsilon V  & 0\\
                                    0 & 0       & 1\\
                                    0 & 1       & 0}}
&
\partial E_{t,q}=F^*\oplus F^*\oplus G^*
}
\end{eqnarray*}

\begin{eqnarray*}
\omega'_0 = \mat{   0 & 0       & 0\\
            0 & -\gamma & 0\\
            0 & 0       & 0\\
            0 & -1      &0} &\colon& {\partial E'}^q= F \oplus G \oplus F 
\\&&\quad\ra \partial E'_{q+1}=G\oplus (F\oplus F^*) \oplus G 
\\
\omega'_0 = \mat{   0 & 0 & 0 & 0\\
            0 & 0 & 0 & 0\\
            0 & -1& 0 & 0} &\colon& {\partial E'}^{q+1}=G^*\oplus F^*\oplus F\oplus G^* 
\\
&&\quad\ra \partial E'_{q}= F^* \oplus G^* \oplus F^*
\end{eqnarray*}
\begin{eqnarray*}
\omega'_1 = \mat{   0 & 0           & 0\\
            0 & -\epsilon\theta & 0\\
            0 & 0           & 0} &\colon& {\partial E'}^q = F \oplus G \oplus F \ra \partial E'_{q}= F^* \oplus G^* \oplus F^*
\\
\sigma_0 = \mat{    -V^*    & 0 & 0}
    &\colon& \partial E_t^q=F\oplus F\oplus G \ra \partial E_{t,q+2}= G
\\
\sigma_0 = \mat{    \eta_0  & 0 & 0 & 0\\
            0   & 0 & 0 & 0\\
            0   & 0 & 0 & 0\\
            -R^*    & 0 & 0 & 0}
    &\colon& \partial E_t^{q+1}=G^*\oplus G^*\oplus (F^*\oplus F) 
\\&&\quad\ra \partial E_{t,q}=F^*\oplus F^*\oplus G^*
\\
\sigma_1 = \mat{    \eta_1  & 0 & 0\\
            0   & 0 & 0\\
            0   & 0 & 0\\
            0   & 0 & 0}
    &\colon& \partial E_t^q=F\oplus F\oplus G 
\\
&&\quad\ra \partial E_{t,q+1}=G\oplus G\oplus (F\oplus F^*)
\\
\sigma_1 = \mat{    \eta_1          & 0 & 0 & 0\\
            0           & 0 & 0     & 0\\
            \epsilon\gamma^* R^*    & 0 & 0 & 0}
    &\colon& \partial E_t^{q+1}=G^*\oplus G^*\oplus (F^*\oplus F) 
\\&&\quad\ra \partial E_{t,q}=F^*\oplus F^*\oplus G^*
\\
\sigma_2 = \mat{    \eta_2          & 0 & 0\\
            0           & 0 & 0\\
            \epsilon\theta^* V^*    & 0 & 0}
    &\colon& \partial E_t^q=F\oplus F\oplus G 
\\
&&\quad\ra \partial E_{t,q}=F^*\oplus F^*\oplus G^*
\end{eqnarray*}

Then we change the boundary of the $(2q+2)$-dimensional quadratic Poincar\'e complex
$$x_t=(g_t \colon \partial E_t \ra E, (\delta\omega_t, \omega_t))$$ 
from Section \ref{flipqutwdblsec} using Lemma \ref{deltapairlem} and the above isomorphism.
We obtain
\begin{eqnarray}
\label{wtdefeqn}\\\nonumber
w_t         &=& (g'_t=g_t a \colon \partial E' \ra E, (\delta\omega'=\delta\omega_t+(-)^{2q+1}g_t\sigma g_t^*, \omega'))\\\nonumber
g'_{t,q+1}      &=& \mat{1 & -\epsilon R & -\epsilon S & -\beta} \colon \partial E'_{q+1}=G\oplus (F\oplus F^*) \oplus G  \ra E_{q+1}=G\\\nonumber
\delta\omega'_{t,0}     &=& -\eta_0 \colon E^{q+1}=G^* \ra E_{q+1}=G
\end{eqnarray}

The next step is to glue $w_t$ to an algebraic $h$-cobordism $y=(m\colon \partial E'=D\cup_C D \ra D, (0, \omega'))$
with $m_r=\mat{-1&0&1} \colon \partial E'_r=D_r\oplus C_{r-1}\oplus D_r \ra D_r$. 
Let the result be the $(2q+2)$-dimensional quadratic Poincar\'e complex $(V,\tau)=w_t \cup -y$.
Using the isomorphism
($V_r=E_r\oplus \partial E'_{r-1}\oplus D_r$)
\begin{eqnarray*}
\xymatrix@+30pt
{
V_{q+3}=0\oplus G \oplus 0              \ar[d]_-{-\epsilon\svec{1\\\gamma\\\mu\\1}}
                            \ar[r]^-{-\epsilon}
&G                          \ar[d]^-{\svec{ 1\\0\\0\\0}}
\\
V_{q+2}=0\oplus (G\oplus F\oplus F^* \oplus G) \oplus 0 \ar[d]_-{\svec{ -\epsilon   & R & S         & \epsilon\beta\\
                                    -\epsilon\mu    & 0 & \epsilon      & 0\\
                                    0       & \mu^* & \epsilon\gamma^*  & 0\\
                                    0       & 0 & \epsilon      &-\epsilon\mu\\
                                    \epsilon    & 0 & 0         &-\epsilon}}
                            \ar[r]^-{\svec{ 1   & 0 & 0 & 0\\
                                    -\gamma & 1 & 0 & 0\\
                                    -\mu    & 0 & 1 & 0\\
                                    -1  & 0 & 0 & 1}}
&G\oplus F\oplus F^* \oplus G               \ar[d]^-{\svec{ 0 & R     & 0 & 0\\
                                    0 & 0     & 1 & 0\\
                                    0 & \mu^* & 0 & 0\\
                                    0 & 0     & 0 & 0\\
                                    0 & 0     & 0 & 1}}
\\
V_{q+1}=G\oplus (F^*\oplus G^*\oplus F^*)\oplus G   \ar[d]_-{\svec {0 & -\epsilon  & 0 & \epsilon &-\epsilon\mu}}
                            \ar[r]^-{\svec {1 & -\epsilon S & 0 & 0 & \beta\\
                                    0 & \epsilon    & 0 & 0 & 0\\
                                    0 & -\gamma^*   & 1 & 0 & 0\\
                                    0 & -1      & 0 & 1 & -\mu\\
                                    0 & 0       & 0 & 0 & -\epsilon}}
&G\oplus F^*\oplus G^*\oplus F^*\oplus G        \ar[d]^-{\svec {0 & 0 & 0 & 1 & 0}}
\\
V_q    =0\oplus 0 \oplus F^*                \ar[r]^-{\epsilon}
&F^*
}
\end{eqnarray*}
we can shrink the chain complex $V$ to a smaller chain complex $V'$ via the chain equivalence given by
\begin{eqnarray*}
\xymatrix@C+30pt
{
V_{q+2}=G\oplus F\oplus F^* \oplus G \ar[r]_-{\svec{    -\gamma & 1 & 0 & 0}} \ar[d]
& V'_{q+2}=F                        \ar[d]^-{\svec{R\\\mu^*}} 
\\
V_{q+1}=G\oplus F^*\oplus G^*\oplus F^*\oplus G \ar[r]_-{\svec{1 & -\epsilon S  & 0 & 0 & \beta\\0 & -\gamma^*  & 1 & 0 & 0}}
& V'_{q+1}=G\oplus G^*
}
\end{eqnarray*}

The induced quadratic structure on $V'$ is given by
\begin{eqnarray*}
\tau'_0&=&\mat{-\eta_0 & \beta \\ 0 & -\epsilon\theta^*} \colon {V'}^{q+1}=G^*\oplus G \ra V'_{q+1}=G\oplus G^*
\end{eqnarray*}

Using Lemma \ref{trcplxlem} we can define the quadratic signature
\begin{defi}
\label{quadsgndef}
Let $z=\kreckfs$ be a regular $\epsilon$-quadratic split preformation and $t=(\alpha, \beta, \bar\nu)$
a flip-isomorphism rel$\partial$ of $z$.

Choose representatives $\theta$ of $\bar\theta\in Q_{-\epsilon}(G)$ and 
a representative $\nu$ for $\bar\nu \in Q_{-\epsilon}(F^*)$
and a $\kappa\in \Hom_\Lambda(G,G^*)$ such that 
$\beta^*\theta\beta+\theta+\mu^*\nu\mu=\kappa+\epsilon\kappa^*$
and such that there is a homotopy
$(\Delta,\eta)\colon (h_t,\chi_t)\simeq (1,0) \colon (C,\psi) \isora (C,\psi)$ 
of the isomorphism defined in (\ref{htdefeqn}) and (\ref{defceqn}). 
The {\bf quadratic signature $\rho^*(z,t,\nu,\theta,\kappa,\Delta,\eta)=[(M,\xi)] \in L_{2q+2}(\Lambda)$} is given
\index{Quadratic signature} by the non-singular $(-\epsilon)$-quadratic form
\begin{eqnarray*}
\xi&=&\mat{-\eta_0 & \beta & 0 \\ 0 & -\epsilon\theta^* & 0 \\ R^* & \mu & 0}\colon M=G^*\oplus G \oplus F \ra M^*
\end{eqnarray*}
\hfill\qed\end{defi}

\section{Properties of the Quadratic Signatures}
\label{quadelemsec}

First we show that the quadratic signatures detect whether an element
in $\el$ is elementary or not.

\begin{thm}
\label{quadelemthm}
$[z']\in l_{2q+2}(\Lambda)$ is elementary if and only if 
there is a flip-isomorph\-ism rel$\partial$ $t$ of $z$
with $[z]=[z']\in\el$ and
$\Delta$, $\kappa$, $\theta$, $\nu$ and $\eta$ as in Definition \ref{quadsgndef} such that
$\rho^*(z,t,\nu,\theta,\kappa,\Delta,\eta)=0\in L_{2q+2}(\Lambda)$.
\end{thm}

\begin{proof}
Let $z$ be a regular $\epsilon$-quadratic split preformation.
If $z$ is elementary then
use the flip-isomorphism and the choices for $\Delta$, $\nu$, \etc 
made in the proof of Proposition \ref{fliprelbprp}\ref{fliprelbit}.

On the other hand assume that there is $t$, $\Delta$, \etc such that 
$\rho^*(z,t,\nu,\theta,\kappa,\Delta,\eta)=0\in L_{2q+2}(\Lambda)$. Then $(V,\tau)= y\cup -w_t$ constructed
in the previous Section is null-cobordant.
Hence the Poincar\'e pairs $y$ and $w_t$ from the previous section are cobordant rel$\partial$.
$y$ is an algebraic $h$-cobordism.  The Poincar\'e pair $x$
constructed in Section \ref{flipqutwdblsec} is cobordant rel$\partial$
to an algebraic $h$-cobordism. By Theorem \ref{algmainthm}
$[z]\in l_{2q+2}(\Lambda)$ is elementary.
\end{proof}

Quadratic signature depend on a lot of choices. We can restrict that dependency to
a certain degree: 
\begin{lem}
\label{quadamblem}
Let $z$, $t$, $\nu$, $\theta$, $\kappa$, $\Delta$ and $\eta$ as in Definition \ref{quadsgndef}.

Let $\tilde\theta \in\Hom_\Lambda(G,G^*)$, $\tilde\nu\in\Hom_\Lambda(F^*, F)$,
$\tilde\kappa \in \Hom_\Lambda(G,G^*)$ and $\tilde\eta\in W_\%(C,\epsilon)_{2q+3}$.

Define
\begin{eqnarray*}
\widehat\nu     &=& \nu   +\tilde\nu   +\epsilon\tilde\nu^*\in\Hom_\Lambda(F^*, F)\\
\widehat\theta  &=& \theta+\tilde\theta+\epsilon\tilde\theta^*\in\Hom_\Lambda(G,G^*)\\
\widehat\kappa  &=& \kappa+\tilde\kappa-\epsilon\tilde\kappa^* + \beta^*\tilde\theta\beta + \tilde\theta + \mu^*\tilde\nu\mu \in\Hom_\Lambda(G,G^*)\\
\widehat\eta    &=& \eta  +\tilde\chi + d(\tilde\eta)
\end{eqnarray*}
with $\tilde\chi$ defined in Remark \ref{welldef2rem}.

Then $\rho^*(z,t,\widehat\nu,\widehat\theta,\widehat\kappa,\Delta,\widehat\eta)$ exists and is equal
to  $\rho^*(z,t,\nu,\theta,\kappa,\Delta,\eta)$.
\end{lem}
\begin{proof}
Straightforward. See also Remarks \ref{welldef1rem}, \ref{welldef2rem} and \ref{welldef3rem}.  
\end{proof}

\begin{lem}
\label{quadaddlem}
Let $z$, $t$, $\nu$, $\theta$, $\kappa$, $\Delta$ and $\eta$ as in Definition \ref{quadsgndef}.
Let $z'$, $t'$, $\nu'$, $\theta'$, $\kappa'$, $\Delta'$ and $\eta'$ another set of data
as in Definition \ref{quadsgndef}.
Then
\begin{eqnarray*}
\rho^*(z,t,\nu,\theta,\kappa,\Delta,\eta) &+& \rho^*(z',t',\nu',\theta',\kappa',\Delta',\eta')\\
&=&\rho^*(z\oplus z', t\oplus t',\nu\oplus\nu',\theta\oplus\theta',\kappa\oplus\kappa',\Delta\oplus\Delta',\eta\oplus\eta')
\end{eqnarray*}
\end{lem}

In a special case the quadratic signatures only depend on $z$, $t$ and $\Delta$.
Later it will be shown that the quadratic signature does only depend on
the preformation and flip-isomorphism rel$\partial$ if
$\Lambda=\bZ$ and $\epsilon=-1$ (see Corollary \ref{coolintcor}).
\begin{lem}
\label{quaddeplem}
Let $z$, $t$, $\nu$, $\theta$, $\kappa$, $\Delta$ and $\eta$ as in Definition \ref{quadsgndef}.
Let $z$, $t$, $\nu'$, $\theta'$, $\kappa'$, $\Delta$ and $\eta'$ another set of data
as in Definition \ref{quadsgndef}.
Assume that the map $Q_{-\epsilon}(G^*) \ra Q_{-\epsilon}(F)$, $\theta \mt \mu\theta\mu^*$ is injective.
Then $\rho^*(z,t,\nu,\theta,\kappa,\Delta,\eta)=\rho^*(z,t,\nu',\theta',\kappa',\Delta,\eta')$.
Hence the quadratic signature only depends on $z$, $t$ and $\Delta$.

If $\Lambda=\bZ$ this is the case if $\mu$ is injective and, for $\epsilon=1$, if additionally 
$|\Tor\coker\mu|$ is odd.
\end{lem}
\begin{proof}
Using Lemma \ref{quadamblem} we can reduce the problem to the case where $\theta=\theta'$,
$\nu=\nu'$. Then we have two homotopies
\begin{eqnarray*}
(\Delta,\eta)&\colon& (1,0) \simeq (h_t,\chi_t) \colon (C,\psi) \ra (C,\psi)\\
(\Delta,\eta')&\colon& (1,0) \simeq (h_t,\chi'_t) \colon (C,\psi) \ra (C,\psi)
\end{eqnarray*}
Hence $\chi'_t-\chi_t=d(\eta'-\eta)$. It follows that $d(\tilde\eta)_0=0\colon C^{q+1}\ra C_q$ and $d(\tilde\eta)_1=0\colon C^q\ra C_q$
for $\tilde\eta=\eta'-\eta$. Combining the equations 
yields $d\tilde\eta_0d^*=(1+T_\epsilon)(-\epsilon d\tilde\eta_1 -\epsilon\tilde\eta_2) \colon C^q\ra C_q$. Hence
$\mu\tilde\eta_0\mu^*=0\in Q_{-\epsilon}(F)$. Hence $\eta'=\eta\in Q_{-\epsilon}(G)$ and the claim follows.

Now let $\Lambda=\bZ$. We need to prove that $\mu$ induces an injection between the $Q$-groups.
By the Smith Normal Form Theorem (\cite{Newman72})
we can assume that 
$$\mu=\mat{d_1 & 0\cdots 0 & 0 \\ 0 & \ddots & 0 \\ 0 & 0\cdots 0 & d_n \\ 0 & 0\cdots 0 & 0} \colon G=\bZ^n \ra F^*=\bZ^m$$
with all $d_i > 0$ and $m\ge n$. By assumption, all $d_i$ are odd if $\epsilon=1$.
Let $\theta\in Mat(n\times n, \bZ)$ such that $\mu\theta\mu^*=\kappa+\epsilon\kappa^*$
for some $\kappa\in Mat(m\times m, \bZ)$. It follows that $\theta_{ij}=\epsilon\theta_{ji}$ and
$d_i^2\theta_{ii}=(1+\epsilon)\kappa_{ii}$ for all $i,j\in\{1,\dots,n\}$.
If $\epsilon=-1$ then $\theta_{ii}=0$ for all $i$. Hence $
\kappa'_{ij}=\begin{cases}
         \theta_{ij} & \text{if $i < j$}\\
             0       & \text{else}
         \end{cases}$
fulfils $\theta=\kappa'+\epsilon{\kappa'}^*$ and $[\theta]=0\in Q_{-\epsilon}(G^*)$.

If $\epsilon=1$ then $d_i^2\theta_{ii}=2\kappa_{ii}$. By assumption $d_i^2 | \kappa_{ii}$
and hence we can define $\kappa'\in Mat(n\times n, \bZ)$ by 
$\kappa'_{ij} =
\begin{cases}
\theta_{ij}         & \text{if $i < j$}\\
\kappa_{ii}/d_{ii}^2    & \text{if $i=j$}\\
0           & \text{else}
\end{cases}$. It follows that $\theta=\kappa'+\epsilon{\kappa'}^*$ and $[\theta]=0\in Q_{-\epsilon}(G^*)$.
\end{proof}

\section{Quadratic and Asymmetric Signatures}
\label{quadasymsec}

There is a close relationship of asymmetric and quadratic signatures which is not that surprising
if we re-examine their construction. For simplicity let $(W,M,M)$ be a $(2q+2)$-dimensional cobordism
with $\partial M=\emptyset$. By glueing $W$ together along $M$ and after surgery below the middle dimension
we obtain a $(2q+2)$-dimensional closed manifold $V$. Poincar\'e duality induces
a non-singular $(-)^{q}$-symmetric form $\lambda\colon H^{q+1}(V) \isora H_{q+1}(V)$. The asymmetric 
signature of $(W,M,M)$ is the image of that form in $LAsy^0(\Lambda)$ 
(compare Section \ref{asymgeosec}).

Similarly, in the algebraic surgery world we can glue an $(2q+2)$-dimensional quadratic Poincar\'e pair
$d=(D\oplus D \ra E, (\delta\nu, \nu\oplus-\nu))$ together along $D$ and obtain a $(2q+2)$-dimensional
quadratic Poincar\'e pair $(V,\tau)$ (which is basically the quadratic signature in this simple 
situation). By Lemma \ref{asymcor} the image of its symmetrization $(V, (1+T)\tau)$
in $LAsy^0(\Lambda)$ is the asymmetric signature of $d$.
This fact generalizes to the case of all quadratic signatures.

\begin{thm}
\label{asyquadthm}
Let $z$, $t$, $\nu$, $\theta$, $\kappa$, $\Delta$ and $\eta$ as in Definition \ref{quadsgndef}.
The image of the quadratic signature $\rho^*(z,t,\nu,\theta,\kappa, \Delta,\eta)$ under the map
\begin{eqnarray*}
L_{2q+2}(\Lambda) &\ra& LAsy^0(\Lambda)\\
(K,\psi) &\mt& (K,\psi_0-\epsilon\psi^*)
\end{eqnarray*}
is the asymmetric signature $\sigma^*(z,t)$.
\end{thm}
\begin{proof}
Let $(V,\tau)=w_t \cup -y$ be the union of $(2q+2)$-dimensional quadratic Poincar\'e pairs
defined in Section \ref{quadconstrsec}. (Using standard algebraic surgery theory - \eg Lemma \ref{trcplxlem} -
$(V,\tau)$ is a $(2q+2)$-dimensional quadratic Poincar\'e complex representing the quadratic signature
$\rho^*(z,t,\nu,\theta,\kappa, \Delta,\eta)\in L_{2q+2}(\Lambda)$).

By Lemma \ref{unionsymlem}, $(V,(1+T)\tau)\cong(1+T)w_t \cup - (1+T)y$.
Again by Lemma \ref{unionsymlem} there is an isomorphism 
$$(1,\tilde\sigma)\colon(\partial E',\theta')=(1+T)c\cup-(1+T)c \isora (\partial E', (1+T)\omega')=(1+T)(c\cup-c)$$
with $c=(f\colon C \ra D, (\delta\psi,\psi))$ (compare (\ref{defsymceqn}) and (\ref{defceqn}))
which - applied to $(1+T)w_t$ - yields a new $(2q+2)$-dimensional symmetric Poincar\'e pair
\begin{eqnarray*}
w^t &=& (g'_t \colon \partial E' \ra E, (\delta\theta'=(1+T)\delta\omega'+(-)^{2q-1}g'_t\tilde\sigma {g'}^*_t,\theta'))\\
\delta\theta'_0 &=& -((1+T_\epsilon)\eta_0+\epsilon\Delta_{q+1}T\psi_0\Delta_{q+1}^*)\colon E^{q+1}=G^* \ra E_{q+1}=G
\end{eqnarray*}
(we use Lemma \ref{deltapairlem} here).

By Lemma \ref{unionbdrylem} $(1+T)w_t \cup - (1+T)y\cong w^t \cup - y'$ with some Poincar\'e pair $y'$ which is - like $y$ and $(1+T)y$
- an algebraic $h$-cobordism. 

All in all this means that $(V,(1+T)\tau)\cong w^t \cup - y'$.  Note that 
$w^t$ is a Poincar\'e pair with a (trivial) twisted double on its boundary. Hence by Proposition \ref{asymgenprop},
\ref{hcobasymzeroprp} and Corollary \ref{asymcor} 
\begin{eqnarray*}
(1+T_{-\epsilon})\rho^*(z,t,\nu,\theta,\kappa, \Delta,\eta)&=&(V,(1+T)\tau)=\sigma^*(V,(1+T)\tau)\\
&=&\sigma^*(w^t)-\sigma^*(y')=\sigma^*(w^t)\in LAsy^0(\Lambda)
\end{eqnarray*}
Now we verify that $(\Delta, (1+T)\eta+\xi)\colon (h_t,0) \simeq (1,0) \colon (C,\phi)\isora (C,\phi)$ 
is a well-defined homotopy with $\xi$ as in Lemma \ref{chhomlem} which transforms
$w^t$ to $x^t$ using Lemma \ref{homsymtwdblem}.
Then Lemma \ref{asysgnhomlem} shows that $\sigma^*(w^t)=\sigma^*(x^t)=\sigma^*(x,t)\in LAsy^0(\Lambda)$.
\end{proof}

\begin{rem}
\label{asyquadrem}
The relationship between quadratic and asymmetric signatures
can be made more precise. By \cite{Ran98} 30.29 there is an exact sequence
\begin{eqnarray*}
\xymatrix@R-20pt
{
0 \ar[r]
&
DBL_{2q+1}(\Lambda) \ar[r]^{t}
&
L_{2q+2}(\Lambda) \ar[r]^{(1+T_{-\epsilon})}
& 
LAsy^0(\Lambda)\\
&& 
\rho^*(z,t,\nu,\theta,\kappa,\Delta,\eta) \ar@{|->}[r]
&
\sigma^*(z,t)
}
\end{eqnarray*}
with $DBL_{2q+1}(\Lambda)$ the kernel of the map
\begin{eqnarray*}
LAut_{2q+1}(\Lambda) &\ra& L^p_{2q+1}(\Lambda)\\
\left[(h,\chi)\colon (C,\psi) \stackrel{\simeq}{\ra} (C,\psi)\right] &\mt& \left[(C,\psi)\right]
\end{eqnarray*}

The map $t$ is induced by the {\bf algebraic mapping torus}: Let $(h,\chi)\colon (C,\psi) \stackrel{\simeq}{\ra} (C,\psi)$
be a self-equivalence of an $n$-dimensional quadratic Poincar\'e complex. Then the algebraic mapping torus
is the union of the fundamental $(n+1)$-dimensional quadratic Poincar\'e pair
$$((h,1) \colon C\oplus C \ra C, ((-)^n\chi, \psi\oplus-\psi))$$ as described in Definition \ref{unionfunddefi}.

In the case of $\Lambda=\bZ$ and $q=2k$, the map $L_{4k+2}(\bZ) \ra LAsy^0(\bZ)$ factors through $L^{4k+2}(\bZ)=0$
hence the above sequence boils down to
\begin{eqnarray*}
\xymatrix@R-20pt
{
0 \ar[r]
&
DBL_{4k+1}(\bZ)=\bZ/2\bZ \ar[r]^-{\cong}
&
L_{4k+2}(\bZ)=\bZ/2\bZ \ar[r]^-{0}
& 
LAsy^0(\bZ)
}
\end{eqnarray*}

In the case $\Lambda=\bZ$ and $q=2k-1$, the composition $L_{4k}(\bZ)=\bZ \ra LAsy^0(\bZ) \ra LAsy^0(\bC)$
is an injection (because composition with the asymmetric multi-signatures $LAsy^0(\bC)\cong \bZ[S^1]$ 
from \cite{Ran98} Proposition 40.6 and the projection $\bZ[S^1]\ra\bZ$, $\sum_{g\in S^1} n_g \mt n_1$ 
gives the signature on $L_{4k}(\bZ)=\bZ$). The exact sequence becomes
\begin{eqnarray*}
\xymatrix@R-5pt
{
0 \ar[r]
&
DBL_{4k-1}(\bZ)=0 \ar[r]^-{0}
&
L_{4k}(\bZ)=\bZ \ar[r]^-{(1+T)}
& 
LAsy^0(\bZ)
}
\end{eqnarray*}
\end{rem}

\begin{cor}
\label{coolintcor}
Let $q=2m-1$ \ie  $\epsilon=-1$. Let $z=\kreckfs$ be a regular skew-quadratic split preformation over $\bZ$ 
\begin{enumerate}
\item $[z]\in l_{4m}(\bZ)$ is elementary if and only if there is a flip-isomorphism rel$\partial$ $t$ such that $\sigma^*(z,t)=0\in LAsy^0(\bZ)$.
\item The quadratic signature $\rho^*(z,t,\nu,\theta,\kappa,\Delta,\eta)\in L_{4m}(\bZ)$ only depends on $z$ and $t$.
\end{enumerate}
\end{cor}

\chapter{Non-Singular Formations}
\label{nonsingchap}

{\bf Throughout this chapter ``formation'' will mean a non-singular $\epsilon$-quadra\-tic
split formation for $\epsilon=(-)^q$, $q \ge 2$. Let $\Lambda$ be a weakly finite ring with involution and $1$.}

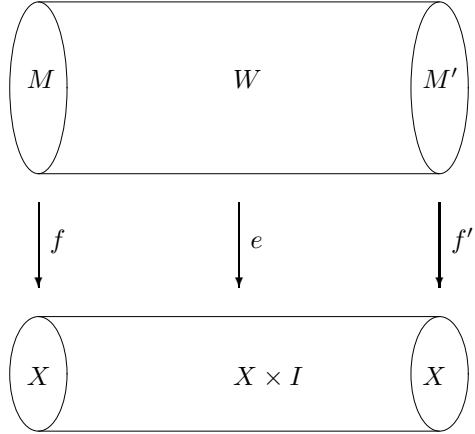
\begin{center}
\setlength{\unitlength}{0.006in}
\begin{picture}(750,400)(0,0)
\put(200,300){\ellipse{50}{150}}
\put(550,300){\ellipse{50}{150}}
\path(200,375)(550,375)
\path(200,225)(550,225)
\put(190,300){$M$}
\put(370,300){$W$}
\put(535,300){${M'}$}

\put(200,200){\vector(0,-1){75}}
\put(375,200){\vector(0,-1){75}}
\put(550,200){\vector(0,-1){75}}
\put(210,160){$f$}
\put(385,160){$e$}
\put(560,160){$f'$}

\put(200,50){\ellipse{50}{100}}
\put(550,50){\ellipse{50}{100}}
\path(200,0)(550,0)
\path(200,100)(550,100)
\put(190,40){$X$}
\put(370,40){$X\times I$}
\put(535,40){$X$}
\end{picture}
\end{center}
Let $(e,f,f')\colon (W, M,M') \ra X\times(I,0,1)$
be a special kind of Kreck surgery situation (\ie $e$, $f$ and $f'$ are highly-connected):
all maps are normal maps and $(X, \partial X)$ is a finite geometric Poincar\'e pair such that $\partial M\homra\partial X$
is a homotopy equivalence.
Such a normal cobordism is called a {\bf presentation} of $f$. Presentations are also used to
define obstructions to odd-dimensional traditional surgery problems (for details see Section \ref{formgeosec}).
Hence the Kreck surgery obstruction $z=\kreckfs$ with $F=K_{q+1}(W,M)$ and $G=K_{q+1}(W)$ 
in this case is a formation \ie $\svec{\gamma\\\mu}\colon G \ra H_{\epsilon}(F)$ is
an inclusion of a lagrangian. 

It is possible to prove much stronger results about formations.
In {\bf Section \ref{nonsingflipsec}} we deal with some useful technical observations about formations and their 
flip-isomorphism.

If one applies the construction of Section \ref{constrpairsec} to a formation
one obtains a $(2q+2)$-dimensional quadratic Poincar\'e pair 
$x=(g \colon \partial E=D'\cup_C D \ra E, (\delta\omega=0, \omega))$ for which $C$ is contractible.

The definition of quadratic signatures was rather awkward because we had to make sure
that a flip-isomorphism induces a self-equivalence $(h_t, \chi_t)$ on $(C,\psi)$ which
is homotopic to the identity. For formations $C$ is contractible and,
hence, every flip-isomorphism is a flip-isomorphism rel$\partial$. This also
leads to a simpler version of the quadratic signatures in {\bf Section \ref{quadsgnsec}}. 

We investigate the behaviour of the asymmetric signatures of
formations in {\bf Section \ref{nonsingasysec}}. They turn out to be
independent from the choice of flip-isomorphism. As an application we
will construct non-elementary preformations for which all asymmetric
signatures are vanishing. We will also show how the asymmetric signatures 
relate to traditional even-dimensional surgery theory \ie
how they behave for boundaries of non-singular forms.

\section{Flip-Isomorphisms}
\label{nonsingflipsec}
We will need to discuss some very technical properties of non-singular formations in order to make the computations 
in the following sections.

A formation $\kreckfs$ is an $\epsilon$-quadratic split preformation such that the map
$$\mat{\gamma\\\mu} \colon G \ra H_\epsilon(F)$$
is an inclusion of a lagrangian. By  \cite{Ran80a} Proposition 2.2., this map can be extended to an isomorphism of 
hyperbolic $\epsilon$-quadratic forms
$$\left(f=\mat{\gamma & \tilde\gamma\\\mu&\tilde\mu},\mat{\theta&0\\\tilde\gamma^*\mu&\tilde\theta}\right)\colon H_\epsilon(G) 
\stackrel{\cong}{\ra} H_\epsilon(F)$$

\begin{rem}
\label{extrem}
For any $\tau\colon G^*\ra G$ the maps $\tilde\gamma'=\tilde\gamma+\gamma (\tau-\epsilon\tau^*)$,
$\tilde\mu'=\tilde\mu+\mu (\tau-\epsilon\tau^*)$, $\tilde\theta'=\tilde\theta+(\tau-\epsilon\tau^*)^*\theta (\tau-\epsilon\tau^*) 
+\tilde\gamma^*\mu (\tau-\epsilon\tau^*)^* -\epsilon\tau$ define another extension to an isomorphism of hyperbolic forms.
Conversely any such extension can be derived from $\tilde\gamma$, $\tilde\mu$, $\tilde\theta$.
\end{rem}

The relationship between those maps and a flip-isomorphism of the formation can be described as follows:
\begin{lem}
Let $\kreckfs$ be a formation and $(\alpha, \beta, \nu)$ a flip isomorphism. Let $f$, $\tilde\gamma$, $\tilde\mu$
and $\tilde\theta$ as explained before.

\begin{enumerate}

\item $f^{-1}=\mat{\tilde\mu^* & \epsilon\tilde\gamma^* \\ \epsilon\mu^* & \gamma^*}$.

\item There is an isomorphism 
\begin{eqnarray*}
&&\left( \mat{\tilde\mu^* & \tilde\gamma^*\\\mu^*&\gamma^*}, \mat{1 & -1 \\ \epsilon\tilde\gamma^*\mu & \tilde\mu^*\gamma}, \mat{0&0\\-1&0}\right)
\colon\\
&&\quad\kreckfs + \kreckfb{F^*}{-\epsilon\mu}{G}{\gamma}{F}{\theta} \stackrel{\cong}{\ra} \partial(G,\theta)+(G^*,G)
\end{eqnarray*}

\item 
\begin{enumerate}
\item $\alpha(\gamma + (\nu -\epsilon\nu^*)^*\mu)=\epsilon\mu\beta$
\item $\alpha^{-*}\mu=\gamma\beta$
\item $\beta^*\theta\beta + \theta + \mu^*\nu\mu = 0 \in Q_{-\epsilon}(G)$
\end{enumerate}

\item There is a $\xi \in Q_{-\epsilon}(G^*)$ and $Y=\xi-\epsilon\xi^*$ such that
\begin{enumerate}
\item $\epsilon\alpha(\tilde\gamma + (\nu -\epsilon\nu^*)^*\tilde\mu)\beta^* = \mu Y+\tilde\mu$
\item $\alpha^{-*}\tilde\mu \beta^* = \gamma Y+\tilde\gamma$
\end{enumerate}
\end{enumerate}
\end{lem}
\begin{proof}
\begin{enumerate}
\item From $f^*H_\epsilon f=H_\epsilon$ it follows that $f^{-1}=H_\epsilon^{-1}f^*H_\epsilon$.
\item Follows straight from the Definition \ref{fisodef}.
\item One can compute the composition of isometries of hyperbolic forms
$$f^{-1} H_\epsilon \mat{\alpha & \alpha(\nu-\epsilon\nu^*)^*\\0&\alpha^{-*}} f \mat{\beta^{-1}&0\\0&\beta^*}=
\mat{1 & Y\\0&1}$$
with 
\begin{eqnarray*}
Y   &=& \tilde\mu^*\alpha^{-*}\tilde\mu\beta^*+\tilde\gamma^*\alpha\tilde\gamma\beta^*+\tilde\gamma^*\alpha(\nu-\epsilon\nu^*)^*\tilde\mu\beta^*\\
\xi     &=& -\beta\tilde\theta^*\beta^*-\beta\tilde\mu^*\nu^*\tilde\mu\beta^*+\epsilon Y^*\theta Y +\epsilon\tilde\gamma^*\mu Y + \epsilon\tilde\theta
\end{eqnarray*}
\end{enumerate}
\end{proof}

\section{Quadratic Signatures}
\label{quadsgnsec}

Let $z=\kreckfs$ be a formation and $t=(\alpha,\beta,\nu)$ be a flip-isomorphism. 
there is a representative $\theta$ of $\bar\theta\in Q_{-\epsilon}(G)$
and a representative $\nu$ for $\bar\nu \in Q_{-\epsilon}(F^*)$
and a $\kappa\in \Hom_\Lambda(G,G^*)$ such that 
$\beta^*\theta\beta+\theta+\mu^*\nu\mu=\kappa+\epsilon\kappa^*$.

As described in the previous section we can extend $\svec{\gamma\\\mu}$ to an
isomorphism $\svec{\gamma&\tilde\gamma\\\mu&\tilde\mu}$ of hyperbolic forms.
A choice of $\tilde\gamma$, $\tilde\mu$, $\tilde\theta$ (compare previous section) 
defines a homotopy $\Delta_C\colon 1 \simeq 0 \colon C \ra C$ with 
\begin{eqnarray*}
\Delta_{C,q+1} &=&  \mat{\epsilon\tilde\mu^* & \tilde\gamma^*} \colon C_q=F\oplus F^*\ra C_{q+1}=G\\
\Delta_{C,q}   &=& -\epsilon\mat{\tilde\gamma\\ \tilde\mu}     \colon C_{q-1}=G^* \ra C_q=F\oplus F^*
\end{eqnarray*}

Then all flip-isomorphisms of a $z$ are flip-isomorphisms rel$\partial$
with a homotopy
\begin{eqnarray*}
(\Delta=\Delta_C(1-h_t), \eta=\Delta_{C,\%}(\chi_t-\Delta_{\%}\psi))\colon (1,0) \simeq (h_t,\chi_t) \colon (C,\psi) \ra (C,\psi)
\end{eqnarray*}

The non-singular $(-\epsilon)$-quadratic form $(M,\xi)$ from the definition of quadratic signatures
(\ref{quadsgndef}) is given by
\begin{eqnarray*}
\xi&=&\mat{-\eta_0 & \beta & 0 \\ 0 & -\epsilon\theta^* & 0 \\ R^* & \mu & 0}\colon M=G^*\oplus G \oplus F^* \ra M^*
\end{eqnarray*}

Using the isomorphism
\begin{eqnarray*}
f&=&\mat{1&0&-\tilde\gamma^*\alpha(\nu^*-\epsilon\nu)-\tilde\mu^*\alpha^{-*}\\0&1&0\\0&0&1} \colon M^* \isora M^*
\end{eqnarray*}
we obtain a prettier non-singular $(-\epsilon)$-quadratic form $(M,\xi'=-f\xi f^*)$
\begin{eqnarray*}
\xi'&=&\mat{    \tilde\gamma^*\tilde\mu+\tilde\gamma^*\alpha\nu\alpha^*\tilde\gamma     & -\tilde\gamma^*\alpha\gamma & 0 \\ 
        0                                   & \epsilon\theta^* & 0 \\ 
        \epsilon(\alpha^*\tilde\gamma-\tilde\mu)                & -\mu & 0}
\colon M=G^*\oplus G \oplus F^* \ra M^*
\end{eqnarray*}

\begin{defi}
\label{quasgndef}\index{Quadratic signature!for non-singular formations}
The {\bf quadratic signature $\tilde\rho^*(z,t,\tilde\gamma,\tilde\mu,\tilde\theta)$} is the element
$(M,\xi')\linebreak \in L_{2q+2}(\Lambda)$.
\hfill\qed\end{defi}

\begin{thm}
\label{formquadthm}
Let $z'$ be a formation.
$[z']\in l_{2q+2}(\Lambda)$ is elementary if and only if there is a stably strongly isomorphic $z=\kreckfs$,
a flip-isomorphism $t$ and $\tilde\gamma$, $\tilde\mu$ and $\tilde\theta$ such that
$$\left(f=\mat{\gamma & \tilde\gamma\\\mu&\tilde\mu},\mat{\theta&0\\\tilde\gamma^*\mu&\tilde\theta}\right)\colon H_\epsilon(G) 
\stackrel{\cong}{\ra} H_\epsilon(F)$$
is an isometry of hyperbolic $\epsilon$-quadratic forms and $\tilde\rho^*(z,t,\tilde\gamma, \tilde\mu, \tilde\theta)=0\in L_{2q+2}(\Lambda)$.
\end{thm}

\begin{proof}
If there exist $t$,$\tilde\gamma$, $\tilde\mu$ and $\tilde\theta$ as above then by construction
$$0=\tilde\rho^*(z,t,\tilde\gamma,\tilde\mu)=-\rho^*(z,t,\Delta,\rho)$$
By Theorem \ref{quadelemthm} $[z]\in l_{2q+2}(\Lambda)$ is elementary.
Assume now that $z=\kreckfs$ is elementary and let it have the special form described in Proposition \ref{elemprop}\ref{elemprop3}.
Clearly $\svec{\sigma\\\tau} \colon  U \ra U\oplus U^*$ is a split injection and even a lagrangian. As in Section \ref{nonsingflipsec}
the map can be extended to an isometry 
$$\left(\mat{\sigma&\tilde\sigma\\\tau&\tilde\tau},\mat{\theta'&0\\\tilde\sigma^*\tau & \tilde\theta'} \right) \colon H_\epsilon(R) \ra H_\epsilon(U^*)$$
of hyperbolic forms. Then the maps
\begin{eqnarray*}
\tilde\gamma &=& \mat{0&0\\0&\tilde\sigma}\colon G^*=U^*\oplus R^* \ra F=U\oplus U^*\\
\tilde\mu    &=& \mat{1&-\epsilon\tilde\sigma\\0&\tilde\tau} \colon G^*=U^*\oplus R^* \ra F^*=U^*\oplus U\\
\tilde\theta &=& \mat{0&0\\0&\tilde\theta'} \colon  G^*=U^*\oplus R^* \ra G=U\oplus R
\end{eqnarray*}
are completing $\svec{\gamma\\\mu}$ to an isometry of hyperbolic forms (compare Section \ref{nonsingflipsec}).

Define a strong flip-isomorphism $t=(\alpha,\beta, \nu=0)$ of $z$ by
\begin{eqnarray*}
\alpha &=& \mat{0&-1\\-\epsilon&0}\colon F=U\oplus U^* \ra F^*=U^*\oplus U\\
\beta  &=& \mat{-1&-\tau\\0&1}\colon    G=U\oplus R \ra  G=U\oplus R
\end{eqnarray*}

Then $\tilde\rho^*(z,t,\tilde\mu,\tilde\theta)$ is 
represented by a non-singular $(-\epsilon)$-quadratic form $(M,\xi')$ (as in Definition \ref{quasgndef})
which has a lagrangian
\begin{eqnarray*}
i=\mat{ 1&0&0\\
    0&1&0\\
    0&0&1\\
    0&0&0\\
    0&0&0\\
    0&\epsilon\tilde\sigma&0}\colon U^*\oplus R^*\oplus U \ra M=U^*\oplus R^*\oplus U\oplus R\oplus U^*\oplus U
\end{eqnarray*}
\end{proof}

\section{Asymmetric Signatures}
\label{nonsingasysec}

The asymmetric signature of formations has one surprising property: it
is independent of the choice of flip-isomorphism (although the
existence of a flip-isomorphism is still a necessity to define it). We
illustrate this fact by showing an analogy in the world of
manifolds. Let $(W,M,M')$ be a manifold with $\partial M=\emptyset$
and let $H\colon M \isora M'$ be a diffeomorphism. Glueing the
cobordism along $H$ yields a closed manifold $W_H$. Different choices
of $H$ lead to different manifolds which however are in the same {\it
Schneiden-und-Kleben-cobordism class} (= cut-and-paste-cobordism
class). These cobordism groups were \eg studied in \cite{KreckKNO73}
(see also \cite{Ran98} Remark 30.30) and are quotients of the ordinary
cobordism groups using the equivalence relation 
$P\cup_f N \sim P\cup_g N$ for manifolds with boundary $(P,\partial P)$ and
$(N,\partial N)$ and homeomorphisms 
$f,g\colon \partial N \isora\partial P$. 
The SKL-cobordism group of an $(n+1)$-dimensional
manifold $V$ with $n>5$ is discovered by an asymmetric signature
similar to the one used for twisted doubles in Section
\ref{asymgeosec}: one takes a singular chain complex $C=C(V)$ with a
chain equivalence $\lambda=[M]\cap - \colon C^{n+1-*} \homra C$ and
looks at the image of $(C,\lambda)\in LAsy^{n+1}(\bZ)$. The SKL-bordism group
is isomorphic of the image of $L^{n+1}(\bZ) \ra LAsy^{n+1}(\bZ)$, hence $\bZ$ if
$n\equiv 0\pmod{4}$ and zero else. 
Our proof will use equivalent facts about symmetric Poincar\'e pairs
(see \cite{Ran98} 30.30).

\begin{thm}
\label{asymformthm}
Let $z$ be a formation. Let $t$ and $t'$ be two flip-isomorphisms. Then $\sigma^*(z,t)=\sigma^*(z,t')\in LAsy^0(\Lambda)$.
Hence the map $\sigma^*$ defined in Theorem \ref{flasymthm} induces a map
\begin{eqnarray*}
\sigma^*\colon \tilde l_{2q+2}(\Lambda) \ra LAsy^0(\Lambda)
\end{eqnarray*}
with $\tilde l_{2q+2}(\Lambda)=\{[z]\in l_{2q+2}(\Lambda) : \text{$z$ allows stable flip-isomorphisms}\}$
a submonoid of $l_{2q+2}(\Lambda)$.
\end{thm}

This theorem can be applied to the boundary of non-singular forms.
They are the obstructions of Wall's surgery theory interpreted as a special case of Kreck's surgery theory
\ie they live in the image of the inclusion
\begin{eqnarray*}
L_{2q+2}(\Lambda) &\hookrightarrow& l_{2q+2}(\Lambda)\\
(K,\theta) &\mt& \partial(K,\theta)=\kreckf{K}{1_K}{K}{\theta-\epsilon\theta^*}{\theta}
\end{eqnarray*}

\begin{cor}
\label{asymbdrycor}
Let $(K,\theta)$ be a $(-\epsilon)$-quadratic form. Then $z=\partial(K,\theta)$ is a formation.
\begin{enumerate}
\item $z$ has a (stable) flip-isomorphism (\ie $[z]\in\tilde l_{2q+2}(\Lambda)$) if and only if $(K,\theta)$ is non-singular.
\item $z$ has a stable strong flip-isomorphism if and only if $(K,\theta)$ is non-singular and $2\cdot [(K,\theta)]=0\in L_{2q+2}(\Lambda)$. 
\item $[z]\in\el$ is elementary if and only if $(K,\theta)$ is non-singular and $[(K,\theta)]=0\in L_{2q+2}(\Lambda)$.
\item If $(K,\theta)$ is non-singular, $\sigma^*([z,t])=[(K,\theta-\epsilon\theta^*)]\in LAsy^0(\Lambda)$ for any stable flip-isomorphism $t$.
\item Assume that either $\Lambda$ is a field of characteristic different to $2$ or that $\Lambda=\bZ$ and $\epsilon=-1$.
      $z$ is elementary if and only if the asymmetric signature vanishes.
\end{enumerate}
\end{cor}
\begin{proof}
\begin{enumerate}
\item Obviously $(K,\theta)$ must be non-singular if $z$ allows a stable flip-isomorphism. If $(K,\theta)$ is non-singular then
      $t=(\lambda^*,1,\epsilon\lambda^{-1})$ is a flip-isomorphism of $z$. \label{abc1}
\item If $(K,\theta)$ has a strong flip-isomorphism $(\alpha,\beta)$ then $\beta^*\theta\beta=-\theta\in Q_{-\epsilon}(K)$.
      Hence $[(K,\theta)]=[(K,-\theta)]\in L_{2q+2}(\Lambda)$. On the other hand if $2\cdot [(K,\theta)]=0\in L_{2q+2}(\Lambda)$
      then, after stabilization, there is an isomorphism $\beta^*\theta\beta=-\theta\in Q_{-\epsilon}(K)$ and
      $(\epsilon\lambda\beta, \beta)$ is a strong flip-isomorphism.
\item Easy.
\item Let $t=(\alpha,\beta,\nu)$ be the flip-isomorphism of \ref{abc1}. Let $\lambda=\theta-\epsilon\theta^*$ and let
$$
    \rho=\mat{0&0&\alpha\\1&0&-\epsilon\\0&1&\epsilon\alpha(\nu^*-\epsilon\nu)\alpha^*}\colon M=F\oplus F^*\oplus F \ra M^*
$$
as in Definition \ref{asymsignflipdef}. Then $\rho\oplus-\lambda$ has the
lagrangian $$\mat{\epsilon & 1\\0&-\epsilon\lambda\\1&0\\-\epsilon&1}$$
Hence the asymmetric signature $\sigma^*(z,t)=(K,\lambda)\in LAsy^0(\Lambda)$. By Theorem \ref{asymformthm}
the asymmetric signature is independent of the choice of stable flip-isomorphism.
\item In these cases the maps
$(1+T_{-\epsilon})\colon L_{2q+2}(\Lambda) \rightarrow L^{2q+2}(\Lambda)$ and
$L^{2q+2}(\Lambda) \rightarrow LAsy^0(\Lambda)$
are injective (see e.g. \cite{Ran98} Chapter 39D and Remark \ref{asyquadrem}).
\end{enumerate}
\end{proof}

Here is a counter-example for the converse of Theorem \ref{flasymthm}.
\begin{ex}
\label{asymcountex}
Let $z$ be the boundary of any non-singular skew-quadratic form over $\bZ$ or $\bZt$ with
non-trivial Arf-invariant. 
By Corollary
\ref{asymbdrycor} $z$ has stable flip-isomorphisms and all asymmetric
signatures vanish but it is not stably elementary.
\end{ex}

Back to the proof of Theorem \ref{asymformthm}.
We recall that in Section \ref{asymsignsec} the asymmetric signature $\sigma^*(z,t)\in LAsy^0(\Lambda)$ 
was defined as the asymmetric signature of the $(2q+2)$-dimensional symmetric Poincar\'e pair
$x_t$. In our case $C$ is contractible and $D\oplus D$ and $\partial E_t$ are chain equivalent.
The following two lemmas treat this situation in general.
\begin{lem}
\label{contrtwdblem}
Let $(f\colon C \ra D,(\delta\phi,\phi))$ be an $n$-dimensional symmetric
Poincar\'e pair and $(h,\chi)\colon (C,\phi) \homra (C,\phi)$ a self-equivalence.
Assume that $C$ is contractible with $\Delta \colon 1 \simeq 0 \colon C \ra C$.

Define $\nu=\delta\phi+(-)^{n-1}f\Delta^\%\phi f^*$
and $\bar\rho=\Delta^\%(\Delta^\%\phi-\chi-h\Delta^\%\phi h^*)$. 

There is an equivalence
\begin{eqnarray*}
(a,\sigma)&\colon& (D\oplus D, \nu\oplus -\nu) \stackrel{\simeq}{\ra} (D\cup_h D, \delta\phi\cup_\chi -\delta\phi)
\\
a&=&\mat{1&0\\0&1\\0&0} \colon D_r\oplus D_r \ra (D\cup_h D)_r=D_r\oplus D_r \oplus C_{r-1}
\\
\sigma_s&=&\mat{    (-)^n f\bar\rho_s f^*   & 0 & 0 \\
            0           & 0 & (-)^{s-1}f\Delta^\%\phi_s\\
            (-)^{n+1-r}\Delta^\%\phi_s h^*f^*& 0 & (-)^{n+1-r+s}T_\epsilon\Delta^\%\phi_{s-1}} \colon
\\
&& 
(D\cup_h D)^{n+1-r+s} = D^{n+1-r+s}\oplus D^{n+1-r+s} \oplus C^{n-r+s}\\
&&\quad\ra(D\cup_h D)_r=D_r\oplus D_r \oplus C_{r-1}
\end{eqnarray*}
of $n$-dimensional $\epsilon$-symmetric Poincar\'e complexes.
\end{lem}
\begin{lem}
\label{contrasymlem}
Let $x=(g\colon \partial E \ra E,(\theta,\partial\theta))$ be an $(n+1)$-dimensional symmetric
Poincar\'e pair such that the boundary $(\partial E, \partial\theta)$ is a twisted double of an $n$-dimensional symmetric
Poincar\'e pair $(f\colon C \ra D,(\delta\phi,\phi))$ with respect to a self-equivalence
$(h,\chi)\colon (C,\phi) \homra (C,\phi)$. 
We write 
$$
    g=\mat{j_0 & j_1 & k}\colon \partial E_r = D_r\oplus D_r \oplus C_{r-1} \ra E_r
$$
Assume that $C$ is contractible with $\Delta \colon 1 \simeq 0 \colon C \ra C$.
Hence there is a chain equivalence $\Delta \colon 0 \simeq 1 \colon C \ra C$ (\ie $d\Delta+\Delta d=1_C$).
Applying Lemmas \ref{contrtwdblem}, \ref{deltapairlem} to $x$ yields an $(n+1)$ -dimensional symmetric
Poincar\'e pair
\begin{eqnarray*}
x'&=&(\mat{j_0& j_1} \colon D\oplus D \ra E, (\theta'=\theta+(-)^n g\sigma g^*, \nu\oplus -\nu))
\end{eqnarray*}

Let $(B,\lambda)$ be the asymmetric complex of $x$ and $(B',\lambda')$ the asymmetric complex of $x'$
as given in Proposition \ref{cplxprop}. Then there is an equivalence 
\begin{eqnarray*}
(b,\xi) &\colon& (B, T\lambda) \ra (B', T\lambda')
\\
b&=&\mat{1&(-)^{r}j_0 f\Delta & 0\\0&0&1} \colon
B_r=E_r\oplus C_{r-1}\oplus D_r \ra B'_r=E_r\oplus D_r
\\
\xi&=&\mat{ (-k\Delta\phi_0 h^*+(-)^{r+1} j_0 f \bar\rho_0)\Delta^*f^*j_0^*     & 0\\
        f\Delta\phi_0((-)^{n-r} k^*-h^*\Delta^*f^*j_0^*)        & (-)^{n-r}f\Delta\phi_0 f^*}\colon
\\
&&{B'}^{n+2-r}=E^{n+2-r}\oplus D^{n+1-r} \ra B'_r=E_r\oplus D_r
\end{eqnarray*}
of $(n+1)$-dimensional asymmetric complexes with $\rho$ and $\bar\rho$ as defined in Lemma \ref{contrtwdblem}.
\end{lem}

\begin{proof}[Proof of Theorem \ref{asymformthm}]
By Lemma \ref{contrasymlem}, $\sigma^*(z,t)\in LAsy^0(\Lambda)$ is the asymmetric signature
of the $(2q+2)$-dimensional symmetric Poin\-car\'e pair
\begin{eqnarray*}
{x'}^t&=&({g'}^t \colon D\oplus D \ra E, (\delta\theta', \nu\oplus -\nu))\\
{g'}^t_{q+1}&=&\mat{1&-\beta} \colon D_{q+1}\oplus D_{q+1}=G\oplus G \ra E_{q+1}=G\\
\delta\theta'_0&=&-\epsilon Y \colon E^{q+1}=G^* \ra E_{q+1}=G\\
\nu_0&=&-\tilde\mu^*\colon  D^q=F \ra D_{q+1}=G\\
\nu_0&=&-\tilde\mu\colon D^{q+1}=G^* \ra D_q=F^*
\end{eqnarray*}
By Corollary \ref{asymcor} $\sigma^*({x'}^t)$ is the image of the union 
of ${x'}^t$ in $LAsy^0(\Lambda)$. But there is another way to construct ${x'}^t$: 
there is a $(2q+2)$-dimensional quadratic Poincar\'e pair
\begin{eqnarray*}
\tilde x&=& (\tilde g \colon D\oplus D' \ra E, (0,\nu\oplus-\nu'))\\
\tilde g&=& \mat{1&-1} \colon D_{q+1}\oplus D_{q+1}=G\oplus G \ra E_{q+1}=G\\
\nu'_0  &=& -\tilde\gamma^* \colon {D'}^q=F^* \ra D'_{q+1}=G\\
\nu'_0  &=& -\tilde\gamma   \colon {D'}^{q+1}=G^* \ra D'_q=F
\end{eqnarray*}
and an isomorphism $(\bar e_t,\bar \chi_t) \colon (D,\nu) \isora (D',\nu')$ given by
\begin{eqnarray*}
\bar e_{t,q+1}=\beta    \colon D_{q+1}=G \ra D'_{q+1}=G\\
\bar e_{t,q}  =\alpha^{-*}\colon D_q=F^*   \ra D'_{q}=F\\
\bar\chi_{t,0}  =-\epsilon Y\colon {D'}^{q+1}=G^*\ra D'_{q+1}=G
\end{eqnarray*}
Lemma \ref{deltapairlem} and the isomorphism can be used to replace
the ``boundary component'' $(D', \nu')$ by $(D,\nu)$. The result will be 
${x'}^t$. Glueing both ends (\ie $D$ and $D'$) of $\tilde x$ together
using $(\bar e_t,\bar \chi_t)$ yields the union of ${x'}^t$. 
Hence all unions of ${x'}^t$ for different choices of $t$ are in the same 
{\bf algebraic Schneiden-und-Kleben-cobordism} class. 
By \cite{Ran98} 30.30(ii) their images in $LAsy^{2q+2}(\Lambda)$ coincide.
Those images are precisely the asymmetric signatures $\sigma^*({x'}^t)=\sigma^*(z,t)$.
\end{proof}

\chapter{Preformations with Linking Forms}
\label{linkformchap}

{\bf For the whole chapter let $q\ge 2$, $\epsilon=(-)^q$ and let $\Lambda$ be a weakly finite ring
with $1$ and involution.}

We consider a special group of preformations $\kreckfs$ 
namely those for which $\mu$ becomes an isomorphism after localization.
The most important examples
are probably preformations over $\bZ$ with injective $\mu$
and $\rk G=\rk F$. For those classes of preformations one can
use the theory of linking forms developed for formations in \cite{Ran81} 3.4.
and improve our results for asymmetric signatures.

In {\bf Section \ref{localsec}} we quickly repeat the concept of localization
and define linking forms following \cite{Ran81}. {\bf Section \ref{linkflipsec}} defines linking forms
on preformations and discusses the relationship between isometries of those linking forms 
and flip-isomorphisms. It turns out that every flip-isomorphism induces
an isometry of linking forms and in turn every isometry of linking forms gives rise
to a stable flip-isomorphism.

Similar to the flip-$l$-monoids, in {\bf Section \ref{linkmodsec}} 
we define linking-$l$-monoids of
preformations with a choice of isometry of their linking forms 

In {\bf Section \ref{linkasymsec}} we show that
the asymmetric signatures we defined in Section \ref{asymsignsec}
only depend on the effect of the flip-isomorphism on the linking forms
of the preformation. If a preformation is stably elementary then
all those asymmetric signatures vanish (see Theorem \ref{elemlinkasymthm}). 
This theorem is an improvement to Theorem \ref{elemasymsignthm} 
because isometries of linking forms are easier to handle then flip-isomorphisms. 
For $\bZ$ there are only finitely many isometries of a given
linking form. Also, it is enough to look at one representative of 
a class in $\el$. 

We will use these results to define asymmetric signatures
for certain simply-connect\-ed Kreck surgery problems using the
topological linking forms of the manifolds involved
(see {\bf Section \ref{linkapplsec}}).

\section{Localization}
\label{localsec}
Although we could generalize our results for Ore-localization, we will focus on
the easier case of localization away from a central multiplicative subset.

We repeat some definitions from \cite{Ran81} Chapter 3.1. and 3.4.

\begin{defi}
\label{multsetdef}\index{Central and multiplicative subset}
A subset $S\subset\Lambda$ is called {\bf central and multiplicative} if
\begin{enumerate}
\item $st \in S$ for all $s,t \in S$,
\item $\bar s\in S$ for all $s \in S$,
\item if $sa=0\in\Lambda$ for some $s\in S$ and $a\in\Lambda$ then $a=0 \in \Lambda$,
\item $sa=as\in\Lambda$ for all $s\in S$ and $a\in\Lambda$. \hfill\qed
\end{enumerate}
\end{defi} 

\begin{defi}
Let $S\subset\Lambda$ be a central and multiplicative subset.
The {\bf localization $S^{-1}\Lambda$ of $\Lambda$ away from $S$}\index{Localization}
is the ring with involution and $1$ defined by the equivalence classes 
of pairs $(a,s)\in \Lambda\times S$ under the relation 
$$
(a,s) \sim (a',s') \Longleftrightarrow as'=a's\in\Lambda
$$
with
\begin{eqnarray*}
(a,s)+(a',s')       &=& (as'+a's, ss')\\
(a,s)\cdot(a',s')   &=& (aa', ss')\\
\overline{(a,s)}&=&(\bar a, \bar s)
\end{eqnarray*}
\hfill\qed\end{defi}

\begin{ex}
\label{grprngex}
Let $\pi$ be a group and $w\colon \pi \ra \bZ/2\bZ$ be a group morphism.
Let $\Lambda=\bZ[\pi]$ be its group ring endowed with the $w$-twisted involution
$n\cdot 1_g \mt w(g)n\cdot 1_{g^{-1}}$. Then $S=\bZ \setminus\{0\}$ is
a central multiplicative subset of $\Lambda$. The localization
of $\Lambda$ away from $S$ is canonically isomorphic to the group ring $\bQ[\pi]$
with the obvious involution.
\end{ex}

\begin{defi}
Let $S\subset\Lambda$ be a central and multiplicative subset.
A morphism $f\colon M \ra N$ of $\Lambda$-modules is called an
{\bf $S$-isomorphism}\index{S-isomorphism@$S$-isomorphism} if the induced   $S^{-1}\Lambda$-module morphism
\begin{eqnarray*}
S^{-1} f \colon S^{-1}M &\ra& S^{-1}N\\
\frac{x}{s} &\mt& \frac{f(x)}{s}
\end{eqnarray*}
is an isomorphism.
\hfill\qed\end{defi}

\begin{defi}
Let $S\subset\Lambda$ be a central and multiplicative subset.
A chain complex $C$ over $\Lambda$ is {\bf $S$-acyclic}\index{S-acyclic@$S$-acyclic}
if the chain complex $S^{-1} C=C\otimes_\Lambda S^{-1}\Lambda$ is acyclic.
\hfill\qed\end{defi}

\begin{defi}
Let $S\subset\Lambda$ be a central and multiplicative subset.
A {\bf $(\Lambda, S)$-module $M$}\index{L@$(\Lambda, S)$-module} is an $\Lambda$-module such that there is
an exact sequence of $\Lambda$-modules
\begin{eqnarray*}
0\ra P \stackrel{d}{\ra} Q \ra M \ra 0
\end{eqnarray*}
with $P$ and $Q$ \fg projective and $d$ an $S$-isomorphism.
\hfill\qed\end{defi}

\begin{defi}
Let $S\subset\Lambda$ be a central and multiplicative subset.
\begin{enumerate}
\item\index{QLS1@$Q^{\epsilon}(\Lambda, S)$}\index{QLS2@$Q_{\epsilon}(\Lambda, S)$}
Define the {\bf relative $Q$-groups}
\begin{eqnarray*}
Q^{\epsilon}(\Lambda, S)&=&\{b\in S^{-1}\Lambda| b-\epsilon \bar b=a-\epsilon \bar  a, a\in\Lambda\}/\Lambda
\\
&\subset&Q^{\epsilon}(S^{-1}\Lambda/\Lambda)=\{b\in S^{-1}\Lambda| b-\epsilon \bar b\in\Lambda\}/\Lambda
\\
Q_{\epsilon}(\Lambda, S)&=&\{b\in S^{-1}\Lambda | b=\epsilon\bar b\}/\{a+\epsilon\bar a\colon a\in\Lambda\}
\\
&\subset&Q_{\epsilon}(S^{-1}\Lambda/\Lambda)=S^{-1}\Lambda/\{a+b-\epsilon\bar b| a\in\Lambda, b\in S^{-1}\Lambda\}
\end{eqnarray*}

\item \index{Linking form!symmetric}
An {\bf $\epsilon$-symmetric linking form $(M,\lambda)$ over $(\Lambda,S)$} is an $(\Lambda,S)$-module $M$
together with a pairing $\lambda \colon M \times M \ra S^{-1}\Lambda/\Lambda$ such that
$\lambda(x,-)\colon M \rightarrow S^{-1}\Lambda/\Lambda$ is $\Lambda$-linear for all $x\in M$ and 
$\lambda(x,y)=\epsilon\overline{\lambda(y,x)}$ for all $x,y\in M$.

\item \index{Linking form!split quadratic}
A {\bf split $\epsilon$-quadratic linking form $(M,\lambda,\nu)$ over $(\Lambda,S)$} is an
$\epsilon$-symmetric linking form $(M,\lambda)$ over $(\Lambda,S)$ together with a map
$\nu \colon M \ra Q_{\epsilon}(S^{-1}\Lambda/\Lambda)$
such that for all $x,y\in M$ and $a\in\Lambda$
\begin{enumerate}
\item $\nu(ax)=a\nu(x)\bar a \in Q_\epsilon(S^{-1}\Lambda/\Lambda)$
\item $\nu(x+y)-\nu(x)-\nu(y)=\lambda(x,y)\in Q_\epsilon(S^{-1}\Lambda/\Lambda)$
\item $(1+T_\epsilon)\nu(x) = \lambda(x,x)\in Q^\epsilon(\Lambda, S)$
\end{enumerate}
\end{enumerate}
\hfill\qed\end{defi}

\begin{ex}
Let $\Lambda=\bZ$ and $S=\bZ\setminus\{0\}$. Then $S^{-1}\Lambda=\bQ/\bZ$.
A $(\Lambda,S)$-module is nothing but a finite abelian group.

An $\epsilon$-symmetric linking form $(M,\lambda)$ over $(\Lambda,S)$
is a finite abelian group $M$ together with a bilinear $\epsilon$-symmetric pairing
$\lambda \colon M\times M \ra \bQ/\bZ$ on it.

For $\epsilon=1$ a  split quadratic linking form $(M,\lambda,\nu)$ over $(\Lambda,S)$ is nothing but
symmetric linking form $(M, \lambda)$ with $\nu \colon M \ra Q_{1}(S^{-1}\Lambda/\Lambda)=\bQ/\bZ$
given by $\nu(x)=\frac{1}{2}\lambda(x,x)$.

For $\epsilon=-1$, $Q_{-1}(\bQ/\bZ)=0$ and a split skew-quadratic linking form $(M,\lambda,\nu)$ over $(\Lambda, S)$ is
a skew-symmetric linking form $(M,\lambda)$ with $\lambda(x,x)=0$ for all $x\in M$.
\end{ex}

\section{Flip-Isomorphisms and Linking Forms of Preformations}
\label{linkflipsec}

Let $S\subset\Lambda$ be a central multiplicative subset of $\Lambda$.

As in the proof of \cite{Ran81} p. 242ff
we define the linking forms of preformations for which $\gamma$ or $\mu$ are $S$-isomorphisms.
\begin{defi}
Let $x=\kreckfs$ be a regular $\epsilon$-quadratic split preformation.
\begin{enumerate}
\item If $\mu$ is an $S$-isomorphism
there is a split $(-\epsilon)$-quadratic
linking form $L_\mu=(\coker\mu, \lambda_\mu, \nu_\mu)$ over $(\Lambda, S)$ given by
\begin{align*}
\lambda_\mu \colon \coker\mu \times \coker\mu &\ra S^{-1}\Lambda/\Lambda, &
(x,y) &\mt \dfrac{1}{s}\gamma^*(x)(g)\\
\nu_\mu \colon \coker\mu &\ra Q_{-\epsilon}(S^{-1}\Lambda/\Lambda), &
y &\mt \dfrac{1}{s}\theta(g)(g)\overline{\dfrac{1}{s}}
\end{align*}
for $x,y \in F^*$, $g\in G$, $s\in S$ such that $sy=\mu(g)$.

\item If $\gamma$ is an $S$-isomorphism
there is a split $(-\epsilon)$-quadratic
linking form $L_\gamma=(\coker\gamma, \lambda_\gamma, \nu_\gamma)$ over $(\Lambda, S)$ given by
\begin{align*}
\lambda_\gamma \colon \coker\gamma \times \coker\gamma &\ra S^{-1}\Lambda/\Lambda,&
(x,y) &\mt \epsilon\dfrac{1}{s}\mu^*(x)(g)\\
\nu_\gamma \colon \coker\gamma &\ra Q_{-\epsilon}(S^{-1}\Lambda/\Lambda),&
y &\mt -\dfrac{1}{s}\theta(g)(g)\overline{\dfrac{1}{s}}
\end{align*}
for $x,y \in F$, $g\in G$, $s\in S$ such that $sy=\gamma(g)$.

\item If $\gamma^*\mu$ is an $S$-isomorphism
there is a split $(-\epsilon)$-quadratic linking form 
$L_{\gamma^*\mu}=(\coker\gamma^*\mu, \lambda_{\gamma^*\mu}, \nu_{\gamma^*\mu})$ over $(\Lambda, S)$ given by
\begin{align*}
\lambda_{\gamma^*\mu} \colon \coker\gamma^*\mu \times \coker\gamma^*\mu &\ra S^{-1}\Lambda/\Lambda,&
(x,y) &\mt \dfrac{1}{s} x(g)\\
\nu_{\gamma^*\mu} \colon \coker\gamma^*\mu &\ra Q_{-\epsilon}(S^{-1}\Lambda/\Lambda), &
y &\mt \dfrac{1}{s}\theta(g)(g)\overline{\dfrac{1}{s}}
\end{align*}
for $x,y \in G^*$, $g\in G$, $s\in S$ such that $sy=\gamma^*\mu(g)$.
\end{enumerate}
Similarly, there exist $(-\epsilon)$-symmetric linking forms $L^\gamma$, $L^\mu$ and $L^{\gamma^*\mu}$
for $\epsilon$-quadratic preformations with $\gamma$ respectively $\mu$ $S$-isomorphisms.
\hfill\qed\end{defi}

\begin{rem}
The definitions are taken from the proof of \cite{Ran81} Proposition 3.4.3 which
establishes a correspondence between {\bf $S$-formations} and linking forms.
It is easy to verify that a regular split $\epsilon$-quadratic preformation $z=\kreckfs$ and its flip $z'$
are $S$-formations if $\mu$ or $\gamma$ are $S$-isomorphisms.
Under that correspondence $z$ is mapped to $L_\mu$ and
$z'$ is mapped to $L_\mu$. We will exploit this fact in the proof
of Proposition \ref{sisoprop}.
\end{rem}

Linking forms can tell us something about elementariness and the existence of
flip-isomorphisms.
\begin{prop}
\label{sisoprop}
Let $z=\kreckfs$ be a regular split $\epsilon$-quadratic  preformation
with either $\mu$ or $\gamma$ an $S$-isomorphism.
\begin{enumerate}

\item If $z$ allows a flip-isomorphism then both $\gamma$ and $\mu$ are 
$S$-isomorphisms. Every flip-isomorphism $t=(\alpha,\beta,\chi)$ induces
an isomorphism of split $(-\epsilon)$-quadratic linking forms 
$[\alpha^{-*}]\colon L_\mu \stackrel{\cong}{\ra} L_\gamma$.

\item Assume $\gamma$ and $\mu$ are both $S$-isomorphisms and $L_\gamma$ and $L_\mu$ are isomorphic.
Every isomorphism $l\colon L_\mu \stackrel{\cong}{\ra} L_\gamma$ induces a 
stable flip-isomorphism $(\alpha,\beta,\chi)$ of $z$
such that $[\alpha^{-*}]=l\colon L_\mu \stackrel{\cong}{\ra} L_\gamma$.

\item If $z$ is stably elementary then
\begin{enumerate}
\item $\gamma$ and $\mu$ are $S$-isomorphism,
\item $L_\mu\cong L_\gamma$,
\item $(G, \gamma^*\mu, \theta)$ is $S$-hyperbolic \ie
$L_{\gamma^*\mu}=0 \in \tilde L_\epsilon(\Lambda, S)$ (see \cite{Ran81} p.271). 
\end{enumerate}

\end{enumerate}
Similar for the non-split case.
\end{prop}

\begin{proof}
For the proof it is necessary to remember the definition
of the $(2q+1)$-dimensional quadratic complexes $(N,\zeta)$ and $(N',\zeta')$ associated with $z$ and its flip $z'$
((\ref{Ndefeqn}) and (\ref{Npdefeqn}) on p.~\pageref{Ndefeqn}).
\begin{enumerate}
\item Direct calculation.

\item 
$z$ and its flip $z'$ are split $\epsilon$-quadratic $S$-formations in the sense of \cite{Ran81} p.240.
Hence we can apply \cite{Ran81} Proposition 3.4.3. The proof shows that there exists a
stable isomorphism of split $\epsilon$-quadratic $S$-formations between $z$ and $z'$.
Using the isomorphism in Remark \ref{weakisorem} \ref{isorem2} it is not difficult to show that
this leads to a stable (weak) flip-isomorphism of $z$.

\item Obvious from Proposition~\ref{elemprop}~\ref{elemprop3} and \cite{Ran81} Proposition 3.4.6ii).
\end{enumerate}
\end{proof}

\section{The Linking-$l$-Monoid $ll^S_{2q}(\Lambda)$}
\label{linkmodsec}

Let $S\subset\Lambda$ be a central multiplicative subset.

Proposition \ref{sisoprop} shows that there is a strong relationship between flip-isomorphisms
of preformations which allow linking forms and the isomorphisms of those linking forms.
Similar to Section \ref{asymflmonoidelemsec} we define a monoid of preformations
with linking forms and a choice of isomorphism between them.

\begin{defi}
\begin{eqnarray*}
l^S_{2q+2}(\Lambda)  &=& \left\{[\kreckfs]\in \el | \text{$\mu$ and $\gamma$ are $S$-isomorphisms}\right\}\\
fl^S_{2q+2}(\Lambda) &=& \left\{[(z,t)] : [z]\in l^S_{2q+2}(\Lambda)\right\}
\end{eqnarray*}
are sub-monoids of $l_{2q+2}(\Lambda)$ and $fl_{2q+2}(\Lambda)$. Similarly one defines
$l_S^{2q+2}(\Lambda)$ and $fl_S^{2q+2}(\Lambda)$ in the non-split case.
\hfill\qed\end{defi}

\begin{defi}
We consider tuples $(z,l)$ with $z=\kreckfs$ a regular $\epsilon$-quadratic preformation
such that $\gamma$ and $\mu$ are $S$-isomorphisms and $l\colon L_\mu \isora L_\gamma$
an isomorphism of linking forms.

An isomorphism $(\eta,\zeta)$ of such tuples $(z,l)$ and $(z',l')$ is a strong isomorphism 
$(\eta,\zeta)\colon z \isora z'$ of preformations such that $l'=\eta^* l \eta$.

Define the hyperbolic elements $(y_k,0)$ with $y_k=\partial H_{-\epsilon}(\Lambda^{k})$.

The stable isomorphism classes of such tuples form the {\bf linking-$l$-monoid $ll^S_{2q+2}(\Lambda)$}
Similarly we can define {\bf $ll_S^{2q+2}(\Lambda)$} for the non-split case.
\index{Linking-$l$-monoid}\index{ll1@$ll^S_{2q+2}(\Lambda)$}\index{ll2@$ll_S^{2q+2}(\Lambda)$}
\hfill\qed\end{defi}

\begin{prop}
\label{linkdiagprp}
There is a commuting diagram of surjective morphisms of abelian monoids with zero
\begin{eqnarray*}
\xymatrix
{
fl^S_{2q+2}(\Lambda)    \ar@{>>}[d]_{L} 
            \ar[r]_{\pi_f}
&
l^S_{2q+2}(\Lambda)
\\
ll^S_{2q+2}(\Lambda)    \ar[ru]_{\pi_l}
}
\end{eqnarray*}
with $L(z,t=(\alpha,\beta,\nu))=(z,[\alpha^{-*}])$, 
$\pi_f(z,t)=z$ and $\pi_l(z,l)=z$. 
The fibre $\pi^{-1}([z])$ of a preformation $z=\kreckfs$ is either empty or
the set of all isometries $L_\mu \isora L_\gamma$.
Similar for the non-split case.
\end{prop}
\begin{proof}
Use Proposition \ref{sisoprop}
\end{proof}

\section{Asymmetric Signatures}
\label{linkasymsec}

Let $S\subset\Lambda$ be a central multiplicative subset.
Let $z=\kreckfs$ be a regular $\epsilon$-quadratic split preformation
such that $\gamma$ and $\mu$ are $S$-isomorphisms. In other words $[z]\in l_{2q+2}^S(\Lambda)$.
Proposition \ref{sisoprop} shows that any 
flip-isomorphism gives rise to an isomorphisms of the linking forms $L_\gamma$ and $L_\mu$
and vice versa. The relationship between flip-isomorphism and linking form isomorphism
goes even further. We will prove that the asymmetric signature of two flip-isomorphisms
which induce the same isomorphism $L^\mu \stackrel{\cong}{\ra} L^\gamma$ are identical.

In the case of $\Lambda=\bZ$, this drastically reduces the amount of work one has to put into
checking all the asymmetric signatures (see Theorem
\ref{elemlinkasymthm}). Instead of going through all flip-isomorphisms
of all the preformations $z+\partial H_{-\epsilon}(\Lambda^n)$, we only have to compute them for the finite
number of isomorphisms of two given linking forms of one arbitrary
representative of $[z]\in\el$.

First we need to check that the asymmetric signature of an isomorphism  $L^\mu \stackrel{\cong}{\ra} L^\gamma$
is well-defined. 

\begin{prop}
\label{asywellprp}
Let $t=(\alpha,\beta,\chi)$ and $t'=(\alpha',\beta',\chi')$ two flip-isomorphisms of $z$
which induce the same isomorphism of linking forms 
$$[\alpha^{-*}]=[{\alpha'}^{-*}]\colon L^\mu \stackrel{\cong}{\ra} L^\gamma$$
Then $\sigma^*(z,t)=\sigma^*(z,t')\in LAsy^0(\Lambda)$.
Similar in the non-split case.
\end{prop}

We need a little lemma which shows that for $1$-dimensional $S$-acyclic complexes
quasi-isomorphisms and chain-equivalences are the same. This is not true for
arbitrary (even free) chain complexes.
\begin{lem}
\label{superlem}
Let $N$ and $N'$ be free $1$-dimensional $S$-acyclic chain complexes
and $f\colon N\ra N'$ a chain map which induces the zero map in homology.
Then there is a chain homotopy $\Delta \colon f\simeq 0$.
\end{lem}

\begin{proof}[Proof of Proposition \ref{asywellprp}]
The asymmetric signature of a flip-isomorph\-ism is constructed in Section \ref{asymsignsec}
as the asymmetric signature of the $(2q+2)$-dimensional symmetric Poincar\'e pair
$x^t=(g^t \colon \partial E_t \ra E, (0, \theta_t))$. Its boundary is
a twisted double of the symmetric Poincar\'e pair $(f\colon C\ra D, (0,\phi))$ 
in respect to the self-equivalence $(h_t,0)\colon (C,\phi) \isora (C,\phi)$.
We will show that the two flip-isomorphisms $t$ and $t'$ lead to homotopic 
equivalences $(h_t,0) \simeq (h_{t'},0)$ and that therefore the twisted doubles
$(\partial E_t,\theta_t)$ and$(\partial E_{t'},\theta_{t'})$ are equivalent
and, finally, that
the asymmetric signatures of $x^t$ and $x^{t'}$ are the same.

As described in Section \ref{fliptranslsec},
 $t$ and $t'$ induce two isomorphisms
$$(e,\rho), (e', \rho') \colon (N, \zeta) \stackrel{\cong}{\ra} (N', \zeta')$$
of the $(2q+1)$-dimensional quadratic complexes defined in (\ref{Ndefeqn}) and (\ref{Npdefeqn}) on p.~\pageref{Ndefeqn}ff.
The fact that $t$ and $t'$ induce the same linking form isomorphisms translates into
$e^*={e'}^* \colon H^*(N') \stackrel{\cong}{\ra} H^*(N)$. 

By Lemma \ref{superlem}, $e$ and $e'$ are chain homotopic.
Let $\Delta \colon e \simeq e' \colon N \ra N'$ be a chain homotopy.

By the proof of \cite{Ran80a} Proposition 3.4. (see also Lemma \ref{bdrmaplem}) $(e,\rho)$ and $(e',\rho')$ induce isomorphisms 
$$(\partial e, 0), (\partial e', 0) \colon (C,\phi)=(\partial N, (1+T)\partial\zeta)
\stackrel{\cong}{\ra} (C',\phi')=(\partial N', (1+T)\partial\zeta')$$
of $2q$-dimensional symmetric Poincar\'e complexes given by
\begin{eqnarray*}
\partial e=\mat{e & (-)^{r+1}(1+T)\rho \\ 0 & e^{-*}} \colon 
\partial N_r=N_{r+1}\oplus N^{2q+1-r} \rightarrow \partial N'_r=N'_{r+1}\oplus {N'}^{2q+1-r}
\end{eqnarray*}

Using the fact that $N$ and $N'$ are short and $S$-acyclic 
one can show that there is a chain equivalence 
$(\partial\Delta,0) \colon (\partial e,0) \simeq (\partial e',0) \colon (C,\phi) \isora (C',\phi')$
given by
\begin{eqnarray*}
\partial\Delta_{q+1}&=&\mat{0 & \epsilon {\beta'}\Delta^*\alpha^{-*}} 
\colon C_q=F\oplus F^* \ra C'_{q+1}=G^*
\\
\partial\Delta_q&=&\mat{\Delta\\0} \colon C_{q-1}=G^* \ra C'_q=F^*\oplus F
\end{eqnarray*}

As explained in Section \ref{fliptranslsec} we compose $\partial e$ 
with the inverse $(h,0)\colon (C,\phi) \stackrel{\cong}{\ra} (C',\phi')$ from (\ref{hdefeqn})
on page~\pageref{hdefeqn} in order to get the self-equivalence
$(h_t,0) \colon (C,\phi) \stackrel{\cong}{\ra} (C,\phi)$. Using Lemma \ref{chhomlem} or by
direct calculation one finds a homotopy of the chain maps
$$(h^{-1}\partial\Delta,0) \colon (h_t,0) \simeq (h_{t'},0) \colon (C,\phi) \isora (C,\phi)$$
which can be fed into Lemma \ref{asysgnhomlem}. Hence $\sigma^*(z,t')=\sigma^*(x^{t'})=\sigma^*(x'')\in LAsy^0(\Lambda)$
with $x''=(g''\colon \partial E_t \ra E, (0,\theta_t))$ given by
$$p=\mat{1 & {\beta'} & 0 -\beta'\Delta^*\alpha^{-*}}\colon 
\partial E_{t,q+1}=G\oplus G\oplus F\oplus F^* \ra E_{q+1}=G$$

Finally, there is a homotopy equivalence $(1,1;l)\colon x'' \ra x^t$ given by 
$$l=\mat{0&\epsilon\beta'\Delta^*\alpha^{-*}&0}\colon \partial E_{t,q}=F^*\oplus F^*\oplus G^* \ra E_{q+1}=G$$
Hence $x''$ and $x^t$ are cobordant rel$\partial$ by Lemma \ref{cobhomlem}. 
By Proposition \ref{asymgenprop} their asymmetric signatures coincide.
\end{proof}

This proposition shows that the asymmetric signature of a flip-isomorphism in $l^S_{2q+2}(\Lambda)$
only depends on its induced isomorphism of linking forms. 

\begin{thm}
\label{elemlinkasymthm}
There is a lift of the asymmetric signature map of Theorem \ref{flasymthm}
\begin{eqnarray*}
\xymatrix
{
fl_S^{2q+2}(\Lambda)    \ar@{>>}[d]_{L} 
            \ar[r]_{\sigma^*}
&
LAsy^0(\Lambda)
\\
ll_S^{2q+2}(\Lambda)    \ar[ru]_{\sigma^*}
}
\end{eqnarray*}
Let $\pi_l$ be as in Proposition \ref{linkdiagprp}. 
If $[z]\in l_S^{2q+2}(\Lambda)$ is elementary then
$\sigma^*(\pi_l^{-1}[z])=\{0\}$ \ie there exist isomorphisms $l\colon
L^\mu \stackrel{\cong}{\ra} L^\gamma$ and the asymmetric signatures
$\sigma^*(z,l)\in LAsy^0(\Lambda)$ vanish for all of them.
Similar in the split case.
\end{thm}

\section{Asymmetric Signatures of Simply-Connected Manifolds}
\label{linkapplsec}
We will now concentrate our efforts to simply-connected manifold \ie the case of
$\Lambda=\bZ$ and $S=\bZ \setminus \{0\}$.
First we observe that results of Section \ref{linkasymsec} can be generalized -
they remain true even if $\gamma$ and $\mu$ are not injective and if $G$ is not free (\ie if
the preformation is not regular). We can also show that the algebraic linking forms of an
obstruction preformation are induced by the topological linking forms in 
certain simply-connected Kreck surgery problems.

Let $z=\skreckfs$ be a  $\epsilon$-quadratic preformation over $\bZ$ such that
$\coker\gamma$ and $\coker\mu$ are finite. ($G$ may have torsion and $\gamma$ and $\mu$
might not be injective).  We will show that the asymmetric signature
for flip-isomorphisms of $z$ also only depend on their behaviour on the linking forms
$L^\gamma$ and $L^\mu$ by constructing a new preformation $z'$ which is closely
linked to $z$ but fulfils all the requirements of Section \ref{linkasymsec}
(\ie it is regular and its maps are $S$-isomorphisms).

\begin{lem}
\label{intfliplem}
\begin{enumerate}
\item $\ker\gamma=\ker\mu$
\item $G/\ker\gamma$ is \fg free and of rank $\rk F$.
\item There is another regular split $\epsilon$-quadratic preformation $z'$ over $\bZ$ defined by
      $z'=\kreckf{F}{[\gamma]}{G/\ker\gamma}{[\mu]}{\psi}$ for which the following diagram commutes
        \begin{eqnarray*}
    \xymatrix
    {
    F&
    G \ar[l]^{\gamma} \ar[r]_{\mu} \ar@{>>}[d]_{\pi}&
    F^*\\
    & G/\ker\gamma \ar[ul]^{[\gamma]} \ar[ur]_{[\mu]}
    }
    \end{eqnarray*}
\item $z$ is elementary if and only if $z'$ is elementary
\item Every flip-isomorphism $t=(\alpha, \beta, \sigma)$ of the $\epsilon$-quadratic
      formation $\skreckfs$ induces a flip-isomorphism $t'=(\alpha, \beta', \sigma)$
      of $z'$. Then $\sigma^*(z,t)=\sigma^*(z',t')\in LAsy^0(\bZ)$.
\end{enumerate}
\end{lem}
\begin{proof}
There is a free \fg submodule $G'$ such that $G=G'\oplus \ker\gamma$. There can be
no torsion in $G'$  because $\gamma$ is a homomorphism into a free module. We write
\begin{eqnarray*}
\gamma &=& \mat{\gamma_1 & 0}  \colon G=G'\oplus \ker\gamma \ra F\\
\mu    &=& \mat{\mu_1 & \mu_2} \colon G=G'\oplus \ker\gamma \ra F^*
\end{eqnarray*}
$\gamma_1$ is obviously injective and hence induces an isomorphism over $\bQ$.
Therefore $G'$ has the same rank as $F$. As $\gamma^*\mu$ is $(-\epsilon)$-symmetric,
$\mu_2$ must vanish. Then $\mu_1$ must be injective as well.
Hence $z'=\kreckf{F}{\gamma_1}{G'}{\mu_1}{\theta|G'}$ is a well-defined preformation
with all the claimed properties.

By Lemma \ref{pushelemlem} $z$ is elementary if and only if $z'$ is elementary.

Let $t=(\alpha,\beta,\chi)$ be a flip-isomorphism of $z$. Obviously
$\beta(\ker\gamma)=\ker\mu=\ker\gamma$. So we can write 
$$\beta=\mat{\beta_1&0\\\beta_2&\beta_3}\colon G=G'\oplus\ker\gamma \ra G=G'\oplus\ker\gamma$$
It follows that $t'=(\alpha, \beta_1, \sigma)$ is a flip-isomorphism of $z'$. The Definition \ref{asymsignflipdef}
of the asymmetric signature is independent of $\beta$ and $G$, hence the signatures of $t$ and $t'$
are the same.
\end{proof}

The lemma justifies the following generalization of asymmetric signatures:

\begin{defi}
Let $l\colon L^\mu \stackrel{\cong}{\ra} L^\gamma$ be an isomorphism of split $(-\epsilon)$-quadratic linking forms.
Let $z'$ be the regular $\epsilon$-quadratic split preformation constructed in Lemma \ref{intfliplem}.
Then $\tilde\sigma^*(z,l)=\sigma^*(z',l)\in LAsy^0(\bZ)$ is the {\bf asymmetric signature of $z$ and the
isomorphism $l\colon L^\mu \stackrel{\cong}{\ra} L^\gamma$ of linking forms}\index{Asymmetric signature!over $\bZ$}.
\hfill\qed\end{defi}

\begin{thm}
\label{simplelinkasymthm}
If $[z]\in {l'}^{2q+2}(\bZ)$ is elementary then the asymmetric signature $\tilde\sigma^*(z,l)\in LAsy^0(\bZ)$ vanishes for all
isomorphisms $l\colon L^\mu \stackrel{\cong}{\ra} L^\gamma$ of linking forms.
\end{thm}
\begin{proof}
Follows from Lemma \ref{intfliplem}, Lemma \ref{pushelemlem} and Theorem \ref{elemlinkasymthm}.
\end{proof}

\begin{defi}[\cite{Ran02} Example 12.44]
\begin{enumerate}
\item   Let $M$ be a $(2q+1)$-dimens\-ion\-al manifold. The {\bf linking form on M}\index{Linking form!topological} is the
    bilinear form on the torsion submodule of $H_q(M)$ and $H_q(M,\partial M)$ given by
    \begin{eqnarray*}
    l_M\colon TH_q(M) \times TH_q(M,\partial M) &\ra& \bQ/\bZ\\
    (x,y) &\mt& \frac{1}{s}<z,y>
    \end{eqnarray*}
    with $z\in C^q(M,\partial M)$ and $s\in\bZ\setminus\{0\}$ such that $sx=d(z\cap[M])\in C_q(M)$.
\item   Let $M \ra B$ be a map of $(2q+1)$-dimensional manifold in a topological space.
        The {\bf B-linking form on M}\index{Linking form! $B$-} is the $(-\epsilon)$-symmetric form on the torsion submodule
        of $H_{q+1}(B,M)$ given by
    \begin{eqnarray*}
    l_M^B\colon TH_{q+1}(B,M) \times TH_{q+1}(B,M) &\ra& \bQ/\bZ\\
    (x,y) &\mt& l_M(p(x),p(y))
    \end{eqnarray*}
    with $p\colon H_{q+1}(B,M) \ra H_q(M)$.
    \hfill\qed
\end{enumerate}
\end{defi}

\begin{rem}
If $\partial M=\emptyset$ then $l_M$ is a non-singular $(-\epsilon)$-symmetric linking form on $TH_q(M)$.
\end{rem}

\begin{prop}
\label{mfdlinkprp}
We repeat the assumptions of Kreck's surgery theory in the simply-connected case:
Let $p\colon B \rightarrow BO$ be a fibration with $\pi_1(B)=0$.  Let $M_i$ be
$(2q+1)$-dimensional manifolds with a {\bf $(q-1)$-smoothings in $B$}  \ie a
lift of the stable normal bundle over $p$ which is $q$-connected.
Let $f\colon \partial M_0 \rightarrow \partial M_1$ be a diffeomorphism compatible
with the smoothings. Let $W$ be a cobordism of $M_0 \cup_f M_1$ with a
compatible $q$-smoothing over $B$. 
As in Corollary \ref{niceobstrcor} we define an obstruction
\begin{eqnarray*}
x(W)&=&\kreckfs\\
&=&(H_{q+1}(W,M_0)\longleftarrow H_{q+2}(B,W) \ra H_{q+1}(W,M_1), \theta)\\
  &\in& l'_{2q+2}(\bZ)
\end{eqnarray*}

If $\coker\gamma=H_{q+1}(B,M_0)$ is finite then $L^\gamma=-l_{M_0}^B$.

If $\coker\mu   =H_{q+1}(B,M_1)$ is finite then $L^\mu = -l_{M_1}^B$.

Assume both cokernels are finite.
If $W$ is bordant rel$\partial$ to an $h$-cobordism then there exist
isomorphisms $l\colon L^\mu=-l^B_{M_1}\isora L^\gamma=\epsilon -l^B_{M_0}$ 
and their asymmetric signatures $\tilde\sigma^*(x(W),l)\in LAsy^0(\bZ)$
will all vanish.
\end{prop}
\begin{proof}
The complex $\tilde C_{q+2}=H_{q+2}(B,W) \stackrel{\gamma}{\ra} \tilde C_{q+1}=H_{q+1}(W,{M_0})$
has homology $H_*(\tilde C)=H_*(B,{M_0})$. There is a homotopy equivalence $m\colon \tilde C \ra C(B,{M_0})$
and there is a chain map $C(B,{M_0}) \ra C_{*-1}({M_0})$ which induces the connecting homomorphism
$\partial_*\colon H_*(B,{M_0}) \ra H_{*-1}({M_0})$. Both maps together yield a chain map
\begin{eqnarray*}
\xymatrix
{
C_{q+1}({M_0}) \ar[r]^{d}&
C_q({M_0})\\
\tilde C_{q+2} \ar[u]^{p} \ar[r]_{\gamma}&
\tilde C_{q+1} \ar[u]^{p}
}
\end{eqnarray*}
which induces the connecting map $p\colon H_{q+1}(B,{M_0}) \ra H_q({M_0})$.

Let $a,b \in \coker\gamma=H_{q+1}(B,{M_0})=H_{q+1}(\tilde C)$. Represent both
homology classes by chains $\bar a, \bar b \in \tilde C_{q+1}$. Then there
is a $g\in\tilde C_{q+2}$ and an $s\in\bZ\setminus\{0\}$ such that 
$s\bar a=\gamma(g)$. Let $z\in C^q({M_0},\partial {M_0})$ such that $p(g)=z\cap[{M_0}]$.
Then $sp(\bar a)=d(z\cap[{M_0}])$. 
Hence $l_{M_0}^B(a,b)=\frac{1}{s}\langle z, p(\bar b)\rangle $. Let $b' \in H^{q+1}(W,{M_0}')$
such that $b'\cap[W]=\bar b$.
Then $l_{M_0}^B(a,b)=\frac{1}{s}\langle p^*z, b'\cap[W]\rangle =-\epsilon\frac{1}{s}\langle b', p^*(z)\cap[W]\rangle $.
Since $p$ is a connecting homomorphism $p^*(z)\cap [W]=-\epsilon i(z \cap [{M_0}])=-\epsilon ip(g)=-\epsilon\mu(g)$.
Hence $l_{M_0}^B(a,b)=\frac{1}{s}\langle b',\mu(g)\rangle =-\epsilon\frac{1}{s}\mu^*(b)(g)=-L^\gamma(a,b)$.
\end{proof}

\appendix
\chapter{A crash course in algebraic surgery theory}
\label{algsurchap}

This chapter is a compilation of the main theorems and constructions 
of algebraic surgery theory, taken from \cite{Ran80a} or \cite{Ran81}.

{\bf Throughout this chapter
let $\Lambda$ be a weakly finite ring with $1$ and an involution
and let $\epsilon\in\Lambda$ such that $\bar\epsilon=\epsilon^{-1}$ (\eg $\epsilon=\pm 1$).}

\section{Quadratic and symmetric complexes}
\label{cplxsec}

\begin{defi}
A {\bf chain complex $C$ (over $\Lambda$)}\index{Chain map} is a collection of homomorphisms of 
\fg free $\Lambda$-modules $\{d_r\colon C_r \ra C_{r-1}|r\in\bZ\}$ such that 
$d_r d_{r+1}=0 \colon C_{r+1}\ra C_{r-1}$ for all $r$.
$C$ is called $n$-dimensional if $C_r=0$ for $r<0$ and $r>n$.

Its {\bf homology $\Lambda$-modules $H_*(C)$}\index{Homology} are defined by
$$
    H_r(C)=\ker(d\colon C_r \ra C_{r-1})/\im(d\colon C_{r+1} \ra C_r).
$$
Its {\bf cohomology $\Lambda$-modules $H^*(C)$}\index{Cohomology} are defined by
$$
    H^r(C)=\ker(d^*\colon C^r \ra C^{r+1})/\im(d^*\colon C^{r-1} \ra C^r).
$$

A {\bf chain map $f\colon C \ra D$}\index{Chain map} of chain complexes over $\Lambda$ is 
a collection of $\Lambda$-module morphisms $\{f_r \colon C_r \ra D_r| r\in\bZ\}$
such that $d_D f_r = f_{r-1} d_C \colon C_r \ra D_{r-1}$ for all $r$.

A {\bf chain homotopy $g\colon f \simeq f' \colon C\ra D$ of two chain maps $f$ and $f'$}\index{Chain homotopy}
is a collection of $\Lambda$-module morphisms $\{g_r \colon C_{r-1} \ra D_r| r\in\bZ\}$
such that $f'_r-f_r= d_D g_{r+1} + g_r d_C \colon C_r \ra D_r$.

A chain map is a {\bf chain equivalence}\index{Chain equivalence} if it has a chain homotopy inverse.
It is an {\bf isomorphism} if it consists of isomorphism of modules only.

The {\bf mapping cone $\cone(f)$ of a chain map $f\colon C \ra D$}\index{Mapping cone} is the chain complex given by
\begin{eqnarray*}
d_\cone &=& \mat{d_D & (-)^{r-1}f\\0& d_C} \colon \cone(f)_r=D_r\oplus C_{r-1} \ra \cone(f)_{r-1}=D_{r-1}\oplus C_{r-2}
\end{eqnarray*}
\hfill\qed\end{defi}

\begin{defi}
\label{qgrpdef}
Let $C$ be a chain complex.
The {\bf $\epsilon$-duality involution $T_\epsilon$}\index{Duality involution map}\index{Tepsilon@$T_\epsilon$} 
is defined by
\begin{eqnarray*}
    T_\epsilon \colon \Hom_\Lambda(C^p,C_q) &\ra& \Hom_\Lambda(C^q, C_p)\\
    \psi &\mt& (-)^{pq}\epsilon\psi^*
\end{eqnarray*}

We define new chain complexes $W^\%(C,\epsilon)$ and $W_\%(C,\epsilon)$ by
\index{WC1@$W^\%(C,\epsilon)$}\index{WC2@$W_\%(C,\epsilon)$}
\begin{eqnarray*}
W^\%(C,\epsilon)_n &=& \{\phi_s \colon C^{n-r+s} \ra C_r | r\in\bZ, s \ge 0\}\\
d^\% \colon W^\%(C,\epsilon)_n &\ra& W^\%(C,\epsilon)_{n-1}\\
\{\phi_s\} &\mt& \{d\phi_s + (-)^r\phi_s d^* + (-)^{n+s-1}(\phi_{s-1} + (-)^s T_\epsilon\phi_{s-1}) \colon\\ 
&&C^{n-r+s-1} \ra C_r | r\in\bZ, s \ge 0\}\\
&&\text{where we set $\phi_{-1}=0$.}
\\
W_\%(C,\epsilon)_n &=& \{\psi_s \colon C^{n-r-s} \ra C_r | r\in\bZ, s \ge 0\}\\
d_\% \colon W_\%(C,\epsilon)_n &\ra& W_\%(C,\epsilon)_{n-1}\\
\{\psi_s\} &\mt& \{d\psi_s + (-)^r\psi_s d^* + (-)^{n-s-1}(\psi_{s+1} + (-)^{s+1} T_\epsilon\psi_{s+1}) \colon \\
&&C^{n-r-s-1} \ra C_r | r\in\bZ, s \ge 0\}
\end{eqnarray*}

\index{QC1@$Q^n(C,\epsilon)$} \index{QC2@$Q_n(C,\epsilon)$}
Their homology groups are the {\bf $\epsilon$-symmetric $Q$-groups $Q^n(C,\epsilon)=H_n(W^\%(C,\epsilon))$} and
the {\bf $\epsilon$-quadratic $Q$-groups $Q_n(C,\epsilon)=H_n(W_\%(C,\epsilon))$}.
They are related by the {\bf $\epsilon$-symmetrization map}\index{Symmetrization map}
\begin{eqnarray*}
Q_n(C,\epsilon)     &\ra& Q^n(C,\epsilon)\\
\{\psi_s\} &\mt& \begin{cases} \{(1+T_\epsilon)\psi_0\}     &: \text{if $s=0$,}\\
                 0              &: \text{if $s\neq 0$.} 
           \end{cases}
\end{eqnarray*}
\hfill\qed\end{defi}

\begin{rem}
In the case of $\epsilon=1$ we omit the $\epsilon$ and we will simply speak of 
symmetric complexes, $Q_n(C)$, $W^\%(C)$, $T$ \etc
\end{rem}

\begin{defi}
\label{pcplxdef}
Let $C$ be a chain complex and $n\in\bN$. Define the chain complex $C^{n-*}$ by
$$
    d_{C^{n-*}}=(-)^r d^*_C \colon (C^{n-*})_r=C^{n-r}=C_{n-r}^* \ra  (C^{n-*})_{r-1}.
$$

An {\bf $\epsilon$-symmetric $n$-dimensional complex $(C,\phi)$}\index{Symmetric complex} is a chain complex $C$
together with a cycle $\phi\in W^\%(C,\epsilon)_n$. It is called {\bf Poincar\'e}\index{Symmetric complex!Poincar\'e}
if the {\bf Poincar\'e duality map}\index{Symmetric complex!Poincar\'e duality map}
$$
    \phi_0 : C^{n-*} \ra C
$$
is a chain equivalence.

An {\bf $\epsilon$-quadratic $n$-dimensional complex $(C,\psi)$}\index{Quadratic complex} is a chain complex $C$
together with a cycle $\psi\in W_\%(C,\epsilon)_n$. 
It is called {\bf Poincar\'e}\index{Quadratic complex!Poincar\'e} 
if the {\bf Poincar\'e duality map}\index{Quadratic complex!Poincar\'e duality map}
$$
    (1+T_\epsilon)\psi_0 : C^{n-*} \ra C
$$
is a chain equivalence.

A {\bf morphism of $\epsilon$-symmetric $n$-dimensional complexes $f=(f,\rho)\colon (C,\phi) \ra (C',\phi')$} 
\index{Symmetric complex!morphism}
is a chain map $f\colon C \ra C'$ and a $\rho \in W^\%(C',\epsilon)_{n+1}$ such that $\phi'-f\phi f^* = d^\%(\rho)$
\ie
\begin{eqnarray*}
\phi'_s - f\phi_s f^*= d\rho_s + (-)^r\rho_s d^* + (-)^{n+s}(\rho_{s-1}+(-)^s T_\epsilon \rho_{s-1}) \colon C^{n-r+s}\rightarrow C_r
\end{eqnarray*}
\ie $f\phi f^*=\phi' \in Q^n(C)$. It is an {\bf equivalence} if $f\colon C\ra C'$ is a chain equivalence.\index{Symmetric complex!equivalence}
It is an {\bf isomorphism} if $f\colon C\ra C'$ is an isomorphism.\index{Symmetric complex!isomorphism}

A {\bf map of $\epsilon$-quadratic $n$-dimensional complexes $f=(f,\sigma)\colon (C,\psi) \ra (C',\psi')$}
\index{Quadratic complex!morphism} 
is a chain map $f\colon C \ra C'$ and a $\sigma \in W_\%(C',\epsilon)_{n+1}$ such that $\psi'-f\psi f^* = d_\%(\sigma)$
\ie
\begin{eqnarray*}
\psi'_s - f\psi_s f^*= d\sigma_s + (-)^r\sigma_s d^* + (-)^{n-s}(\sigma_{s+1}+(-)^{s+1} T_\epsilon \sigma_{s+1}) \colon C^{n-r-s}\rightarrow C_r
\end{eqnarray*}
\ie $ f\psi f^*=\psi' \in Q_n(C)$. It is an {\bf equivalence} if $f\colon C\ra C'$ is a chain equivalence.\index{Quadratic complex!equivalence}
It is an {\bf isomorphism} if $f\colon C\ra C'$ is an isomorphism.\index{Quadratic complex!isomorphism}
\hfill\qed\end{defi}

We can define compositions and inverses of morphisms. One can also define inverses for homotopy equivalences
but we will not need such a construction in this treatise.
\begin{defi}
\label{compinvdef}\index{Symmetric complex!morphism!composition}
The {\bf composition of two morphisms of $\epsilon$-symmetric $n$-dimensional complexes} 
$(f,\rho)\colon (C,\phi) \rightarrow (C',\phi')$ and $(f',\rho')\colon (C',\phi') \rightarrow (C'',\phi'')$
is the morphism $(f'f, \rho'+f'\rho {f'}^*)\colon (C,\phi) \rightarrow (C'',\phi'')$.

\index{Quadratic complex!morphism!composition}
The {\bf composition of two morphisms of $\epsilon$-quadratic $n$-dimensional complexes} 
$(f,\sigma)\colon (C,\psi) \rightarrow (C',\psi')$ and $(f',\sigma')\colon (C',\psi') \rightarrow (C'',\psi'')$
is the morphism $(f'f, \sigma'+f'\sigma {f'}^*)\colon (C,\psi) \rightarrow (C'',\psi'')$.

\index{Symmetric complex!isomorphism!inverse}
The {\bf inverse of an isomorphism $(f,\rho)\colon (C,\phi) \isora (C',\phi')$ of
$\epsilon$-symmetric $n$-dimensional complexes} is the isomorphism 
$(f,\rho)^{-1}=(f^{-1},-f^{-1}\rho f^{-*})\colon (C',\phi')\linebreak \isora (C,\phi)$.

\index{Quadratic complex!isomorphism!inverse}
The {\bf inverse of an isomorphism $(f,\sigma)\colon (C,\psi) \isora (C',\psi')$ of
$\epsilon$-quadratic $n$-dimensional complexes} is the isomorphism 
$(f,\sigma)^{-1}=(f^{-1},-f^{-1}\sigma f^{-*})\colon (C',\psi')\linebreak \isora (C,\psi)$.
\hfill\qed\end{defi}

\section{Quadratic and symmetric pairs}
\label{pairsec}

Whereas the algebraic equivalent of closed manifolds (respectively normal maps of closed manifolds) are symmetric Poincar\'e complexes
(respectively quadratic Poincar\'e complexes), the analogues of manifolds with boundaries or normal maps are
symmetric and quadratic pairs.

\begin{defi}
Let $f \colon C \ra D$ be a chain map.\index{Wf1@$W^\%(f,\epsilon)$}\index{Wf2@$W_\%(f,\epsilon)$}
We define chain complexes $W^\%(f,\epsilon)$ and $W_\%(f,\epsilon)$ by
\begin{eqnarray*}
W^\%(f,\epsilon)_{n+1} &=& \{(\delta\phi_s \colon D^{n-p+s+1}\rightarrow D_p, \phi_s \colon C^{n-r+s} \rightarrow C_r) | p,r\in\bZ, s \ge 0\}\\
d^\% \colon W^\%(f,\epsilon)_{n+1} &\rightarrow& W^\%(f,\epsilon)_{n}\\
\{(\delta\phi_s, \phi_s)\} &\mapsto& \{(d(\delta\phi_s) + (-)^r(\delta\phi_s)d^* + (-)^{n+s}(\delta\phi_{s-1}+(-)^s T_\epsilon(\delta\phi_{s-1})\\
&&+(-)^n f\phi_s f^*\colon D^{n-r+s} \rightarrow D_r,\\
&&d\phi_s + (-)^r\phi_s d^* + (-)^{n+s-1}(\phi_{s-1} + (-)^s T_\epsilon\phi_{s-1}) \colon\\ 
&&C^{n-r+s-1} \rightarrow C_r) | r\in\bZ, s \ge 0\}\\
&&\text{where we set $\phi_{-1}=0$ and $\delta\phi_{-1}=0$.}
\\
W_\%(f,\epsilon)_{n+1} &=& \{(\delta\psi_s \colon D^{n-p-s+1}\rightarrow D_p, \psi_s \colon C^{n-r-s} \rightarrow C_r) | p,r\in\bZ, s \ge 0\}\\
d_\% \colon W_\%(f,\epsilon)_{n+1} &\rightarrow& W_\%(f,\epsilon)_{n}\\
\{(\delta\psi_s, \psi_s)\} &\mapsto& \{(d(\delta\psi_s) + (-)^r(\delta\psi_s)d^* + (-)^{n-s}(\delta\psi_{s+1}+(-)^{s+1}T_\epsilon(\delta\psi_{s+1}))\\
&&+(-)^n f\psi_s f^* \colon D^{n-r-s} \rightarrow D_r,\\
&&d\psi_s + (-)^r\psi_s d^* + (-)^{n-s-1}(\psi_{s+1} + (-)^{s+1} T_\epsilon\psi_{s+1}) \colon\\ 
&&C^{n-r-s-1} \rightarrow C_r) | r\in\bZ, s \ge 0\}
\end{eqnarray*}

\index{Qf1@$Q^n(f,\epsilon)$}\index{Qf2@$Q_n(f,\epsilon)$}
Their homology groups are the {\bf $\epsilon$-symmetric $Q$-groups $Q^n(f,\epsilon)=H_n(W^\%(f,\epsilon))$} and
the {\bf $\epsilon$-quadratic $Q$-groups $Q_n(f,\epsilon)=H_n(W_\%(f,\epsilon))$}.
They are related by the {\bf $\epsilon$-symmetrization map}\index{Duality involution map}\index{Tepsilon@$T_\epsilon$}
\begin{eqnarray*}
Q_n(f,\epsilon)     &\ra& Q^n(f,\epsilon)\\
(\delta\psi_s, \psi_s) &\mt& \begin{cases} ((1+T_\epsilon)\delta\psi_0, (1+T_\epsilon)\psi_0) &: \text{if $s=0$,}\\
                        0                         &: \text{if $s\neq 0$.} 
                 \end{cases}
\end{eqnarray*}
\hfill\qed\end{defi}

\begin{defi}
\label{pairdef}
\index{Poincar\'e pair!symmetric}\index{Symmetric pair}\index{Symmetric pair!Poincar\'e}\index{Symmetric pair!Poincar\'e duality map}
An {\bf $(n+1)$-dimensional $\epsilon$-symmetric pair $(f\colon C \ra D, (\delta\phi,\linebreak\phi))$} is a chain map $f\colon C\ra D$
together with a cycle $ (\delta\phi,\phi) \in W^\%(f,\epsilon)_{n+1}$. It is called 
{\bf Poincar\'e} if the {\bf Poincar\'e duality map $D^{n+1-*} \ra \cone(f)$} given by
$$
    \mat{\delta\phi_0\\ (-)^{n+1-r}\phi_0 f^*} \colon D^{n+1-r} \ra \cone(f)_r
$$
is a chain equivalence.

\index{Poincar\'e pair!quadratic}\index{Quadratic pair}\index{Quadratic pair!Poincar\'e}\index{Quadratic pair!Poincar\'e duality map}
An {\bf $(n+1)$-dimensional $\epsilon$-quadratic pair $(f\colon C \ra D, (\delta\psi,\psi))$} is a chain map $f\colon C\ra D$
together with a cycle $(\delta\psi,\psi) \in W_\%(f,\epsilon)_{n+1}$. It is called 
{\bf Poincar\'e} if the {\bf Poincar\'e duality map $D^{n+1-*} \ra \cone(f)$} given by
$$
    \mat{(1+T_\epsilon)\delta\psi_0\\ (-)^{n+1-r}(1+T_\epsilon)\psi_0 f^*} \colon D^{n+1-r} \ra \cone(f)_r
$$
is a chain equivalence.
\hfill\qed\end{defi}

\begin{rem}
In the above definitions the Poincar\'e duality maps can be replaced by the chain maps
$$
    \mat{\delta\phi_0, f\phi_0 } \colon \cone(f)^{n+1-*} \ra D
$$
in the symmetric case and by
$$
    \mat{(1+T_\epsilon)\delta\psi_0, f(1+T_\epsilon)\psi_0 } \colon \cone(f)^{n+1-*} \ra D
$$
in the quadratic case.
\end{rem}

\begin{defi}
\label{homeqpairdef}\index{Quadratic pair!homotopy equivalence}\index{Homotopy equivalence}
A {\bf homotopy equivalence of $(n+1)$-dimensional $\epsilon$-quadra\-tic pairs}
$$
    (g,h;k) \colon (f\colon C \ra D, (\delta\psi,\psi)) \ra (f'\colon C' \ra D', (\delta\psi',\psi'))
$$
is a triple $(g,h;k)$ consisting of chain equivalences
$$
    g \colon C \ra C', \quad h \colon D \ra D'
$$
and a chain homotopy
$$
    k\colon f'g \simeq hf \colon C \ra D'
$$
such that
$$
    (g,h;k)_\%(\delta\psi,\psi)=(\delta\psi',\psi')\in Q_{n+1}(f',\epsilon)
$$
with
\begin{eqnarray*}
(g,h;k)_\%(\delta\psi,\psi)_s&=&(h\delta\psi_s h^* + (-)^n k\psi_s (hf)^* + (-)^{r+1}k T_\epsilon\psi_{s+1} k^* +\\
&&(-)^{n-r} f'g\psi_s k^*
\colon {D'}^{n+1-s-r} \ra D'_r,\\
&&g\psi_s g^* \colon {C'}^{n-s-q} \ra C'_q)
\end{eqnarray*}
\hfill\qed\end{defi}

Here are some useful lemmas about changing the boundary of a pair and examples for homotopy equivalences of pairs.
\begin{lem}
\label{deltapairlem}
Let $c=(f \colon C \ra D, (\delta\psi, \psi))$ be an $(n+1)$-dimensional $\epsilon$-quadratic pair
and $(g,\sigma) \colon (C',\psi') \ra (C,\psi)$ be a map of $n$-dimensional $\epsilon$-quadratic complexes.
Then 
$$
     c'=(f'=fg \colon C' \ra D, (\delta\psi'=\delta\psi+(-)^n f\sigma f^*, \psi'))
$$
is an $(n+1)$-dimensional $\epsilon$-quadratic pair. If $g$ is a chain equivalence, then $c$ and $c'$ are
homotopy equivalent pairs.
Same for the symmetric case.
\end{lem}

Complexes and pairs are in a one-to-one correspondence.
\begin{defi}
\label{thickthomdef}\index{Connected}\index{Symmetric complex!connected}
An $n$-dimensional $\epsilon$-symmetric complex $(C,\phi)$ is {\bf connected} if $H_0(\phi_0\colon C^{n-*}\ra C)=0$.
An $n$-dimensional $\epsilon$-quadratic complex $(C,\psi)$ is {\bf connected} if $H_0((1+T\epsilon)\psi_0\colon C^{n-*}\ra C)=0$.
\index{Quadratic complex!connected}

\index{Boundary!of a symmetric complex}\index{Symmetric complex!boundary}\index{$\partial(C,\phi)$}
The {\bf boundary $(\partial C,\partial\phi)$ of a connected $n$-dimensional $\epsilon$-symmetric complex $(C,\phi)$}
is the $(n-1)$-dimensional $\epsilon$-symmetric Poincar\'e complex defined by
\begin{eqnarray*}
d_{\partial C}&=&\mat{d_C & (-)^r\phi_0 \\ 0 & (-)^r d_C^*} \colon \partial C_{r}=C_{r+1}\oplus C^{n-r} \ra \partial C_{r-1}=C_r\oplus C^{n-r+1}
\\
\partial\phi_0&=&\mat{(-)^{n-r-1}T_\epsilon\phi_1 & (-)^{rn}\epsilon\\1 & 0} \colon 
\\
&&\quad\partial C^{n-r-1}=C^{n-r}\oplus C_{r+1} \ra\partial C_r=C_{r+1}\oplus C^{n-r}
\\
\partial\phi_s&=&\mat{(-)^{n-r+s-1}T_\epsilon\phi_{s+1}&0\\0&0} \colon 
\\
&&\quad\partial C^{n-r+s-1}=C^{n-r+s}\oplus C_{r-s+1} \ra
\partial C_r=C_{r+1}\oplus C^{n-r}
\quad (s>0)
\end{eqnarray*}

\index{Boundary!of a quadratic complex}\index{Quadratic complex!boundary}\index{$\partial(C,\psi)$}
The {\bf boundary $(\partial C,\partial\psi)$ of a connected $n$-dimensional $\epsilon$-quadratic complex $(C,\psi)$}
is the $(n-1)$-dimensional $\epsilon$-quadratic Poincar\'e complex defined by
\begin{eqnarray*}
d_{\partial C}&=&\mat{d_C & (-)^r(1+T_\epsilon)\psi_0 \\ 0 & (-)^r d_C^*} \colon 
\partial C_{r}=C_{r+1}\oplus C^{n-r} \ra \partial C_{r-1}=C_r\oplus C^{n-r+1}
\\
\partial\psi_0&=&\mat{0&0\\1 & 0} \colon 
\\
&&\partial C^{n-r-1}=C^{n-r}\oplus C_{r+1} \ra \partial C_r=C_{r+1}\oplus C^{n-r}
\\
\partial\psi_s&=&\mat{(-)^{n-r-s-1}T_\epsilon\psi_{s-1}&0\\0&0} \colon 
\\
&&\partial C^{n-r-s-1}=C^{n-r-s}\oplus C_{r+s+1} \ra \partial C_r=C_{r+1}\oplus C^{n-r}\quad (s>0)
\end{eqnarray*}

\index{Thickening!of a symmetric complex}\index{Symmetric complex!thickening}
The {\bf thickening of a connected $n$-dimensional $\epsilon$-symmetric complex $(C,\phi)$} is the $\epsilon$-symmetric
$n$-dimensional Poincar\'e pair $(i_C \colon \partial C \ra C^{n-*}, (0,\partial\phi))$ with 
$i_C=\mat{0&1}\colon \partial C_r = C_{r+1}\oplus C^{n-r} \ra (C^{n-*})_r=C^{n-r}$.

\index{Thickening!of a quadratic complex}\index{Quadratic complex!thickening}
The {\bf thickening of a connected $n$-dimensional $\epsilon$-quadratic complex $(C,\psi)$} is the $\epsilon$-quadratic
$n$-dimensional Poincar\'e pair $(i_C \colon \partial C \ra C^{n-*}, (0,\partial\psi))$ with 
$i_C=\mat{0&1}\colon \partial C_r = C_{r+1}\oplus C^{n-r} \ra (C^{n-*})_r=C^{n-r}$.

\index{Thom complex!of a symmetric pair}\index{Symmetric pair!Thom complex}
The {\bf Thom complex  of an $(n+1)$-dimensional $\epsilon$-symmetric Poincar\'e pair
$(f\colon C \ra D, (\delta\phi,\phi))$} is the connected $(n+1)$-dimensional $\epsilon$-symmetric complex
\linebreak
$(\cone(f), \delta\phi/\phi)$ given by
\begin{eqnarray*}
(\delta\phi/\phi)_s &=& \mat{\delta\phi_s & 0\\ (-)^{n+1-r}\phi_s f^* & (-)^{n-r+s+1}T_\epsilon\phi_{s-1}}\colon\\
&&\cone(f)^{n+1-r+s}=D^{n+1-r+s}\oplus C^{n-r+s} \ra \cone(f)_r=D_r\oplus C_{r-1}
\end{eqnarray*}

\index{Thom complex!of a quadratic pair}\index{Quadratic pair!Thom complex}
The {\bf Thom complex  of an $(n+1)$-dimensional $\epsilon$-quadratic Poincar\'e pair
$(f\colon C \ra D, (\delta\psi,\psi))$} is the connected $n$-dimensional $\epsilon$-quadratic complex
\linebreak
$(\cone(f), \delta\psi/\psi)$ given by
\begin{eqnarray*}
(\delta\psi/\psi)_s &=& \mat{\delta\psi_s & 0\\ (-)^{n+1-r}\psi_s f^* & (-)^{n-r-s}T_\epsilon\psi_{s+1}}\colon\\
&&\cone(f)^{n+1-r-s}=D^{n+1-r-s}\oplus C^{n-r-s} \ra \cone(f)_r=D_r\oplus C_{r-1}
\end{eqnarray*}
\hfill\qed\end{defi}

\begin{lem}[\cite{Ran80a} Proposition 3.4.]
\label{bdrmaplem}
Let $(f,\chi)\colon (C,\psi) \isora (C',\psi')$ be an isomorphism of $n$-dimensional $\epsilon$-quadratic
complexes. Then there is an isomorphism
\begin{eqnarray*}
(\partial f, \partial\chi)&\colon& (\partial C,\partial\psi) \isora (\partial C',\partial\psi')
\\
\partial f  = \mat{f & (-)^{r-1}(1+T_\epsilon)\chi_0 f^{-*}\\0&f^{-*}} &\colon& 
\partial C_r=C_{r+1}\oplus C^{n-r}\\&&\quad\isora \partial C'_r=C'_{r+1}\oplus {C'}^{n-r}
\\
\partial\chi_0  = 0 &\colon& \partial {C'}^{n-r} \ra \partial C'_r
\\
\partial\chi_s  = \mat{(-)^{n-r-s}T\chi_{s-1} & 0 \\ 0 & 0} &\colon& 
\partial {C'}^{n-r-s}={C'}^{n-r-s+1}\oplus {C'}_{r+s}\\ 
&&\quad\ra \partial C'_r=C'_{r+1}\oplus {C'}^{n-r}
\quad (s > 0)
\end{eqnarray*}
and a homotopy equivalence $(\partial f, f^{-*};0)$ between the thickening-ups of $(C,\psi)$ and $(C',\psi')$.
\end{lem}

\begin{prop}[\cite{Ran80a} Proposition 3.4.]
\label{thomthickprop}
The Thom complex operation induces an
natural one-to-one correspondence between the homotopy equivalence classes
of $n$-dimensional $\epsilon$-symmetric Poincar\'e pairs and 
the homotopy equivalence classes of connected $n$-dimensional $\epsilon$-symmetric complexes.
Poincar\'e pairs with contractible boundaries correspond to Poincar\'e complexes.
Thickening is the inverse operation.
Similar for the quadratic case.
\end{prop}

\section{Unions of pairs}
\label{unionsec}

The union-construction is an algebraic analogue of glueing two $(n+1)$-dimensional cobordisms $(W,M, M')$ and
$(W', M', M'')$ together at $M'$.
\begin{defi}[\cite{Ran80a} p.135]
\label{uniondef}
\index{Union!of two symmetric pairs}\index{Symmetric pair!union}
The {\bf union of two adjoining $\epsilon$-symmetric $(n+1)$-dimensional cobordisms}
\begin{eqnarray*}
c&=&(\mat{f_C    & f_{C'}}  \colon C  \oplus C'  \ra D,  (\delta\phi,  \phi  \oplus -\phi'))\\
c'&=&(\mat{f'_{C'}& f'_{C''}}\colon C' \oplus C'' \ra D', (\delta\phi', \phi' \oplus -\phi''))
\end{eqnarray*}
is the $\epsilon$-symmetric $(n+1)$-dimensional cobordism
\begin{eqnarray*}
c \cup c' &=& (\mat{f''_{C'}& f''_{C''}}\colon C \oplus C'' \ra D'', (\delta\phi'', \phi \oplus -\phi''))
\end{eqnarray*}
given  by
\begin{eqnarray*}
d_{D''}&=&\mat{d_D & (-)^{r-1}f_{C'} & 0 \\ 0 & d_{C'} & 0 \\ 0 & (-)^{r-1} f'_{C'} & d_{D'}}\colon
\\&&\quad
D''_r=D_r\oplus C'_{r-1}\oplus D'_r \ra D''_{r-1}=D_{r-1}\oplus C'_{r-2}\oplus D'_{r-1}\\
f''_{C}   &=& \mat{f_C \\ 0 \\ 0} \colon C_r \ra D''_r\\
f''_{C''} &=& \mat{0 \\ 0 \\ f'_{C''}}\colon C''_r \ra D''_r\\
\delta\phi''_s&=&\mat{  \delta\phi_s            & 0                 & 0\\
            (-)^{n-r}\phi'_s f^*_{C'}   & (-)^{n-r+s}T_\epsilon \phi'_{s-1} & 0\\
            0               & (-)^s f'_{C'}\phi'_s          & \delta\phi'_s}\colon\\
&&{D''}^{n-r+s+1}=D^{n-r+s+1}\oplus {C'}^{n-r+s} \oplus {D'}^{n-r+s+1} 
\\&&\quad\ra D''_r=D_r \oplus C'_{r-1} \oplus D'_r
\end{eqnarray*}
We shall normally write
\begin{eqnarray*}
 D''=D\cup_{C'} D', \quad \delta\phi''=\delta\phi \cup_{\phi'} \delta\phi'
\end{eqnarray*}

\index{Union!of two quadratic pairs}\index{Quadratic pair!union}
The {\bf union of two adjoining $\epsilon$-quadratic $(n+1)$-dimensional cobordisms}
\begin{eqnarray*}
c &=&(\mat{f_C    & f_{C'}}  \colon C  \oplus C'  \ra D,  (\delta\psi,  \psi  \oplus -\psi'))\\
c'&=&(\mat{f'_{C'}& f'_{C''}}\colon C' \oplus C'' \ra D', (\delta\psi', \psi' \oplus -\psi''))
\end{eqnarray*}
is the $\epsilon$-quadratic $(n+1)$-dimensional cobordism
\begin{eqnarray*}
c \cup c' &=& (\mat{f''_{C'}& f''_{C''}}\colon C \oplus C'' \ra D'', (\delta\psi'', \psi \oplus -\psi''))
\end{eqnarray*}
given  by the same complex $D''$ and the same maps $f''_{C}$, $f''_{C''}$ as in the symmetric case and 
\begin{eqnarray*}
\delta\psi''_s&=&\mat{  \delta\psi_s            & 0                 & 0\\
            (-)^{n-r}\psi'_s f^*_{C'}   & (-)^{n-r-s+1}T_\epsilon \psi'_{s+1}   & 0\\
            0               & (-)^s f'_{C'}\psi'_s          & \delta\psi'_s}\colon\\
&&{D''}^{n-r-s+1}=D^{n-r-s+1}\oplus {C'}^{n-r-s} \oplus {D'}^{n-r-s+1} 
\\&&\quad\ra
D''_r=D_r \oplus C'_{r-1} \oplus D'_r
\end{eqnarray*}
We shall normally write
\begin{eqnarray*}
 D''=D\cup_{C'} D', \quad \delta\psi''=\delta\psi \cup_{\psi'} \delta\psi'
\end{eqnarray*}
\hfill\qed\end{defi}

Glueing and symmetrizing cobordisms are commutative operations as the following example illustrates.
\begin{lem}
\label{unionsymlem}
Let
\begin{eqnarray*}
c &=&(f_{C'}  \colon C' \ra D,  (\delta\psi, -\psi'))\\
c'&=&(f'_{C'} \colon C' \ra D', (\delta\psi', \psi'))
\end{eqnarray*}
be two $\epsilon$-quadratic $(n+1)$-dimensional Poincar\'e pairs. Then
\begin{eqnarray*}
(1, \chi) &\colon& (D''=D\cup_{C'} D', (1+T_\epsilon)(\delta\psi \cup_{\psi'} \delta\psi'))\\
&&\ra (D'', (1+T_\epsilon)(\delta\psi) \cup_{(1+T_\epsilon)\psi'} (1+T_\epsilon)(\delta\psi'))
\\
\chi_0=\mat{0&0&0\\0&(-)^{r-1}T\psi_0&0\\0&0&0}&\colon& {D''}^{n+2-r}=D^{n+2-r}\oplus {C'}^{n+1-r} \oplus {D'}^{n+2-r}
\\
&&\quad\ra D''_r=D_r\oplus C'_{r-1}\oplus D'_r
\end{eqnarray*}
is an isomorphism of $(n+2)$-dimensional $\epsilon$-symmetric Poincar\'e complexes.
\end{lem}

Next we show that changing the common boundary of two pairs doesn't change their union
\begin{lem}
\label{unionbdrylem}
Let
\begin{eqnarray*}
c &=&(f  \colon C \ra D,  (\delta\psi,  \psi))\\
c'&=&(f' \colon C \ra D', (\delta\psi', \psi))
\end{eqnarray*}
be two $\epsilon$-quadratic $(n+1)$-dimensional Poincar\'e pairs.
Let $(h,\chi) \colon (\widehat{C}, \widehat{\psi}) \stackrel{\simeq}{\ra} (C, \psi)$ be an equivalence of quadratic complexes.
Using Lemma \ref{deltapairlem}, define the $(n+1)$-dimensional $\epsilon$-quadratic Poincar\'e pairs
\begin{eqnarray*}
\widehat{c}&=& (\widehat{f} = fh \colon \widehat{C}\ra D, (\widehat{\delta\psi}=\delta\psi+(-)^n f\chi f^*, \widehat{\psi}))\\
\widehat{c}'&=& (\widehat{f}' = f'h \colon \widehat{C}\ra D', (\widehat{\delta\psi}'=\delta\psi'+(-)^n f'\chi {f'}^*, \widehat{\psi}))
\end{eqnarray*}
Then there is an chain equivalence of $(n+2)$-dimensional $\epsilon$-quadratic Poincar\'e complexes
\begin{eqnarray*}
\left(\mat{1&0&0\\0&h&0\\0&0&1},\sigma\right) &\colon& 
\widehat{c}\cup -\widehat{c}'=(D \cup_{\widehat{C}} D',  \widehat{\delta\psi}\cup_{\widehat{\psi}} -\widehat{\delta\psi}')
\\
&&\quad\stackrel{\simeq}{\ra} c\cup -c = (D\cup_C D', \delta\psi \cup_\psi -\delta\psi')
\\
\sigma_s &=& \mat{  0           & 0                 & 0\\
            (-)^{n-r}\chi_s f^* & (-)^{n+1-r-s}T_\epsilon\chi_{s+1} & 0\\
            0           & (-)^{s-1} f' \chi_s           & 0} \colon
\\ 
&&(D\cup_C D')^{n+2-r-s}=D^{n+2-r-s}\oplus C^{n+1-r-s}\oplus {D'}^{n+2-r-s}
\\
&&\quad \ra (D\cup_C D')_r = D_r\oplus C_{r-1}\oplus D'_{r-1}
\end{eqnarray*}
\end{lem}
 
Another construction we will use is the union of a {\bf fundamental pair}.
The geometrical analogue can be described as such: Let $(W,M,M)$ be
an $(n+1)$-dimensional cobordism and glue it together along its boundaries.
Using Mayer-Vietoris one sees that the resulting manifold $V$ has the chain complex $C(f-g)$ with
$(f,g)\colon M+M \ra W$ being the inclusion of the boundary. 
\begin{deflem}[\cite{Ran98} Definition 24.1.]
\label{unionfunddefi}\index{Fundamental pair}\index{Symmetric pair!fundamental}\index{Quadratic pair!fundamental}
An $(n+1)$-dimensional $\epsilon$-sym\-metric pair is called {\bf fundamental}
if it is of the form $((f,g)\colon C\oplus C \ra D, (\delta\phi, \phi\oplus -\phi))$.
Similar for the quadratic case.

\index{Symmetric pair!fundamental!union}\index{Fundamental pair!union}\index{Union!of a symmetric fundamental pair}
The {\bf union of a fundamental $(n+1)$-dimensional $\epsilon$-symmetric Poincar\'e pair
$((f,g)\colon C\oplus C \ra D, (\delta\phi, \phi\oplus -\phi))$ over $\Lambda$} is the 
$(n+1)$-dimensional $\epsilon$-symmetric Poincar\'e $(U,\rho)$ complex over $\Lambda$
given by
\begin{eqnarray*}
U&=&\cone(f-g \colon C \ra D)
\\
\rho_s&=&\mat{\delta\phi_s & (-)^s g\phi_s\\ (-)^{n-r+1}\phi_s f^* & (-)^{n-r+s+1}T_\epsilon\phi_{s-1}}\colon
\\
&&U^{n+1-r+s}= D^{n+1-r+s}\oplus C^{n-r+s} \ra U_r=D_{r}\oplus C_{r-1}
\end{eqnarray*}

\index{Quadratic pair!fundamental!union}\index{Fundamental pair!union}\index{Union!of a quadratic fundamental pair}
The {\bf union of a fundamental $(n+1)$-dimensional $\epsilon$-quadratic Poincar\'e pair
$((f,g)\colon C\oplus C \ra D, (\delta\psi, \psi\oplus -\psi))$ over $\Lambda$} is the 
$(n+1)$-dimensional $\epsilon$-symmetric Poincar\'e $(U,\sigma)$ complex over $\Lambda$
given by
\begin{eqnarray*}
U&=&\cone(f-g \colon C \ra D)
\\
\sigma_s&=&\mat{\delta\psi_s & (-)^s g\psi_s\\ (-)^{n-r+1}\psi_s f^* & (-)^{n-r-s}T_\epsilon\psi_{s+1}}\colon
\\
&&U^{n+1-r-s}= D^{n+1-r-s}\oplus C^{n-r-s} \ra U_r=D_{r}\oplus C_{r-1}
\end{eqnarray*}
\end{deflem}
\begin{proof}
Compute the union $(D'',\phi'')$ of the two cobordism over $\Lambda$
\begin{eqnarray*}
c  &=& ((f,g) \colon C\oplus C \ra D, (\delta\phi, \phi\oplus -\phi))\\
c' &=& ((1,1) \colon C\oplus C \ra C, (0, -\phi \oplus \phi))
\end{eqnarray*}
(or the quadratic analogue). Then use the isomorphism of chain complexes
\begin{eqnarray*}
a_r&=&\mat{ 1 & 0 & 0 & -g\\
        0 & 1 & 0 & 0\\
        0 & 1 & 1 & (-)^{r+1}d\\
        0 & 0 & 0 & 1}
\end{eqnarray*}
\begin{eqnarray*}
\xymatrix@H+30pt@C-10pt
{
D''_r=D_r\oplus (C_{r-1}\oplus C_{r-1}) \oplus C_r
\ar[r]^-{a_r} 
\ar[d]_{\svec{  d & (-)^{r-1}f  & (-)^{r-1}g    & 0\\
        0 & d           & 0             & 0\\
        0 & 0       & d     & 0\\
        0 & (-)^{r-1}   & (-)^{r-1} & d}}
&
D''_r
\ar[d]^{\svec{  d & (-)^{r-1}(f-g)  & 0 & 0\\
        0 & d           & 0 & 0\\
        0 & 0           & 0 & 0\\
        0 & 0           & (-)^{r-1} & 0}}
\\
D''_{r-1}=D_{r-1}\oplus (C_{r-2}\oplus C_{r-2}) \oplus C_{r-1}
\ar[r]^-{a_{r-1}}&
D''_{r-1}
}
\end{eqnarray*}
\end{proof}

\section{Surgery on complexes}
\label{surgcplxsec}

\begin{defi}
\label{surgerydef}

Let $(C,\psi)$ be a connected $n$-dimensional $\epsilon$-quadratic complex
and $c=(f \colon C \ra D, (\delta\psi,\psi))$ an $\epsilon$-quadratic $(n+1)$-dimensional
pair.

\index{Connected}\index{Quadratic complex!connected}
$c$ is {\bf connected} if the zeroth homology of its Poincar\'e duality map vanishes.

\index{Surgery on quadratic complexes}\index{Quadratic complex!surgery}
The {\bf result of an $\epsilon$-quadratic surgery on a connected pair $c$} is
the connected $n$-dimensional $\epsilon$-quadratic complex $(C',\psi')$ given by
\begin{eqnarray*}
d_{C'}&=&\mat{  d_C     & 0     & (-)^{n+1}(1+T_\epsilon)\psi_0 f^*\\
        (-)^r f & d_D   & (-)^r(1+T_\epsilon)\delta\psi_0\\
        0   & 0 & (-)^r d_D^*}\colon
\\
&&C'_r=C_r\oplus D_{r+1}\oplus D^{n-r+1} \ra C'_{r-1}=C_{r-1}\oplus D_r\oplus D^{n-r+2}
\\
\psi'_0&=&\mat{ \psi_0              & 0                 & 0\\
        0               & 0                 & 0\\
        0               & 1                 & 0}\colon
\\
&&{C'}^{n-r}=C^{n-r}\oplus D^{n-r+1}\oplus D_{r+1} \ra C'_r = C_r\oplus D_{r+1}\oplus D^{n-r+1}
\\
\psi'_s&=&\mat{ \psi_s              & (-)^{s}T_\epsilon\psi_{s-1}f^*    & 0\\
        0               & (-)^{n-r-s}T_\epsilon\delta\psi_{s-1}& 0\\
        0               & 0                 & 0}\colon
\\
&&{C'}^{n-r-s}=C^{n-r-s}\oplus D^{n-r-s+1}\oplus D_{r+s+1} 
\\
&&\quad\ra C'_r = C_r\oplus D_{r+1}\oplus D^{n-r+1}
\quad(s > 0)
\end{eqnarray*}
Similar for the symmetric case.
\hfill\qed\end{defi}

In an obvious way we can introduce the notion of an cobordism of complexes.
It turns out to be an equivalence relation and its equivalence classes are the
surgery $L$-groups.
\begin{defi}
\label{cobodef}\index{Cobordism of complexes}\index{Quadratic complexes!cobordism}
A {\bf cobordism of two $n$-dimensional $\epsilon$-quadratic Poincar\'e complexes $(C,\psi)$ 
and $(C',\psi')$} is an $(n+1)$-dimensional $\epsilon$-symmetric Poincar\'e pair
$(f \colon C\oplus C' \ra D, (\delta\psi, \phi\oplus-\psi'))$.
Similar for the symmetric case.
\hfill\qed\end{defi}

The well-known relations between surgery and cobordism hold also in the algebraic world
\begin{prop}[\cite{Ran80a} Proposition 4.1.]
\label{surgbordprop}
\begin{enumerate}
\item Algebraic surgery preserves the homotopy type of the boundary, sending Poincar\'e 
complexes to Poincar\'e complexes.
\item Two Poincar\'e complexes are cobordant if and only if the one can obtained from
the other by surgery and homotopy equivalence.
\end{enumerate}
\end{prop}

\begin{prop}[\cite{Ran80a} Proposition 3.2]
\label{cobordprop}
Cobordism is an equivalence relation on $n$-dimensional Poincar\'e complexes.
Homotopy equivalent Poincar\'e complexes are cobordant. 
\end{prop}

The cobordism classes of Poincar\'e complexes are groups and are the preferred definition for the $L$-groups
amongst algebraic surgeons because they are related to the $L$-groups defined as Witt-groups of forms and
formations. See \cite{Ran80a} Chapter 4 and 5 for more details. 
For our purposes the only fact we really need is the following lemma:
\begin{lem}
\label{trcplxlem}
Let $(C,\psi)$ be an $2m$-dimensional quadratic Poincar\'e complex with $C_i=0$ for $i\not\in\{m+1,m\}$.
Then  $$\left(M=C^{m}\oplus C_{m+1}, \theta=\mat{\psi_0&0\\d^*&0}\in Q_{(-)^m}(M)\right)$$
is a non-singular $(-)^m$-quadratic form.

$(C,\psi)$ is null-cobordant if and only if $(M,\theta)=0 \in L_{2m+2}(\Lambda)$.
Homotopic or cobordant $\epsilon$-quadratic complexes lead to the same element in $L_{2m+2}(\Lambda)$.
\end{lem}
\begin{rem}
This is a special case of the instant surgery obstruction given in \cite{Ran80a}.
A similar result does not hold in general for symmetric complexes.
\end{rem}
\begin{proof}
Define the connected $(2m+1)$-dimensional quadratic pair $(f\colon C \ra D, (0,\psi))$
with $f=1\colon C_{m+1}\ra D_{m+1}=C_{m+1}$ and $D_i=0$ for $i\neq m+1$.
We simplify the result $C'$ of the surgery on $C$ using the homotopy equivalence
\begin{eqnarray*}
\xymatrix@+30pt
{
C'_{m+1}=C_{m+1} \ar[d]_-{\svec{d\\(-)^{m+1}\\0}}
&
\\
C'_m=C_m\oplus D_{m+1} \oplus D^{m+1} \ar[r]_-{\svec{(-)^m&d&0\\0&0&1}}
&
M^*
}
\end{eqnarray*}
\cite{Ran80a} Proposition 4.3. and Proposition 5.1. finish the proof
\end{proof}

\backmatter
\bibliographystyle{alpha}

\printindex
\end{document}